\documentclass{amsproc}
\usepackage{breqn,float,amssymb,pdflscape,rotating,hyperref,xcolor,setspace,graphicx}
\makeatletter
\@namedef{subjclassname@2020}{%
  \textup{2020} Mathematics Subject Classification}
\makeatother
\newtheorem{theorem}{Theorem}[section]

\theoremstyle{definition}

\newtheorem{example}[theorem]{Example}

\theoremstyle{remark}

\numberwithin{equation}{section}
\allowdisplaybreaks
\begin{document}
\onehalfspacing
%\doublespacing
%
\title[Definite integrals involving Bessel functions]{Definite integrals involving Bessel functions expressed as a series of special functions}

%    Remove any unused author tags.

%    author one information
\author{Robert Reynolds}
\address[Robert Reynolds]{Department of Mathematics and Statistics, York University, Toronto, ON, Canada, M3J1P3}
\email[Corresponding author]{milver73@gmail.com}
\thanks{}

%%    author two information
%\author{ Allan Stauffer}
%\address[Allan Stauffer]{Department of Mathematics and Statistics, York University, Toronto, ON, Canada, M3J1P3}
%\email{stauffer@yorku.ca}
%\thanks{This research is supported by NSERC Canada under Grant 504070}

\subjclass[2020]{Primary  30E20, 33-01, 33-03, 33-04}

\keywords{Bessel function, contour integration, Hurwitz-Lerch zeta function}

\date{}

\dedicatory{}

\begin{abstract}
Series involving hypergeometric functions are used to derive, extend and evaluate integrals involving the product of two Bessel functions of the first kind $J_{u}(a z)$ $J_{v}(b z)$ with order $u,v$, studied by Landau et al. The method used in this work is contour integration.
\end{abstract}

\maketitle
\section{Preliminaries}
Throughout this work we will use the Hurwitz-Lerch zeta function $\Phi(z,s,a)$ given in [DLMF,\href{https://dlmf.nist.gov/25.14}{25.14}], and its special cases; the Hurwitz zeta function $\zeta(s,a)$, [DLMF,\href{https://dlmf.nist.gov/25.11}{25.11}] and the polylogarithm $Li_{s}(z)$, [DLMF,\href{https://dlmf.nist.gov/25.12#ii}{25.12(ii)}], the log-gamma function given in [DLMF,\href{https://dlmf.nist.gov/25.11#vi.info}{25.11.18}], Catalan's constant $C$, [DLMF,\href{https://dlmf.nist.gov/26.6.E12}{25.11.40}], Glaisher's constant $A$, [DLMF,\href{https://dlmf.nist.gov/5.17.E5}{5.17.5}], Apery's constant $\zeta(3)$, [Wolfram,\href{https://mathworld.wolfram.com/AperysConstant.html}{1}], the Pochhammer's symbol $(a)_{n}$, [DLMF,\href{https://dlmf.nist.gov/5.2.iii}{5.2.5}], Stieltjes constant $\gamma_{n}$, [DLMF,\href{https://dlmf.nist.gov/25.2.E5}{25.2.5}]
\begin{equation}
\Phi'(i,0,u)=\log \left(\frac{\Gamma \left(\frac{u}{4}\right)}{2 \Gamma
   \left(\frac{u+2}{4}\right)}\right)+i \log \left(\frac{\Gamma \left(\frac{u+1}{4}\right)}{2 \Gamma
   \left(\frac{u+3}{4}\right)}\right)
\end{equation}
where $u\in\mathbb{C}$.
\begin{equation}
\Phi'(-i,0,a)=\log \left(\frac{\Gamma \left(\frac{a}{4}\right)}{2 \Gamma
   \left(\frac{a+2}{4}\right)}\right)-i \log \left(\frac{\Gamma \left(\frac{a+1}{4}\right)}{2 \Gamma
   \left(\frac{a+3}{4}\right)}\right)
\end{equation}
where $u\in\mathbb{C}$.
\begin{multline}
\text{Li}_s(z)=(2 \pi )^{s-1} \Gamma (1-s) \left(i^{1-s} \zeta \left(1-s,\frac{1}{2}-\frac{i \log (-z)}{2 \pi}\right)\right. \\ \left.
+i^{s-1} \zeta \left(1-s,\frac{i \log (-z)}{2 \pi }+\frac{1}{2}\right)\right)
\end{multline}
where $Re(s)>0$.
\begin{equation}\label{eq:hurwitz}
(2 i \pi )^{-s} \Gamma (s) \left((-1)^s \text{Li}_s\left(\frac{1}{z}\right)+\text{Li}_s(z)\right)=\zeta
   \left(1-s,\frac{\pi -i \log (-z)}{2 \pi }\right)
\end{equation}
where $Re(s)>0$.
\subsection{The Incomplete Gamma~Function}
The incomplete gamma functions are given in equation [DLMF,\href{https://dlmf.nist.gov/8.4.E13}{8.4.13}], $\gamma(a,z)$ and $\Gamma(a,z)$, are defined by
\begin{equation}
\gamma(a,z)=\int_{0}^{z}t^{a-1}e^{-t}dt
\end{equation}
and
\begin{equation}
\Gamma(a,z)=\int_{z}^{\infty}t^{a-1}e^{-t}dt
\end{equation}
where $Re(a)>0$. The~incomplete gamma function has a recurrence relation given by
\begin{equation}
\gamma(a,z)+\Gamma(a,z)=\Gamma(a)
\end{equation}
where $a\neq 0,-1,-2,..$. The~incomplete gamma function is continued analytically by
\begin{equation}
\gamma(a,ze^{2m\pi i})=e^{2\pi mia}\gamma(a,z)
\end{equation}
and
\begin{equation}\label{eq:7}
\Gamma(a,ze^{2m\pi i})=e^{2\pi mia}\Gamma(a,z)+(1-e^{2\pi m i a})\Gamma(a)
\end{equation}
where $m\in\mathbb{Z}$, $\gamma^{*}(a,z)=\frac{z^{-a}}{\Gamma(a)}\gamma(a,z)$ is entire in $z$ and $a$. When $z\neq 0$, $\Gamma(a,z)$ is an entire function of $a$ and $\gamma(a,z)$ is meromorphic with simple poles at $a=-n$ for $n=0,1,2,...$ with residue $\frac{(-1)^n}{n!}$. These definitions are listed in [DLMF,\href{https://dlmf.nist.gov/8.2.i}{8.2(i)}] and [DLMF,\href{https://dlmf.nist.gov/8.2.ii}{8.2(ii)}].
The incomplete gamma functions are particular cases of the more general hypergeometric and Meijer G functions see section (5.6) and equation (6.9.2) in \cite{erd}. 
Some Meijer G representations we will use in this work are given by;
\begin{equation}\label{g1}
\Gamma (a,z)=\Gamma (a)-G_{1,2}^{1,1}\left(z\left|
\begin{array}{c}
 1 \\
 a,0 \\
\end{array}
\right.\right)
\end{equation}
and
\begin{equation}\label{g2}
\Gamma (a,z)=G_{1,2}^{2,0}\left(z\left|
\begin{array}{c}
 1 \\
 0,a \\
\end{array}
\right.\right)
\end{equation}
from equations (2.4) and (2.6a) in \cite{milgram}.  We will also use the derivative notation given by;
\begin{equation}\label{g3}
\frac{\partial \Gamma (a,z)}{\partial a}=\Gamma (a,z) \log (z)+G_{2,3}^{3,0}\left(z\left|
\begin{array}{c}
 1,1 \\
 0,0,a \\
\end{array}
\right.\right)
\end{equation}
from equations (2.19a) in \cite{milgram}, (9.31.3) in \cite{grad} and equations (5.11.1), (6.2.11.1) and (6.2.11.2) in \cite{luke}, and (6.36) in \cite{aslam}.
\section{Introduction}
Definite integrals over finite intervals involving the product of Bessel functions are used in the field of diffraction theory. According to literature, in 1897, Carslaw \cite{carslaw}, corresponded with Professor Sommerfeld on his work of obtaining multiform solutions of differential equations of physical mathematics. Carslaw's work involved the solutions of second order partial differential equations applied to the two-dimensional problem of the diffraction and reflection of plane waves of light incident on an opaque semi-infinite plane bounded by a straight edge. 
Definite integrals over finite intervals involving the product of Bessel functions are used in the field of diffraction theory. Work showcasing such integrals are from authors Heins \cite{heins}, where his work shed further light on representation theorems for partial differential equations of a special form, in particular the two-dimensional wave equation. Naylor \cite{naylor}, where an integral transform adapted to the solution of certain boundary value problems connected with the Helmholtz equation in cylindrical or spherical polar coordinates when the radial variable r varies over some infinite interval. Miller \cite{miller}, worked on an integral, involving the functions $J_{n+½}$ and $J_{n+1/2}(x)$ and $J_{-n-1/2}(x)$ where n is a positive integer, arises in a problem in electrical network theory and represents a measure of the overall error in the diagonal Pad\'{e} approximant to the function $e^{-s}$. Ross \cite{ross}, worked on an integral equation arising in diffraction theory and many mixed boundary value problems in diffraction theory which could be reduced to the solution of dual integral equations. A modified form of the analysis of Schmelter and Lewin is given and a relatively simple form for the solution was also obtained.
\\\\
The contour integral representation for the generalized Landau-Luswili definite integral listed in \cite{landau}, where the definition of the Bessel function of the first kind is given in [DLMF,\href{https://dlmf.nist.gov/10.2.E2}{10.2.2}] is given by;
\begin{multline}\label{eq:def_int_cont}
\frac{1}{2\pi i}\int_{C}\int_{0}^{b}a^w w^{-1-k} x^{m+v+w+\mu } \left(1+c x^p\right)^{\beta } \lambda ^v (x \lambda )^{-v} \rho ^{\mu } (x \rho
   )^{-\mu } J_v(x \lambda ) J_{\mu }(x \rho )dxdw\\
=\frac{1}{2\pi i}
   \int_{C}\sum _{j=0}^{\infty } \sum _{l=0}^{\infty } \sum _{h=0}^{\infty }\frac{(-1)^{j+l} 2^{-2 j-2 l-v-\mu } a^w b^{1+2 j+2 l+m+h p+v+w+\mu } c^h w^{-1-k} }{(1+2 j+2 l+m+h p+v+w+\mu ) j! l! }\\
   \frac{\lambda ^v \left(\lambda
   ^2\right)^j \rho ^{\mu } \left(\rho ^2\right)^l \binom{\beta }{h}}{\Gamma(1+j+v) \Gamma (1+l+\mu )}dw
\end{multline}
where $Re(b)>0$.
where $|Re(m)|<1,w\in\mathbb{C},Re(b)>0,k\in\mathbb{C}$.\\\\
In this paper we derive the definite integral given by
\begin{multline}
\int_0^b x^m \left(1+c x^p\right)^{\beta } J_v(x \lambda ) J_{\mu }(x \rho ) \log ^k\left(\frac{1}{a x}\right)
   \, dx
=\frac{\left(\lambda ^v \rho ^{\mu }\right) }{2^{v+\mu } a^{1+m+v+\mu }}\\ \times\sum _{j=0}^{\infty } \sum _{l=0}^{\infty } \sum _{h=0}^{\infty }
   \frac{(-1)^{j+l} c^h \lambda ^{2 j} \rho ^{2 l} \binom{\beta }{h} }{2^{2 j+2 l} a^{2 j+2 l+h p} (1+2 j+2 l+m+h p+v+\mu )^{k+1} }\\ \times
   \frac{\Gamma \left(1+k,(1+2 j+2 l+m+h p+v+\mu ) \log
   \left(\frac{1}{a b}\right)\right)}{j! l! \Gamma (1+j+v)
   \Gamma (1+l+\mu )}
\end{multline}
where $Re(\beta)>0,Re(b)>0$.
where $Re(m)>0$. and the parameters $k,b,c$ are general complex numbers. This definite integral has a wider range of evaluation for the parameters relative to previous results in the latter. This integral can be evaluated in terms of the Hurwitz-Lerch zeta and incomplete gamma function. The derivation involves two definite integrals and follows the method used by us in~\cite{reyn4}. This method involves using a form of the generalized Cauchy's integral formula given by
\begin{equation}\label{intro:cauchy}
\frac{y^k}{\Gamma(k+1)}=\frac{1}{2\pi i}\int_{C}\frac{e^{wy}}{w^{k+1}}dw.
\end{equation}
where $C$ is in general an open contour in the complex plane where the bilinear concomitant has the same value at the end points of the contour. We then multiply both sides by a function of $x$, then take a definite integral of both sides. This yields a definite integral in terms of a contour integral. Then we multiply both sides of Equation~(\ref{intro:cauchy})  by another function of $y$ and take the infinite sum of both sides such that the contour integral of both equations are the same.
\section{Derivation of the contour integral representations}
\subsection{Left-hand side contour integral representation}
Using a generalization of Cauchy's integral formula \ref{intro:cauchy}, we form the definite integral by replacing $y$ by $\log{ax}$ and multiply both sides by $\rho ^{\mu } \lambda ^v (\lambda  x)^{-v} (\rho  x)^{-\mu } \left(cx^p+1\right)^{\beta } x^{\mu +m+v} J_v(x \lambda ) J_{\mu }(x \rho )$ and take the definite integral over $x\in[0,b]$ and simplify to get;
\begin{multline}\label{eq:lhs}
\int_{0}^{b}\frac{x^{m+v+\mu } \left(1+c x^p\right)^{\beta } \lambda ^v (x \lambda )^{-v} \rho ^{\mu } (x \rho )^{-\mu }
   J_v(x \lambda ) J_{\mu }(x \rho ) \log ^k(a x)}{k!}dx\\
=\frac{1}{2\pi i}\int_{C}\int_{0}^{b}w^{-1-k} x^{m+v+\mu } (a x)^w \left(1+c x^p\right)^{\beta }
   \lambda ^v (x \lambda )^{-v} \rho ^{\mu } (x \rho )^{-\mu } J_v(x \lambda ) J_{\mu }(x \rho )dxdw
\end{multline}
where $Re(b)>0$.
We are able to switch the order of integration over $x$ and $w$ using Fubini's theorem for multiple integrals see page 178 in \cite{gelca}, since the integrand is of bounded measure over the space $\mathbb{C} \times [0,b]$.
\subsection{Right-hand side contour integral representation}
Using a generalization of Cauchy's integral formula \ref{intro:cauchy}, we replace $y\to y+xi+\log(a)$ and multiply both sides by $e^{mxi}$. Then form a second equation by replacing $a\to -a$ and take their difference to get;
 \begin{multline}\label{eq:4.3}
\frac{i e^{-i m x} \left((\log (a)-i x+y)^k-e^{2 i m x} (\log (a)+i x+y)^k\right)}{2 k!}\\
=\frac{1}{2\pi i}\int_{C}a^w w^{-k-1} e^{w y}
   \sin (x (m+w))dw
\end{multline}
Next we take the indefinite integral with respect to $y$ to get;
\begin{multline}\label{eq:4.4}
-\frac{i (-1)^{k+1} a^{-m} m^{-k-1} (\Gamma (k+1,-m (-i x+y+\log (a)))-\Gamma (k+1,-m (i x+y+\log (a))))}{2 k!}\\
=\frac{1}{2\pi i}\int_{C}\frac{a^w w^{-k-1} e^{y (m+w)} \sin (x (m+w))}{m+w}dw
\end{multline}
from equation equation (2.325.6(12i)) in \cite{grad}. Finally we replace $y\to -\gamma, x\to \beta$, multiply both sides by $2\pi$ and simplify to get;
\begin{multline}\label{eq:4.5}
\frac{i \pi  (-1)^{k+1} a^{-m} m^{-k-1} (\Gamma (k+1,m (-i \beta +\gamma -\log (a)))-\Gamma (k+1,m (i \beta
   +\gamma -\log (a))))}{k!}\\
=\frac{1}{2\pi i}\int_{C}\frac{2 \pi  a^w w^{-k-1} e^{-(\gamma  (m+w))} \sin (\beta  (m+w))}{m+w}dw
\end{multline}
Use equation (\ref{eq:4.4}) and set $y=0$, replace $m\to h p+2 j+2 l+\mu +m+v+1$ and multiply both sides by 
\begin{equation}
\frac{c^h \left(\lambda ^2\right)^j (-1)^{j+l} \rho ^{\mu } \left(\rho
   ^2\right)^l \lambda ^v \binom{\beta }{h} 2^{-2 j-2 l-\mu -v} b^{h p+2 j+2
   l+\mu +m+v+1}}{j! l! \Gamma (j+v+1) \Gamma (l+\mu +1)}
\end{equation}
then take the infinite series over $j\in[0,\infty),h\in[0,\infty),l\in[0,\infty)$ to get;
\begin{multline}\label{eq:rhs}
\sum _{j=0}^{\infty } \sum _{l=0}^{\infty } \sum _{h=0}^{\infty } \frac{(-1)^{j+l} 2^{-2 j-2 l-v-\mu }
   c^h \lambda ^v \left(\lambda ^2\right)^j \rho ^{\mu
   } \left(\rho ^2\right)^l \binom{\beta }{h}a^{-1-2 j-2 l-m-h p-v-\mu }  }{j! k! l! \Gamma
   (1+j+v) \Gamma (1+l+\mu )(-1-2 j-2 l-m-h p-v-\mu )^{k+1} }\\ \times
\Gamma (1+k,-((1+2 j+2 l+m+h p+v+\mu ) \log (a b)))\\
=-\frac{1}{2\pi i}\sum _{j=0}^{\infty } \sum _{l=0}^{\infty } \sum _{h=0}^{\infty }\int_{C}\frac{c^h \left(\lambda ^2\right)^j (-1)^{j+l} w^{-k-1} \rho ^{\mu }
   \left(\rho ^2\right)^l \lambda ^v (a b)^w }{j! l! \Gamma (j+v+1) \Gamma (l+\mu +1) }\\ \times
\frac{\binom{\beta }{h} 2^{-2 j-2 l-\mu
   -v} b^{h p+2 j+2 l+\mu +m+v+1}}{(h
   p+2 j+2 l+\mu +m+v+w+1)}dw\\
=-\frac{1}{2\pi i}\int_{C}\sum _{j=0}^{\infty } \sum _{l=0}^{\infty } \sum _{h=0}^{\infty }\frac{c^h \left(\lambda ^2\right)^j (-1)^{j+l} w^{-k-1} \rho ^{\mu }
   \left(\rho ^2\right)^l \lambda ^v (a b)^w }{j! l! \Gamma (j+v+1) \Gamma (l+\mu +1) }\\ \times
\frac{\binom{\beta }{h} 2^{-2 j-2 l-\mu
   -v} b^{h p+2 j+2 l+\mu +m+v+1}}{(h
   p+2 j+2 l+\mu +m+v+w+1)}dw
\end{multline}
from equation [Wolfram, \href{http://functions.wolfram.com/07.27.02.0002.01}{07.27.02.0002.01}] where $|Re(b)|<1$. We are able to switch the order of integration and summation over $w$ using Tonellii's theorem for  integrals and sums see page 177 in \cite{gelca}, since the summand is of bounded measure over the space $\mathbb{C} \times [0,\infty)\times [0,\infty)\times [0,\infty)$, where $Re(b)>0$.
%}
\section{Derivations and evaluations}
\begin{theorem}
Extended form of the Landau-Luswili definite integral \cite{landau}.
\begin{multline}\label{eq:theorem}
\int_0^b x^m \left(1+c x^p\right)^{\beta } J_v(x \lambda ) J_{\mu }(x \rho ) \log ^k\left(\frac{1}{a x}\right)
   \, dx
=\frac{\left(\lambda ^v \rho ^{\mu }\right) }{2^{v+\mu } a^{1+m+v+\mu }}\\ \times\sum _{j=0}^{\infty } \sum _{l=0}^{\infty } \sum _{h=0}^{\infty }
   \frac{(-1)^{j+l} c^h \lambda ^{2 j} \rho ^{2 l} \binom{\beta }{h} }{2^{2 j+2 l} a^{2 j+2 l+h p}  j! l! \Gamma (1+j+v)\Gamma (1+l+\mu )}\\ \times
\frac{\Gamma \left(1+k,(1+2 j+2 l+m+h p+v+\mu ) \log\left(\frac{1}{a b}\right)\right)}{(1+2 j+2 l+m+h p+v+\mu )^{k+1}}
\end{multline}
where $Re(\beta)>0,Re(b)>0$.
\end{theorem}
\begin{proof}
Since the right-hand sides of equations (\ref{eq:lhs}) and (\ref{eq:rhs}) are equal relative to equation (\ref{eq:def_int_cont}), we may equate the left-hand sides and simplify the gamma function to yield the stated result.
\end{proof}
\begin{example}
Definite integral involving the product of Bessel functions with trigonometric functions. This is a special case of equation (10.46.1) in \cite{hansen}. In this example we use equation (\ref{eq:theorem}) and set $k\to 0,a\to 1,b\to 1,\rho \to \lambda ,c\to -1,p\to 2,\beta \to
   \frac{\beta -1}{2},x\to \sin (\theta )$ amd simplify the series on the right-hand side.
\begin{multline}\label{eq:lan1}
\int_0^{\frac{\pi }{2}} J_v(\lambda  \sin (\theta )) J_{\mu }(\lambda  \sin (\theta )) \cos ^{\beta }(\theta )
   \sin ^m(\theta ) \, d\theta \\
=\sum _{j=0}^{\infty } \frac{(-1)^j 2^{-1-2 j-v-\mu } \lambda ^{2 j+v+\mu } \Gamma
   \left(\frac{1+\beta }{2}\right) \Gamma \left(\frac{1}{2} (1+2 j+m+v+\mu )\right) }{\Gamma (1+j) \Gamma (1+j+v) \Gamma (1+\mu ) \Gamma
   \left(\frac{1}{2} (2+2 j+m+v+\beta +\mu )\right)}\\ \times
\,_1F_2\left(\frac{1}{2}+j+\frac{m}{2}+\frac{v}{2}+\frac{\mu }{2};1+j+\frac{m}{2}+\frac{v}{2}+\frac{\beta
   }{2}+\frac{\mu }{2},1+\mu ;-\frac{\lambda ^2}{4}\right)
\end{multline}
where $Re(\beta)>0,Re(m)>0$.
\end{example}
\begin{example}
Derivation of equation (4) in \cite{landau}. Extended Stoyanov-Farrel integral \cite{stoyanov}.  In this example we use equation (\ref{eq:lan1}) and set $m\to 0,\beta \to 0,v\to 0,\mu \to 0$ and simplify;
\begin{multline}\label{eq:lan2}
\int_0^{\frac{\pi }{2}} J_0(\lambda  \sin (\theta )){}^2 \, d\theta =\sqrt{\pi } \sum _{j=0}^{\infty }
   \frac{\left(-\frac{\lambda ^2}{4}\right)^j \Gamma \left(\frac{1}{2} (3+2 j)\right) \,
   _1F_2\left(\frac{1}{2}+j;1,1+j;-\frac{\lambda ^2}{4}\right)}{(1+2 j) \Gamma (1+j)^3}
\end{multline}
where $\lambda\in\mathbb{C}$.
\end{example}
\begin{example}
Here we use equation (\ref{eq:lan2}) and set $\lambda \to 2 i$ and simplify.
\begin{equation}
\int_0^{\frac{\pi }{2}} I_0(2 \sin (\theta )){}^2 \, d\theta =\frac{ \sqrt{\pi }}{2} \sum _{j=0}^{\infty }
   \frac{\Gamma \left(\frac{1}{2}+j\right) }{\Gamma (1+j)^3}\, _1F_2\left(\frac{1}{2}+j;1,1+j;1\right)
\end{equation}
\end{example}
\begin{example}
Alternate form for equation (10.45.1) in \cite{hansen}. Here we use equation (\ref{eq:lan1}) and set $v\to c-1,\lambda \to 2 x,\theta \to t,\mu \to d-1,m\to 2 a-c-d+1,\beta \to -2
   a+2 b-1$ and equate the right-hand sides.
\begin{multline}\label{eq:lan3}
\sum _{l=0}^{\infty } \frac{(-1)^l x^{2 l} \Gamma (a+l) \, _1F_2\left(a+l;c,b+l;-x^2\right)}{\Gamma (1+l)
   \Gamma (b+l) \Gamma (d+l)}\\
   =\frac{\Gamma (a) }{\Gamma (b) \Gamma (d)}\sum _{k=0}^{\infty } \frac{\left((a)_k (c+d-1)_{2 k}\right) (-1)^k
   x^{2 k}}{k! (b)_k (c)_k (d)_k (c+d-1)_k}
\end{multline}
where $Re(a)>0,Re(b-a)>0$.
\end{example}
\begin{example}
Extended form of Equation (7) in \cite{landau}. Here we repeat the derivation used in equation (\ref{eq:lan1}) with $x\to \cos (\theta )$ and simplify.
\begin{multline}\label{eq:lan4}
\int_0^{\frac{\pi }{2}} J_v(\lambda  \cos (\theta )) J_{\mu }(\lambda  \cos (\theta )) \cos ^m(\theta ) \sin
   ^{\beta }(\theta ) \, d\theta \\
=\sum _{j=0}^{\infty } \frac{(-1)^j 2^{-1-2 j-v-\mu } \lambda ^{2 j+v+\mu } \Gamma
   \left(\frac{1+\beta }{2}\right) \Gamma \left(\frac{1}{2} (1+2 j+m+v+\mu )\right) }{\Gamma (1+j) \Gamma (1+j+v) \Gamma (1+\mu ) \Gamma
   \left(\frac{1}{2} (2+2 j+m+v+\beta +\mu )\right)}\\ \times
\,_1F_2\left(\frac{1}{2}+j+\frac{m}{2}+\frac{v}{2}+\frac{\mu }{2};1+j+\frac{m}{2}+\frac{v}{2}+\frac{\beta
   }{2}+\frac{\mu }{2},1+\mu ;-\frac{\lambda ^2}{4}\right)
\end{multline}
where $Re(v)>-1/2,Re(\mu)<-1/2$.
\end{example}
\begin{example}
Derivation of a series form of equation (6) in \cite{landau}. Extended Wong integral \cite{wong}. Here we use equation (\ref{eq:lan4}) and set $\mu \to v,m\to 0,\beta \to 0$ and simplify;
\begin{multline}
\int_0^{\frac{\pi }{2}} J_v(\lambda  \cos (\theta )){}^2 \, d\theta \\
=\sum _{j=0}^{\infty } \frac{(-1)^j
   2^{-1-2 j-2 v} \sqrt{\pi } \lambda ^{2 j+2 v} \Gamma \left(\frac{1}{2} (1+2 j+2 v)\right) 
  }{\Gamma (1+j) \Gamma (1+v) \Gamma(1+j+v)^2}\\   \times \,_1F_2\left(\frac{1}{2}+j+v;1+v,1+j+v;-\frac{\lambda ^2}{4}\right)
  \\
\approx \frac{\log (\lambda )+2 \log (2)-\psi \left(\frac{1}{2}+v\right)}{\pi  \lambda }+\frac{\sin \left(2
   \lambda -v \pi -\frac{\pi }{4}\right)}{2 \sqrt{\pi } \lambda ^{3/2}}+\frac{\left(16 v^2-9\right) \cos \left(2
   \lambda -v \pi -\frac{\pi }{4}\right)}{32 \sqrt{\pi } \lambda ^{5/2}}\\
+\frac{\left(4 v^2-1\right) \left(\log
   (\lambda )+2 \log (2)-\psi \left(\frac{1}{2}+v\right)\right)-\left(4 v^2-3\right)}{16 \pi  \lambda ^3}
\end{multline}
where $\lambda\in\mathbb{C}$.
\end{example}
\begin{example}
Derivation of Equation (8) in \cite{landau}.  In this example we use equation (\ref{eq:lan1}) and set $m\to 0,\beta \to 0$ and simplify;
\begin{multline}
\int_0^{\frac{\pi }{2}} J_v(\lambda  \sin (\theta )) J_{\mu }(\lambda  \sin (\theta )) \, d\theta
   \\
=\frac{\left(\sqrt{\pi } \lambda ^{+v+\mu }\right) }{2^{1+v+\mu } \Gamma (1+\mu )}\sum _{j=0}^{\infty } \frac{(-1)^j 2^{-2 j} \lambda ^{2 j}
   \Gamma \left(\frac{1}{2} (1+2 j+v+\mu )\right) }{\Gamma (1+j) \Gamma (1+j+v) \Gamma
   \left(\frac{1}{2} (2+2 j+v+\mu )\right)}\\ \times
\, _1F_2\left(\frac{1}{2}+j+\frac{v}{2}+\frac{\mu
   }{2};1+j+\frac{v}{2}+\frac{\mu }{2},1+\mu ;-\frac{\lambda ^2}{4}\right)
\end{multline}
where $Re(v)>-1/2$.
\end{example}
\begin{example}
The Hurwitz-Lerch zeta function. In this example we use equation (\ref{eq:theorem}) and set $a\to 1,b\to 1,\beta \to -1,c\to z,m\to m-1$ and simplify using [DLMF,\href{https://dlmf.nist.gov/25.14.i}{25.14.1}].
\begin{multline}\label{eq:lerch}
\int_0^1 \frac{x^{-1+m} J_v(x \lambda ) J_{\mu }(x \rho ) \log ^k\left(\frac{1}{x}\right)}{1+x^p z} \,
   dx\\
=\frac{\lambda ^v \rho ^{\mu } }{2^{v+\mu } p^{1+k}}\sum _{j=0}^{\infty }\sum _{l=0}^{\infty } \frac{(-1)^{j+l} 2^{-2
   j-2 l} \lambda ^{2 j} \rho ^{2 l} \Gamma (1+k) \Phi \left(-z,1+k,\frac{2 j+2 l+m+v+\mu }{p}\right)}{j! l! \Gamma
   (1+j+v) \Gamma (1+l+\mu )}
\end{multline}where $Re(p)>1$.
\end{example}
\begin{example}
Definite integral in terms of the log-gamma function. In this example we use equation (\ref{eq:lerch}) and apply l'Hopital's rule as $k\to -1$ and set $z=i$ and repeat for $z=-i$ and add the two equations and simplify.
\begin{multline}
\int_0^1 \frac{\left(x^s-x^m\right) J_v(x \lambda ) J_{\mu }(x \rho )}{\left(1+x^p\right) \log
   \left(\frac{1}{x}\right)} \, dx
=\frac{\lambda ^v \rho ^{\mu } }{2^{v+\mu }}\sum _{j=0}^{\infty } \sum
   _{l=0}^{\infty } \frac{(-1)^{j+l} \lambda ^{2 j} \rho ^{2 l} }{2^{2
   j+2 l} j! l! \Gamma (1+j+v) \Gamma (1+l+\mu )}\\
\times \log \left(\frac{\Gamma
   \left(\frac{1+2 j+2 l+m+p+v+\mu }{2 p}\right) \Gamma \left(\frac{1+2 j+2 l+s+v+\mu }{2 p}\right)}{\Gamma
   \left(\frac{1+2 j+2 l+m+v+\mu }{2 p}\right) \Gamma \left(\frac{1+2 j+2 l+p+s+v+\mu }{2 p}\right)}\right)
\end{multline}
where $|Re(m)|<1,|Re(s)|<1$.
\end{example}
\begin{example}
Hurwitz zeta functions. In this example we use equation (\ref{eq:lerch}) and set $z=1$ and simplify;
\begin{multline}\label{eq:hurwitz_funs}
\int_0^1 \frac{x^{-1+m} J_v(x \lambda ) J_{\mu }(x \rho ) \log ^k\left(\frac{1}{x}\right)}{1+x^p} \, dx\\
=\sum
   _{j=0}^{\infty } \sum _{l=0}^{\infty } \frac{(-1)^{j+l} 2^{-2 j-2 l-v-\mu } p^{-1-k} \lambda ^{2 j+v} \rho ^{2
   l+\mu } \Gamma (1+k)}{j! l! \Gamma (1+j+v) \Gamma
   (1+l+\mu )}\\ \times
 \left(2^{-1-k} \zeta \left(1+k,\frac{2 j+2 l+m+v+\mu }{2 p}\right)\right. \\ \left.
-2^{-1-k} \zeta
   \left(1+k,\frac{1}{2} \left(1+\frac{2 j+2 l+m+v+\mu }{p}\right)\right)\right)
\end{multline}
where $|Re(m)|<1,Re(p)>0$.
\end{example}
\begin{example}
Malmsten logarithmic integral form \cite{malmsten}. In this example we use equation (\ref{eq:hurwitz_funs}) and take the first derivative with respect to $k$ and set $k=0$ and simplify;
\begin{multline}
\int_0^1 \frac{x^{-1+m} J_v(x \lambda ) J_{\mu }(x \rho ) \log \left(\log
   \left(\frac{1}{x}\right)\right)}{1+x^p} \, dx\\
=\sum _{j=0}^{\infty } \sum _{l=0}^{\infty } \frac{(-1)^{j+l} 2^{-1-2
   j-2 l-v-\mu } \lambda ^{2 j+v} \rho ^{2 l+\mu } }{p \Gamma (1+j) \Gamma (1+l)
   \Gamma (1+j+v) \Gamma (1+l+\mu )}\\ \times
\left((\gamma +\log (2 p)) \left(\psi ^{(0)}\left(\frac{2 j+2
   l+m+v+\mu }{2 p}\right)-\psi ^{(0)}\left(\frac{2 j+2 l+m+p+v+\mu }{2 p}\right)\right)\right. \\ \left.
-\gamma _1\left(\frac{2 j+2
   l+m+v+\mu }{2 p}\right)+\gamma _1\left(\frac{2 j+2 l+m+p+v+\mu }{2 p}\right)\right)
\end{multline}
where $|Re(m)|<1,Re(p)>0$.
\end{example}
\begin{example}
The Hurwitz zeta function. In this example we use equation (\ref{eq:lerch}) and set $z=-1$ and simplify;
\begin{multline}
\int_0^1 \frac{x^{-1+m} J_v(x \lambda ) J_{\mu }(x \rho ) \log ^k\left(\frac{1}{x}\right)}{1-x^p} \, dx\\
=\sum
   _{j=0}^{\infty } \sum _{l=0}^{\infty } \frac{(-1)^{j+l} 2^{-2 j-2 l-v-\mu } p^{-1-k} \lambda ^{2 j+v} \rho ^{2
   l+\mu } \Gamma (1+k) \zeta \left(1+k,\frac{2 j+2 l+m+v+\mu }{p}\right)}{j! l! \Gamma (1+j+v) \Gamma (1+l+\mu
   )}
\end{multline}
where $|Re(m)|<1,Re(p)>0$.
\end{example}
\begin{example}
Derivation of equation on Page 381 question (19) in \cite{watson}. Here we use equation (\ref{eq:theorem}) and set $k\to 0,a\to 1,b\to 1,\beta \to 0,p\to 0,c\to 1,\mu \to v,\rho \to \lambda$ and simplify;
\begin{multline}
\frac{4^{-v} \lambda ^{2 v} \,
   _2F_3\left(\frac{1}{2}+v,\frac{1}{2}+\frac{m}{2}+v;1+v,\frac{3}{2}+\frac{m}{2}+v,1+2 v;-\lambda ^2\right)}{(1+m+2
   v) \Gamma (1+v)^2}\\
=\sum _{l=0}^{\infty } \frac{(-1)^l 4^{-l-v} \lambda ^{2 l+2 v} \,
   _1F_2\left(\frac{1}{2}+l+\frac{m}{2}+v;1+v,\frac{3}{2}+l+\frac{m}{2}+v;-\frac{\lambda ^2}{4}\right)}{(1+2 l+m+2 v)
   \Gamma (1+l) \Gamma (1+v) \Gamma (1+l+v)}
\end{multline}
where $Re(m)>0,Re(v)>0,Re(\lambda)>0$.
\end{example}
\begin{example}
Bessel and logarithm function. Here we use equation (\ref{eq:theorem}) and set $\beta \to 0,p\to 0,c\to 1,\mu \to v,\rho \to \lambda ,m\to 0,a\to 1,b\to 1$ and simplify;
\begin{multline}
\int_0^1 J_v(x \lambda ){}^2 \log ^k\left(\frac{1}{x}\right) \, dx\\
=\left(\frac{\lambda }{2}\right)^{2 v} \Gamma (1+k) \sum _{j=0}^{\infty } \sum _{l=0}^{\infty }
   \frac{\left(\frac{\lambda }{2}\right)^{2 j+2 l} (-1)^{j+l}}{(1+2 j+2 l+2 v)^{k+1} (j! l! \Gamma (1+j+v) \Gamma (1+l+v))}
\end{multline}
where $Re(v)>-1/2$.
\end{example}
\section{Integrals involving Bessel functions with trigonometric arguments \cite{luke} pp. 298-306}
In section 13.3.2 in \cite{luke} pages 298-306, the definite integrals are derived by expanding the Bessel function(s) in the integrand using equations 1.4.1(1) or 1.4.1(17-19) appropriately. Term wise integration is used with the aid of equations 1.2(14-16). The series expansions are then expressed in terms of the hypergeometric function or as a series of hypergeometric functions. Special case results are obtained by specializing the parameters and or variables along with using equations 1.3.5, 1.3.6 and 1.4.1. The aim of this section of this article is to derive alternate forms for some the definite integrals in this section of the book by Luke in terms of the series of the hypergeometric function which is not present in the book by Luke. These derivations are achieved by substituting values for the parameters in equation (\ref{eq:theorem}).
\begin{example}
Derivation of equations 13.3.2(19,20,21) in \cite{luke}. Here we use equation (\ref{eq:theorem}) and set $k\to 0,a\to 1,b\to 1,\lambda \to z,\rho \to z,\beta \to -\frac{1}{2},c\to -1,p\to 2,m\to 2 (a+b-1)$ and simplify;
\begin{multline}
\int_0^{\frac{\pi }{2}} J_v(z \sin (t)) J_{\mu }(z \sin (t)) \sin ^{-2+2 a+2 b}(t) \, dt=\int_0^1 \frac{x^{2 (-1+a+b)} J_v(x z) J_{\mu }(x z)}{\sqrt{1-x^2}} \, dx\\
=\sum
   _{j=0}^{\infty } \frac{(-1)^j 2^{-1-2 j-v-\mu } \sqrt{\pi } z^{2 j+v+\mu } \Gamma \left(\frac{1}{2} (-1+2 a+2 b+2 j+v+\mu )\right) }{\Gamma (1+j) \Gamma (1+j+v) \Gamma (1+\mu ) \Gamma
   \left(\frac{1}{2} (2 a+2 b+2 j+v+\mu )\right)}\\ \times \,
   _1F_2\left(-\frac{1}{2}+a+b+j+\frac{v}{2}+\frac{\mu }{2};a+b+j+\frac{v}{2}+\frac{\mu }{2},1+\mu ;-\frac{z^2}{4}\right)
\end{multline}
where $Re(\mu+v+2a)>0,Re(b)>0$.
\end{example}
\begin{example}
Derivation of equation 13.3.2(22) in \cite{luke}. Here we use equation (\ref{eq:theorem}) and set $k\to 0,a\to 1,b\to 1,\mu \to v-1,c\to -1,\beta \to -v-\frac{1}{2},p\to 2,m\to 2 v,\rho \to z,\lambda \to z$ and simplify;
\begin{multline}
\int_0^{\frac{\pi }{2}} J_{v-1}(z \sin (t)) J_v(z \sin (t)) \tan ^{2 v}(t) \, dt=\int_0^1 x^{2 v} \left(1-x^2\right)^{-\frac{1}{2}-v} J_{-1+v}(x z) J_v(x z) \, dx\\
=\sum
   _{j=0}^{\infty } \frac{(-1)^j 2^{-2 j-2 v} z^{-1+2 j+2 v} \Gamma \left(\frac{1}{2} (1-2 v)\right) \Gamma (j+2 v) \, _1F_2\left(j+2 v;v,\frac{1}{2}+j+v;-\frac{z^2}{4}\right)}{\Gamma
   (1+j) \Gamma (v) \Gamma (1+j+v) \Gamma \left(\frac{1}{2} (1+2 j+2 v)\right)}
\end{multline}
where $0< Re(v)<1/2$.
\end{example}
\begin{example}
Derivation of equation 13.3.2(23) in \cite{luke}. Here we use equation (\ref{eq:theorem}) and set $k\to 0,a\to 1,b\to 1,c\to -1,\beta \to -u-\frac{1}{2},p\to 2,m\to 2 u+1,\rho \to z,\lambda \to z$ and simplify;
\begin{multline}
\int_0^{\frac{\pi }{2}} J_v(z \sin (t)) J_{\mu }(z \sin (t)) \sin (t) \tan ^{2 \mu }(t) \, dt\\
=\int_0^1 x^{1+2 \mu } \left(1-x^2\right)^{-\frac{1}{2}-\mu } J_v(x z) J_{\mu }(x z)
   \, dx\\
=\sum _{j=0}^{\infty } \frac{(-1)^j 2^{-1-2 j-v-\mu } z^{2 j+v+\mu } \Gamma \left(\frac{1}{2} (1-2 \mu )\right) \Gamma \left(\frac{1}{2} (2+2 j+v+3 \mu )\right) }{\Gamma (1+j) \Gamma (1+j+v) \Gamma (1+\mu ) \Gamma \left(\frac{1}{2}
   (3+2 j+v+\mu )\right)}\\ \times
\,
   _1F_2\left(1+j+\frac{v}{2}+\frac{3 \mu }{2};\frac{3}{2}+j+\frac{v}{2}+\frac{\mu }{2},1+\mu ;-\frac{z^2}{4}\right)
\end{multline}
where $Re(3\mu+v)>-2$.
\end{example}
\begin{example}
Derivation of equation 13.3.2(24) in \cite{luke}. Here we use equation (\ref{eq:theorem}) and set $k\to 0,a\to 1,b\to 1,v\to 1,\mu \to 1,\lambda \to z,\rho \to z,m\to -1,p\to 2,c\to -1,\beta \to -\frac{1}{2}$ and simplify;
\begin{multline}
\int_0^{\frac{\pi }{2}} J_1(z \sin (t)){}^2 \csc (t) \, dt=\int_0^1 \frac{J_1(x z){}^2}{x \sqrt{1-x^2}} \, dx\\
=\sum _{j=0}^{\infty } \frac{(-1)^j 2^{-3-2 j} \sqrt{\pi } z^{2+2 j}
   \, _1F_2\left(1+j;2,\frac{3}{2}+j;-\frac{z^2}{4}\right)}{\Gamma (2+j) \Gamma \left(\frac{1}{2} (3+2 j)\right)}
\end{multline}
where $z\in\mathbb{C}$.
\end{example}
\begin{example}
Derivation of equation 13.3.2(33) in \cite{luke}. Here we use equation (\ref{eq:theorem}) and set $k\to 0,a\to 1,b\to 1,\lambda \to z,\rho \to w,x\to \sin ^2(t),c\to -1,p\to 1$ and simplify;
\begin{multline}\label{eq:luke33}
\int_0^{\frac{\pi }{2}} J_v\left(w \sin ^2(t)\right) J_{\mu }\left(z \sin ^2(t)\right) \cos ^{-1+2 b}(t) \sin ^{-1+2 a}(t) \, dt\\
=\sum _{l=0}^{\infty } \frac{(-1)^l 2^{-1-2 l-v-\mu
   } w^{2 l+v} z^{\mu } \Gamma (b) \Gamma (a+2 l+v+\mu )}{\Gamma (1+l) \Gamma (1+l+v)
   \Gamma (1+\mu ) \Gamma (a+b+2 l+v+\mu )}\\ \times 
 \, _2F_3\left(\frac{a}{2}+l+\frac{v}{2}+\frac{\mu }{2},\frac{1}{2}+\frac{a}{2}+l+\frac{v}{2}+\frac{\mu
   }{2};\right. \\ \left.
   \frac{a}{2}+\frac{b}{2}+l+\frac{v}{2}+\frac{\mu }{2},\frac{1}{2}+\frac{a}{2}+\frac{b}{2}+l+\frac{v}{2}+\frac{\mu }{2},1+\mu ;-\frac{z^2}{4}\right)
\end{multline}
where $Re(b)>Re(a)$.
\end{example}
\begin{example}
Derivation of equation 13.3.2(34) in \cite{luke}. Errata. Here we use equation (\ref{eq:luke33}) and set $w\to z$;
\begin{multline}
\int_0^{\frac{\pi }{2}} J_v\left(z \sin ^2(t)\right) J_{\mu }\left(z \sin ^2(t)\right) \cos ^{-1+2 b}(t) \sin
   ^{-1+2 a}(t) \, dt\\
=\sum _{j=0}^{\infty } \frac{(-1)^j 2^{-1-2 j-v-\mu } z^{2 j+v+\mu } \Gamma (b) \Gamma (a+2
   j+v+\mu ) }{\Gamma (1+j) \Gamma
   (1+v) \Gamma (1+j+\mu ) \Gamma (a+b+2 j+v+\mu )}\\ \times\, _2F_3\left(\frac{a}{2}+j+\frac{v}{2}+\frac{\mu }{2},\frac{1}{2}+\frac{a}{2}+j+\frac{v}{2}+\frac{\mu
   }{2};1+v,\frac{a}{2}\right. \\ \left.
+\frac{b}{2}+j+\frac{v}{2}+\frac{\mu
   }{2},\frac{1}{2}+\frac{a}{2}+\frac{b}{2}+j+\frac{v}{2}+\frac{\mu }{2};-\frac{z^2}{4}\right)\\
\neq\frac{\left(2^{-1-v-\mu } z^{v+\mu } \Gamma (b) \Gamma (a+v+\mu
   )\right) }{\Gamma (1+v)
   \Gamma (1+\mu ) \Gamma (a+b+v+\mu )}\\ \times
\, _4F_5\left(\frac{1+v}{2},\frac{2+v}{2},\frac{1}{2} (a+v+\mu ),\frac{1}{2} (1+a+v+\mu
   );1+v,\frac{1}{2} (a+b+v+\mu ),\right. \\ \left.
\frac{1}{2} (1+a+b+v+\mu ),1+\mu ,1+v+\mu ;-\frac{z^2}{4}\right)
\end{multline}
where $Re(v)>0,Re(z)>0,Re(\mu)>0$.
\end{example}
\begin{example}
Derivation of equation 13.3.2(51) in \cite{luke}. Here we use equation (\ref{eq:theorem}) and set $k\to 0,a\to 1,b\to 1,\beta \to -\frac{1}{2},p\to 2,c\to -1,\mu \to \frac{1}{2},\rho \to b z,\lambda \to z,m\to \frac{1}{2}$ and simplify;
\begin{multline}
\int_0^{\frac{\pi }{2}} J_v(z \sin (t)) \sin (b z \sin (t)) \, dt\\
=\sum _{l=0}^{\infty } \frac{(-1)^l 2^{-2-2 l-v} b \pi  z^{1+v} (b z)^{2 l} \Gamma \left(\frac{1}{2} (2+2 l+v)\right)
   \, _1F_2\left(1+l+\frac{v}{2};\frac{3}{2}+l+\frac{v}{2},1+v;-\frac{z^2}{4}\right)}{\Gamma (1+l) \Gamma \left(\frac{1}{2} (3+2 l)\right) \Gamma (1+v) \Gamma \left(\frac{1}{2} (3+2
   l+v)\right)}\\
=\pi  \sum _{k=0}^{\infty } (-1)^k J_{2 k+1}(b z) J_{\frac{v}{2}+k+\frac{1}{2}}\left(\frac{z}{2}\right) J_{\frac{v}{2}-k-\frac{1}{2}}\left(\frac{z}{2}\right)\\
=\pi  J_1(b z)
   J_{\frac{1}{2} (-1+v)}\left(\frac{z}{2}\right) J_{\frac{1+v}{2}}\left(\frac{z}{2}\right)\\
-\pi  J_3(b z) J_{\frac{1}{2} (-3+v)}\left(\frac{z}{2}\right)
   J_{\frac{3+v}{2}}\left(\frac{z}{2}\right)
\end{multline}
where $Re(v)>-1$.
\end{example}
\begin{example}
Derivation of equation 13.3.2(52) in \cite{luke}. Here we use equation (\ref{eq:theorem}) and set $k\to 0,a\to 1,b\to 1,\beta \to -\frac{1}{2},p\to 2,c\to -1,\mu \to -\frac{1}{2},\rho \to b z,\lambda \to z,m\to \frac{1}{2}$ and simplify;
\begin{multline}
\int_0^{\frac{\pi }{2}} J_v(z \sin (t)) \cos (b z \sin (t)) \, dt\\
=\sum _{l=0}^{\infty } \frac{(-1)^l 2^{-1-2 l-v} \pi  z^v (b z)^{2 l} \Gamma \left(\frac{1}{2} (1+2 l+v)\right) \,
   _1F_2\left(\frac{1}{2}+l+\frac{v}{2};1+l+\frac{v}{2},1+v;-\frac{z^2}{4}\right)}{\Gamma \left(\frac{1}{2}+l\right) \Gamma (1+l) \Gamma \left(1+l+\frac{v}{2}\right) \Gamma
   (1+v)}\\
=\frac{1}{2} \pi  J_0(b z) J_{\frac{v}{2}}\left(\frac{z}{2}\right){}^2+\pi  \sum _{k=1}^{\infty } (-1)^k J_{2 k}(b z) J_{\frac{v}{2}+k}\left(\frac{z}{2}\right)
   J_{\frac{v}{2}-k}\left(\frac{z}{2}\right)
\end{multline}
where $Re(v)>-1$.
\end{example}
\begin{example}
Definite integral involving the reciprocal logarithm function. Here we use equation (\ref{eq:theorem}) and set $a\to \frac{1}{b},\mu \to v,\rho \to \lambda ,\beta \to 0,p\to 0,c\to 1,v\to 0,b\to 1, k\to s$ then form a second equation by replacing $m\to r$ and take their difference.
\begin{multline}
\int_0^1 \frac{\left(x^m-x^r\right) J_0(x \lambda ){}^2}{\log (x)} \, dx\\
=\sum _{j=0}^{\infty } \sum
   _{l=0}^{\infty }\left(-\frac{\lambda ^2}{4}\right)^{j+l} \frac{1}{\Gamma (1+j)^2 \Gamma (1+l)^2} \log \left(\frac{1+2 j+2 l+m}{1+2 j+2
   l+r}\right)
\end{multline}
where $|Re(m)|<1,|Re(r)|<1$
\end{example}
\begin{example}
Logarithmic singularity and trigonometric function integrand. Here we use equation (\ref{eq:theorem}) and set $\beta \to 0,p\to 0,c\to 1,v\to \frac{1}{2},\mu \to -\frac{1}{2},b\to 1,\rho \to \lambda ,m\to 1$, then take the first partial derivative with respect to $k$ and simplify;
\begin{multline}
\int_0^1 \log ^k\left(\frac{1}{a x}\right) \log \left(\log \left(\frac{1}{a x}\right)\right) \sin (x \lambda
   ) \, dx\\
=\sum _{j=0}^{\infty } \sum _{l=0}^{\infty } \frac{(-1)^{j+l} 2^{-1-2 j-k-2 l} a^{-2 (1+j+l)}
   (1+j+l)^{-1-k} \lambda ^{1+2 j+2 l} }{\Gamma (2+2 j) \Gamma (1+2 l)}\\ \times \left(\Gamma \left(1+k,2 (1+j+l) \log \left(\frac{1}{a}\right)\right)
   \left(-\log (1+j+l)+\log \left((1+j+l) \log \left(\frac{1}{a}\right)\right)\right)\right. \\ \left.
+G_{2,3}^{3,0}\left(2 (1+j+l)
   \log \left(\frac{1}{a}\right)\right |\left.
\begin{array}{c}
 1,1 \\
 0,0,1+k \\
\end{array}
\right)\right)
\end{multline}
where $Re(k)>0,Re(\lambda)>0$.
\end{example}
\begin{figure}[H]
\caption{Plot of $Re\left(\sqrt{\log \left(\frac{1}{2 x}\right)} \log \left(\log \left(\frac{1}{2 x}\right)\right) \sin (3 x)\right)$}
\includegraphics[width=8cm]{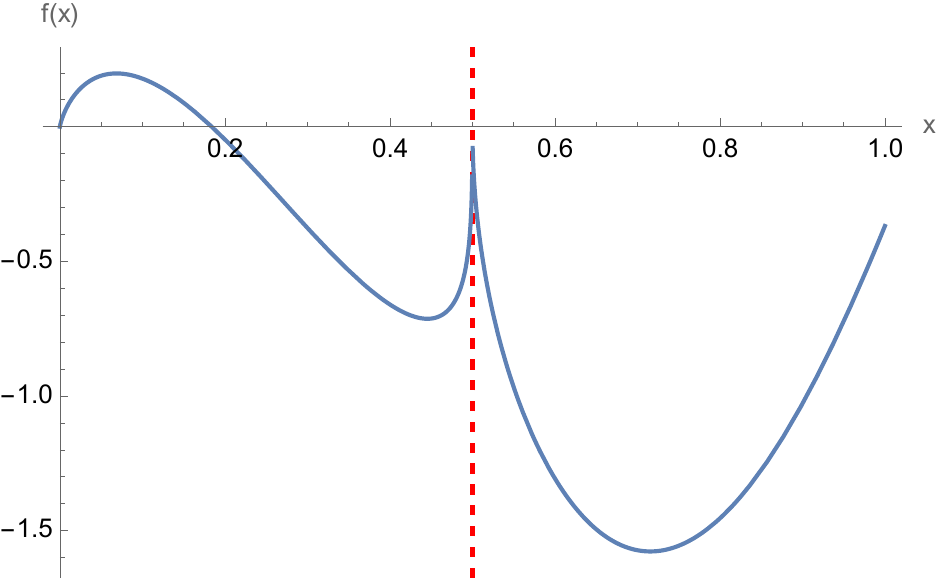}
\end{figure}
\begin{figure}[H]
\caption{Plot of $Re \left(\sin (2 x) \sqrt{\log \left(\frac{1}{a x}\right)} \log \left(\log \left(\frac{1}{a x}\right)\right) \right)$}
\includegraphics[width=8cm]{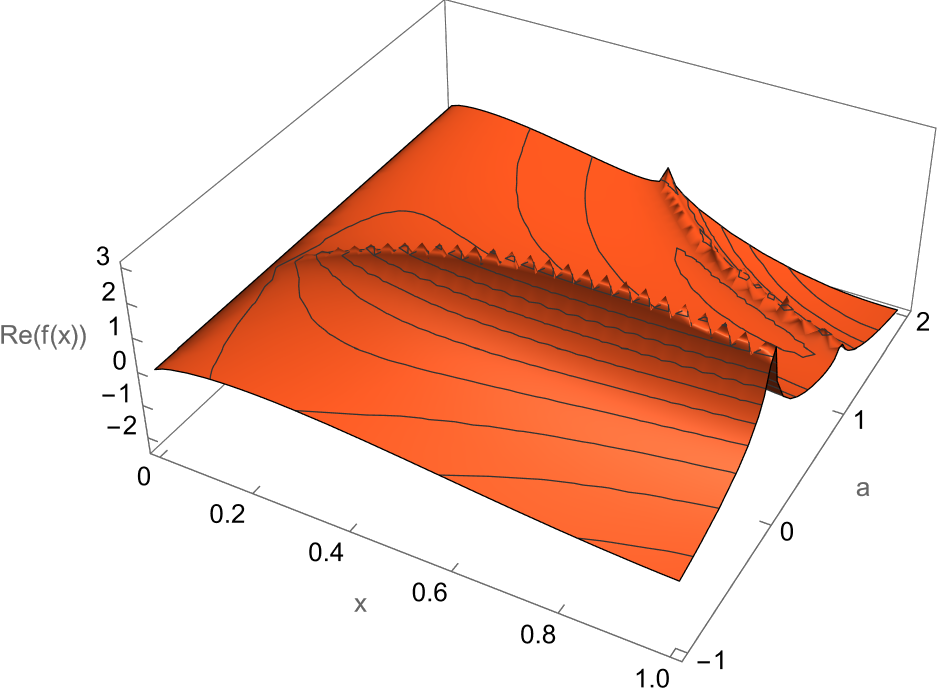}
\end{figure}
\begin{figure}[H]
\caption{Plot of $Im \left(\sin (2 x) \sqrt{\log \left(\frac{1}{a x}\right)} \log \left(\log \left(\frac{1}{a x}\right)\right) \right)$}
\includegraphics[width=8cm]{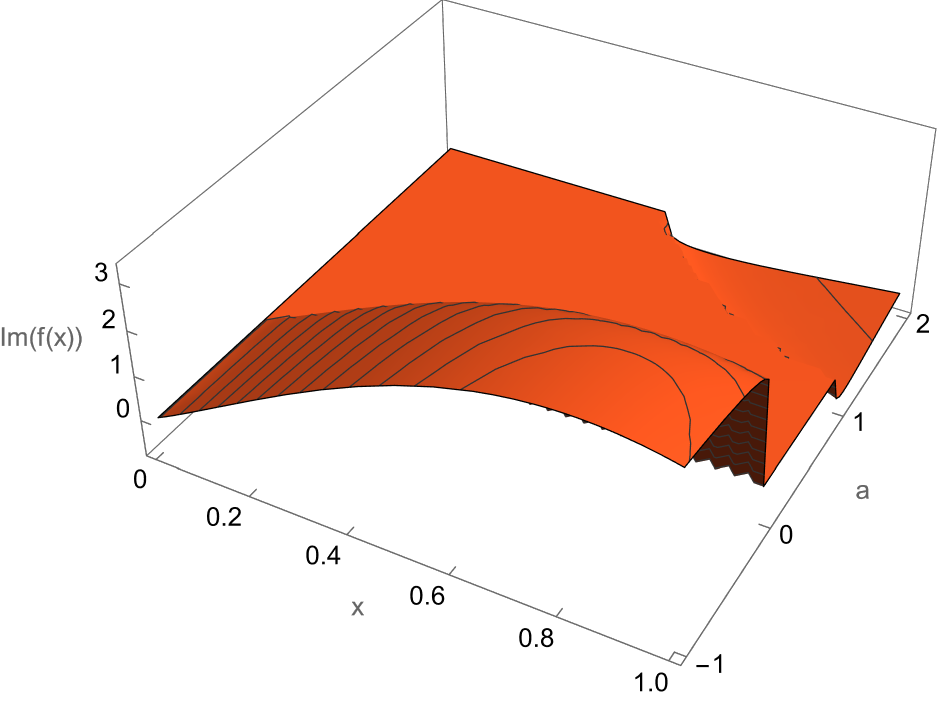}
\end{figure}
\begin{example}
Orthogonality of the Bessel function where $j_{v,l}$ represents the $k^th$ zero of the Bessel function $J_n(x)$. Here we use equation (\ref{eq:theorem}) and set $k\to 0,a\to 1,\beta \to 0,p\to 0,c\to 1,\mu \to v,b\to 1,m\to 1$ and simplify using [DLMF,\href{https://dlmf.nist.gov/10.22.E37}{10.22.37}].
\begin{multline}
\int_0^1 x J_v\left(x j_{v,l}\right) J_v\left(x j_{v,m}\right) \, dx\\
=\sum _{j=0}^{\infty } \frac{(-1)^j
   2^{-1-2 j-2 v} \left(j_{v,l}\right){}^{2 j+v} \left(j_{v,m}\right){}^v \, _1F_2\left(1+j+v;1+v,2+j+v;-\frac{1}{4}
   \left(j_{v,m}\right){}^2\right)}{\Gamma (1+j) \Gamma (1+v) \Gamma (2+j+v)}
\end{multline}
where $Re(v)>0,Re(l)>0,Re(m)>0$.
\end{example}
\begin{example}
Extended Wolfram series. Here we use equation (\ref{eq:theorem}) and set  $k\to 0,a\to 1,\beta \to 0,p\to 0,c\to 1,m\to 0,\rho \to \lambda $ and simplify using equation [Wolfram, \href{http://functions.wolfram.com/03.01.21.0048.01}{01}] ;
\begin{multline}
\frac{\, _3F_4\left(\frac{1}{2}+\frac{v}{2}+\frac{\mu }{2},\frac{1}{2}+\frac{v}{2}+\frac{\mu
   }{2},1+\frac{v}{2}+\frac{\mu }{2};1+v,\frac{3}{2}+\frac{v}{2}+\frac{\mu }{2},1+\mu ,1+v+\mu ;-b^2 \lambda
   ^2\right)}{(1+v+\mu ) \Gamma (1+v) \Gamma (1+\mu )}\\
=\sum _{j=0}^{\infty } \frac{(-1)^j 2^{-2 j} b^{2 j} \lambda
   ^{2 j} \, _1F_2\left(\frac{1}{2}+j+\frac{v}{2}+\frac{\mu }{2};\frac{3}{2}+j+\frac{v}{2}+\frac{\mu }{2},1+\mu
   ;-\frac{1}{4} b^2 \lambda ^2\right)}{(1+2 j+v+\mu ) \Gamma (1+j) \Gamma (1+j+v) \Gamma (1+\mu )}
\end{multline}
where $Re(v)>0,Re(\lambda)>0$.
\end{example}
\section{Derivation and extension of some finite integrals in Prudnikov et al.}
In this section we derive series forms and errata for some finite integrals of the product of two Bessel functions listed in Prudnikov et al. The general formula (\ref{eq:theorem}) used to derive these integrals offer a wider range of evaluation of the parameters involved. This allows for more finite integral formulae derivations. A consequence of the formulae in this section is the derivation of the generalized hypergeometric function in terms of the series of the generalized hypergeometric function. 
\begin{example}
Derivation of equation (2.12.32.2) in \cite{prud2}. Here we use equation (\ref{eq:theorem}) and set $k\to 0,a\to 1,c\to -\frac{1}{a},\beta \to \beta -1,\lambda \to c,\rho \to c,b\to a,m\to \alpha -1,p\to 1$ and simplify;
\begin{multline}
\int_0^a (a-x)^{-1+\beta } x^{-1+\alpha } J_v(c x) J_{\mu }(c x) \, dx\\
=\frac{\left(\Gamma (\beta ) 2^{-v-\mu } a^{-1+v+\alpha +\beta +\mu } c^{v+\mu }\right) }{\Gamma (1+\mu )}\sum _{j=0}^{\infty
   } \frac{(-1)^j 2^{-2 j} a^{2 j} c^{2 j} \Gamma (2 j+v+\alpha +\mu ) }{\Gamma (1+j) \Gamma (1+j+v) \Gamma (2 j+v+\alpha +\beta +\mu )}\\ \times \, _2F_3\left(j+\frac{v}{2}+\frac{\alpha }{2}+\frac{\mu }{2},\frac{1}{2}+j+\frac{v}{2}+\frac{\alpha }{2}+\frac{\mu
   }{2};\right. \\ \left.
j+\frac{v}{2}+\frac{\alpha }{2}+\frac{\beta }{2}+\frac{\mu }{2},\frac{1}{2}+j+\frac{v}{2}+\frac{\alpha }{2}+\frac{\beta }{2}+\frac{\mu }{2},1+\mu ;-\frac{1}{4} a^2
   c^2\right)\\
=\frac{a^{\alpha +\beta +\mu +v-1} \left(\frac{c}{2}\right)^{\mu +v} (\Gamma (\beta )
   \Gamma (\alpha +\mu +v)) }{\Gamma (\mu +1) \Gamma (v+1) \Gamma (\alpha +\beta +\mu +v)}\\ \times
\, _4F_5\left(\frac{1}{2} (\mu +v+1),\frac{1}{2} (\mu +v+2),\frac{1}{2} (\mu +v+\alpha ),\frac{1}{2} (\mu +v+\alpha +1);\right. \\ \left.
\mu +1,v+1,\mu +v+1,\frac{1}{2} (\mu
   +v+\alpha +\beta ),\frac{1}{2} (\mu +v+\alpha +\beta +1);-(a c)^2\right)
\end{multline}
where $Re(a)>0,Re(\beta)>0,Re(\mu+\alpha+v)>0$.
\end{example}
\begin{example}
Derivation of equation (2.12.32.3) in \cite{prud2}. Here we use equation (\ref{eq:theorem}) and set $k\to 0,a\to 1,p\to 2,\beta \to \beta -1,c\to -\frac{1}{a^2},\lambda \to c,\rho \to c,m\to \alpha -1,b\to a$ and simplify;
\begin{multline}
\int_0^a x^{-1+\alpha } \left(a^2-x^2\right)^{-1+\beta } J_v(c x) J_{\mu }(c x) \, dx
=\frac{\left(2^{-1-v-\mu } a^{v+\alpha +2 (-1+\beta )+\mu } c^{v+\mu } \Gamma (\beta )\right)
  }{\Gamma (1+\mu )} \\
  \sum _{j=0}^{\infty } \frac{(-1)^j 2^{-2 j} a^{2 j} c^{2 j} \Gamma \left(\frac{1}{2} (2 j+v+\alpha +\mu )\right) }{\Gamma (1+j) \Gamma (1+j+v) \Gamma \left(\frac{1}{2} (2 j+v+\alpha +2 \beta +\mu)\right)}\\ \times
\, _1F_2\left(j+\frac{v}{2}+\frac{\alpha }{2}+\frac{\mu
   }{2};j+\frac{v}{2}+\frac{\alpha }{2}+\beta +\frac{\mu }{2},1+\mu ;-\frac{1}{4} a^2 c^2\right)\\
=\frac{a^{\alpha +2 \beta +\mu +v-2} \left(\frac{c}{2}\right)^{\mu +v} \left(\Gamma (\beta ) \Gamma \left(\frac{1}{2} (\alpha +\mu +v)\right)\right) }{2 \left(\Gamma \left(\beta+\frac{1}{2} (\alpha +\mu +v)\right) \Gamma (\mu +1) \Gamma (v+1)\right)}\\ \times
\,_3F_4\left(\frac{1}{2} (\alpha +\mu +v),\frac{1}{2} (\mu +v+1),\frac{\mu +v}{2}+1;\right. \\ \left.
\beta +\frac{1}{2} (\alpha +\mu +v),\mu +1,\mu +v+1,v+1;-(a c)^2\right)
\end{multline}
where $Re(a)>0,Re(\beta)>0,Re(\mu+\alpha+v)>0$.
\end{example}
\begin{example}
Derivation of equation (2.12.32.4a) in \cite{prud2}. Here we use equation (\ref{eq:theorem}) and set $k\to 0,a\to 1,\mu \to -v,\lambda \to c,\rho \to c,c\to -\frac{1}{a^2},p\to 2,\beta \to -v-\frac{3}{2},m\to 2 v+1,b\to a$ and simplify;
\begin{multline}
\int_0^a \frac{x^{1+2 v} J_{-v}(c x) J_v(c x)}{\left(a^2-x^2\right)^{\frac{3}{2}+v}} \, dx\\
=\sum _{j=0}^{\infty } \frac{(-1)^j 2^{-1-2 j} a^{-1+2 j} c^{2 j} \Gamma
   \left(-\frac{1}{2}-v\right) \, _1\tilde{F}_2\left(1+j+v;\frac{1}{2}+j,1-v;-\frac{1}{4} a^2 c^2\right)}{\Gamma (1+j)}\\
=\frac{(a c)^v \Gamma \left(-v-\frac{1}{2}\right) J_{-v}(2 a c)}{2
   a \sqrt{\pi }}
\end{multline}
where $Re(a)>0,-1< \pm Re(v)<1/2$.
\end{example}
\begin{example}
Derivation of equation (2.12.32.4b) in \cite{prud2}. Errata. Here we use equation (\ref{eq:theorem}) and set $k\to 0,a\to 1,\mu \to -v,\lambda \to c,\rho \to c,c\to -\frac{1}{a^2},p\to 2,\beta \to v-\frac{3}{2},m\to 1-2 v,b\to a$ and simplify;
\begin{multline}
\int_0^a x^{1-2 v} \left(a^2-x^2\right)^{-\frac{3}{2}+v} J_{-v}(c x) J_v(c x) \, dx\\
=\sum _{j=0}^{\infty } \frac{(-1)^j 2^{1-2 v} a^{-1+2 j} c^{2 j} \Gamma (1+j-v) \Gamma (-1+2 v)
    \sin (\pi  v)}{\pi  \Gamma (1+2 j) \Gamma (1+j+v)}\\ \times \, _1F_2\left(1+j-v;\frac{1}{2}+j,1-v;-\frac{1}{4} a^2 c^2\right)\\
=\frac{(a c)^{-v} \Gamma \left(v-\frac{1}{2}\right) J_v(2 a c)}{2 a
   \sqrt{\pi }}
\end{multline}
where $Re(a)>0,-1< Re(v)<1/2$.
\end{example}
\begin{example}
Derivation of equation (2.12.32.5) in \cite{prud2}. Here we use equation (\ref{eq:theorem}) and set $k\to 0,a\to 1,\mu \to v+1,\lambda \to c,\rho \to c,c\to -\frac{1}{a^2},p\to 2,\beta \to -v-\frac{3}{2},m\to 2 v+2,b\to a$ and simplify;
\begin{multline}
\int_0^a x^{2+2 v} \left(a^2-x^2\right)^{-\frac{3}{2}-v} J_v(c x) J_{1+v}(c x) \, dx\\
=\sum _{j=0}^{\infty } \frac{(-1)^j 4^{1+v} a^{1+2 j+2 v} c^{1+2 j+2 v} \Gamma (-2 (1+v))
   \Gamma (2+j+2 v) \sin (\pi  v)}{\pi  \Gamma (1+j) \Gamma (2 (1+j+v))}\\ \times
\, _1F_2\left(2+j+2 v;2+v,\frac{3}{2}+j+v;-\frac{1}{4} a^2 c^2\right) \\
=\frac{(a c)^v \Gamma
   \left(-v-\frac{1}{2}\right) J_{v+1}(2 a c)}{2 \sqrt{\pi }}
\end{multline}
where $Re(a)>0,-1< Re(v)<1/2$.
\end{example}
\begin{example}
Derivation of equation (2.12.32.6) in \cite{prud2}. Here we use equation (\ref{eq:theorem}) and set $k\to 0,a\to 1,\mu \to 1,v\to 1,\lambda \to c,\rho \to c,c\to -\frac{1}{a^2},p\to 2,\beta \to -\frac{1}{2},m\to -1,b\to a$ and simplify;
\begin{multline}
\int_0^a \frac{J_1(c x){}^2}{x \sqrt{a^2-x^2}} \, dx\\
=\sum _{j=0}^{\infty } \frac{(-1)^j 2^{-3-2 j} a^{1+2 j} c^{2+2 j} \sqrt{\pi } \, _1F_2\left(1+j;2,\frac{3}{2}+j;-\frac{1}{4}
   a^2 c^2\right)}{\Gamma (2+j) \Gamma \left(\frac{1}{2} (3+2 j)\right)}\\=\frac{1}{2 a}-\frac{J_1(2 a c)}{2 a^2 c}
\end{multline}
where $Re(a)>0$.
\end{example}
\begin{example}
Derivation of equation (2.12.32.7) in \cite{prud2}. Here we use equation (\ref{eq:theorem}) and set $k\to 0,a\to 1,\mu \to 2,v\to 2,\lambda \to c,\rho \to c,c\to -\frac{1}{a^2},p\to 2,\beta \to -\frac{1}{2},m\to -3,b\to a$ and simplify;
\begin{multline}
\int_0^a \frac{J_2(c x){}^2}{x^3 \sqrt{a^2-x^2}} \, dx\\
=\sum _{j=0}^{\infty } \frac{(-1)^j 2^{-6-2 j} a^{1+2 j} c^{4+2 j} \sqrt{\pi } \, _1F_2\left(1+j;3,\frac{3}{2}+j;-\frac{1}{4}
   a^2 c^2\right)}{\Gamma (3+j) \Gamma \left(\frac{1}{2} (3+2 j)\right)}\\
=\frac{\frac{1}{2} (a c)^2+\frac{1}{6} (a c)^4-J_2(2 a c)-2 a c J_3(2 a c)}{4 a^5 c^2}
\end{multline}
where $Re(a)>0$.
\end{example}
\section{Generalized Lipschitz-Hankel integral with two Bessel functions}
%
%\section{Introduction}
%
Definite integrals involving the product of Bessel functions of the first kind $J_{\mu }(x \alpha ) J_{\nu }(x \beta )$ according to Eason et al. \cite{eason}, were first studied and applied to physical problems by (Smythe, 1939), \cite{smythe} in his work on the potential in cylindrical co-ordinates $(\rho,z,\phi)$ due to a ring of radius $a$ and carrying a charge $Q$, (Bateman 1944, pp. 417), \cite{batemen} in his work on the steam function for a ring, and (Terazawa 1916), \cite{terazawa} in his work on the components of stress in a semi-infinite solid. These problems are related to hydrodynamics and the theory of steady-state conduction of heat and the solutions of more difficult problems involving series solutions. Other applications involving the product of Bessel functions were studied by (Carlson, 1980), \cite{carlson} in his work on the Laplace of the product of Bessel functions. (Benton, 1973) \cite{benton} studied infinite integrals involving the product of Bessel functions and expressed these integrals in terms of complete elliptic integrals $K$ and $E$. (Oliver, 2021) \cite{oliver}, extnded a useful identity relating the infinite sum of two Bessel functions to their infinite integral. In [DLMF,\href{https://dlmf.nist.gov/10.22.vi}{10.22(vi)}], for collections of integrals of Bessel functions, including integrals with respect to the order, see Andrews et al. (\cite{andrews}, 1999, pp. 216-225), Apelblat (\cite{apelblat},1983, §12), Erdélyi et al. (\cite{erd2}, 1953b, §§7.7.1-7.7.7 and 7.14-7.14.2), Erdélyi et al. (\cite{erdt1}, 1954a, b), Gradshteyn and Ryzhik (\cite{grad}, 2015, §§5.5 and 6.5-6.7), Gröbner and Hofreiter (\cite{grobner}, 1950, pp. 196-204), Luke (\cite{luke}, 1962), Magnus et al. (\cite{magnus}, 1966, §3.8), Marichev (\cite{marichev}, 1983, pp. 191-216), Oberhettinger (\cite{oberm}, 1974, §§1.10 and 2.7), Oberhettinger (\cite{oberf}, 1990, §§1.13-1.16 and 2.13-2.16), Oberhettinger and Badii (\cite{oberl}, 1973, §§1.14 and 2.12), Okui (\cite{okui_1974,okui_1975}, 1974, 1975), Prudnikov et al. (\cite{prud2}, 1986b, §§1.8-1.10, 2.12-2.14, 3.2.4-3.2.7, 3.3.2, and 3.4.1), Prudnikov et al. (\cite{prud4}, 1992a, §§3.12-3.14), Prudnikov et al. (\cite{prud5}, 1992b, §§3.12-3.14), Watson (\cite{watson}, 1944, Chapters 5, 12, 13, and 14), and Wheelon (\cite{wheelon}, 1968). This compendium is intended for applied mathematicians , physicists and engineers studying problems involving finite and infinite integrals involving Bessel functions.
\\\\
In this work using contour integration, we will derive and evaluate the finite integral given by;
\begin{multline}
\int_0^b e^{-x^r \gamma } x^m \left(1+c x^s\right)^{\delta } J_{\mu }(x \alpha ) J_{\nu }(x \beta ) \log
   ^k\left(\frac{1}{a x}\right) \, dx\\
=\sum _{j=0}^{\infty } \sum _{h=0}^{\infty } \sum _{l=0}^{\infty } \sum
   _{p=0}^{\infty } \frac{(-1)^{h+j+l} 2^{-2 h-2 j-\mu -\nu } a^{-1-2 h-2 j-m-l r-p s-\mu -\nu } c^p \alpha ^{2
   j+\mu } \beta ^{2 h+\nu } \gamma ^l \binom{\delta }{p} }{h! j! l! \Gamma (1+j+\mu ) \Gamma (1+h+\nu ) (1+2 h+2 j+m+l r+p s+\mu +\nu
   )^{k+1}}\\ \times 
\Gamma \left(1+k,(1+2 h+2 j+m+l r+p s+\mu +\nu ) \log
   \left(\frac{1}{a b}\right)\right)
\end{multline}
where $Re(\mu+\nu+m)>0,Re(r)>0$.
The derivation involves a definite integral over a finite interval and follows the method used by us in~\cite{reyn4}. This method involves using a form of the generalized Cauchy's integral formula given by
\begin{equation}\label{intro:cauchy}
\frac{y^k}{\Gamma(k+1)}=\frac{1}{2\pi i}\int_{C}\frac{e^{wy}}{w^{k+1}}dw.
\end{equation}
where $C$ is in general an open contour in the complex plane where the bilinear concomitant has the same value at the end points of the contour. We then multiply both sides by a function of $x$, then take a definite integral of both sides. This yields a definite integral in terms of a contour integral. Then we multiply both sides of Equation~(\ref{intro:cauchy})  by another function of $y$ and take the infinite sum of both sides such that the contour integral of both equations are the same.
The contour integral representation form involving the generalized Lipschitz-Hankel form is given by;
\begin{multline}\label{eq:p1}
\frac{1}{2\pi i}\int_{C}\int_{0}^{b}a^w e^{-x^r \gamma } w^{-1-k} x^{m+w} \left(1+c x^s\right)^{\delta } J_{\mu }(x \alpha ) J_{\nu }(x \beta)dxdw\\
=\frac{1}{2\pi i}\int_{C}\sum _{j=0}^{\infty } \sum_{h=0}^{\infty }\sum _{l=0}^{\infty }\sum _{p=0}^{\infty } \frac{(-1)^{h+j} 2^{-2 h-2 j-\mu -\nu } a^w b^{1+2h+2 j+m+l r+p s+w+\mu +\nu } c^p w^{-1-k}}{(1+2 h+2 j+m+l r+p s+w+\mu +\nu ) h! j! l! }\\ \times\frac{ \alpha ^{2 j+\mu } \beta ^{2 h+\nu } (-\gamma )^l \binom{\delta}{p}}{\Gamma (1+j+\mu ) \Gamma (1+h+\nu )}dw
\end{multline}
where $Re(\mu+\nu+m+w)>0, Re(b)>0$.
\subsection{Left-hand side contour integral representation}
Using a generalization of Cauchy's integral formula \ref{intro:cauchy}, we form the definite integral by replacing $y$ by $\log{ax}$ and multiply both sides by 
\begin{equation}
\frac{c^p (-1)^{h+j} (-\gamma )^l \binom{\delta }{p} \beta ^{2 h+\nu } \alpha ^{2 j+\mu } 2^{-2 h-2 j-\mu -\nu } x^{2
   h+2 j+l r+\mu +m+\nu +p s}}{h! j! l! \Gamma (h+\nu +1) \Gamma (j+\mu +1)}\end{equation} 
then take the infinite sum of both sides over $h \in [0,\infty),j\in [0,\infty),l \in [0,\infty),p \in [0,\infty)$ to get;
\begin{multline}\label{eq:p2}
\int_0^b e^{-x^r \gamma } x^m \left(1+c x^s\right)^{\delta } J_{\mu }(x \alpha ) J_{\nu }(x \beta ) \log
   ^k\left(\frac{1}{a x}\right) \, dx\\
=\frac{1}{2\pi i}\sum _{h=0}^{\infty } \sum _{j=0}^{\infty } \sum _{l=0}^{\infty }\sum _{p=0}^{\infty }\int_{C}\int_{0}^{b}\frac{\alpha ^{\mu } a^w \beta ^{\nu } c^p \left(\beta ^2\right)^h (-1)^{h+j} \left(\alpha ^2\right)^j w^{-k-1}(-\gamma )^l }{h! j! l! \Gamma (h+\nu+1) \Gamma (j+\mu +1)}\\ \times
\binom{\delta }{p} 2^{-2 h-2 j-\mu -\nu } x^{2 h+2 j+l r+\mu +m+\nu +p s+w}dxdw\\
=\frac{1}{2\pi i}\int_{C}\int_{0}^{b}\sum _{h=0}^{\infty } \sum _{j=0}^{\infty } \sum _{l=0}^{\infty }\sum _{p=0}^{\infty }\frac{\alpha ^{\mu } a^w \beta ^{\nu } c^p \left(\beta ^2\right)^h (-1)^{h+j} \left(\alpha ^2\right)^j w^{-k-1}(-\gamma )^l }{h! j! l! \Gamma (h+\nu+1) \Gamma (j+\mu +1)}\\ \times
\binom{\delta }{p} 2^{-2 h-2 j-\mu -\nu } x^{2 h+2 j+l r+\mu +m+\nu +p s+w}dxdw\\
=\frac{1}{2\pi i}\int_{C}\int_{0}^{b}w^{-k-1} x^m (a x)^w e^{\gamma  \left(-x^r\right)} \left(c x^s+1\right)^{\delta } J_{\mu }(x
   \alpha ) J_{\nu }(x \beta )dxdw
\end{multline}
where $Re(\alpha)>0,Re(\beta)>0,Re(\gamma)>0,Re(z)<1$.
We are able to switch the order of integration over $t$ and $w$ using Fubini's theorem for multiple integrals see page 178 in \cite{gelca}, since the integrand is of bounded measure over the space $\mathbb{C} \times [0,b]$.
\subsection{Right-hand side contour integral}
Using a generalization of Cauchy's integral formula \ref{intro:cauchy}, we replace $y\to y+xi+\log(a)$ and multiply both sides by $e^{mxi}$. Then form a second equation by replacing $a\to -a$ and take their difference to get;
 \begin{multline}\label{eq:4.3}
\frac{i e^{-i m x} \left((\log (a)-i x+y)^k-e^{2 i m x} (\log (a)+i x+y)^k\right)}{2 k!}\\
=\frac{1}{2\pi i}\int_{C}a^w w^{-k-1} e^{w y}
   \sin (x (m+w))dw
\end{multline}
Next we take the indefinite integral with respect to $y$ to get;
\begin{multline}\label{eq:4.4}
-\frac{i (-1)^{k+1} a^{-m} m^{-k-1} (\Gamma (k+1,-m (-i x+y+\log (a)))-\Gamma (k+1,-m (i x+y+\log (a))))}{2 k!}\\
=\frac{1}{2\pi i}\int_{C}\frac{a^w w^{-k-1} e^{y (m+w)} \sin (x (m+w))}{m+w}dw
\end{multline}
from equation equation (2.325.6(12i)) in \cite{grad}. Finally we replace $y\to -\gamma, x\to \beta$, multiply both sides by $2\pi$ and simplify to get;
\begin{multline}\label{eq:4.5}
\frac{i \pi  (-1)^{k+1} a^{-m} m^{-k-1} (\Gamma (k+1,m (-i \beta +\gamma -\log (a)))-\Gamma (k+1,m (i \beta
   +\gamma -\log (a))))}{k!}\\
=\frac{1}{2\pi i}\int_{C}\frac{2 \pi  a^w w^{-k-1} e^{-(\gamma  (m+w))} \sin (\beta  (m+w))}{m+w}dw
\end{multline}
Use equation (\ref{eq:4.4}) and set $y=0$, replace $m\to 2 h+2 j+l r+\mu +m+\nu +p s+1,a\to a b$ and multiply both sides by 
\begin{equation}
\frac{\alpha ^{\mu } \beta ^{\nu } c^p \left(\beta ^2\right)^h (-1)^{h+j} \left(\alpha ^2\right)^j (-\gamma )^l
   \binom{\delta }{p} 2^{-2 h-2 j-\mu -\nu } b^{2 h+2 j+l r+\mu +m+\nu +p s+1}}{h! j! l! \Gamma (h+\nu +1) \Gamma
   (j+\mu +1)}.
\end{equation}

Next take the infinite sum of both sides over $h \in [0,\infty),j \in [0,\infty),l \in [0,\infty),p \in [0,\infty)$ and simplify to get;
\begin{multline}\label{eq:p3}
\sum _{h=0}^{\infty } \sum _{j=0}^{\infty } \sum _{l=0}^{\infty }\sum_{p=0}^{\infty} \frac{(-1)^{h+j} 2^{-2 h-2 j-\mu -\nu } a^{-1-2 h-2 j-m-l r-p s-\mu -\nu } c^p \alpha ^{2 j+\mu } \beta ^{2 h+\nu } (-\gamma )^l \binom{\delta }{p} }{h! j! k! l! \Gamma (1+j+\mu ) \Gamma (1+h+\nu ) (-1-2 h-2 j-m-l r-p s-\mu -\nu )^{k+1}}\\ \times
\Gamma (1+k,-((1+2 h+2 j+m+l r+p s+\mu +\nu ) \log (a b)))\\
=-\frac{1}{2\pi i}\sum _{h=0}^{\infty } \sum _{j=0}^{\infty } \sum
   _{l=0}^{\infty }\sum _{l=0}^{\infty } \int_{C}\frac{(-1)^{h+j} 2^{-2 h-2 j-\mu -\nu } b^{1+2 h+2 j+m+l r+p
   s+\mu +\nu } (a b)^w c^p w^{-1-k} \alpha ^{\mu } \alpha ^{2 j} }{(1+2 h+2 j+m+l r+p s+w+\mu +\nu ) h! j! l! }\\ \times
   \frac{\beta ^{\nu } \beta ^{2 h} (-\gamma )^l \binom{\delta }{p}}{\Gamma (1+j+\mu
   ) \Gamma (1+h+\nu )}dw\\
   =-\frac{1}{2\pi i} \int_{C}\sum _{h=0}^{\infty } \sum _{j=0}^{\infty } \sum
   _{l=0}^{\infty }\sum _{l=0}^{\infty }\frac{(-1)^{h+j} 2^{-2 h-2 j-\mu -\nu } b^{1+2 h+2 j+m+l r+p
   s+\mu +\nu } (a b)^w c^p w^{-1-k} \alpha ^{\mu } \alpha ^{2 j} }{(1+2 h+2 j+m+l r+p s+w+\mu +\nu ) h! j! l! }\\ \times
   \frac{\beta ^{\nu } \beta ^{2 h} (-\gamma )^l \binom{\delta }{p}}{\Gamma (1+j+\mu
   ) \Gamma (1+h+\nu )}dw\\ 
\end{multline}
where $Re(\alpha)>0,Re(\beta)>0,Re(\gamma)>0,Re(z)<1$.
from equation [Wolfram, \href{http://functions.wolfram.com/07.27.02.0002.01}{07.27.02.0002.01}] where $|Re(b)|<1$. We are able to switch the order of integration and summation over $w$ using Tonellii's theorem for  integrals and sums see page 177 in \cite{gelca}, since the summand is of bounded measure over the space $\mathbb{C} \times [0,\infty)$
\begin{theorem}
Generalized Lipschitz-Hankel integral with two Bessel functions.
\begin{multline}\label{eq:thm}
\int_0^b e^{-x^r \gamma } x^m \left(1+c x^s\right)^{\delta } J_{\mu }(x \alpha ) J_{\nu }(x \beta ) \log
   ^k\left(\frac{1}{a x}\right) \, dx\\
=\sum _{j=0}^{\infty } \sum _{h=0}^{\infty } \sum _{l=0}^{\infty } \sum
   _{p=0}^{\infty } \frac{(-1)^{h+j+l} 2^{-2 h-2 j-\mu -\nu } a^{-1-2 h-2 j-m-l r-p s-\mu -\nu } c^p \alpha ^{2
   j+\mu } \beta ^{2 h+\nu } \gamma ^l \binom{\delta }{p} }{h! j! l! \Gamma (1+j+\mu ) \Gamma (1+h+\nu ) (1+2 h+2 j+m+l r+p s+\mu +\nu
   )^{k+1}}\\ \times 
\Gamma \left(1+k,(1+2 h+2 j+m+l r+p s+\mu +\nu ) \log
   \left(\frac{1}{a b}\right)\right)
\end{multline}
where $Re(\mu+\nu+m)>0,Re(r)>0$.
\end{theorem}
\begin{proof}
Since the right-hand sides of equations (\ref{eq:p2}) and (\ref{eq:p3}) are equal relative to equation (\ref{eq:p1}), we may equate the left-hand sides and simplify the gamma function to yield the stated result.
\end{proof}
\section{Evaluations derivations involving finite and infinite integrals}
In the following Tables we will evaluate equation (\ref{eq:thm}) in terms of an infinite series involving the generalized hypergeometric function and derive alternate forms in terms of special functions present in current literature. We will evaluate finite and infinite definite integrals where the restrictions on the parameters are stated for each example. We also give alternate forms for Errata found in some examples.
\begin{example}
Mellin transform. This integral form can be used derive entries (3.12.15.1-64), (3.12.16.1-19) in \cite{prud4}. Here we use equation (\ref{eq:thm}) and set $k\to 0,a\to 1,\delta \to 0,s\to 0,c\to 1$ and simplify;
\begin{multline}\label{eq:thm1}
\int_0^{\infty } e^{-x^r \gamma } x^m J_{\mu }(x \alpha ) J_{\nu }(x \beta ) \, dx\\
=\sum _{h=0}^{\infty } \sum
   _{j=0}^{\infty } \frac{(-1)^{h+j} 2^{-2 h-2 j-\mu -\nu } \alpha ^{2 j+\mu } \beta ^{2 h+\nu } \gamma ^{-\frac{1+2
   h+2 j+m+\mu +\nu }{r}} \Gamma \left(\frac{1+2 h+2 j+m+\mu +\nu }{r}\right)}{r h! j! \Gamma (1+j+\mu ) \Gamma
   (1+h+\nu )}
\end{multline}
where $Re(\mu+m+v)>0,Re(r)>1$.
\end{example}
\begin{example}
Derivation of equation 7.7.3.20 in \cite{erd2}. Here we use equation (\ref{eq:thm1}) and set $r\to 2,\gamma \to \gamma ^2,m\to \lambda -1$ and simplify;
\begin{multline}\label{eq:thm2}
\int_0^{\infty } e^{-x^2 \gamma ^2} x^{-1+\lambda } J_{\mu }(x \alpha ) J_{\nu }(x \beta ) \, dx\\
=\sum
   _{j=0}^{\infty } \frac{(-1)^j 2^{-1-2 j-\mu -\nu } \alpha ^{2 j+\mu } \beta ^{\nu } \gamma ^{-2 j-\lambda -\mu
   -\nu } \Gamma \left(\frac{1}{2} (2 j+\lambda +\mu +\nu )\right) }{\Gamma (1+j) \Gamma (1+j+\mu ) \Gamma (1+\nu
   )}\\ \times
\, _1F_1\left(j+\frac{\lambda }{2}+\frac{\mu}{2}+\frac{\nu }{2};1+\nu ;-\frac{\beta ^2}{4 \gamma ^2}\right)\\
=\frac{1}{2^{\nu +\mu +1} \alpha ^{-\mu } \beta ^{-\nu } \gamma ^{\nu +\mu +\lambda } \Gamma (\nu
   +1)}\sum _{m=0}^{\infty } \frac{\Gamma \left(m+\frac{\nu }{2}+\frac{\mu }{2}+\frac{\lambda }{2}\right) }{m!\Gamma (m+\mu +1)}\\ \times
\,
   _2F_1\left(-m,-\mu -m;\nu +1;\frac{\beta ^2}{\alpha ^2}\right) \left(-\frac{\alpha ^2}{4 \gamma ^2}\right)^m
\end{multline}
where $Re(\mu+\nu+\lambda)>0,Re(r)>0$.
\end{example}
\begin{example}
Series form for 7.7.3.20 in \cite{erd2}. Here we use equation (\ref{eq:thm2}) and simplify using [DLMF,\href{https://dlmf.nist.gov/15.4.E20}{15.4.20}].
\begin{multline}
\sum _{j=0}^{\infty } \frac{(-1)^j \left(\frac{\alpha }{2 \gamma }\right)^{2 j} \Gamma \left(\frac{1}{2} (2
   j+\lambda +\mu +\nu )\right) \, _1F_1\left(j+\frac{\lambda }{2}+\frac{\mu }{2}+\frac{\nu }{2};1+\nu ;-\frac{\beta
   ^2}{4 \gamma ^2}\right)}{\Gamma (1+j) \Gamma (1+j+\mu )}\\
=\sum _{m=0}^{\infty } \frac{\Gamma \left(m+\frac{\nu
   }{2}+\frac{\mu }{2}+\frac{\lambda }{2}\right) \, _2F_1\left(-m,-\mu -m;\nu +1;\frac{\beta ^2}{\alpha ^2}\right)
   \left(-\frac{\alpha ^2}{4 \gamma ^2}\right)^m}{m! \Gamma (m+\mu +1)}
\end{multline}
where $Re(\mu+\nu+\lambda)>0,Re(r)>0$.
\end{example}
\begin{example}
Derivation of equation 7.7.3.21 in \cite{erd2}. 
\begin{multline}\label{eq:thm3}
\sum _{j=0}^{\infty } \frac{\left(-\frac{1}{4}\right)^j \left(\frac{\alpha }{\gamma }\right)^{2 j} \Gamma
   \left(\frac{1}{2} (2 j+\lambda +\mu +\nu )\right) \, _1F_1\left(\frac{1}{2} (2 j+\lambda +\mu +\nu );1+\nu
   ;-\frac{\alpha ^2}{4 \gamma ^2}\right)}{\Gamma (1+j) \Gamma (1+j+\mu )}\\
=\frac{\Gamma \left(\frac{1}{2} (\lambda
   +\mu +\nu )\right) \, _3F_3\left(\frac{1}{2}+\frac{\mu }{2}+\frac{\nu }{2},1+\frac{\mu }{2}+\frac{\nu
   }{2},\frac{\lambda }{2}+\frac{\mu }{2}+\frac{\nu }{2};1+\mu ,1+\nu ,1+\mu +\nu ;-\frac{\alpha ^2}{\gamma
   ^2}\right)}{\Gamma (1+\mu )}
\end{multline}
where $Re(\mu+\lambda+v)>0,Re(s)>0$.
\end{example}
\section{Derivations of some Bessel integrals in Gradshteyn and Ryzhik Table 6.612}
\begin{example}
Derivation of equation 6.612.3 in \cite{grad}. Here we use equation (\ref{eq:thm1}) and set $r\to 1,m\to 0,\mu \to \nu ,\gamma \to \alpha ,\alpha \to b,\beta \to \gamma $ and simplify;
\begin{multline}\label{eq:thm4}
\int_0^{\infty } e^{-x \alpha } J_{\nu }(b x) J_{\nu }(x \gamma ) \, dx\\
=\sum _{h=0}^{\infty } \frac{(-1)^h
   b^{\nu } \alpha ^{-1-2 h-2 \nu } \gamma ^{2 h+\nu } \Gamma \left(\frac{1}{2} (1+2 h+2 \nu )\right) \,
   _2F_1\left(\frac{1}{2}+h+\nu ,1+h+\nu ;1+\nu ;-\frac{b^2}{\alpha ^2}\right)}{\sqrt{\pi } \Gamma (1+h) \Gamma
   (1+\nu )}\\
\neq\frac{Q_{\nu -\frac{1}{2}}\left(\frac{\alpha ^2+b^2+\gamma ^2}{2 b \gamma }\right)}{\pi  \sqrt{\gamma 
   b}}
\end{multline}
where $Re(\alpha\pm ib\pm i\gamma)>0,Re(\gamma)>0,Re(\nu)>-1/2$.
\end{example}
\begin{example}
Derivation of equation (6.612.4) in \cite{grad}. Here we use equation (\ref{eq:thm1}) and set $r\to 1,m\to 0,\mu \to 0,\nu \to 0,\alpha \to b,\beta \to b$ and simplify;
\begin{multline}\label{eq:thm5}
\int_0^{\infty } e^{-x \gamma } J_0(b x){}^2 \, dx\\
=\sum _{h=0}^{\infty } \frac{(-1)^h b^{2 h} \gamma ^{-1-2
   h} \Gamma \left(\frac{1}{2} (1+2 h)\right) \, _2F_1\left(\frac{1}{2}+h,1+h;1;-\frac{b^2}{\gamma
   ^2}\right)}{\sqrt{\pi } \Gamma (1+h)}\\
=\frac{2 K\left(\frac{4 b^2}{4 b^2+\gamma ^2}\right)}{\pi  \sqrt{\gamma ^2+4
   b^2}}
\end{multline}
where $Re(\gamma)>0$.
\end{example}
\begin{example}
Derivation of equation (6.612.5) in \cite{grad}. Here we use equation (\ref{eq:thm1}) and set $r\to 1,m\to 0,\mu \to 1,\nu \to 1,\alpha \to b,\beta \to b,\gamma \to 2 \alpha$ and simplify;
\begin{multline}\label{eq:thm6}
\int_0^{\infty } e^{-2 x \alpha } J_1(b x){}^2 \, dx\\
=\sum _{h=0}^{\infty } \frac{(-1)^h 2^{-3-2 h} b^{2+2 h}
   \alpha ^{-3-2 h} \Gamma \left(\frac{1}{2} (3+2 h)\right) \, _2F_1\left(\frac{3}{2}+h,2+h;2;-\frac{b^2}{4 \alpha
   ^2}\right)}{\sqrt{\pi } \Gamma (1+h)}\\
=\frac{\left(2 \alpha ^2+b^2\right) K\left(\frac{b^2}{b^2+\alpha
   ^2}\right)-2 \left(\alpha ^2+b^2\right) E\left(\frac{b^2}{b^2+\alpha ^2}\right)}{\pi  b^2 \sqrt{\alpha
   ^2+b^2}}
\end{multline}
where $Re(\alpha)>0$.
\end{example}
\begin{example}
Derivation of equation (6.6.5) in \cite{grad}. Here we use equation (\ref{eq:thm1}) and set $\gamma \to \alpha ,\mu \to \nu ,\beta \to 2 \beta ,\alpha \to 2 \gamma ,x\to \sqrt{t}$ then $t\to x,r\to 2,m\to 1$ and simplify;
\begin{multline}\label{eq:thm7}
\int_0^{\infty } e^{-x \alpha } J_{\nu }\left(2 \sqrt{x} \beta \right) J_{\nu }\left(2 \sqrt{x} \gamma
   \right) \, dx\\
=\sum _{h=0}^{\infty } \frac{(-1)^h \alpha ^{-1-h-\nu } \beta ^{2 h+\nu } \gamma ^{\nu } \,
   _1F_1\left(1+h+\nu ;1+\nu ;-\frac{\gamma ^2}{\alpha }\right)}{\Gamma (1+h) \Gamma (1+\nu )}\\
=\frac{I_{\nu
   }\left(\frac{2 \beta  \gamma }{\alpha }\right) \exp \left(-\frac{\beta ^2+\gamma ^2}{\alpha }\right)}{\alpha
   }
\end{multline}
where $Re(\nu)>-1$.
\end{example}
\begin{example}
 Derivation of equation (6.618.5) in \cite{grad}. Here we use equation (\ref{eq:thm1}) and set $r\to 2,\gamma \to \alpha ,m\to 0,\alpha \to b,\beta \to b$ and simplify;
 \begin{multline}\label{eq:thm8}
\int_0^{\infty } e^{-x^2 \alpha } J_{\mu }(b x) J_{\nu }(b x) \, dx\\
=\sum _{h=0}^{\infty } \frac{(-1)^h
   2^{-1-2 h-\mu -\nu } b^{2 h+\mu +\nu } \alpha ^{-\frac{1}{2}-h-\frac{\mu }{2}-\frac{\nu }{2}} \Gamma
   \left(\frac{1}{2} (1+2 h+\mu +\nu )\right) }{\Gamma (1+h) \Gamma (1+\mu ) \Gamma (1+h+\nu )}\\ \times 
\, _1F_1\left(\frac{1}{2}+h+\frac{\mu }{2}+\frac{\nu }{2};1+\mu
   ;-\frac{b^2}{4 \alpha }\right)\\
=\frac{2^{-\nu -\mu -1} \alpha
   ^{-\frac{1}{2} (\nu +\mu +1)} b^{\nu +\mu } \Gamma \left(\frac{1}{2} (\mu +\nu +1)\right) }{\Gamma (\mu +1) \Gamma (\nu +1)}\\ \times
\,_3F_3\left(\frac{1}{2} (\mu +\nu +1),\frac{1}{2} (\mu +\nu +2),\frac{1}{2} (\mu +\nu +1);\mu +1,\nu +1,\nu +\mu
   +1;-\frac{b^2}{\alpha }\right)
\end{multline}
where $Re(\nu+\mu)>-1,Re(\alpha)>0$.
\end{example}
\begin{example}
 Derivation of equation (6.618.5) in \cite{grad} series representation. Here we use equation (\ref{eq:thm8}) and simplify;
 \begin{multline}\label{eq:thm9}
 \sum _{h=0}^{\infty } \frac{\left(-\frac{b^2}{4 \alpha }\right)^h \Gamma \left(\frac{1}{2} (1+2 h+\mu +\nu
   )\right) \, _1F_1\left(\frac{1}{2}+h+\frac{\mu }{2}+\frac{\nu }{2};1+\mu ;-\frac{b^2}{4 \alpha }\right)}{\Gamma
   (1+h) \Gamma (1+h+\nu )}\\
   =\frac{\Gamma \left(\frac{1}{2} (\mu +\nu +1)\right) \, _3F_3\left(\frac{1}{2} (\mu +\nu
   +1),\frac{1}{2} (\mu +\nu +2),\frac{1}{2} (\mu +\nu +1);\mu +1,\nu +1,\nu +\mu +1;-\frac{b^2}{\alpha
   }\right)}{\Gamma (\nu +1)}
 \end{multline}
\end{example}
\begin{example}
Derivation of Table entries (6.621.5(10)-19) in \cite{grad} generalized series form.  Here we use equation (\ref{eq:thm1}) and set $r\to 1,\gamma \to z,\alpha \to a,\beta \to b$ and simplify;
\begin{multline}\label{eq:thm10}
\int_0^{\infty } e^{-x z} x^m J_{\mu }(a x) J_{\nu }(b x) \, dx\\
=\sum _{h=0}^{\infty } \frac{(-1)^h 2^{-2
   h-\mu -\nu } a^{\mu } b^{2 h+\nu } z^{-1-2 h-m-\mu -\nu } \Gamma (1+2 h+m+\mu +\nu ) }{\Gamma (1+h) \Gamma (1+\mu ) \Gamma (1+h+\nu )}\\ \times
\,
   _2F_1\left(\frac{1}{2}+h+\frac{m}{2}+\frac{\mu }{2}+\frac{\nu }{2},1+h+\frac{m}{2}+\frac{\mu }{2}+\frac{\nu
   }{2};1+\mu ;-\frac{a^2}{z^2}\right)
\end{multline}
where $\arg(a)>0,\arg(b)>0,\arg(z)>0$.
\end{example}
\begin{example}
Derivation of Table (6.626.1.-3,5) in \cite{grad} Laplace transform in terms of generalized series. Here we use equation (\ref{eq:thm1}) and set $r\to 1,\gamma \to \alpha ,m\to \lambda -1,\alpha \to b,\beta \to c$ and simplify;
\begin{multline}\label{eq:thm11}
\int_0^{\infty } e^{-x \alpha } x^{-1+\lambda } J_{\mu }(b x) J_{\nu }(c x) \, dx\\
=\sum _{h=0}^{\infty }
   \frac{(-1)^h 2^{-2 h-\mu -\nu } b^{\mu } c^{2 h+\nu } \alpha ^{-2 h-\lambda -\mu -\nu } \Gamma (2 h+\lambda +\mu
   +\nu ) }{\Gamma (1+h) \Gamma (1+\mu ) \Gamma(1+h+\nu )}\\ \times
\, _2F_1\left(h+\frac{\lambda }{2}+\frac{\mu }{2}+\frac{\nu }{2},\frac{1}{2}+h+\frac{\lambda
   }{2}+\frac{\mu }{2}+\frac{\nu }{2};1+\mu ;-\frac{b^2}{\alpha ^2}\right)
\end{multline}
where $Re(\lambda+\mu+\nu)>0$.
\end{example}
\begin{example}
Derivation of equation (1.14.14.16) in \cite{oberl}. Here we use equation (\ref{eq:thm11}) and set $\mu \to \nu ,b\to a,c\to a,\lambda \to 2 \nu +1$ and simplify;
\begin{multline}\label{eq:thm12}
\int_0^{\infty } e^{-x \alpha } x^{2 \nu } J_{\nu }(a x){}^2 \, dx\\
=\sum _{h=0}^{\infty } \frac{(-1)^h
   4^{-h-\nu } a^{2 (h+\nu )} \alpha ^{-1-2 h-4 \nu } \Gamma (1+2 h+4 \nu ) }{\Gamma (1+h) \Gamma (1+\nu ) \Gamma (1+h+\nu
   )}\\ \times
   \, _2F_1\left(\frac{1}{2}+h+2 \nu ,1+h+2
   \nu ;1+\nu ;-\frac{a^2}{\alpha ^2}\right)\\
=\frac{\left(\frac{4 a}{\alpha ^2}\right)^{2 \nu } \Gamma \left(\frac{1}{2}+\nu \right) \Gamma
   \left(\frac{1}{2}+2 \nu \right) \, _2F_1\left(\frac{1}{2}+\nu ,\frac{1}{2}+2 \nu ;1+\nu ;-\frac{4
   a^2}{\alpha ^2}\right)}{\pi  \alpha  \Gamma (\nu +1)}
\end{multline}
where $Re(\alpha)>2 Im(a)$.
\end{example}
\begin{example}
 Derivation of equation (1.14.14.16) in \cite{oberl} in terms of series. 
 \begin{multline}
\sum _{h=0}^{\infty } \frac{\left(-\frac{a^2}{4 \alpha ^2}\right)^h \Gamma (1+2 h+4 \nu ) \,
   _2F_1\left(\frac{1}{2}+h+2 \nu ,1+h+2 \nu ;1+\nu ;-\frac{a^2}{\alpha ^2}\right)}{\Gamma (1+h) \Gamma (1+h+\nu
   )}\\
=\frac{4^{3 \nu } }{\pi }\Gamma \left(\frac{1}{2}+\nu \right) \Gamma \left(\frac{1}{2}+2 \nu \right) \,
   _2F_1\left(\frac{1}{2}+\nu ,\frac{1}{2}+2 \nu ;1+\nu ;-\frac{4 a^2}{\alpha ^2}\right)
\end{multline}
where $Re(\alpha)>2 Im(a)$.
\end{example}
\begin{example}
 Derivation of equation (1.14.14.17) in \cite{oberl}, series form for a definite integral. Here we use equation (\ref{eq:thm11}) and set $\alpha \to p,\lambda \to \mu +\nu +1,b\to a,c\to a,x\to t$ and simplify;
 \begin{multline}
\int_0^{\infty } e^{-p t} t^{\mu +\nu } J_{\mu }(a t) J_{\nu }(a t) \, dt\\
=\sum _{h=0}^{\infty } \frac{(-1)^h
   2^{-2 h-\mu -\nu } a^{2 h+\mu +\nu } p^{-1-2 h-2 \mu -2 \nu } \Gamma (1+2 h+2 \mu +2 \nu ) 
}{\Gamma (1+h) \Gamma (1+\mu ) \Gamma (1+h+\nu )}\\ \times
\,_2F_1\left(h+\frac{\mu }{2}+\frac{\nu }{2}+\frac{1}{2} (1+\mu +\nu ),\frac{1}{2}+h+\frac{\mu }{2}+\frac{\nu}{2}+\frac{1}{2} (1+\mu +\nu );1+\mu ;-\frac{a^2}{p^2}\right)\\
=2
   \pi ^{-3/2} (4 a)^{\mu +\nu } \Gamma \left(\frac{1}{2}+\nu +\mu \right) \int_0^{\frac{\pi }{2}} \left(p^2+4 a^2
   \cos ^2(t)\right)^{-\mu -\nu -\frac{1}{2}} \cos ^{\mu +\nu }(t) \cos ((\mu -\nu ) t) \, dt
\end{multline}
where $Re(p)>2 Im(a)$.
\end{example}
\section{Table entries (8.11.1-30) in Erd\'{e}yli et al. (1954)}
\begin{example}
Derivation of equation (8.11.18) in \cite{erdt2} as a series representation for a definite integral. Here we use equation (\ref{eq:thm11}) and set $\lambda \to \mu ,\mu \to \nu ,b\to \beta ,c\to y$ and simplify;
\begin{multline}
\int_0^{\infty } e^{-x \alpha } x^{-1+\mu } \sqrt{y} J_{\nu }(x y) J_{\nu }(x \beta ) \, dx\\
=\sum_{h=0}^{\infty } \frac{(-1)^h 2^{-2 h-2 \nu } y^{\frac{1}{2}+2 h+\nu } \alpha ^{-2 h-\mu -2 \nu } \beta ^{\nu }\Gamma (2 h+\mu +2 \nu )}{\Gamma (1+h) \Gamma (1+\nu ) \Gamma (1+h+\nu )}\\ \times 
 \, _2F_1\left(h+\frac{\mu }{2}+\nu ,\frac{1}{2}+h+\frac{\mu }{2}+\nu ;1+\nu;-\frac{\beta ^2}{\alpha ^2}\right)\\
=\frac{\beta ^{\nu }y^{\nu +\frac{1}{2}} \Gamma (\mu +2 \nu )  }{\pi  \alpha ^{\mu +2 \nu } \Gamma (2 \nu +1)}\int_0^{\pi } \, _2F_1\left(\frac{\mu }{2}+\nu ,\frac{\mu+1}{2}+\nu ;\nu +1;-\left(\frac{\sqrt{\beta ^2+y^2-2 \beta  y \cos (\phi )}}{\alpha }\right)^2\right)\\
 \sin ^{2 \nu }(\phi ) \, d\phi
\end{multline}
where $Re(\alpha)>|Im(\beta)|, Re(\mu+2\nu)>0$.
\end{example}
\begin{example}
 Derivation of equation (8.11.21) in \cite{erdt2} as a series representation for a definite integral. Here we use equation (\ref{eq:thm1}) and set $r\to 1,\gamma \to \alpha ,m\to \mu -\nu ,\beta \to y,\alpha \to \beta$ and simplify;
 \begin{multline}
\int_0^{\infty } e^{-x \alpha } x^{\mu -\nu } \sqrt{y} J_{\mu }(x \beta ) J_{\nu }(x y) \, dx\\
=\sum _{j=0}^{\infty } \frac{(-1)^j 2^{\mu -\nu } y^{\frac{1}{2}+\nu } \alpha ^{-1-2
   j-2 \mu } \beta ^{2 j+\mu } \Gamma \left(\frac{1}{2} (1+2 j+2 \mu )\right) }{\sqrt{\pi } \Gamma (1+j)
   \Gamma (1+\nu )}\\ \times
\, _2F_1\left(\frac{1}{2}+j+\mu ,1+j+\mu ;1+\nu ;-\frac{y^2}{\alpha ^2}\right)\\
=\frac{\beta ^{\mu } y^{\nu +\frac{1}{2}} \Gamma \left(\mu +\frac{1}{2}\right)  }{2^{\nu -\mu } \pi  \Gamma \left(\nu +\frac{1}{2}\right)}\int_0^{\pi } \sin ^{2 \nu }(\phi ) \left((\alpha +i y \cos (\phi))^2+\beta ^2\right)^{-\mu -\frac{1}{2}} \, d\phi
\end{multline}
where $Re(\alpha)>|Im(\beta)|$.
\end{example}
\begin{example}
Derivation of equation (8.11.24) in \cite{erdt2} as a series representation. Errata. Here we use equation (\ref{eq:thm1}) and set $r\to 2,\gamma \to \alpha ,m\to \lambda +1,\alpha \to \beta ,\beta \to y$ and simplify;

\begin{multline}
\int_0^{\infty } e^{-\alpha  x^2} x^{1+\lambda } \sqrt{y} J_{\mu }(x \beta ) J_{\nu }(x y) \, dx\\
=\sum _{j=0}^{\infty } \frac{(-1)^j 2^{-1-2 j-\mu -\nu } y^{\frac{1}{2}+\nu }
   \alpha ^{-1-j-\frac{\lambda }{2}-\frac{\mu }{2}-\frac{\nu }{2}} \beta ^{2 j+\mu } \Gamma \left(\frac{1}{2} (2+2 j+\lambda +\mu +\nu )\right) }{\Gamma (1+j) \Gamma (1+j+\mu ) \Gamma (1+\nu )}\\ \times
\, _1F_1\left(1+j+\frac{\lambda
   }{2}+\frac{\mu }{2}+\frac{\nu }{2};1+\nu ;-\frac{y^2}{4 \alpha }\right)\\
=\sqrt{y} \sum _{m=0}^{\infty } \frac{\Gamma \left(m+\frac{\nu
   }{2}+\frac{\mu }{2}+\frac{\lambda }{2}\right) \left(-\frac{\beta ^2}{4 \alpha }\right)^m \, _2F_1\left(-m,-\mu -m;\nu +1;\left(\frac{y}{\beta }\right)^2\right)}{m! \Gamma (m+\mu
   +1)}
\end{multline}
where $Re(s)>0$.
\end{example}
\begin{example}
Derivation of equation (2.6) in \cite{eason} as a series representation for a definite integral. Errata. Here we use equation (\ref{eq:thm1}) and set $r\to 1,m\to \lambda ,\gamma \to c,\alpha \to a,\beta \to b,x\to t$ and simplify;
\begin{multline}
\int_0^{\infty } e^{-c t} t^{\lambda } J_{\mu }(a t) J_{\nu }(b t) \, dt\\
=\sum _{h=0}^{\infty } \frac{(-1)^h
   2^{-2 h-\mu -\nu } a^{\mu } b^{2 h+\nu } c^{-1-2 h-\lambda -\mu -\nu } \Gamma (1+2 h+\lambda +\mu +\nu ) }{\Gamma (1+h) \Gamma (1+\mu ) \Gamma (1+h+\nu )}\\ \times
\,_2F_1\left(\frac{1}{2} (1+2 h+\lambda +\mu +\nu ),\frac{1}{2} (2+2 h+\lambda +\mu +\nu );1+\mu
   ;-\frac{a^2}{c^2}\right)\\
\neq\frac{\Gamma (\mu -\nu +\lambda +1)
  }{2^{\mu -\nu } \pi  c^{\mu
   -\nu +\lambda +1} \Gamma (\mu -\nu +1)} \int_0^{\pi } (a-b \exp (-i \theta ))^{\mu -\nu } \exp (-i \nu  \theta )\\
 \, _2F_1\left(\frac{\mu }{2}-\frac{\nu}{2}+\frac{\lambda }{2}+\frac{1}{2},\frac{\mu }{2}-\frac{\nu }{2}+\frac{\lambda }{2}+1;\mu -\nu
   +1;-\left(\frac{\sqrt{a^2+b^2-2 a b \cos (\theta )}}{c}\right)^2\right) \, d\theta 
\end{multline}
where $Re(\mu+\nu+c)>0$.
\end{example}
\section{Table (3.12.15) in Prudnikov et al. (1992a) Laplace transform involving the product of two Bessel functions.}
\begin{example}
Derivation of entry (3.12.15.2) in \cite{prud4} series representation. Here we use equation (\ref{eq:thm1}) and set $m\to 0,r\to 1,\gamma \to p,\alpha \to a,\beta \to a$ and simplify;
\begin{multline}
\int_0^{\infty } e^{-p x} J_{\mu }(a x) J_{\nu }(a x) \, dx
=\sum _{j=0}^{\infty } \frac{(-1)^j 2^{-2 j-\mu
   -\nu } a^{2 j+\mu +\nu } p^{-1-2 j-\mu -\nu } \Gamma (1+2 j+\mu +\nu ) }{\Gamma (1+j) \Gamma(1+j+\mu ) \Gamma (1+\nu )}\\ \times
\, _2F_1\left(\frac{1}{2}+j+\frac{\mu}{2}+\frac{\nu }{2},1+j+\frac{\mu }{2}+\frac{\nu }{2};1+\nu ;-\frac{a^2}{p^2}\right)\\
=\frac{\left(\frac{a}{2}\right)^{\mu +\nu } p^{-\mu -\nu -1} \Gamma (\mu +\nu +1) }{\Gamma (\mu +1) \Gamma (\nu +1)}\\ \times
\,
   _4F_3\left(\frac{1}{2} (\mu +\nu +1),\frac{1}{2} (\mu +\nu +1),\frac{\mu +\nu }{2}+1,\frac{\mu +\nu }{2}+1;\right. \\ \left.
\mu+\nu +1,\mu +1,\nu +1;-\left(\frac{2 a}{p}\right)^2\right)
\end{multline}
where $Re(\mu+\nu)>-1, Re(p)>2|Im(a)|$.
\end{example}
\section{Table of Kausel entries in terms of infinite series of the hypergeometric function}
\begin{example}
Derivation of equation ENS-4.3 in \cite{kausel} series representation. Here we use equation (\ref{eq:thm1}) and set $m\to 1,\mu \to 0,\nu \to 0,r\to 1,\gamma \to s,\alpha \to a,\beta \to b$ and simplify;
\begin{multline}
\int_0^{\infty } e^{-s x} x J_0(a x) J_0(b x) \, dx\\
=\sum _{j=0}^{\infty } \frac{\left(-\frac{1}{4}\right)^j a^{2 j} s^{-2-2 j} \Gamma (2 (1+j)) \,
   _2F_1\left(1+j,\frac{3}{2}+j;1;-\frac{b^2}{s^2}\right)}{\Gamma (1+j)^2}\\
=\frac{\left(s \left(\frac{2 \sqrt{a b}}{\sqrt{(a+b)^2+s^2}}\right)^3\right) E\left(\left(\frac{2 \sqrt{a
   b}}{\sqrt{(a+b)^2+s^2}}\right)^2\right)}{4 \pi  \left(1-\left(\frac{2 \sqrt{a b}}{\sqrt{(a+b)^2+s^2}}\right)^2\right) a b \sqrt{a b}}
\end{multline}
where $Re(s)>0$.
\end{example}
\begin{example}
Derivation of equation ENS-4.8 in \cite{kausel} series representation. Here we use equation (\ref{eq:thm1}) and set $m\to 1,\mu \to 1,\nu \to 0,r\to 1,\gamma \to s,\alpha \to a,\beta \to b$ and simplify;
\begin{multline}
\int_0^{\infty } e^{-s x} x J_0(b x) J_1(a x) \, dx\\
=\sum _{j=0}^{\infty } \frac{(-1)^j 2^{-1-2 j} a^{1+2 j} s^{-3-2 j} \Gamma (3+2 j) \,
   _2F_1\left(\frac{3}{2}+j,2+j;1;-\frac{b^2}{s^2}\right)}{\Gamma (1+j) \Gamma (2+j)}\\
=\frac{\frac{\left(a^2-b^2-s^2\right) E\left(\frac{4 a
   b}{(a+b)^2+s^2}\right)}{(a-b)^2+s^2}+K\left(\frac{4 a b}{(a+b)^2+s^2}\right)}{a \pi  \sqrt{(a+b)^2+s^2}}
\end{multline}
where $Re(s)>0$.
\end{example}
\begin{example}
Derivation of equation ENS-4.4 in \cite{kausel} series representation. Here we use equation (\ref{eq:thm1}) and set $m\to 1,\mu \to 1,\nu \to 1,r\to 1,\gamma \to s,\alpha \to a,\beta \to b$ and simplify;
\begin{multline}
\int_0^{\infty } e^{-s x} x J_1(a x) J_1(b x) \, dx\\
=\sum _{j=0}^{\infty } \frac{2 (-1)^j a^{1+2 j} b s^{-4-2 j} \Gamma \left(\frac{1}{2} (5+2 j)\right) \,
   _2F_1\left(2+j,\frac{5}{2}+j;2;-\frac{b^2}{s^2}\right)}{\sqrt{\pi } \Gamma (1+j)}\\
=\frac{s \left(\frac{\left(a^2+b^2+s^2\right) E\left(\frac{4 a
   b}{(a+b)^2+s^2}\right)}{(a-b)^2+s^2}-K\left(\frac{4 a b}{(a+b)^2+s^2}\right)\right)}{a b \pi  \sqrt{(a+b)^2+s^2}}
\end{multline}
where $Re(s)>0$.
\end{example}
\begin{example}
Derivation of equation ENS-4.1 in \cite{kausel} series representation. Here we use equation (\ref{eq:thm1}) and set $m\to 0,\mu \to 0,\nu \to 0,r\to 1,\gamma \to s,\alpha \to a,\beta \to b$ and simplify;
\begin{multline}
\int_0^{\infty } e^{-s x} J_0(a x) J_0(b x) \, dx\\
=\sum _{j=0}^{\infty } \frac{\left(-\frac{1}{4}\right)^j a^{2 j} s^{-1-2 j} \Gamma (1+2 j) \,
   _2F_1\left(\frac{1}{2}+j,1+j;1;-\frac{b^2}{s^2}\right)}{\Gamma (1+j)^2}\\
=\frac{2 K\left(\frac{4 a b}{(a+b)^2+s^2}\right)}{\pi  \sqrt{(a+b)^2+s^2}}
\end{multline}
where $Re(s)>0$.
\end{example}
\begin{example}
Derivation of equation ENS-4.7 in \cite{kausel} series representation. Here we use equation (\ref{eq:thm1}) and set $m\to 0,\mu \to 1,\nu \to 0,r\to 1,\gamma \to s,\alpha \to a,\beta \to b$ and simplify;
\begin{multline}
\int_0^{\infty } e^{-s x} J_0(b x) J_1(a x) \, dx\\
=\sum _{j=0}^{\infty } \frac{(-1)^j 2^{-1-2 j} a^{1+2 j} s^{-2-2 j} \Gamma (2 (1+j)) \,
   _2F_1\left(1+j,\frac{3}{2}+j;1;-\frac{b^2}{s^2}\right)}{\Gamma (1+j) \Gamma (2+j)}\\
=-\frac{s K\left(\frac{4 a b}{(a+b)^2+s^2}\right)}{a \pi  \sqrt{(a+b)^2+s^2}}+\frac{\theta
   (a-b)}{a}-\frac{s | a-b|  \Pi \left(\frac{4 a b}{(a+b)^2}|\frac{4 a b}{(a+b)^2+s^2}\right) \text{sgn}(a-b)}{a (a+b) \pi  \sqrt{(a+b)^2+s^2}}
\end{multline}
where $Re(s)>0$.
\end{example}
\begin{example}
Derivation of equation ENS-4.2 in \cite{kausel} series representation. Here we use equation (\ref{eq:thm1}) and set $m\to 0,\mu \to 1,\nu \to 1,r\to 1,\gamma \to s,\alpha \to a,\beta \to b$ and simplify;
\begin{multline}
\int_0^{\infty } e^{-s x} J_1(a x) J_1(b x) \, dx\\
=\sum _{j=0}^{\infty } \frac{(-1)^j a^{1+2 j} b s^{-3-2 j} \Gamma \left(\frac{1}{2} (3+2 j)\right) \,
   _2F_1\left(\frac{3}{2}+j,2+j;2;-\frac{b^2}{s^2}\right)}{\sqrt{\pi } \Gamma (1+j)}\\
=\frac{\sqrt{(a+b)^2+s^2} \left(-E\left(\frac{4 a b}{(a+b)^2+s^2}\right)+\left(1-\frac{2 a
   b}{(a+b)^2+s^2}\right) K\left(\frac{4 a b}{(a+b)^2+s^2}\right)\right)}{a b \pi }
\end{multline}
where $Re(s)>0$.
\end{example}
\begin{example}
Derivation of equation ENS-4.6 in \cite{kausel} series representation. Here we use equation (\ref{eq:thm1}) and set $m\to -1,\mu \to 1,\nu \to 0,r\to 1,\gamma \to s,\alpha \to a,\beta \to b$ and simplify;
\begin{multline}
\int_0^{\infty } \frac{e^{-s x} J_0(b x) J_1(a x)}{x} \, dx\\
=\sum _{j=0}^{\infty } \frac{(-1)^j 2^{-1-2 j} a^{1+2 j} s^{-1-2 j} \Gamma (1+2 j) \,
   _2F_1\left(\frac{1}{2}+j,1+j;1;-\frac{b^2}{s^2}\right)}{\Gamma (1+j) \Gamma (2+j)}\\
=\frac{\sqrt{(a+b)^2+s^2} E\left(\frac{4 a b}{(a+b)^2+s^2}\right)+\frac{\left(a^2-b^2\right)
   K\left(\frac{4 a b}{(a+b)^2+s^2}\right)}{\sqrt{(a+b)^2+s^2}}}{a \pi }\\
-\frac{s \theta (a-b)}{a}+\frac{s^2 | a-b|  \Pi \left(\frac{4 a b}{(a+b)^2}|\frac{4 a b}{(a+b)^2+s^2}\right)
   \text{sgn}(a-b)}{a (a+b) \pi  \sqrt{(a+b)^2+s^2}}
\end{multline}
where $Re(s)>0$.
\end{example}
\begin{example}
Derivation of equation ENS-4.9 in \cite{kausel} series representation. Here we use equation (\ref{eq:thm1}) and set $m\to -1,\mu \to 1,\nu \to 1,r\to 1,\gamma \to s,\alpha \to a,\beta \to b$ and simplify;
\begin{multline}
\int_0^{\infty } \frac{e^{-s x} J_1(a x) J_1(b x)}{x} \, dx\\
=\sum _{j=0}^{\infty } \frac{(-1)^j a^{1+2 j} b s^{-2-2 j} \Gamma \left(\frac{1}{2} (3+2 j)\right) \,
   _2F_1\left(1+j,\frac{3}{2}+j;2;-\frac{b^2}{s^2}\right)}{2 \sqrt{\pi } \Gamma (2+j)}\\
=\frac{1}{4 a b}\left(a^2+b^2+\frac{2 s \sqrt{(a+b)^2+s^2} E\left(\frac{4 a b}{(a+b)^2+s^2}\right)}{\pi }-\frac{2
   s \left(2 a^2+2 b^2+s^2\right) K\left(\frac{4 a b}{(a+b)^2+s^2}\right)}{\pi  \sqrt{(a+b)^2+s^2}}\right. \\ \left.+\left(a^2-b^2\right) \left(-1+\frac{2 s | a-b|  \Pi \left(\frac{4 a
   b}{(a+b)^2}|\frac{4 a b}{(a+b)^2+s^2}\right)}{(a+b) \pi  \sqrt{(a+b)^2+s^2}}\right) \text{sgn}(a-b)\right)
\end{multline}
where $Re(s)>0$.
\end{example}
\begin{example}
Mellin transform involving polynomials, exponential and two Bessel functions. Here we use equation (\ref{eq:thm}) and set $k\to 0,a\to 1,c\to -z$ and simplify using equation [Wolfram,\href{http://functions.wolfram.com/06.07.03.0005.01}{01}];
\begin{multline}\label{eq:mellin_poly}
\int_0^{\infty } e^{-x^r \gamma } x^m \left(1-x^s z\right)^{\delta }
   J_{\mu }(x \alpha ) J_{\nu }(x \beta ) \, dx\\
=\sum _{h=0}^{\infty } \sum
   _{j=0}^{\infty } \sum _{p=0}^{\infty } \frac{(-1)^{h+j} 2^{-2 h-2 j-\mu -\nu
   } (-z)^p \alpha ^{2 j+\mu } \beta ^{2 h+\nu } \gamma ^{-\frac{1+2 h+2 j+m+p
   s+\mu +\nu }{r}}}{r h! j! \Gamma (1+j+\mu ) \Gamma (1+h+\nu )}\\ \times
    \binom{\delta }{p} \Gamma \left(\frac{1+2 h+2 j+m+p s+\mu
   +\nu }{r}\right)
\end{multline}
where $Re(\mu+\nu+m)>0, Re(\gamma)>2|Im(a)|$.
\end{example}
\section{Tables in Prudnikov et al. (1992)}
In this section we derive generalized forms for entries (2.12.38), (2.12.39), (2.12.40.1-16), (2.12.41) in \cite{prud2}.
\begin{example}
Derivation of equation (2.12.38.1) in \cite{prud2}. In this example we use equation (\ref{eq:mellin_poly}) and set $r\to 1,\gamma \to q,m\to 1,\delta \to 0,z\to -1,s\to 0,\alpha \to b,\beta \to c$ and simplify;
\begin{multline}
\int_0^{\infty } e^{-q x} x J_{\mu }(b x) J_{\nu }(c x) \, dx\\
=\sum _{j=0}^{\infty } \frac{(-1)^j 2^{-2 j-\mu -\nu } b^{2 j+\mu } c^{\nu } q^{-2-2 j-\mu -\nu } \Gamma (2+2 j+\mu +\nu ) }{\Gamma (1+j) \Gamma (1+j+\mu ) \Gamma (1+\nu )}\\ \times
\, _2F_1\left(1+j+\frac{\mu }{2}+\frac{\nu }{2},\frac{3}{2}+j+\frac{\mu }{2}+\frac{\nu }{2};1+\nu;-\frac{c^2}{q^2}\right)
\end{multline}
where $Re(p)>|Im(b)|+|Im(c)|,Re(\mu+\nu)>-1$.
\end{example}
\begin{example}
Derivation of equation (2.12.38.2) in \cite{prud2}. In this example we use equation (\ref{eq:mellin_poly}) and set $r\to 1,\gamma \to q,m\to \alpha -1,\delta \to 0,z\to -1,s\to 0,\alpha \to b,\beta \to c$ and simplify;
\begin{multline}
\int_0^{\infty } e^{-q x} x^{-1+\alpha } J_{\mu }(b x) J_{\nu }(c x) \, dx\\
=\sum _{j=0}^{\infty } \frac{(-1)^j 2^{-2 j-\mu -\nu } b^{2 j+\mu } c^{\nu } q^{-2 j-\alpha -\mu -\nu } \Gamma (2 j+\alpha +\mu +\nu ) }{\Gamma (1+j) \Gamma (1+j+\mu ) \Gamma (1+\nu )}\\ \times
\, _2F_1\left(j+\frac{\alpha }{2}+\frac{\mu }{2}+\frac{\nu
   }{2},\frac{1}{2}+j+\frac{\alpha }{2}+\frac{\mu }{2}+\frac{\nu }{2};1+\nu ;-\frac{c^2}{q^2}\right)
\end{multline}
where $Re(p)>|Im(b)|+|Im(c)|,Re(\mu+\nu)>-1$.
\end{example}
\begin{example}
Derivation of equation (2.12.39.5) in \cite{prud2} pp. 223. Errata. In this example we use equation (\ref{eq:mellin_poly}) and set $r\to 2,\gamma \to p,\delta \to 0,s\to 0,z\to -1,m\to 0,\mu \to \nu ,\alpha \to e^{-\frac{i \pi }{4}} c,\beta \to e^{\frac{i \pi }{4}} c$ and simplify;
\begin{multline}
\int_0^{\infty } e^{-p x^2} J_{\nu }\left(c e^{-\frac{1}{4} (i \pi )} x\right) J_{\nu }\left(c e^{\frac{i \pi }{4}} x\right) \, dx\\
=\sum _{h=0}^{\infty } \frac{(-1)^h 2^{-1-2 h-2 \nu } \left(\sqrt[4]{-1} c\right)^{\nu } \left(c e^{-\frac{1}{4} (i \pi )}\right)^{\nu } \left(c
   e^{\frac{i \pi }{4}}\right)^{2 h} p^{-\frac{1}{2}-h-\nu } \Gamma \left(\frac{1}{2} (1+2 h+2 \nu )\right) }{\Gamma (1+h) \Gamma (1+\nu ) \Gamma (1+h+\nu )}\\ \times
\, _1F_1\left(\frac{1}{2}+h+\nu ;1+\nu ;\frac{i c^2}{4 p}\right)\\
=\frac{2^{-\nu -1} \Gamma \left(\nu +\frac{1}{2}\right)
   M_{\frac{1}{4},\frac{\nu }{2}}\left(\frac{c^2}{2 p}\right) M_{-\frac{1}{4},\frac{\nu }{2}}\left(\frac{c^2}{2 p}\right)}{\left(c^2 \sqrt{p}\right) (\Gamma (\nu +1) \Gamma (\nu +1))}
\end{multline}
where $Re(p)>|Im(b)|+|Im(c)|,Re(\mu+\nu)>-1$.
\end{example}
\begin{example}
 Derivation of equation (2.12.39.1) in \cite{prud2}. In this example we use equation (\ref{eq:mellin_poly}) and set $r\to 2,\gamma \to p,\delta \to 0,s\to 0,z\to -1,m\to 0,\mu \to 0,\nu \to 1,\alpha \to b,\beta \to c$ and simplify;
 \begin{multline}
\int_0^{\infty } e^{-p x^2} J_0(b x) J_1(c x) \, dx\\
=\sum _{j=0}^{\infty } \frac{(-1)^j 2^{-2-2 j} b^{2 j} c p^{-1-j} \, _1F_1\left(1+j;2;-\frac{c^2}{4 p}\right)}{\Gamma (1+j)}\\
=\frac{1}{c}\int_0^{\frac{c^2}{4 p}} e^{-\frac{b^2}{4 p}-t} I_0\left(\sqrt{\frac{b^2 t}{p}}\right) \,
   dt
\end{multline}
where $Re(p)>|Im(b)|+|Im(c)|,Re(\mu+\nu)>-1$.
\end{example}
\begin{example}
Derivation of equation (2.12.39.2) in \cite{prud2}. Errata. In this example we use equation (\ref{eq:mellin_poly}) and set $r\to 2,\gamma \to p,\delta \to 0,s\to 0,z\to -1,m\to \alpha -1,\alpha \to b,\beta \to c$ and simplify;
\begin{multline}
\int_0^{\infty } e^{-p x^2} x^{-1+\alpha } J_{\mu }(b x) J_{\nu }(c x) \, dx\\
=\sum _{j=0}^{\infty } \frac{(-1)^j 2^{-1-2 j-\mu -\nu } b^{2 j+\mu } c^{\nu } p^{-j-\frac{\alpha }{2}-\frac{\mu }{2}-\frac{\nu }{2}} \Gamma \left(\frac{1}{2} (2 j+\alpha +\mu +\nu )\right) }{\Gamma (1+j) \Gamma (1+j+\mu ) \Gamma (1+\nu )}\\ \times
\,
   _1F_1\left(j+\frac{\alpha }{2}+\frac{\mu }{2}+\frac{\nu }{2};1+\nu ;-\frac{c^2}{4 p}\right)\\
=\frac{\left(b^{\mu } c^{\nu } p^{-\frac{1}{2} (\alpha +\mu +\nu )}\right) }{2^{\mu +\nu +1} \Gamma (\nu +1)}\sum _{k=0}^{\infty } \frac{\Gamma \left(k+\left(\alpha +\mu +\frac{\nu
   }{2}\right)\right) \left(-\frac{b^2}{4 p}\right)^k \, _2F_1\left(-k,-\mu -k;\nu +1;\left(\frac{c}{b}\right)^2\right)}{\Gamma (\mu +k+1) k!}
\end{multline}
where $Re(p)>|Im(b)|+|Im(c)|,Re(\mu+\nu)>-1$.
\end{example}
\section{Finite Hankel transforms}
According to (Poularikas, 2018) section (11.1.4) in \cite{poularikas}, the finite Hankel transform is used in determining the temperature distribution in a long circular cylinder, finding the unsteady viscous flow in a rotating long circular cylinder, calculating the vibrations of a circular membrane, and evaluating the temperature distribution of cooling of a circular cylinder.
According to (Patra, ) in section 7 in \cite{patra}, finite Hankel transform arise in work involving the solution of certain special types of boundary value problems. Sneddon (Phil. Mag., 37, 1946) \cite{sneddon_f}, was the first to introduce this transform. Later applications of this transform was found in the works of several other authors in discussing solutions of axisymmetric physical problems in long circular cylinders and membranes. In this section we look at alternative forms in the derivations of known formulae involving finite Hankel transforms.
\begin{example}
Derivation of equation (7.4) in \cite{patra} series representation using Wolfram Mathematica. In this example we use equation (\ref{eq:thm}) and set $k\to 0,a\to 1,r\to 0,\gamma \to 1,\delta \to 0,s\to 0,c\to 1,m\to 1,x\to r,\mu \to \nu ,b\to a$ and simplify;
\begin{multline}
\sum _{j=0}^{\infty } \frac{(-1)^j 2^{-1-2 j-2 \nu } a^{1+2 j+2 \nu } \alpha ^{2 j+\nu } \beta ^{\nu } \, _1F_2\left(1+j+\nu ;1+\nu ,2+j+\nu ;-\frac{1}{4} a^2 \beta ^2\right)}{(1+j+\nu ) \Gamma (1+j) \Gamma (1+\nu ) \Gamma (1+j+\nu )}\\
=\frac{\beta  J_{-1+\nu }(a \beta ) J_{\nu }(a
   \alpha )-\alpha  J_{-1+\nu }(a \alpha ) J_{\nu }(a \beta )}{\alpha ^2-\beta ^2}
\end{multline}
where $Re(\nu)>-1$.
\end{example}
\begin{example}
Derivation of equation (1.8.3.1) in \cite{prud2}, series representation. Table (1.8.3) can also be derived using equation (\ref{eq:thm}). In this example we use equation (\ref{eq:thm}) and set $k\to 0,a\to 1,r\to 0,\gamma \to 1,\delta \to 0,s\to 0,c\to 1,\alpha \to 1,\beta \to 1,m\to \lambda ,x\to r,b\to x$ and simplify;
\begin{multline}
\int_0^x r^{\lambda } J_{\mu }(r) J_{\nu }(r) \, dr\\
=\sum _{j=0}^{\infty } \frac{(-1)^j 2^{-2 j-\mu -\nu } x^{1+2 j+\lambda +\mu +\nu } }{(1+2 j+\lambda +\mu +\nu ) \Gamma (1+j) \Gamma (1+j+\mu ) \Gamma (1+\nu )}\\ \times
\, _1F_2\left(\frac{1}{2}+j+\frac{\lambda }{2}+\frac{\mu}{2}+\frac{\nu }{2};\frac{3}{2}+j+\frac{\lambda }{2}+\frac{\mu }{2}+\frac{\nu }{2},1+\nu ;-\frac{x^2}{4}\right)\\
=\frac{2^{-\mu-\nu } x^{1+\lambda +\mu +\nu } }{(1+\lambda +\mu +\nu ) \Gamma (1+\mu ) \Gamma (1+\nu )}\, _3F_4\left(\frac{1}{2}+\frac{\mu }{2}+\frac{\nu }{2},1+\frac{\mu }{2}+\frac{\nu }{2},\frac{1}{2}+\frac{\lambda }{2}+\frac{\mu }{2}+\frac{\nu }{2};\right. \\ \left.
1+\mu,\frac{3}{2}+\frac{\lambda }{2}+\frac{\mu }{2}+\frac{\nu }{2},1+\nu ,1+\mu +\nu ;-x^2\right)
\end{multline}
where $Re(x)>0$.
\end{example}
\begin{example}
Derivation of equation (3.3) in \cite{bailey} series representation of the product of Legendre polynomials. In this example we use equation (\ref{eq:mellin_poly}) and set $r\to 1,\delta \to 0,s\to 0,z\to -1,x\to t$ then $\alpha \to \sin (\phi ),\beta \to \sin (\Phi ),\gamma \to \cos (\Phi ) \cos (\phi ),m\to -\frac{1}{2}$ and simplify;
\begin{multline}
\int_0^{\infty } \frac{e^{-t \cos (\phi ) \cos (\Phi )} J_{\mu }(t \sin (\phi )) J_{\nu }(t \sin (\Phi
   ))}{\sqrt{t}} \, dt\\
=\sum _{j=0}^{\infty } \frac{(-1)^j 2^{-2 j-\mu -\nu } (\cos (\phi ) \cos (\Phi
   ))^{-\frac{1}{2}-2 j-\mu -\nu } \Gamma \left(\frac{1}{2} (1+4 j+2 \mu +2 \nu )\right) 
}{\Gamma (1+j) \Gamma (1+j+\mu ) \Gamma(1+\nu )}\\ \times
\,_2F_1\left(\frac{1}{4}+j+\frac{\mu }{2}+\frac{\nu }{2},\frac{3}{4}+j+\frac{\mu }{2}+\frac{\nu }{2};1+\nu ;-\sec^2(\phi ) \tan ^2(\Phi )\right) \sin ^{2 j+\mu }(\phi ) \sin ^{\nu }(\Phi )\\
=\Gamma \left(\mu +\nu +\frac{1}{2}\right) P_{\nu -\frac{1}{2}}^{-\mu }(\cos (\phi )) P_{\mu
   -\frac{1}{2}}^{-\nu }(\cos (\Phi ))
\end{multline}
where $|Im(\phi)|<\pi/2,|Im(\Phi)|<\pi/2, Re(\mu+\nu)>-1/2$
\end{example}
\begin{example}
Derivation of equation (2.12.29) in \cite{erdt1} series representation for the ratio of gamma functions.  In this example we use equation (\ref{eq:mellin_poly}) and set $r\to 1,\delta \to 0,s\to 0,z\to -1$ next form a second equation by replacing $\gamma\to-\gamma$ and take their difference and simplify;
\begin{multline}
\int_0^{\infty } x^{-2-\mu +\nu } J_{\mu }(a x) J_{\nu }(b x) \sin (x y) \, dx\\
=-\sum _{j=0}^{\infty }
   \frac{(-1)^j 2^{-1-2 j-\mu -\nu } a^{2 j+\mu } b^{\nu } \left((-i y)^{2 j+2 \nu }+(i y)^{2 j+2 \nu }\right) (-i
   y)^{-2 j-2 \nu } (i y)^{-2 j-2 \nu } y }{\Gamma (1+j) \Gamma (1+j+\mu ) \Gamma (1+\nu )}\\ 
\times \Gamma (-1+2 j+2 \nu ) \, _2F_1\left(-\frac{1}{2}+j+\nu ,j+\nu ;1+\nu
   ;\frac{b^2}{y^2}\right)\\
\neq\frac{2^{\nu -\mu -1} a^{\mu } b^{-\nu }
   (\Gamma (\nu ) y)}{\Gamma (\mu +1)}
\end{multline}
where $Re(\mu+\nu)>0,Re(y)>1$.
\end{example}
\begin{example}
Derivation of equation (4.14.39) in \cite{erdt1}, series representation for the modified Bessel function of the first kind.  In this example we use equation (\ref{eq:mellin_poly}) and set $\delta \to 0,s\to 0,z\to -1,x\to \sqrt{t},\gamma \to p$ followed by $r\to 2,m\to 1,\mu \to \nu$ and simplify;
 \begin{multline}
\int_0^{\infty } e^{-p t} J_{\nu }\left(\sqrt{t} \alpha \right) J_{\nu }\left(\sqrt{t} \beta \right) \, dt\\
=\sum
   _{j=0}^{\infty } \frac{(-1)^j 2^{-2 j-2 \nu } p^{-1-j-\nu } \alpha ^{2 j+\nu } \beta ^{\nu } \, _1F_1\left(1+j+\nu
   ;1+\nu ;-\frac{\beta ^2}{4 p}\right)}{\Gamma (1+j) \Gamma (1+\nu )}\\
=\frac{\exp \left(-\frac{\alpha ^2+\beta ^2}{4
   p}\right) I_{\nu }\left(\frac{\alpha  \beta }{2 p}\right)}{p}
\end{multline}
where $Re(p)>0$.
\end{example}
\begin{example}
Derivation of equation (4.14.42) in \cite{erdt1}, series representation for the generalized hypergeometric function. Errata. In this example we use equation (\ref{eq:mellin_poly}) and set $\delta \to 0,s\to 0,z\to -1,x\to \sqrt{t},\gamma \to p,\alpha \to 2 \sqrt{\alpha },\beta \to 2 \sqrt{\alpha
   },\mu \to 2 \mu ,\nu \to 2 \nu ,r\to 2,m\to 2 \lambda -2$ and $\lambda \to \lambda +\frac{1}{2}$ then simplify;
   \begin{multline}
\int_0^{\infty } e^{-p t} t^{-1+\lambda } J_{2 \mu }\left(2 \sqrt{t} \sqrt{\alpha }\right) J_{2 \nu }\left(2
   \sqrt{t} \sqrt{\alpha }\right) \, dt\\
=\sum _{j=0}^{\infty } \frac{(-1)^j p^{-j-\lambda -\mu -\nu } \alpha ^{j+\mu
   +\nu } \Gamma \left(\frac{1}{2} \left(-1+2 j+2 \left(\frac{1}{2}+\lambda \right)+2 \mu +2 \nu \right)\right) }{\Gamma (1+j) \Gamma (1+j+2 \mu ) \Gamma (1+2\nu )}\\ \times
\,_1F_1\left(j+\lambda +\mu +\nu ;1+2 \nu ;-\frac{\alpha }{p}\right)\\
=\frac{\left(\Gamma (\lambda +\mu +\nu ) \alpha ^{\mu +\nu }\right) }{\Gamma (2 \mu +1) \Gamma
   (2 \nu +1) p^{\lambda +\mu +\nu }}\\ \times
\, _3F_3\left(\mu +\nu +\frac{1}{2},\mu+\nu +1,\lambda +\mu +\nu ;2 \mu +1,2 \nu +1,2 \mu +2 \nu +1;-\frac{4 \alpha }{p}\right)
\end{multline}
where $Re(p)>0$.
\end{example}
\section{Table of Okui entries involving Bessel functions in terms of infinite series of the hypergeometric function}
\begin{example}
Generalized Okui integral form for Table (2.2). In this example we use equation (\ref{eq:mellin_poly}) and set $\delta \to 0,s\to 0,z\to -1,x\to t^2,\gamma \to p,r\to 1,m\to -\frac{1}{2}$ and simplify;
\begin{multline}\label{eq:okui}
\int_0^{\infty } e^{-p x^2} J_{\mu }\left(x^2 \alpha \right) J_{\nu }\left(x^2 \beta \right) \, dx\\
=\sum
   _{j=0}^{\infty } \frac{(-1)^j 2^{-1-2 j-\mu -\nu } p^{-\frac{1}{2}-2 j-\mu -\nu } \alpha ^{2 j+\mu } \beta ^{\nu }
   \Gamma \left(\frac{1}{2} (1+4 j+2 \mu +2 \nu )\right) }{\Gamma (1+j) \Gamma (1+j+\mu
   ) \Gamma (1+\nu )}\\ \times
\, _2F_1\left(\frac{1}{4}+j+\frac{\mu }{2}+\frac{\nu
   }{2},\frac{3}{4}+j+\frac{\mu }{2}+\frac{\nu }{2};1+\nu ;-\frac{\beta ^2}{p^2}\right)
\end{multline}
where $Re(p)>0$.
\end{example}
\begin{example}
Derivation of equation (2.2.1) in \cite{okui_1974}. In this example we use equation (\ref{eq:okui}) and set $\mu \to 0,\nu \to 0,\alpha \to a,\beta \to b$ and simplify;
\begin{multline}
\int_0^{\infty } e^{-p x^2} J_0\left(a x^2\right) J_0\left(b x^2\right) \, dx\\
=\sum _{j=0}^{\infty }
   \frac{(-1)^j 2^{-1-2 j} a^{2 j} p^{-\frac{1}{2}-2 j} \Gamma \left(\frac{1}{2} (1+4 j)\right) \,
   _2F_1\left(\frac{1}{4}+j,\frac{3}{4}+j;1;-\frac{b^2}{p^2}\right)}{\Gamma (1+j)^2}\\
=\frac{1}{\pi ^{3/2}}2 \sqrt[4]{2}
   \sqrt[4]{\frac{1}{a^2+b^2+p^2+\sqrt{-4 a^2 b^2+\left(a^2+b^2+p^2\right)^2}}}\\
 K\left(\frac{1}{2}
   \left(1-\sqrt{1-\frac{2 a^2}{a^2+b^2+p^2+\sqrt{-4 a^2 b^2+\left(a^2+b^2+p^2\right)^2}}}\right)\right)\\
   K\left(\frac{1}{2} \left(1-\sqrt{1-\frac{2 b^2}{a^2+b^2+p^2+\sqrt{-4 a^2
   b^2+\left(a^2+b^2+p^2\right)^2}}}\right)\right)
\end{multline}
where $Re(p)>0$.
\end{example}
\begin{example}
Derivation of equation (2.3.1) in \cite{okui_1974} . In this example we use equation (\ref{eq:mellin_poly}) and set $\delta \to 0,s\to 0,z\to -1,\gamma \to 2 p,r\to 1,m\to 2,\mu \to 0,\nu \to 0,\alpha \to a,\beta \to a$ and simplify;
\begin{multline}
\int_0^{\infty } e^{-2 p x} x^2 J_0(a x){}^2 \, dx\\
=\sum _{j=0}^{\infty } \frac{(-1)^j 2^{-3-4 j} a^{2 j}
   p^{-3-2 j} \Gamma (3+2 j) \, _2F_1\left(\frac{3}{2}+j,2+j;1;-\frac{a^2}{4 p^2}\right)}{\Gamma
   (1+j)^2}\\
=\frac{\sqrt{\frac{a^2}{a^2+p^2}} \left(\left(3-\frac{2 a^2}{a^2+p^2}\right)
   E\left(\frac{a^2}{a^2+p^2}\right)-\left(1-\frac{a^2}{a^2+p^2}\right) K\left(\frac{a^2}{a^2+p^2}\right)\right)}{4 a
   p^2 \pi }
\end{multline}
where $Re(p)>0$.
\end{example}
\begin{example}
Derivation of equation (2.3.2) in \cite{okui_1974} . In this example we use equation (\ref{eq:mellin_poly}) and set $\delta \to 0,s\to 0,z\to -1,\gamma \to 2 p,r\to 1,m\to 2,\mu \to 1,\nu \to 1,\alpha \to a,\beta \to a$ and simplify;
\begin{multline}
\int_0^{\infty } e^{-2 p x} x^2 J_1(a x){}^2 \, dx\\
=\sum _{j=0}^{\infty } \frac{(-1)^j 2^{-7-4 j} a^{2+2 j}
   p^{-5-2 j} \Gamma (5+2 j) \, _2F_1\left(\frac{5}{2}+j,3+j;2;-\frac{a^2}{4 p^2}\right)}{\Gamma (1+j) \Gamma
   (2+j)}\\
=\frac{\sqrt{\frac{a^2}{a^2+p^2}} \left(-\left(\left(1-\frac{2 a^2}{a^2+p^2}\right)
   E\left(\frac{a^2}{a^2+p^2}\right)\right)+\left(1-\frac{a^2}{a^2+p^2}\right)
   K\left(\frac{a^2}{a^2+p^2}\right)\right)}{4 a p^2 \pi }
\end{multline}
where $Re(p)>0$.
\end{example}
\begin{example}
Derivation of equation (2.3.3) in \cite{okui_1974} . In this example we use equation (\ref{eq:mellin_poly}) and set $\delta \to 0,s\to 0,z\to -1,\gamma \to 2 p,r\to 1,m\to -2,\mu \to 1,\nu \to 1,\alpha \to a,\beta \to a$ and simplify;
\begin{multline}
\int_0^{\infty } \frac{e^{-2 p x} J_1(a x){}^2}{x^2} \, dx\\
=\sum _{j=0}^{\infty } \frac{(-1)^j 2^{-3-2 j} a^{2+2
   j} p^{-1-2 j} \Gamma \left(\frac{1}{2} (1+2 j)\right) \, _2F_1\left(\frac{1}{2}+j,1+j;2;-\frac{a^2}{4
   p^2}\right)}{\sqrt{\pi } \Gamma (2+j)}\\
=-p+\frac{4 a \left(-\left(\left(1-\frac{2 a^2}{a^2+p^2}\right)
   E\left(\frac{a^2}{a^2+p^2}\right)\right)+\left(1-\frac{a^2}{a^2+p^2}\right)
   K\left(\frac{a^2}{a^2+p^2}\right)\right)}{3 \left(\frac{a^2}{a^2+p^2}\right)^{3/2} \pi }
\end{multline}
where $Re(p)>0$.
\end{example}
\begin{example}
Derivation of equation (2.3.4) in \cite{okui_1974} . In this example we use equation (\ref{eq:mellin_poly}) and set $\delta \to 0,s\to 0,z\to -1,\gamma \to 2 p,r\to 1,m\to 1,\mu \to 2,\nu \to 2,\alpha \to a,\beta \to a$ and simplify;
\begin{multline}
\int_0^{\infty } e^{-2 p x} x J_2(a x){}^2 \, dx\\
=\sum _{j=0}^{\infty } \frac{(-1)^j 2^{-6-2 j} a^{4+2 j}
   p^{-6-2 j} \Gamma \left(\frac{1}{2} (7+2 j)\right) \, _2F_1\left(3+j,\frac{7}{2}+j;3;-\frac{a^2}{4
   p^2}\right)}{\sqrt{\pi } \Gamma (1+j)}\\
=\frac{\left(16+\frac{a^4}{\left(a^2+p^2\right)^2}-\frac{16
   a^2}{a^2+p^2}\right) E\left(\frac{a^2}{a^2+p^2}\right)-8 \left(1-\frac{a^2}{a^2+p^2}\right)
   \left(2-\frac{a^2}{a^2+p^2}\right) K\left(\frac{a^2}{a^2+p^2}\right)}{2 a p \left(\frac{a^2}{a^2+p^2}\right)^{3/2}
   \pi }
\end{multline}
where $Re(p)>0$.
\end{example}
\begin{example}
Derivation of equation (2.4.1) in \cite{okui_1974} . In this example we use equation (\ref{eq:mellin_poly}) and set $\delta \to 0,s\to 0,z\to -1,\gamma \to 2 p,r\to 1,m\to -1,\mu \to 1,\nu \to 2,\alpha \to a,\beta \to a$ and simplify;
\begin{multline}
\int_0^{\infty } \frac{e^{-2 p x} J_1(a x) J_2(a x)}{x} \, dx\\
=\sum _{j=0}^{\infty } \frac{(-1)^j 2^{-5-2 j}
   a^{3+2 j} p^{-3-2 j} \Gamma \left(\frac{1}{2} (3+2 j)\right) \, _2F_1\left(\frac{3}{2}+j,2+j;3;-\frac{a^2}{4
   p^2}\right)}{\sqrt{\pi } \Gamma (1+j)}\\
=-\frac{p}{a}+\frac{2 \left(-\left(\left(4-\frac{5 a^2}{a^2+p^2}\right)
   E\left(\frac{a^2}{a^2+p^2}\right)\right)+4 \left(1-\frac{a^2}{a^2+p^2}\right)
   K\left(\frac{a^2}{a^2+p^2}\right)\right)}{3 \left(\frac{a^2}{a^2+p^2}\right)^{3/2} \pi }
\end{multline}
where $Re(p)>0$.
\end{example}
\begin{example}
Derivation of equation (2.5.1) in \cite{okui_1974}.  In this example we use equation (\ref{eq:mellin_poly}) and set $\delta \to 0,s\to 0,z\to -1,\gamma \to p,r\to 1,m\to 0,\mu \to 0,\nu \to 0,\alpha \to a,\beta \to b$ and simplify;
\begin{multline}
\int_0^{\infty } e^{-p x} J_0(a x) J_0(b x) \, dx\\
=\sum _{j=0}^{\infty } \frac{(-1)^j a^{2 j} p^{-1-2 j} \Gamma
   \left(\frac{1}{2} (1+2 j)\right) \, _2F_1\left(\frac{1}{2}+j,1+j;1;-\frac{b^2}{p^2}\right)}{\sqrt{\pi } \Gamma
   (1+j)}\\
=\frac{2 \sqrt{\frac{a b}{(a+b)^2+p^2}} K\left(\frac{4 a b}{(a+b)^2+p^2}\right)}{\sqrt{a b} \pi }
\end{multline}
where $Re(p)>0$.
\end{example}
\begin{example}
Derivation of equation (2.5.1) in \cite{okui_1974}.  In this example we use equation (\ref{eq:mellin_poly}) and set $\delta \to 0,s\to 0,z\to -1,\gamma \to p,r\to 1,m\to 1,\mu \to 0,\nu \to 0,\alpha \to a,\beta \to b$ and simplify;
\begin{multline}
\int_0^{\infty } e^{-p x} x J_0(a x) J_0(b x) \, dx\\
=\sum _{j=0}^{\infty } \frac{2 (-1)^j a^{2 j} p^{-2-2 j}
   \Gamma \left(\frac{1}{2} (3+2 j)\right) \, _2F_1\left(1+j,\frac{3}{2}+j;1;-\frac{b^2}{p^2}\right)}{\sqrt{\pi }
   \Gamma (1+j)}\\
=\frac{2 p \left(\frac{a b}{(a+b)^2+p^2}\right)^{3/2} E\left(\frac{4 a
   b}{(a+b)^2+p^2}\right)}{\sqrt{a^3 b^3} \left(1-\frac{4 a b}{(a+b)^2+p^2}\right) \pi }
\end{multline}
where $Re(p)>0$.
\end{example}
\section{Table of Okui entries}
In this section we derive series representations for entries in Okui (1975) expressed in terms of elliptic functions.
\begin{example}
Derivation of generalized form for equations in Table (2.3) in \cite{okui_1975}. In this example we use equation (\ref{eq:mellin_poly}) and set $\delta \to 0,s\to 0,z\to -1,\gamma \to 2 i b,r\to 1$ next we form a second equation by replacing $b\to -b$ and simplify;
\begin{multline}\label{eq:okui_1975}
\int_0^{\infty } x^{-1+m} J_{\mu }(x \alpha ) J_{\nu }(x \beta ) \sin (2 b x) \, dx\\
=\sum _{j=0}^{\infty }
   \frac{i (-1)^j 2^{-1-4 j-m-2 \mu -2 \nu } \left((-i b)^{2 j+m+\mu +\nu }-(i b)^{2 j+m+\mu +\nu }\right) b^{-2 (2
   j+m+\mu +\nu )} \alpha ^{2 j+\mu } \beta ^{\nu } }{\Gamma (1+j) \Gamma (1+j+\mu ) \Gamma (1+\nu )}\\ \times
\Gamma (2 j+m+\mu +\nu ) \, _2F_1\left(j+\frac{m}{2}+\frac{\mu
   }{2}+\frac{\nu }{2},\frac{1}{2}+j+\frac{m}{2}+\frac{\mu }{2}+\frac{\nu }{2};1+\nu ;\frac{\beta ^2}{4
   b^2}\right)
\end{multline}
where $Re(b)\geq 1$ in order for the series to converge.
\end{example}
\begin{example}
Derivation of equation (2.3.1) in \cite{okui_1975}. In this example we use equation (\ref{eq:okui_1975}) and set $m\to -1,\mu \to 1,\nu \to 1,\alpha \to a,\beta \to a$ and simplify;
\begin{multline}
\int_0^{\infty } \frac{J_1(a x){}^2 \sin (2 b x)}{x^2} \, dx,\\
=\sum _{j=0}^{\infty } \frac{i (-1)^j 2^{-4-4 j}
   a^{2+2 j} \left((-i b)^{1+2 j}-(i b)^{1+2 j}\right) b^{-2 (1+2 j)} \Gamma (1+2 j) \,
   _2F_1\left(\frac{1}{2}+j,1+j;2;\frac{a^2}{4 b^2}\right)}{\Gamma (1+j) \Gamma (2+j)}\\
=b-\frac{2 (a+b)
   \left(\left(a^2+b^2\right) E\left(\frac{4 a b}{(a+b)^2}\right)-(a-b)^2 K\left(\frac{4 a
   b}{(a+b)^2}\right)\right)}{3 a^2 \pi }
\end{multline}
where $Re(b)>Re(a)$.
\end{example}
\begin{example}
Derivation of equation (2.3.2) in \cite{okui_1975}. In this example we use equation (\ref{eq:okui_1975}) and set $m\to -1,\mu \to 2,\nu \to 2,\alpha \to a,\beta \to a$ and simplify;
\begin{multline}
\int_0^{\infty } \frac{J_2(a x){}^2 \sin (2 b x)}{x^2} \, dx\\
=\sum _{j=0}^{\infty } \frac{i (-1)^j 2^{-9-4 j}
   a^{4+2 j} \left((-i b)^{3+2 j}-(i b)^{3+2 j}\right) b^{-2 (3+2 j)} \Gamma (3+2 j) \,
   _2F_1\left(\frac{3}{2}+j,2+j;3;\frac{a^2}{4 b^2}\right)}{\Gamma (1+j) \Gamma (3+j)}\\
=\frac{b}{2}-\frac{2 (a+b)
   \left(\left(a^4+9 a^2 b^2-4 b^4\right) E\left(\frac{4 a b}{(a+b)^2}\right)+(a-b)^2 \left(-a^2+4 b^2\right)
   K\left(\frac{4 a b}{(a+b)^2}\right)\right)}{15 a^4 \pi }
\end{multline}
where $Re(b)>Re(a)$.
\end{example}
\begin{example}
Derivation of equation (2.3.3) in \cite{okui_1975}. In this example we use equation (\ref{eq:okui_1975}) and set $m\to 0,\mu \to 1,\nu \to 2,\alpha \to a,\beta \to a$ and simplify;
\begin{multline}
\int_0^{\infty } \frac{J_1(a x) J_2(a x) \sin (2 b x)}{x} \, dx\\
=\sum _{j=0}^{\infty } \frac{i (-1)^j 2^{-8-4 j}
   a^{3+2 j} \left((-i b)^{3+2 j}-(i b)^{3+2 j}\right) b^{-2 (3+2 j)} \Gamma (3+2 j) \,
   _2F_1\left(\frac{3}{2}+j,2+j;3;\frac{a^2}{4 b^2}\right)}{\Gamma (1+j) \Gamma (2+j)}\\
=\frac{b}{a}+\frac{-\left((a+b)
   \left(a^2+4 b^2\right) E\left(\frac{4 a b}{(a+b)^2}\right)\right)+(-a+b) \left(-a^2+4 b^2\right) K\left(\frac{4 a
   b}{(a+b)^2}\right)}{3 a^3 \pi }
\end{multline}
where $Re(b)>Re(a)$.
\end{example}
\begin{example}
Derivation of equation (2.3.4) in \cite{okui_1975}. Errata. In this example we use equation (\ref{eq:okui_1975}) and set $m\to 1,\mu \to 1,\nu \to 3,\alpha \to a,\beta \to a$ and simplify;
\begin{multline}
\int_0^{\infty } J_1(a x) J_3(a x) \sin (2 b x) \, dx\\
=\sum _{j=0}^{\infty } \frac{i (-1)^j 2^{-11-4 j} a^{4+2
   j} \left((-i b)^{5+2 j}-(i b)^{5+2 j}\right) b^{-2 (5+2 j)} \Gamma (5+2 j) }{3 \Gamma (1+j) \Gamma (2+j)}\\ \times
   \,
   _2F_1\left(\frac{5}{2}+j,3+j;4;\frac{a^2}{4 b^2}\right)\\
\neq\frac{4
   b}{a^2}-\frac{-(a+b)^2 \left(7 a^2+16 b^2\right) E\left(\frac{4 a b}{(a+b)^2}\right)+\left(4 a^4-17 a^2 b^2+16
   b^4\right) K\left(\frac{4 a b}{(a+b)^2}\right)}{3 a^4 (a+b) \pi }
\end{multline}
where $Re(b)>Re(a)$.
\end{example}
\begin{example}
Derivation of equation (2.3.5) in \cite{okui_1975}. Errata. In this example we use equation (\ref{eq:okui_1975}) and set $m\to 0,\mu \to 3,\nu \to 2,\alpha \to a,\beta \to a$ and simplify;
\begin{multline}
\int_0^{\infty } \frac{J_2(a x) J_3(a x) \sin (2 b x)}{x} \, dx\\
=\sum _{j=0}^{\infty } \frac{i (-1)^j 2^{-12-4
   j} a^{5+2 j} \left((-i b)^{5+2 j}-(i b)^{5+2 j}\right) b^{-2 (5+2 j)} \Gamma (5+2 j) \,
   _2F_1\left(\frac{5}{2}+j,3+j;3;\frac{a^2}{4 b^2}\right)}{\Gamma (1+j) \Gamma (4+j)}\\
\neq\frac{b}{a}+\frac{(a+b) \left(3
   a^4+52 a^2 b^2-32 b^4\right) E\left(\frac{4 a b}{(a+b)^2}\right)+(-a+b) \left(3 a^4-20 a^2 b^2+32 b^4\right)
   K\left(\frac{4 a b}{(a+b)^2}\right)}{15 a^5 \pi }
\end{multline}
where $Re(b)>Re(a)$.
\end{example}
\begin{example}
Derivation of equation (2.5.1) in \cite{okui_1975}. In this example we use equation (\ref{eq:okui_1975}) and set $m\to 1,\mu \to 0,\nu \to 0,\alpha \to a,\beta \to b,b\to \frac{c}{2}$ and simplify;
\begin{multline}
\int_0^{\infty } J_0(a x) J_0(b x) \sin (c x) \, dx\\
=\sum _{j=0}^{\infty } \frac{i (-1)^j 2^{-2-4 j+2 (1+2 j)}
   a^{2 j} \left(2^{-1-2 j} (-i c)^{1+2 j}-2^{-1-2 j} (i c)^{1+2 j}\right) c^{-2 (1+2 j)} \Gamma (1+2 j) }{\Gamma (1+j)^2}\\ \times
\,_2F_1\left(\frac{1}{2}+j,1+j;1;\frac{b^2}{c^2}\right)\\
=\frac{2 \sqrt{\frac{a b}{-(a-b)^2+c^2}}
   K\left(\frac{4 a b}{-(a-b)^2+c^2}\right)}{\sqrt{a b} \pi }
\end{multline}
where $Re(b)>Re(a)$.
\end{example}
\begin{example}
Derivation of equation (2.5.2) in \cite{okui_1975}. In this example we use equation (\ref{eq:okui_1975}) and set $m\to 1,\mu \to 1,\nu \to 1,\alpha \to a,\beta \to b,b\to \frac{c}{2}$ and simplify;
\begin{multline}
\int_0^{\infty } J_1(a x) J_1(b x) \sin (c x) \, dx\\
=\sum _{j=0}^{\infty } \frac{i (-1)^j 2^{-6-4 j+2 (3+2 j)}
   a^{1+2 j} b \left(2^{-3-2 j} (-i c)^{3+2 j}-2^{-3-2 j} (i c)^{3+2 j}\right) c^{-2 (3+2 j)} \Gamma (3+2 j) }{\Gamma (1+j) \Gamma (2+j)}\\ \times
\,
   _2F_1\left(\frac{3}{2}+j,2+j;2;\frac{b^2}{c^2}\right)\\
=-\frac{-2 E\left(\frac{4 a
   b}{-(a-b)^2+c^2}\right)+\left(2-\frac{4 a b}{-(a-b)^2+c^2}\right) K\left(\frac{4 a b}{-(a-b)^2+c^2}\right)}{2
   \sqrt{a b} \sqrt{\frac{a b}{-(a-b)^2+c^2}} \pi }
\end{multline}
where $Re(b)>Re(a)$.
\end{example}
\begin{example}
Derivation of equation (2.5.3) in \cite{okui_1975}. In this example we use equation (\ref{eq:okui_1975}) and set $m\to 1,\mu \to 2,\nu \to 2,\alpha \to a,\beta \to b,b\to \frac{c}{2}$ and simplify;
\begin{multline}
\int_0^{\infty } J_2(a x) J_2(b x) \sin (c x) \, dx\\
=\sum _{j=0}^{\infty } \frac{i (-1)^j 2^{-11-4 j+2 (5+2 j)}
   a^{2+2 j} b^2 \left(2^{-5-2 j} (-i c)^{5+2 j}-2^{-5-2 j} (i c)^{5+2 j}\right) c^{-2 (5+2 j)} \Gamma (5+2 j) }{\Gamma (1+j) \Gamma (3+j)}\\ \times
\,
   _2F_1\left(\frac{5}{2}+j,3+j;3;\frac{b^2}{c^2}\right)\\
=\frac{-8 \left(2-\frac{4 a
   b}{-(a-b)^2+c^2}\right) E\left(\frac{4 a b}{-(a-b)^2+c^2}\right)+\left(16+\frac{48 a^2
   b^2}{\left(-(a-b)^2+c^2\right)^2}-\frac{64 a b}{-(a-b)^2+c^2}\right) K\left(\frac{4 a b}{-(a-b)^2+c^2}\right)}{24
   \sqrt{a b} \left(\frac{a b}{-(a-b)^2+c^2}\right)^{3/2} \pi }
\end{multline}
where $Re(b)>Re(a)$.
\end{example}
\begin{example}
Derivation of equation (2.5.4) in \cite{okui_1975}. In this example we use equation (\ref{eq:okui_1975}) and set $m\to 1,\mu \to 3,\nu \to 3,\alpha \to a,\beta \to b,b\to \frac{c}{2}$ and simplify;
\begin{multline}
\int_0^{\infty } J_3(a x) J_3(b x) \sin (c x) \, dx\\
=\sum _{j=0}^{\infty } \frac{i (-1)^j 2^{-15-4 j+2 (7+2 j)}
   a^{3+2 j} b^3 \left(2^{-7-2 j} (-i c)^{7+2 j}-2^{-7-2 j} (i c)^{7+2 j}\right) c^{-2 (7+2 j)} }{3 \Gamma (1+j) \Gamma (4+j)}\\ \times
\Gamma (7+2 j) \,
   _2F_1\left(\frac{7}{2}+j,4+j;4;\frac{b^2}{c^2}\right)\\
=-\frac{1}{480\sqrt{a b} \left(\frac{a b}{-(a-b)^2+c^2}\right)^{5/2} \pi }\left( -2 \left(128+\frac{368
   a^2 b^2}{\left(-(a-b)^2+c^2\right)^2}-\frac{512 a b}{-(a-b)^2+c^2}\right)\right. \\ \left.
 E\left(\frac{4 a
   b}{-(a-b)^2+c^2}\right)
+\left(128+\frac{240 a^2 b^2}{\left(-(a-b)^2+c^2\right)^2}-\frac{512 a
   b}{-(a-b)^2+c^2}\right)\right. \\ \left.
 \left(2-\frac{4 a b}{-(a-b)^2+c^2}\right) K\left(\frac{4 a b}{-(a-b)^2+c^2}\right)\right)
\end{multline}
where $Re(b)>Re(a)$.
\end{example}
\begin{example}
Derivation of equation (2.5.5) in \cite{eq:thm}. In this example we use equation (\ref{eq:okui_1975}) and set $k\to 0,a\to 1,b\to 1,r\to 0,\gamma \to 1,\delta \to 0,s\to 0,c\to 1,\mu \to \frac{1}{4},\nu \to
   \frac{1}{4},\alpha \to a,\beta \to b,x\to t^2$ and simplify;
   \begin{multline}
\int_0^1 x^2 J_{\frac{1}{4}}\left(a x^2\right) J_{\frac{1}{4}}\left(b x^2\right) \, dx=\sum _{j=0}^{\infty }
   \frac{(-1)^j 2^{-\frac{1}{2}-2 j} a^{\frac{1}{4}+2 j} \sqrt[4]{b} \,
   _1F_2\left(1+j;\frac{5}{4},2+j;-\frac{b^2}{4}\right)}{(1+j) \Gamma \left(\frac{1}{4}\right) \Gamma (1+j) \Gamma
   \left(\frac{1}{4} (5+4 j)\right)}
\end{multline}
where $Re(b)>Re(a)$.
\end{example}
\begin{example}
Derivation of equation (2.2.1) in \cite{okui_1975}. In this example we use equation (\ref{eq:okui_1975}) and set $\delta \to 0,s\to 0,z\to -1,\gamma \to p,\mu \to 0,\nu \to 0,\alpha \to a,\beta \to a$ and simplify;
\begin{multline}
\int_0^{\infty } J_0\left(a x^2\right){}^2 \cos \left(2 b x^2\right) \, dx\\
=\sum _{j=0}^{\infty } \frac{(-1)^j
   2^{-\frac{5}{2}-4 j} a^{2 j} \left((-i b)^{\frac{1}{2}+2 j}+(i b)^{\frac{1}{2}+2 j}\right) b^{-1-4 j} \Gamma
   \left(\frac{1}{2} (1+4 j)\right) }{\Gamma
   (1+j)^2}\\ \times
\, _2F_1\left(\frac{1}{4}+j,\frac{3}{4}+j;1;\frac{a^2}{4 b^2}\right)\\
=\frac{2 K\left(\frac{2 \left(-\sqrt{-a+b}+\sqrt{a+b}\right)}{2
   \sqrt{b}-\sqrt{-a+b}+\sqrt{a+b}}\right)^2}{\left(2 \sqrt{b}-\sqrt{-a+b}+\sqrt{a+b}\right) \pi ^{3/2}}
\end{multline}
where $Re(b)>Re(a)$.
\end{example}
\begin{example}
Derivation of equation (13.22.1) in \cite{watson}. Here we use equation (\ref{eq:mellin_poly}) and set $\delta \to 0,s\to 0,z\to -1,\gamma \to a,\mu \to \nu ,\alpha \to b,\beta \to c,r\to 1,m\to \mu -1,x\to t$ and simplify;
\begin{multline}
\int_0^{\infty } e^{-a t} t^{-1+\mu } J_{\nu }(b t) J_{\nu }(c t) \, dt\\
=\sum _{j=0}^{\infty } \frac{(-1)^j
   2^{-2 j-2 \nu } a^{-2 j-\mu -2 \nu } b^{2 j+\nu } c^{\nu } \Gamma (2 j+\mu +2 \nu ) }{\Gamma (1+j) \Gamma (1+\nu ) \Gamma(1+j+\nu )}\\ \times
\, _2F_1\left(j+\frac{\mu}{2}+\nu ,\frac{1}{2}+j+\frac{\mu }{2}+\nu ;1+\nu ;-\frac{c^2}{a^2}\right)\\
=\frac{\left((b c)^{\nu } \Gamma (\mu +2 \nu )\right) }{\pi  a^{\mu +2 \nu } \Gamma (2 \nu +1)}\int_0^{\pi } \, _2F_1\left(\frac{1}{2} (\mu +2
   \nu ),\frac{1}{2} (\mu +2 \nu +1);\nu +1;-\frac{b^2+c^2-2 b c \cos (\phi )}{a^2}\right)\\ \sin ^{2 \nu }(\phi ) \,
   d\phi 
\end{multline}
where $Re(\mu+2\nu)>0$.
\end{example}
\begin{example}
Derivation of equation (13.22.2) in \cite{watson}.  Here we use equation (\ref{eq:mellin_poly}) and set $\mu \to -1,\nu \to 1,c\to 1,b\to 1$ and simplify;
\begin{multline}
\int_0^{\infty } \frac{e^{-a t} J_1(t){}^2}{t^2} \, dt\\=\sum _{j=0}^{\infty } \frac{(-1)^j a^{-1-2 j} \Gamma
   \left(\frac{1}{2}+j\right) \, _2F_1\left(\frac{1}{2}+j,1+j;2;-\frac{1}{a^2}\right)}{4 \sqrt{\pi } \Gamma
   (2+j)}\\
=\frac{a \left(-\pi +2 E\left(-\frac{4}{a^2}\right)+\frac{1}{3} \sqrt{1+\frac{4}{a^2}}
   \left(-\left(\left(2+a^2\right) E\left(\frac{4}{4+a^2}\right)\right)+a^2
   K\left(\frac{4}{4+a^2}\right)\right)\right)}{2 \pi }
\end{multline}
where $Re(a)>1/2$.
\end{example}
\begin{example}
Derivation (3.12.15.25) in \cite{prud4}. In this example we use equation (\ref{eq:mellin_poly}) and set $\delta \to 0,s\to 0,z\to -1,\gamma \to p,r\to 1,m\to 2 \nu ,\mu \to \nu ,\alpha \to a,\beta \to a$ and simplify;
\begin{multline}
\int_0^{\infty } e^{-p x} x^{2 \nu } J_{\nu }(a x){}^2 \, dx\\
=\sum _{j=0}^{\infty } \frac{(-1)^j 2^{-2 j-2 \nu }
   a^{2 j+2 \nu } p^{-1-2 j-4 \nu } \Gamma (1+2 j+4 \nu ) \, _2F_1\left(\frac{1}{2}+j+2 \nu ,1+j+2 \nu ;1+\nu
   ;-\frac{a^2}{p^2}\right)}{\Gamma (1+j) \Gamma (1+\nu ) \Gamma (1+j+\nu )}\\
=\frac{\left(2^{4 \nu } a^{2 \nu }\right)
   \left(\Gamma \left(\nu +\frac{1}{2}\right) \Gamma \left(2 \nu +\frac{1}{2}\right)\right) \, _2F_1\left(\nu
   +\frac{1}{2},2 \nu +\frac{1}{2};\nu +1;-\left(\frac{2 a}{p}\right)^2\right)}{\left(\pi  p^{4 \nu +1}\right) \Gamma
   (\nu +1)}
\end{multline}
where $Re(a)>0, Re(p)>0,Re(\nu)>-1/2$
\end{example}
\begin{example}
Derivation (3.12.16.1) in \cite{prud4}. In this example we use equation (\ref{eq:mellin_poly}) and set $\delta \to 0,s\to 0,z\to -1,\gamma \to p,\mu \to \nu ,\nu \to \nu +\frac{1}{4},\alpha \to a,\beta \to a,x\to t^2$ and simplify;
\begin{multline}
\int_0^{\infty } e^{-p x} x^{2 \nu } J_{\nu }(a x){}^2 \, dx\\
=\sum _{j=0}^{\infty } \frac{(-1)^j 2^{-2 j-2 \nu }
   a^{2 j+2 \nu } p^{-1-2 j-4 \nu } \Gamma (1+2 j+4 \nu ) \, _2F_1\left(\frac{1}{2}+j+2 \nu ,1+j+2 \nu ;1+\nu ;-\frac{a^2}{p^2}\right)}{\Gamma (1+j) \Gamma (1+\nu ) \Gamma (1+j+\nu )}\\
=\frac{\left(2^{4 \nu } a^{2 \nu }\right) \left(\Gamma \left(\nu +\frac{1}{2}\right) \Gamma \left(2 \nu +\frac{1}{2}\right)\right) \, _2F_1\left(\nu+\frac{1}{2},2 \nu +\frac{1}{2};\nu +1;-\left(\frac{2 a}{p}\right)^2\right)}{\left(\pi  p^{4 \nu +1}\right) \Gamma (\nu +1)}
\end{multline}
where $Re(a)>0, Re(p)>0,Re(\nu)>-1/2$
\end{example}
\begin{example}
The Hurwitz-Lerch zeta function. Here we use equation (\ref{eq:thm}) and set $\delta \to -1,b\to 1,r\to 0,\gamma \to 1,a\to 1,c\to -z$ and simplify the series over $l\in[0,\infty), p\in[0,\infty)$;
\begin{multline}\label{eq:lerch_bessel}
\int_0^1 \frac{x^m J_{\mu }(x \alpha ) J_{\nu }(x \beta ) \log ^k\left(\frac{1}{x}\right)}{1-x^s z} \, dx\\
=\sum
   _{h=0}^{\infty } \sum _{j=0}^{\infty } \frac{(-1)^{h+j} 2^{-2 h-2 j-\mu -\nu } s^{-1-k} \alpha ^{2 j+\mu } \beta
   ^{2 h+\nu } \Gamma (1+k) \Phi \left(z,1+k,\frac{1+2 h+2 j+m+\mu +\nu }{s}\right)}{\Gamma (1+h) \Gamma (1+j) \Gamma
   (1+j+\mu ) \Gamma (1+h+\nu )}
\end{multline}
where $Re(s)>1$.
\end{example}
\begin{example}
The  Hurwitz zeta function. Here we use equation (\ref{eq:lerch_bessel}) and set $z\to -1$ and simplify using equation [Wolfram,\href{http://functions.wolfram.com/10.06.03.0074.01}{01}];
\begin{multline}\label{eq:lerch_bessel_1}
\int_0^1 \frac{x^m J_{\mu }(x \alpha ) J_{\nu }(x \beta ) \log ^k\left(\frac{1}{x}\right)}{1+x^s} \, dx\\
=\sum
   _{h=0}^{\infty } \sum _{j=0}^{\infty } \frac{(-1)^{h+j} 2^{-2 h-2 j-\mu -\nu } s^{-1-k} \alpha ^{2 j+\mu } \beta
   ^{2 h+\nu } \Gamma (1+k) }{\Gamma (1+h) \Gamma (1+j)
   \Gamma (1+j+\mu ) \Gamma (1+h+\nu )}\\ \times
\left(2^{-1-k} \zeta \left(1+k,\frac{1+2 h+2 j+m+\mu +\nu }{2 s}\right)\right. \\ \left.-2^{-1-k} \zeta
   \left(1+k,\frac{1}{2} \left(1+\frac{1+2 h+2 j+m+\mu +\nu }{s}\right)\right)\right)
\end{multline}
where $Re(s)>1$.
\end{example}
\begin{example}
The digamma function. In the example we use equation (\ref{eq:lerch_bessel_1}) and apply l'Hopital's rule as $k\to 0$ and simplify using equation [Wolfram, \href{https://mathworld.wolfram.com/HurwitzZetaFunction.html}{19}];
\begin{multline}
\int_0^1 \frac{x^m J_{\mu }(x \alpha ) J_{\nu }(x \beta )}{1+x^s} \, dx\\
=\sum _{h=0}^{\infty } \sum
   _{j=0}^{\infty } \frac{(-1)^{h+j} 2^{-1-2 h-2 j-\mu -\nu } \alpha ^{2 j+\mu } \beta ^{2 h+\nu } }{s \Gamma (1+h) \Gamma (1+j) \Gamma (1+j+\mu ) \Gamma (1+h+\nu )}\\ \times
\left(-\psi
   ^{(0)}\left(\frac{1+2 h+2 j+m+\mu +\nu }{2 s}\right)+\psi ^{(0)}\left(\frac{1+2 h+2 j+m+s+\mu +\nu }{2
   s}\right)\right)
\end{multline}
where $Re(s)>1$.
\end{example}
\begin{example}
Nested logarithm function. Here we use equation (\ref{eq:lerch_bessel_1}) and take the first partial derivative with respect to $k$ and apply l'Hopital's rule as $k\to 0$ and simplify using [Wolfram, \href{http://functions.wolfram.com/10.02.20.0011.01}{01}];
\begin{multline}\label{eq:hurwitz}
\int_0^1 \frac{x^m J_{\mu }(x \alpha ) J_{\nu }(x \beta ) \log \left(\log
   \left(\frac{1}{x}\right)\right)}{1+x^s} \, dx\\
=\sum _{h=0}^{\infty } \sum _{j=0}^{\infty } \frac{(-1)^{h+j} 2^{-1-2
   h-2 j-\mu -\nu } \alpha ^{2 j+\mu } \beta ^{2 h+\nu } }{s \Gamma (1+h) \Gamma (1+j) \Gamma (1+j+\mu ) \Gamma (1+h+\nu )}\\ \times
\left((\gamma +\log (2 s)) \left(\psi ^{(0)}\left(\frac{1+2
   h+2 j+m+\mu +\nu }{2 s}\right)\right.\right. \\ \left.\left.
-\psi ^{(0)}\left(\frac{1+2 h+2 j+m+s+\mu +\nu }{2 s}\right)\right)-\gamma
   _1\left(\frac{1+2 h+2 j+m+\mu +\nu }{2 s}\right)\right. \\ \left.
+\gamma _1\left(\frac{1+2 h+2 j+m+s+\mu +\nu }{2
   s}\right)\right)
\end{multline}
where $Re(s)>1$.
\end{example}
\begin{example}
The log-gamma function. Here we use equation (\ref{eq:hurwitz}) and set $\mu \to 0,\nu \to 0,s\to 2,\alpha \to a,\beta \to a$ then form a second equation by replacing $m\to s$ and take their difference apply l'Hopital's rule as $k\to -1$ and simplify;
\begin{multline}
\int_0^1 \frac{\left(-x^m+x^s\right) J_0(a x){}^2}{\left(1+x^2\right) \log \left(\frac{1}{x}\right)} \, dx
=\sum
   _{h=0}^{\infty } \sum _{j=0}^{\infty } \left(-\left(\frac{a}{2}\right)^2\right)^{h+j}\frac{1}{\Gamma (1+h)^2 \Gamma (1+j)^2}\\ \times
 \log \left(\frac{\Gamma
   \left(\frac{1}{4} (3+2 h+2 j+m)\right) \Gamma \left(\frac{1}{4} (1+2 h+2 j+s)\right)}{\Gamma \left(\frac{1}{4} (1+2h+2 j+m)\right) \Gamma \left(\frac{1}{4} (3+2 h+2 j+s)\right)}\right)
\end{multline}
where $|Re(m)|<1, |Re(s)|<1$.
\end{example}
\section{Integrals of the convolution kind}
In 1930, integrals involving products of Bessel functions of the convolution kind were first studied with applications by Watson \cite{watson}, Bailey \cite{bailey1} where he generalized the results of Bateman \cite{watson} and Kapteyn \cite{watson} sections 12.2 and 12.21, Modell \cite{mordell}.  In 1931, Bailey \cite{bailey2} further studied these integrals where he generalized the trigonometric integrals of Kapteyn \cite{watson} section 12.21. In this section we use contour integration \cite{reyn4} to derive an alternate series representation for convolution integrals listed in Luke \cite{luke}, Prudikov et al. \cite{prud2}and Errata for some formulae which are in error.
\\\\
The contour integral representation form involving the generalized convolution form is given by;
\begin{multline}\label{eq:c1}
\frac{1}{2\pi i}\int_{C}\int_{0}^{b}a^w e^{-x^r \theta } w^{-1-k} x^{m+w} \left(1+x^{\eta } \rho \right)^{\lambda } J_{\mu }(\alpha +x \beta )
   J_{\nu }\left(d x^c\right)dxdw\\
=\frac{1}{2\pi i}\int_{C}\sum _{p=0}^{\infty } \sum _{l=0}^{\infty } \sum _{j=0}^{\infty } \sum _{q=0}^{\infty }
   \frac{(-1)^{h+j} 2^{-2 h-2 j-\mu -\nu } a^w b^{1+2 c j+m+p+l r+w+q \eta +c \nu } d^{2 j+\nu } w^{-1-k} \alpha ^{2
   h-p+\mu } }{(1+2 c j+m+p+l r+w+q \eta +c \nu )
   h! j! l! \Gamma (1+h+\mu ) \Gamma (1+j+\nu )}\\ \times
\beta ^p (-\theta )^l \rho ^q \binom{\lambda }{q} \binom{2 h+\mu }{p}dw
\end{multline}
where $Re(\mu+\nu)>0,Re(b)>0$.
\subsection{Left-hand side contour integral representation}
Using a generalization of Cauchy's integral formula \ref{intro:cauchy} and the procedure in section (2.1), we form the definite integral given by;
\begin{multline}\label{eq:c2}
\frac{1}{2\pi i}\int_{0}^{b}\frac{x^m e^{-\theta  x^r} \log ^k(a x) \left(\rho  x^{\eta }+1\right)^{\lambda } J_{\nu }\left(dx^c\right) J_{\mu }(\alpha +x \beta )}{k!}dx\\
=\frac{1}{2\pi i}\int_{C}\int_{0}^{b}w^{-k-1} x^m (a x)^w e^{-\theta  x^r} \left(\rho  x^{\eta
   }+1\right)^{\lambda } J_{\nu }\left(d x^c\right) J_{\mu }(\alpha +x \beta )dxdw
\end{multline}
where $Re(\mu+\nu)>0,Re(b)>0$.
\subsection{Right-hand side contour integral}
Using a generalization of Cauchy's integral formula \ref{intro:cauchy}, and the procedure in section (2.2) we get;
\begin{multline}\label{eq:c3}
\sum _{p=0}^{\infty } \sum _{l=0}^{\infty } \sum _{j=0}^{\infty } \sum _{q=0}^{\infty } \frac{(-1)^{h+j} 2^{-2
   h-2 j-\mu -\nu } a^{-1-2 c j-m-p-l r-q \eta -c \nu } d^{2 j+\nu } \alpha ^{2 h-p+\mu } }{h! j! k! l! \Gamma (1+h+\mu ) \Gamma (1+j+\nu )(-1-2c j-m-p-l r-q \eta -c \nu )^{1+k} }\\ \times
\beta ^p (-\theta )^l \rho ^q \binom{\lambda }{q} \binom{2 h+\mu }{p} \Gamma (1+k,-((1+2 c j+m+p+l r+q\eta +c \nu ) \log (a b)))\\
=-\frac{1}{2\pi i}\int_{C}\sum _{p=0}^{\infty } \sum
   _{l=0}^{\infty } \sum _{j=0}^{\infty } \sum _{q=0}^{\infty } \frac{(-1)^{h+j} 2^{-2 h-2 j-\mu -\nu } b^{1+2 c
   j+m+p+l r+q \eta +c \nu } (a b)^w d^{2 j+\nu } w^{-1-k} \alpha ^{2 h-p+\mu } }{(1+2 c j+m+p+l r+w+q \eta +c \nu ) h! j! l! \Gamma (1+h+\mu ) \Gamma(1+j+\nu )}\\ \times
\beta ^p (-\theta )^l \rho ^q
   \binom{\lambda }{q} \binom{2 h+\mu }{p}dw
\end{multline}
where $Re(\mu+\nu)>0,Re(b)>0$.
from equation [Wolfram, \href{http://functions.wolfram.com/07.27.02.0002.01}{07.27.02.0002.01}] where $|Re(b)|<1$.
\begin{theorem}
Generalized convolution integral with two Bessel functions.
\begin{multline}\label{eq:thm_c}
\int_0^b e^{-x^r \theta } x^m \left(1+x^{\eta } \rho \right)^{\lambda } J_{\mu }(\alpha +x \beta ) J_{\nu
   }\left(d x^c\right) \log ^k(a x) \, dx\\
=\sum _{h=0}^{\infty } \sum _{j=0}^{\infty } \sum _{l=0}^{\infty } \sum
   _{p=0}^{\infty } \sum _{q=0}^{\infty } \frac{(-1)^{h+j} 2^{-2 h-2 j-\mu -\nu } a^{-1-2 c j-m-p-l r-q \eta -c \nu }d^{2 j+\nu } \alpha ^{2 h-p+\mu } \beta ^p (-\theta )^l  \rho ^q
   \binom{\lambda }{q} }{h! j! l!\Gamma (1+h+\mu ) \Gamma (1+j+\nu )(-1-2 c j-m-p-l r-q \eta -c \nu )^{1+k}}\\ \times
\binom{2 h+\mu }{p} \Gamma (1+k,-((1+2 c j+m+p+l r+q \eta +c \nu ) \log (a b)))
\end{multline}
where $Re(b)>0,Re(m+\mu+\nu)>0$.
\end{theorem}
\begin{proof}
Since the right-hand sides of equations (\ref{eq:c2}) and (\ref{eq:c3}) are equal relative to equation (\ref{eq:c1}), we may equate the left-hand sides and simplify the gamma function to yield the stated result.
\end{proof}
\section{Derivation of Entries (2.12.33.7-13) in \cite{prud2}}
In this section we derive some entries in Prudnikov et al. 1986b.
\begin{example}
Derivation of equation (2.12.33.7) in \cite{prud2}. 
Here we use equation (\ref{eq:thm_c}) and set $k\to 0,a\to 1,b\to a,\theta \to 1,r\to 0,m\to \alpha -1,\lambda \to -1,\eta \to 1,\rho \to -\frac{1}{a},\alpha
   \to a c,\beta \to -c,d\to c,c\to 1$ and simplify using [DLMF,\href{https://dlmf.nist.gov/15.4.E20}{15.4.20}] to get;
   \begin{multline}\label{eq:pr1}
\int_0^a \frac{x^{-1+\alpha } J_{\mu }(c (a-x)) J_{\nu }(c x)}{a-x} \, dx\\
=\sum _{j=0}^{\infty } \frac{(-1)^j2^{-2 j-\mu -\nu } a^{-1+2 j+\alpha +\mu +\nu } c^{2 j+\mu +\nu } \Gamma (\mu ) \Gamma (2 j+\alpha +\nu ) }{\Gamma (1+j) \Gamma(1+\mu ) \Gamma (1+j+\nu ) \Gamma (2 j+\alpha +\mu +\nu )}\\ \times
\,_2F_3\left(\frac{1}{2}+\frac{\mu }{2},\frac{\mu }{2};1+\mu ,j+\frac{\alpha }{2}+\frac{\mu }{2}+\frac{\nu}{2},\frac{1}{2}+j+\frac{\alpha }{2}+\frac{\mu }{2}+\frac{\nu }{2};-\frac{1}{4} a^2 c^2\right)\\
=\frac{\left(\frac{2}{c}\right)^{\alpha } }{a \mu }\sum
   _{k=0}^{\infty } \frac{\left((\alpha )_k (-1)^k\right) \Gamma (\alpha +\nu +k) (\alpha +\mu +\nu +2 k) J_{\alpha
   +\mu +\nu +2 k}(a c)}{k! \Gamma (\nu +k+1)}
\end{multline}
where $Re(\alpha+\mu+\nu)>0,Re(a)>0$.
\end{example}
\begin{example}
Derivation of equation (2.12.33.8) in \cite{prud2}. Here we use equation (\ref{eq:pr1}) and set $\alpha \to 0,k\to 0$ and simplify to get;
\begin{multline}
\int_0^a \frac{J_{\mu }(c (a-x)) J_{\nu }(c x)}{(a-x) x} \, dx\\
=\sum _{j=0}^{\infty } \frac{(-1)^j 2^{-2 j-\mu
   -\nu } a^{-1+2 j+\mu +\nu } c^{2 j+\mu +\nu } \Gamma (\mu ) \Gamma (2 j+\nu ) }{\Gamma (1+j) \Gamma (1+\mu ) \Gamma (1+j+\nu ) \Gamma (2 j+\mu +\nu )}\\ \times
\, _2F_3\left(\frac{1}{2}+\frac{\mu
   }{2},\frac{\mu }{2};1+\mu ,j+\frac{\mu }{2}+\frac{\nu }{2},\frac{1}{2}+j+\frac{\mu }{2}+\frac{\nu }{2};-\frac{1}{4}
   a^2 c^2\right)\\
=\frac{(\mu +\nu ) J_{\mu
   +\nu }(a c) \Gamma (\nu )}{a \mu  \Gamma (1+\nu )}
\end{multline}
where $Re(\alpha+\mu+\nu)>0,Re(a)>0$.
\end{example}
\begin{example}
Derivation of equation (2.12.33.9(i)) in \cite{prud2}. Here we use equation (\ref{eq:thm_c}) and set $k\to 0,a\to 1,r\to 0,\theta \to 1,b\to a,m\to \nu ,\lambda \to \mu +1,\rho \to -\frac{1}{a},\eta \to 1,\alpha
   \to a c,\beta \to -c,d\to c,c\to 1$ and simplify using [DLMF,\href{https://dlmf.nist.gov/15.4.E20}{15.4.20}] to get;
   \begin{multline}
\int_0^a (a-x)^{1+\mu } x^{\nu } J_{\mu }(c (a-x)) J_{\nu }(c x) \, dx\\
=\sum _{j=0}^{\infty } \frac{(-1)^j
   2^{1+\mu +\nu } a^{2+2 j+2 \mu +2 \nu } c^{2 j+\mu +\nu } \Gamma \left(\frac{1}{2} (3+2 \mu )\right) \Gamma
   \left(\frac{1}{2} (1+2 j+2 \nu )\right) }{\pi  \Gamma (1+j) \Gamma (3+2 j+2 \mu +2 \nu )}\\ \times
\, _1F_2\left(\frac{3}{2}+\mu ;\frac{3}{2}+j+\mu +\nu ,2+j+\mu +\nu
   ;-\frac{1}{4} a^2 c^2\right)\\
=\frac{a^{\mu +\nu +\frac{3}{2}}
   }{\sqrt{2 \pi  c}}B\left(\mu +\frac{3}{2},\nu +\frac{1}{2}\right) J_{\mu +\nu +\frac{1}{2}}(a c)
\end{multline}
where $Re(\alpha+\mu+\nu)>0,Re(a)>0$.
\end{example}
\begin{example}
Derivation of equation (2.12.33.9(ii)) in \cite{prud2}. Here we use equation (\ref{eq:thm_c}) and set $k\to 0,a\to 1,r\to 0,\theta \to 1,b\to a,m\to \nu ,\lambda \to \mu ,\rho \to -\frac{1}{a},\eta \to 1,\alpha \to
   a c,\beta \to -c,d\to c,c\to 1$ and simplify using [DLMF,\href{https://dlmf.nist.gov/15.4.E20}{15.4.20}] to get;
   \begin{multline}
\int_0^a (a-x)^{\mu } x^{\nu } J_{\mu }(a c-c x) J_{\nu }(c x) \, dx\\
=\sum _{j=0}^{\infty } \frac{(-1)^j 2^{\mu
   +\nu } a^{1+2 j+2 \mu +2 \nu } c^{2 j+\mu +\nu } \Gamma \left(\frac{1}{2} (1+2 \mu )\right) \Gamma
   \left(\frac{1}{2} (1+2 j+2 \nu )\right) }{\pi  \Gamma (1+j) \Gamma (2 (1+j+\mu +\nu ))}\\ \times
\, _1F_2\left(\frac{1}{2}+\mu ;1+j+\mu +\nu ,\frac{3}{2}+j+\mu +\nu
   ;-\frac{1}{4} a^2 c^2\right)\\
=\frac{a^{\mu +\nu +\frac{1}{2}}B\left(\mu +\frac{1}{2},\nu +\frac{1}{2}\right) J_{\mu +\nu +\frac{1}{2}}(a c)}{\sqrt{2 \pi  c}}
\end{multline}
where $Re(\alpha+\mu+\nu)>0,Re(a)>0$.
\end{example}
\begin{example}
Derivation of equation (2.12.33.10) in \cite{prud2}. Errata. Here we use equation (\ref{eq:thm_c}) and set $k\to 0,a\to 1,r\to 0,\theta \to 1,b\to a,m\to \nu +1,\lambda \to \mu +1,\rho \to -\frac{1}{a},\eta \to 1,\alpha
   \to a c,\beta \to -c,d\to c,c\to 1$ and simplify using [DLMF,\href{https://dlmf.nist.gov/15.4.E20}{15.4.20}] to get;
   \begin{multline}
\int_0^a (a-x)^{1+\mu } x^{1+\nu } J_{\mu }(c (a-x)) J_{\nu }(c x) \, dx\\
=\sum _{j=0}^{\infty } \frac{(-1)^j
   2^{2+\mu +\nu } a^{3+2 j+2 \mu +2 \nu } c^{2 j+\mu +\nu } \Gamma \left(\frac{1}{2} (3+2 \mu )\right) \Gamma
   \left(\frac{1}{2} (3+2 j+2 \nu )\right) }{\pi  \Gamma (1+j) \Gamma (2 (2+j+\mu +\nu ))}\\ \times
\, _1F_2\left(\frac{3}{2}+\mu ;2+j+\mu +\nu ,\frac{5}{2}+j+\mu +\nu
   ;-\frac{1}{4} a^2 c^2\right)\\
=\frac{2 \left(\Gamma \left(\mu+\frac{3}{2}\right) \Gamma \left(\nu +\frac{3}{2}\right)\right) a^{\mu +\nu +\frac{3}{2}} c^{-3/2} }{\sqrt{2 \pi } \Gamma(\mu +\nu +2)}\left(J_{\mu+\nu +\frac{3}{2}}(a c)-\frac{(a c) J_{\mu +\nu +\frac{5}{2}}(a c)}{2 (\mu +\nu +2)}\right)
\end{multline}
where $Re(\alpha+\mu+\nu)>0,Re(a)>0$.
\end{example}
\begin{example}
Derivation of equation (2.12.33.11) in \cite{prud2}. Here we use equation (\ref{eq:thm_c}) and set $k\to 0,a\to 1,r\to 0,\theta \to 1,b\to a,m\to \nu +1,\lambda \to -\nu -2,\rho \to -\frac{1}{a},\eta \to
   1,\alpha \to a c,\beta \to -c,d\to c,c\to 1$ and simplify using [DLMF,\href{https://dlmf.nist.gov/15.4.E20}{15.4.20}] to get;
   \begin{multline}
\int_0^a (a-x)^{-2-\nu } x^{1+\nu } J_{\mu }(c (a-x)) J_{\nu }(c x) \, dx\\
=\sum _{j=0}^{\infty } \frac{(-1)^j
   2^{1-\mu +\nu } a^{2 j+\mu +\nu } c^{2 j+\mu +\nu } \Gamma (-1+\mu -\nu ) \Gamma \left(\frac{1}{2} (3+2 j+2 \nu)\right) }{\sqrt{\pi } \Gamma (1+j) \Gamma (1+\mu ) \Gamma (1+2 j+\mu +\nu )}\\ \times
\, _2F_3\left(-\frac{1}{2}+\frac{\mu }{2}-\frac{\nu }{2},\frac{\mu }{2}-\frac{\nu }{2};1+\mu ,\frac{1}{2}+j+\frac{\mu }{2}+\frac{\nu }{2},1+j+\frac{\mu }{2}+\frac{\nu }{2};-\frac{1}{4} a^2
   c^2\right)\\
=\frac{2^{2 \nu +1} \left(\Gamma
   \left(\nu +\frac{3}{2}\right) \Gamma (\mu -\nu -1)\right) \left(\frac{a c}{2}\right)^{\nu } \left(J_{\mu }(a
   c)-\frac{(a c) J_{\mu +1}(a c)}{\mu +\nu +1}\right)}{\sqrt{\pi } \Gamma (\mu +\nu +1)}
\end{multline}
where $Re(\alpha+\mu+\nu)>0,Re(a)>0,-1< Re(\nu)<0,1/2< Re(\mu)<1$.
\end{example}
\begin{example}
Derivation of equation (2.12.33.12) in \cite{prud2}. Here we use equation (\ref{eq:thm_c}) and set $k\to 0,a\to 1,r\to 0,\theta \to 1,b\to a,m\to \nu ,\lambda \to -\nu -1,\rho \to -\frac{1}{a},\eta \to 1,\alpha
   \to a c,\beta \to -c,d\to c,c\to 1$ and simplify using [DLMF,\href{https://dlmf.nist.gov/15.4.E20}{15.4.20}] to get;
   \begin{multline}
\int_0^a (a-x)^{-1-\nu } x^{\nu } J_{\mu }(c (a-x)) J_{\nu }(c x) \, dx\\
=\sum _{j=0}^{\infty } \frac{(-1)^j2^{-\mu +\nu } a^{2 j+\mu +\nu } c^{2 j+\mu +\nu } \Gamma (\mu -\nu ) \Gamma \left(\frac{1}{2} (1+2 j+2 \nu)\right) }{\sqrt{\pi } \Gamma (1+j) \Gamma (1+\mu ) \Gamma (1+2 j+\mu +\nu )}\\ \times
\, _2F_3\left(\frac{\mu }{2}-\frac{\nu }{2},\frac{1}{2}+\frac{\mu }{2}-\frac{\nu }{2};1+\mu ,\frac{1}{2}+j+\frac{\mu }{2}+\frac{\nu }{2},1+j+\frac{\mu }{2}+\frac{\nu }{2};-\frac{1}{4} a^2c^2\right)\\
=\frac{(2 a c)^{\nu } \left(\Gamma\left(\nu +\frac{1}{2}\right) \Gamma (\mu -\nu )\right) J_{\mu }(a c)}{\sqrt{\pi } \Gamma (\mu +\nu +1)}
\end{multline}
where $Re(\alpha+\mu+\nu)>0,Re(a)>0,Re(\nu)>0,1/2< Re(\mu)<1$.
\end{example}
\begin{example}
Derivation of equation (2.12.33.13) in \cite{prud2}. Here we use equation (\ref{eq:thm_c}) and set $k\to 0,a\to 1,r\to 0,\theta \to 1,b\to a,m\to -2,\lambda \to 1,\rho \to -\frac{1}{a},\eta \to 1,\alpha \to a
   c,\beta \to -c,d\to c,c\to 1,\mu \to 0$ and simplify using [DLMF,\href{https://dlmf.nist.gov/15.4.E20}{15.4.20}] to get;
   \begin{multline}
\int_0^a \frac{(a-x) J_0(c (a-x)) J_{\nu }(c x)}{x^2} \, dx\\
=\sum _{j=0}^{\infty } \frac{(-1)^j 2^{-2 j-\nu }a^{2 j+\nu } c^{2 j+\nu } \Gamma (-1+2 j+\nu ) \, _1F_2\left(\frac{3}{2};\frac{1}{2}+j+\frac{\nu }{2},1+j+\frac{\nu}{2};-\frac{1}{4} a^2 c^2\right)}{\Gamma (1+j) \Gamma (1+j+\nu ) \Gamma (1+2 j+\nu )}\\
=\frac{J_{\nu }(a c)-\frac{(ac) J_{\nu +1}(a c)}{\nu +1}}{\nu  (\nu -1)}
\end{multline}
where $Re(a)>0,Re(\nu)>1$.
\end{example}
\section{Derivation of convolution integrals in Luke (1962)}
In this section we derive entries (13.3.2.35-46) in \cite{luke} in terms of infinite series of the hypergeometric function.
\begin{example}
Derivation of equation (13.3.2.35) in \cite{luke}. Here we use equation (\ref{eq:thm_c}) and set $r\to 0,\theta \to 1,k\to 0,a\to 1,b\to 1,x\to t^2$ then $\alpha \to z,\beta \to -z,t\to x,c\to 1,\eta \to 1,\rho \to -1$ and $d\to z,z\to w,\mu \to \nu ,\nu \to \mu ,\lambda \to \beta -1,m\to \alpha -1$ simplify using [DLMF,\href{https://dlmf.nist.gov/15.4.E20}{15.4.20}] to get;
\begin{multline}\label{eq:eq_trig}
\int_0^{\frac{\pi }{2}} J_{\mu }\left(z \sin ^2(t)\right) J_{\nu }\left(w \cos ^2(t)\right) \cos ^{-1+2 \beta
   }(t) \sin ^{-1+2 \alpha }(t) \, dt\\
=\sum _{j=0}^{\infty } \frac{(-1)^j 2^{-1-2 j-\mu -\nu } w^{\nu } z^{2 j+\mu }
   \Gamma (2 j+\alpha +\mu ) \Gamma (\beta +\nu ) }{\Gamma (1+j) \Gamma (1+j+\mu ) \Gamma (1+\nu ) \Gamma (2 j+\alpha +\beta +\mu +\nu
   )}\\ \times
\, _2F_3\left(\frac{\beta }{2}+\frac{\nu
   }{2},\frac{1}{2}+\frac{\beta }{2}+\frac{\nu }{2};j+\frac{\alpha }{2}+\frac{\beta }{2}+\frac{\mu }{2}+\frac{\nu
   }{2},\frac{1}{2}+j+\frac{\alpha }{2}+\frac{\beta }{2}+\frac{\mu }{2}+\frac{\nu }{2},1+\nu
   ;-\frac{w^2}{4}\right)\\
=\frac{\left(\left(\frac{z}{2}\right)^{\mu } \left(\frac{w}{2}\right)^{\nu } \Gamma (\alpha +\mu ) \Gamma (\beta
   +\nu )\right)}{2 (\Gamma (\mu +1) \Gamma (\nu +1) \Gamma (\alpha +\beta
   +\mu +\nu ))}\\
 \sum _{k=0}^{\infty } \frac{\left((-1)^k \left(\frac{z}{2}\right)^{2 k} \left(\frac{\alpha +\mu
   }{2}\right)_k \left(\frac{1}{2} (\alpha +\mu +1)\right)_k\right) }{k! (\mu +1)_k \left(\frac{1}{2} (\alpha +\beta +\mu +\nu )\right)_k
   \left(\frac{1}{2} (\alpha +\beta +\mu +\nu +1)\right)_k}\\ \times
\, _4F_3\left(-k,-\mu -k,\frac{\beta +\nu
   }{2},\frac{1}{2} (\beta +\nu +1);\frac{1}{2} (1-\alpha -\mu )-k,\frac{1}{2} (2-\alpha -\mu )-k,\nu
   +1;\left(\frac{w}{z}\right)^2\right)
\end{multline}
where $Re(\mu+\alpha)>0,Re(\nu+\beta)>0$.
\end{example}
\begin{example}
Derivation of equation (13.3.2.36) in \cite{luke}. Here we use equation (\ref{eq:eq_trig}) and set $w\to z$ simplify to get;
\begin{multline}
\int_0^{\frac{\pi }{2}} J_{\mu }\left(z \sin ^2(t)\right) J_{\nu }\left(z \cos ^2(t)\right) \cos ^{-1+2 \beta
   }(t) \sin ^{-1+2 \alpha }(t) \, dt\\
=\sum _{j=0}^{\infty } \frac{(-1)^j 2^{-1-2 j-\mu -\nu } z^{2 j+\mu +\nu } \Gamma
   (2 j+\alpha +\mu ) \Gamma (\beta +\nu ) }{\Gamma (1+j) \Gamma (1+j+\mu )
   \Gamma (1+\nu ) \Gamma (2 j+\alpha +\beta +\mu +\nu )}\\ \times
\, _2F_3\left(\frac{\beta }{2}+\frac{\nu }{2},\frac{1}{2}+\frac{\beta
   }{2}+\frac{\nu }{2};j+\frac{\alpha }{2}+\frac{\beta }{2}+\frac{\mu }{2}+\frac{\nu }{2},\frac{1}{2}+j+\frac{\alpha
   }{2}+\frac{\beta }{2}+\frac{\mu }{2}+\frac{\nu }{2},1+\nu ;-\frac{z^2}{4}\right)\\
=\sum _{k=0}^{\infty } \frac{(-1)^k 2^{-1-2 k-\mu -\nu } z^{2
   k+\mu +\nu } \Gamma (\alpha +\mu ) \Gamma (\beta +\nu ) }{k! \Gamma (1+\mu ) \Gamma (1+\nu ) \Gamma (\alpha +\beta +\mu +\nu ) (1+\mu )_k \left(\frac{1}{2}
   (\alpha +\beta +\mu +\nu )\right)_k \left(\frac{1}{2} (1+\alpha +\beta +\mu +\nu )\right)_k}\\ \times
\, _4F_3\left(-k,-k-\mu ,\frac{\beta }{2}+\frac{\nu
   }{2},\frac{1}{2}+\frac{\beta }{2}+\frac{\nu }{2};\frac{1}{2}-k-\frac{\alpha }{2}-\frac{\mu }{2},1-k-\frac{\alpha
   }{2}-\frac{\mu }{2},1+\nu ;1\right) \\
\left(\frac{\alpha +\mu }{2}\right)_k \left(\frac{1}{2} (1+\alpha +\mu
   )\right)_k
\end{multline}
where $Re(\mu+\alpha)>0,Re(\nu+\beta)>0$.
\end{example}
\begin{example}
Derivation of equation (13.3.2.37) in \cite{luke}. Errata. Here we use equation (\ref{eq:eq_trig}) and set $w\to z,\alpha \to \mu +1$ simplify to get;
\begin{multline}
\int_0^{\frac{\pi }{2}} J_{\mu }\left(z \sin ^2(t)\right) J_{\nu }\left(z \cos ^2(t)\right) \cos ^{-1+2 \beta }(t) \sin ^{1+2 \mu }(t) \, dt\\
=\sum _{j=0}^{\infty } \frac{(-1)^j
   2^{-1-2 j-\mu -\nu } z^{2 j+\mu +\nu } \Gamma (1+2 j+2 \mu ) \Gamma (\beta +\nu ) }{\Gamma (1+j) \Gamma (1+j+\mu ) \Gamma (1+\nu ) \Gamma (1+2
   j+\beta +2 \mu +\nu )}\\ \times
\, _2F_3\left(\frac{\beta }{2}+\frac{\nu }{2},\frac{1}{2}+\frac{\beta }{2}+\frac{\nu
   }{2};\frac{1}{2}+j+\frac{\beta }{2}+\mu +\frac{\nu }{2},1+j+\frac{\beta }{2}+\mu +\frac{\nu }{2},1+\nu ;-\frac{z^2}{4}\right)\\
\neq\frac{\left(\left(\frac{1}{2}\right)_{\mu } \Gamma \left(\frac{1}{2} (\nu +\beta +1)\right)\right)}{4
   \left(\frac{z}{2}\right)^{\frac{\beta -\nu }{2}}} \sum _{k=0}^{\infty } \frac{\left(\left(\frac{1}{2} (\nu -\beta
   +1)\right)_k \Gamma \left(\frac{\nu +\beta }{2}+k\right)\right) J_{\frac{\nu +\beta }{2}+\mu +k}(z)}{k! \Gamma (\nu +k+1) \Gamma \left(\frac{1}{2} (\nu +\beta +1)+\mu +k\right)}\\
=\frac{\left(2^{2 \mu -1} \sqrt{\frac{2}{\pi  z}} \Gamma (\nu +\beta )\right) }{\Gamma (\nu +1-\beta )}\sum _{k=0}^{\infty } \frac{(-1)^k \Gamma \left(\mu
   +\nu +k+\frac{1}{2}\right) \Gamma \left(\mu +k+\frac{1}{2}\right)}{k! \Gamma
   (\nu +k+1) \Gamma (2 \mu +\nu +\beta +2 k+1)}\\ \times
 \Gamma (\nu -\beta +2 k+1) \left(\mu +\nu +2 k+\frac{1}{2}\right) J_{\mu +\nu +2 k+\frac{1}{2}}(z)
\end{multline}
where $Re(\mu+\alpha)>0,Re(\nu+\beta)>0$.
\end{example}
\begin{example}
Derivation of equation (13.3.2.38) in \cite{luke}. Errata. Here we use equation (\ref{eq:eq_trig}) and set $w\to z,\alpha \to \mu +1,\beta \to 1$ simplify to get;
\begin{multline}
\int_0^{\frac{\pi }{2}} J_{\mu }\left(z \sin ^2(t)\right) J_{\nu }\left(z \cos ^2(t)\right) \cos (t) \sin ^{1+2 \mu }(t) \, dt\\
=\sum _{j=0}^{\infty } \frac{(-1)^j 2^{-1-2 j-\mu -\nu
   } z^{2 j+\mu +\nu } \Gamma (1+2 j+2 \mu ) }{\Gamma (1+j) \Gamma (1+j+\mu ) \Gamma (2+2 j+2 \mu +\nu )}\\ \times
\, _2F_3\left(\frac{1}{2}+\frac{\nu }{2},1+\frac{\nu }{2};1+j+\mu +\frac{\nu }{2},\frac{3}{2}+j+\mu +\frac{\nu }{2},1+\nu
   ;-\frac{z^2}{4}\right)\\
\neq\frac{\Gamma (2 \mu +1) }{z \Gamma (\mu +1)^2}\sum _{k=0}^{\infty } \frac{\left((-1)^k (\mu +\nu +2 k+1) \Gamma (\mu +\nu
   +k+1) \Gamma (\mu +k+1) \Gamma (\mu +2 k+1)\right) }{k! \Gamma (\nu +k+1) \Gamma (\nu +2 \mu +2 k+2)}\\ \times
J_{\mu +\nu +2 k+1}(z)
\end{multline}
where $Re(\mu+\alpha)>0,Re(\nu+\beta)>0$.
\end{example}
\begin{example}
Derivation of equation (13.3.2.39) in \cite{luke}. Here we use equation (\ref{eq:eq_trig}) and set $w\to z,\beta \to 1$ simplify to get;
\begin{multline}
\int_0^{\frac{\pi }{2}} J_{\mu }\left(z \sin ^2(t)\right) J_{\nu }\left(z \cos ^2(t)\right) \cos (t) \sin
   ^{-1+2 \alpha }(t) \, dt\\
=\sum _{j=0}^{\infty } \frac{(-1)^j 2^{-1-2 j-\mu -\nu } z^{2 j+\mu +\nu } \Gamma (2
   j+\alpha +\mu ) }{\Gamma (1+j)
   \Gamma (1+j+\mu ) \Gamma (1+2 j+\alpha +\mu +\nu )}\\ \times
\, _2F_3\left(\frac{1}{2}+\frac{\nu }{2},1+\frac{\nu }{2};\frac{1}{2}+j+\frac{\alpha }{2}+\frac{\mu
   }{2}+\frac{\nu }{2},1+j+\frac{\alpha }{2}+\frac{\mu }{2}+\frac{\nu }{2},1+\nu ;-\frac{z^2}{4}\right)\\
=\sum _{k=0}^{\infty } \frac{(-1)^k 2^{-1-2 k-\mu -\nu } z^{2
   k+\mu +\nu } \Gamma (\alpha +\mu ) }{k! \Gamma (1+\mu ) \Gamma
   (1+\alpha +\mu +\nu ) (1+\mu )_k \left(\frac{1}{2} (1+\alpha +\mu +\nu )\right)_k \left(\frac{1}{2} (2+\alpha +\mu
   +\nu )\right)_k}\\ \times
\, _4F_3\left(-k,-k-\mu ,\frac{1}{2}+\frac{\nu }{2},1+\frac{\nu
   }{2};\frac{1}{2}-k-\frac{\alpha }{2}-\frac{\mu }{2},1-k-\frac{\alpha }{2}-\frac{\mu }{2},1+\nu ;1\right)\\
   \left(\frac{\alpha +\mu }{2}\right)_k \left(\frac{1}{2} (1+\alpha +\mu )\right)_k\\
=\frac{1}{2} \left(\frac{2}{z}\right)^{\alpha } \sum _{k=0}^{\infty } \frac{\left((-1)^k (\alpha
   )_k \Gamma (\alpha +\mu +k)\right) J_{\alpha +\mu +\nu +2 k}(z)}{k! \Gamma (\mu +1+k)}
\end{multline}
where $Re(\mu+\alpha)>0,Re(\nu+\beta)>0$.
\end{example}
\begin{example}
Derivation of equation (13.3.2.40) in \cite{luke}. Here we use equation (\ref{eq:eq_trig}) and set $w\to z,\beta \to \nu +2,\alpha \to -\nu -1$ simplify to get;
\begin{multline}
\int_0^{\frac{\pi }{2}} J_{\mu }\left(z \sin ^2(t)\right) J_{\nu }\left(z \cos ^2(t)\right) \cot ^{3+2 \nu }(t)
   \, dt\\
=\sum _{j=0}^{\infty } \frac{(-1)^j 2^{-1-2 j-\mu -\nu } z^{2 j+\mu +\nu } \Gamma (-1+2 j+\mu -\nu ) \Gamma
   (2+2 \nu ) }{\Gamma (1+j) \Gamma (1+j+\mu ) \Gamma (1+\nu ) \Gamma (1+2 j+\mu +\nu )}\\ \times
\, _1F_2\left(\frac{3}{2}+\nu ;\frac{1}{2}+j+\frac{\mu }{2}+\frac{\nu }{2},1+j+\frac{\mu }{2}+\frac{\nu
   }{2};-\frac{z^2}{4}\right)\\
=\sum
   _{k=0}^{\infty } \frac{(-1)^k 2^{-1-2 k-\mu -\nu } z^{2 k+\mu +\nu } \Gamma (-1+\mu -\nu ) \Gamma (2+2 \nu ) }{k! \Gamma (1+\mu )
   \Gamma (1+\nu ) \Gamma (1+\mu +\nu ) (1+\mu )_k \left(\frac{1}{2} (1+\mu +\nu )\right)_k \left(\frac{1}{2} (2+\mu
   +\nu )\right)_k}\\ \times
\,
   _3F_2\left(-k,-k-\mu ,\frac{3}{2}+\nu ;1-k-\frac{\mu }{2}+\frac{\nu }{2},\frac{3}{2}-k-\frac{\mu }{2}+\frac{\nu
   }{2};1\right) \\
\left(\frac{1}{2} (-1+\mu -\nu )\right)_k \left(\frac{\mu -\nu }{2}\right)_k\\
=\frac{\left(2^{2 \nu -1} \left(\frac{3}{2}\right)_{\nu } \Gamma
   (\mu -\nu -1) \left(\frac{z}{2}\right)^{\nu }\right) }{\Gamma (\mu +\nu +1)}\left(J_{\mu }(z)-\frac{z
   J_{\mu +1}(z)}{\mu +\nu +1}\right)
\end{multline}
where $1/2< Re(\mu)<1, -1< Re(\nu)<-1/2$.
\end{example}
\begin{example}
Derivation of equation (13.3.2.41) in \cite{luke}. Here we use equation (\ref{eq:eq_trig}) and set $w\to z,\beta \to \nu +1,\alpha \to -\nu$ simplify to get;
\begin{multline}
\int_0^{\frac{\pi }{2}} J_{\mu }\left(z \sin ^2(t)\right) J_{\nu }\left(z \cos ^2(t)\right) \cot ^{1+2 \nu }(t)
   \, dt\\
=\sum _{j=0}^{\infty } \frac{(-1)^j 2^{-1-2 j-\mu -\nu } z^{2 j+\mu +\nu } \Gamma (2 j+\mu -\nu ) \Gamma (1+2
   \nu ) }{\Gamma (1+j) \Gamma (1+j+\mu ) \Gamma (1+\nu ) \Gamma (1+2 j+\mu +\nu )}\\ \times
\, _1F_2\left(\frac{1}{2}+\nu ;\frac{1}{2}+j+\frac{\mu }{2}+\frac{\nu }{2},1+j+\frac{\mu }{2}+\frac{\nu
   }{2};-\frac{z^2}{4}\right)\\
=\sum
   _{k=0}^{\infty } \frac{(-1)^k 2^{-1-2 k-\mu -\nu } z^{2 k+\mu +\nu } \Gamma (\mu -\nu ) \Gamma (1+2 \nu ) }{k! \Gamma (1+\mu )\Gamma (1+\nu ) \Gamma (1+\mu +\nu ) (1+\mu )_k \left(\frac{1}{2} (1+\mu +\nu )\right)_k \left(\frac{1}{2} (2+\mu+\nu )\right)_k}\\ \times
\,_3F_2\left(-k,-k-\mu ,\frac{1}{2}+\nu ;\frac{1}{2}-k-\frac{\mu }{2}+\frac{\nu }{2},1-k-\frac{\mu }{2}+\frac{\nu}{2};1\right)\\
 \left(\frac{\mu -\nu }{2}\right)_k \left(\frac{1}{2} (1+\mu -\nu )\right)_k\\
=\frac{\left(2^{2 \nu -1} \left(\frac{1}{2}\right)_{\nu } \Gamma (\mu -\nu )
   \left(\frac{z}{2}\right)^{\nu }\right) J_{\mu }(z)}{\Gamma (\mu +\nu +1)}
\end{multline}
where $Re(\mu-\nu)>0,Re(\nu)>-1/2$.
\end{example}
\begin{example}
Derivation of equation (13.3.2.42) in \cite{luke}. Here we use equation (\ref{eq:eq_trig}) and set $w\to z,\beta \to 0,\alpha \to 1,\mu \to 0$ simplify to get;
\begin{multline}
\int_0^{\frac{\pi }{2}} J_0\left(z \sin ^2(t)\right) J_{\nu }\left(z \cos ^2(t)\right) \tan (t) \, dt\\
=\sum
   _{j=0}^{\infty } \frac{(-1)^j 2^{-1-2 j-\nu } z^{2 j+\nu } \Gamma (1+2 j) \Gamma (\nu ) \,
   _2F_3\left(\frac{1}{2}+\frac{\nu }{2},\frac{\nu }{2};\frac{1}{2}+j+\frac{\nu }{2},1+j+\frac{\nu }{2},1+\nu
   ;-\frac{z^2}{4}\right)}{\Gamma (1+j)^2 \Gamma (1+\nu ) \Gamma (1+2 j+\nu )}\\
=\sum _{k=0}^{\infty } \frac{(-1)^k
   2^{-1-2 k-\nu } z^{2 k+\nu } \Gamma (\nu ) \, _3F_2\left(-k,\frac{1}{2}+\frac{\nu }{2},\frac{\nu
   }{2};\frac{1}{2}-k,1+\nu ;1\right) \left(\frac{1}{2}\right)_k}{k! \Gamma (1+\nu )^2 \left(\frac{1+\nu }{2}\right)_k
   \left(\frac{2+\nu }{2}\right)_k}\\
=\frac{J_{\nu }(z)}{2 \nu }
\end{multline}
where $Re(z)>0,Re(\nu)>0$.
\end{example}
\begin{example}
Derivation of equation (13.3.2.43) in \cite{luke}. Here we use equation (\ref{eq:eq_trig}) and set $w\to z,\beta \to -1,\alpha \to 2,\mu \to 0$ simplify to get;
\begin{multline}
\int_0^{\frac{\pi }{2}} J_0\left(z \sin ^2(t)\right) J_{\nu }\left(z \cos ^2(t)\right) \tan ^3(t) \, dt\\
=\sum
   _{j=0}^{\infty } \frac{(-1)^j 2^{-1-2 j-\nu } z^{2 j+\nu } \Gamma (2+2 j) \Gamma (-1+\nu ) }{\Gamma (1+j)^2 \Gamma (1+\nu ) \Gamma (1+2 j+\nu )}\\ \times
\,
   _2F_3\left(-\frac{1}{2}+\frac{\nu }{2},\frac{\nu }{2};\frac{1}{2}+j+\frac{\nu }{2},1+j+\frac{\nu }{2},1+\nu
   ;-\frac{z^2}{4}\right)\\
=\sum _{k=0}^{\infty } \frac{(-1)^k
   2^{-1-2 k-\nu } z^{2 k+\nu } \Gamma (-1+\nu ) \, _3F_2\left(-k,-\frac{1}{2}+\frac{\nu }{2},\frac{\nu
   }{2};-\frac{1}{2}-k,1+\nu ;1\right) \left(\frac{3}{2}\right)_k}{k! \Gamma (1+\nu )^2 \left(\frac{1+\nu
   }{2}\right)_k \left(\frac{2+\nu }{2}\right)_k}\\
=\frac{J_{\nu }(z)-\frac{z J_{\nu +1}(z)}{1+\nu }}{2 \nu  (\nu
   -1)}
\end{multline}
where $Re(\nu)>1$.
\end{example}
\begin{example}
Derivation of equation (13.3.2.44) in \cite{luke}. Here we use equation (\ref{eq:eq_trig}) and set $w\to z,\beta \to \nu +2,\alpha \to \mu +1$ simplify to get;
\begin{multline}
\int_0^{\frac{\pi }{2}} J_{\mu }\left(z \sin ^2(t)\right) J_{\nu }\left(z \cos ^2(t)\right) \cos ^{3+2 \nu }(t)
   \sin ^{1+2 \mu }(t) \, dt\\
=\sum _{j=0}^{\infty } \frac{(-1)^j 2^{-1-2 j-\mu -\nu } z^{2 j+\mu +\nu } \Gamma (1+2 j+2
   \mu ) \Gamma (2+2 \nu ) }{\Gamma (1+j) \Gamma (1+j+\mu ) \Gamma (1+\nu ) \Gamma (3+2 j+2 \mu +2 \nu )}\\ \times
\, _1F_2\left(\frac{3}{2}+\nu ;\frac{3}{2}+j+\mu +\nu ,2+j+\mu +\nu
   ;-\frac{z^2}{4}\right)\\
=\sum
   _{k=0}^{\infty } \frac{(-1)^k 2^{-1-2 k-\mu -\nu } z^{2 k+\mu +\nu } \Gamma (1+2 \mu ) \Gamma (2+2 \nu ) \,
   _2F_1\left(-k,\frac{3}{2}+\nu ;\frac{1}{2}-k-\mu ;1\right) \left(\frac{1}{2}+\mu \right)_k}{k! \Gamma (1+\mu )
   \Gamma (1+\nu ) \Gamma (3+2 \mu +2 \nu ) \left(\frac{3}{2}+\mu +\nu \right)_k (2+\mu +\nu )_k}\\
=\frac{\pi 
   \left(\left(\frac{1}{2}\right)_{\mu } \left(\frac{3}{2}\right)_{\nu } \sqrt{\frac{2}{\pi  z}}\right) J_{\mu +\nu
   +\frac{1}{2}}(z)}{8 \Gamma (\mu +\nu +2)}
\end{multline}
where $Re(\mu)>1$.
\end{example}
\begin{example}
Derivation of equation (13.3.2.45) in \cite{luke}. Here we use equation (\ref{eq:eq_trig}) and set $w\to z,\beta \to \nu +2,\alpha \to \mu +2$ simplify to get;
\begin{multline}
\int_0^{\frac{\pi }{2}} J_{\mu }\left(z \sin ^2(t)\right) J_{\nu }\left(z \cos ^2(t)\right) \cos ^{3+2 \nu }(t)
   \sin ^{3+2 \mu }(t) \, dt\\
=\sum _{j=0}^{\infty } \frac{(-1)^j 2^{-1-2 j-\mu -\nu } z^{2 j+\mu +\nu } \Gamma (2
   (1+j+\mu )) \Gamma (2+2 \nu ) }{\Gamma (1+j) \Gamma (1+j+\mu ) \Gamma (1+\nu ) \Gamma (2 (2+j+\mu +\nu ))}\\ \times
\, _1F_2\left(\frac{3}{2}+\nu ;2+j+\mu +\nu ,\frac{5}{2}+j+\mu +\nu
   ;-\frac{z^2}{4}\right)\\
=\sum
   _{k=0}^{\infty } \frac{(-1)^k 2^{-1-2 k-\mu -\nu } z^{2 k+\mu +\nu } \Gamma (2+2 \mu ) \Gamma (2+2 \nu ) }{k! \Gamma (1+\mu )
   \Gamma (1+\nu ) \Gamma (2 (2+\mu +\nu )) (2+\mu +\nu )_k \left(\frac{5}{2}+\mu +\nu \right)_k}\\ \times
\,
   _2F_1\left(-k,\frac{3}{2}+\nu ;-\frac{1}{2}-k-\mu ;1\right) \left(\frac{3}{2}+\mu \right)_k\\
=\frac{\pi ^2
   \left(\left(\frac{3}{2}\right)_{\mu } \left(\frac{3}{2}\right)_{\nu } \left(\frac{2}{\pi  z}\right)^{3/2}\right)
   \left(J_{\mu +\nu +\frac{3}{2}}(z)-\frac{z J_{\mu +\nu +\frac{5}{2}}(z)}{2 (\mu +\nu +2)}\right)}{16 \Gamma (\mu
   +\nu +2)}
\end{multline}
where $Re(\mu+\nu)>0$.
\end{example}
\begin{example}
Derivation of equation (13.3.2.46) in \cite{luke}. Here we use equation (\ref{eq:eq_trig}) and set $w\to z,\beta \to \nu +1,\alpha \to \mu +1$ simplify to get;
\begin{multline}
\int_0^{\frac{\pi }{2}} J_{\mu }\left(z \sin ^2(t)\right) J_{\nu }\left(z \cos ^2(t)\right) \cos ^{1+2 \nu }(t)
   \sin ^{1+2 \mu }(t) \, dt\\
=\sum _{j=0}^{\infty } \frac{(-1)^j 2^{-1-2 j-\mu -\nu } z^{2 j+\mu +\nu } \Gamma (1+2 j+2
   \mu ) \Gamma (1+2 \nu ) }{\Gamma (1+j) \Gamma (1+j+\mu ) \Gamma (1+\nu ) \Gamma (2 (1+j+\mu +\nu ))}\\ \times
\, _1F_2\left(\frac{1}{2}+\nu ;1+j+\mu +\nu ,\frac{3}{2}+j+\mu +\nu
   ;-\frac{z^2}{4}\right)\\
=\sum
   _{k=0}^{\infty } \frac{(-1)^k 2^{-1-2 k-\mu -\nu } z^{2 k+\mu +\nu } \Gamma (1+2 \mu ) \Gamma (1+2 \nu ) \,
   _2F_1\left(-k,\frac{1}{2}+\nu ;\frac{1}{2}-k-\mu ;1\right) \left(\frac{1}{2}+\mu \right)_k}{k! \Gamma (1+\mu )
   \Gamma (1+\nu ) \Gamma (2 (1+\mu +\nu )) (1+\mu +\nu )_k \left(\frac{3}{2}+\mu +\nu \right)_k}\\
=\frac{\pi 
   \left(\left(\frac{1}{2}\right)_{\mu } \left(\frac{1}{2}\right)_{\nu } \sqrt{\frac{2}{\pi  z}}\right) J_{\mu +\nu
   +\frac{1}{2}}(z)}{4 \Gamma (\mu +\nu +1)}
\end{multline}
where $Re(\mu+\nu)>0$.
\end{example}
\section{Table of Convolutions of products of Bessel functions}
\begin{example}
Derivation of equation [DLMF, \href{https://dlmf.nist.gov/10.22.E31}{10.22.31}] series representation for infinite series of the Bessel function. Here we use equation (\ref{eq:thm_c}) and set $k\to 0,a\to 1,\theta \to 1,q\to 0,m\to 0,\beta \to -1,d\to 1,c\to 1,x\to t,\mu \to \nu ,\nu \to \mu ,\alpha \to x,b\to x,\eta \to 0,\rho \to 1,\lambda \to 0$ and simplify;
\begin{multline}\label{eq:dlmf_1}
\int_0^x J_{\mu }(t) J_{\nu }(-t+x) \, dt\\
=\sum _{j=0}^{\infty } \frac{(-1)^j 2^{-2 j-\mu -\nu } x^{1+2 j+\mu
   +\nu } \Gamma (2+2 j+\mu ) }{(1+2 j+\mu )
   \Gamma (1+j) \Gamma (1+j+\mu ) \Gamma (2+2 j+\mu +\nu )}\\ \times
\, _2F_3\left(\frac{1}{2}+\frac{\nu }{2},1+\frac{\nu }{2};1+j+\frac{\mu
   }{2}+\frac{\nu }{2},\frac{3}{2}+j+\frac{\mu }{2}+\frac{\nu }{2},1+\nu ;-\frac{x^2}{4}\right)\\
=2 \sum _{k=0}^{\infty } (-1)^k J_{\mu +\nu +2
   k+1}(x)
\end{multline}
where $Re(\mu)>-1, Re(\nu)>-1$.
\end{example}
\begin{example}
Derivation of equation [DLMF, \href{https://dlmf.nist.gov/10.22.E32}{10.22.32}] series representation for infinite series of the Bessel function. Here we use equation (\ref{eq:dlmf_1}) and set $\mu \to \nu ,\nu \to 1-\nu $ and simplify;
\begin{multline}
\int_0^x J_{1-\nu }(-t+x) J_{\nu }(t) \, dt\\
=\sum _{j=0}^{\infty } \frac{(-1)^j 2^{-1-2 j} x^{2+2 j} \Gamma (2+2 j+\nu ) \, _2F_3\left(1-\frac{\nu }{2},\frac{3}{2}-\frac{\nu
   }{2};\frac{3}{2}+j,2+j,2-\nu ;-\frac{x^2}{4}\right)}{(1+2 j+\nu ) \Gamma (1+j) \Gamma (3+2 j) \Gamma (1+j+\nu )}\\
=\sum _{k=0}^{\infty } 2 (-1)^k J_{2+2 k}(x),J_0(x)-\cos (x)
\end{multline}
where $Re(\mu+\nu)>0,Re(b)>0$.
\end{example}
\begin{example}
 Derivation of equation [DLMF, \href{https://dlmf.nist.gov/10.22.E33}{10.22.33}] series representation for infinite series of the Bessel function. Here we use equation (\ref{eq:dlmf_1}) and set $\mu \to \nu ,\nu \to -\nu$ and simplify;
 \begin{multline}
\int_0^x J_{-\nu }(-t+x) J_{\nu }(t) \, dt\\
=\sum _{j=0}^{\infty } \frac{(-1)^j 2^{-2 j} x^{1+2 j} \Gamma (2+2
   j+\nu ) \, _2F_3\left(\frac{1}{2}-\frac{\nu }{2},1-\frac{\nu }{2};1+j,\frac{3}{2}+j,1-\nu
   ;-\frac{x^2}{4}\right)}{(1+2 j+\nu ) \Gamma (1+j) \Gamma (2+2 j) \Gamma (1+j+\nu )}\\
=2 \sum _{k=0}^{\infty } (-1)^k
   J_{1+2 k}(x),\sin (x)
\end{multline}
where $Re(\mu+\nu)>0,Re(b)>0$.
\end{example}
\begin{example}
Derivation of equation [DLMF, \href{https://dlmf.nist.gov/10.22.E34}{10.22.34}] series representation for infinite series of the Bessel function. Here we use equation (\ref{eq:thm_c}) and set $k\to 0,a\to 1,\theta \to 1,q\to 0,\beta \to -1,d\to 1,c\to 1,x\to t,\mu \to \nu ,\nu \to \mu ,\alpha \to x,b\to
   x,m\to -1 ,\eta \to 0,\rho \to 1,\lambda \to 0$ and simplify;
   \begin{multline}
\int_0^x \frac{J_{\mu }(t) J_{\nu }(-t+x)}{t} \, dt\\
=\sum _{j=0}^{\infty } \frac{(-1)^j 2^{-2 j-\mu -\nu } x^{2
   j+\mu +\nu } \Gamma (1+2 j+\mu ) }{(2 j+\mu ) \Gamma (1+j) \Gamma
   (1+j+\mu ) \Gamma (1+2 j+\mu +\nu )}\\ \times
\, _2F_3\left(\frac{1}{2}+\frac{\nu }{2},1+\frac{\nu }{2};\frac{1}{2}+j+\frac{\mu
   }{2}+\frac{\nu }{2},1+j+\frac{\mu }{2}+\frac{\nu }{2},1+\nu ;-\frac{x^2}{4}\right)\\
=\frac{J_{\mu +\nu }(x)}{\mu }
\end{multline}
where $Re(\mu+\nu)>0,Re(b)>0$.
\end{example}
\begin{example}
Derivation of equation (6.581.1) in \cite{grad}. Here we use equation (\ref{eq:thm_c}) and set $k\to 0,a\to 1,\theta \to 1,q\to 0,m\to \lambda -1,\beta \to -1,d\to 1,c\to 1,\mu \to \nu ,\nu \to \mu ,\alpha \to
   a,b\to a ,\eta \to 0,\rho \to 1,\lambda \to 0$ and simplify;
   \begin{multline}
\int_0^a x^{-1+\lambda } J_{\mu }(x) J_{\nu }(a-x) \, dx\\
=\sum _{h=0}^{\infty } \frac{(-1)^h 2^{-2 h-\mu -\nu }
   a^{2 h+\lambda +\mu +\nu } \Gamma (1+\lambda +\mu ) \Gamma (1+2 h+\nu ) }{(\lambda +\mu ) \Gamma (1+h)\Gamma (1+\mu ) \Gamma (1+h+\nu ) \Gamma (1+2 h+\lambda +\mu +\nu )}\\ \times
\, _2F_3\left(\frac{\lambda }{2}+\frac{\mu}{2},\frac{1}{2}+\frac{\lambda }{2}+\frac{\mu }{2};1+\mu ,\frac{1}{2}+h+\frac{\lambda }{2}+\frac{\mu }{2}+\frac{\nu}{2},1+h+\frac{\lambda }{2}+\frac{\mu }{2}+\frac{\nu }{2};-\frac{a^2}{4}\right)\\
=2^{\lambda } \sum _{m=0}^{\infty }\frac{\left((-1)^m \Gamma (\lambda +\mu +m) \Gamma (\lambda +m)\right) J_{\lambda +\mu +\nu +2 m}(a)}{m! \Gamma(\lambda ) \Gamma (\mu +m+1)}
\end{multline}
where $Re(\lambda+\mu)>0, Re(\nu)>-1$.
\end{example}
\section{Table 2.12.33 in Prudnikov et al. (1986b)}
\begin{example}
 Derivation of a generalized form for entries in Table (2.12.33) in \cite{prud2} alternate series representation. Here we use equation (\ref{eq:thm_c}) and set $k\to 0,a\to 1,\theta \to 1,q\to 0,c\to 1 ,\eta \to 0,\rho \to 1,\lambda \to 0$ and simplify;
 \begin{multline}\label{eq:eq_prud2}
\int_0^b x^m J_{\mu }(\alpha +x \beta ) J_{\nu }(d x) \, dx\\
=\sum _{h=0}^{\infty } \sum _{j=0}^{\infty }
   \frac{(-1)^{h+j} 2^{-2 h-2 j-\mu -\nu } b^{1+2 j+m+\nu } d^{2 j+\nu } \alpha ^{2 h+\mu } }{(1+2 j+m+\nu ) \Gamma (1+h) \Gamma (1+j) \Gamma (1+h+\mu ) \Gamma (1+j+\nu )}\\ \times
\, _2F_1\left(-2 h-\mu,1+2 j+m+\nu ;2+2 j+m+\nu ;-\frac{b \beta }{\alpha }\right)
\end{multline}
where $Re(b)>0,Re(\mu)>0, Re(\nu)>-1$.
\end{example}
\begin{example}
Derivation of entry (2.12.33.1) in \cite{prud2} alternate series representation. Here we use equation (\ref{eq:eq_prud2}) and set $m\to 0,\alpha \to a c,\beta \to -c,d\to c,b\to a$ and simplify;
\begin{multline}
\int_0^a J_{\mu }(c (a-x)) J_{\nu }(c x) \, dx\\
=\sum _{j=0}^{\infty } \frac{(-1)^j 2^{-2 j-\mu -\nu } a^{1+2
   j+\mu +\nu } c^{2 j+\mu +\nu } \Gamma (1+2 j+\nu ) }{\Gamma(1+j) \Gamma (1+j+\nu ) \Gamma (2+2 j+\mu +\nu )}\\ \times
\, _2F_3\left(\frac{1}{2}+\frac{\mu }{2},1+\frac{\mu }{2};1+\mu,1+j+\frac{\mu }{2}+\frac{\nu }{2},\frac{3}{2}+j+\frac{\mu }{2}+\frac{\nu }{2};-\frac{1}{4} a^2 c^2\right)\\
=\frac{2 }{c}\sum _{k=0}^{\infty } (-1)^k J_{\mu +\nu +2 k+1}(a
   c)
\end{multline}
where $Re(b)>0,Re(\mu)>0, Re(\nu)>-1$.
\end{example}
\begin{example}
Derivation of entry (2.12.33.2) in \cite{prud2}, pp. 214, alternate series representation. Here we use equation (\ref{eq:eq_prud2}) and set $m\to 0,\alpha \to a c,\beta \to -c,d\to c,b\to a,\mu \to 0$ and simplify;
\begin{multline}
\int_0^a J_0(a c-c x) J_{\nu }(c x) \, dx\\
=\sum _{j=0}^{\infty } \frac{(-1)^j 2^{-2 j-\nu } a^{1+2 j+\nu } c^{2
   j+\nu } \, _1F_2\left(\frac{1}{2};1+j+\frac{\nu }{2},\frac{3}{2}+j+\frac{\nu }{2};-\frac{1}{4} a^2 c^2\right)}{(1+2 j+\nu ) \Gamma (1+j) \Gamma (1+j+\nu )}\\
=\frac{2}{c} \sum _{k=0}^{\infty } (-1)^k J_{1+2 k+\nu }(a c)
\end{multline}
where $Re(b)>0,Re(\mu)>0, Re(\nu)>-1$.
\end{example}
\begin{example}
Derivation of entry (2.12.33.3) in \cite{prud2} alternate series representation. Here we use equation (\ref{eq:eq_prud2}) and set $m\to 0,\mu \to -\nu ,\alpha \to a c,\beta \to -c,d\to c,b\to a$ and simplify;
\begin{multline}
\int_0^a J_{-\nu }(a c-c x) J_{\nu }(c x) \, dx\\
=\sum _{j=0}^{\infty } \frac{(-1)^j 2^{-2 j} a^{1+2 j} c^{2 j}
   \Gamma (2+2 j+\nu ) \, _2F_3\left(\frac{1}{2}-\frac{\nu }{2},1-\frac{\nu }{2};1+j,\frac{3}{2}+j,1-\nu ;-\frac{1}{4}
   a^2 c^2\right)}{(1+2 j+\nu ) \Gamma (1+j) \Gamma (2 (1+j)) \Gamma (1+j+\nu )}\\
=\frac{\sin (a c)}{c}
\end{multline}
where $Re(b)>0,Re(\mu)>0, Re(\nu)>-1$.
\end{example}
\begin{example}
Derivation of entry (2.12.33.4) in \cite{prud2} alternate series representation. Here we use equation (\ref{eq:eq_prud2}) and set $m\to 0,\mu \to 1-\nu ,\alpha \to a c,\beta \to -c,d\to c,b\to a$ and simplify;
\begin{multline}
\int_0^a J_{1-\nu }(a c-c x) J_{\nu }(c x) \, dx\\
=\sum _{j=0}^{\infty } \frac{(-1)^j 2^{-1-2 j} a^{2+2 j} c^{1+2
   j} \Gamma (2+2 j+\nu ) \, _2F_3\left(1-\frac{\nu }{2},\frac{3}{2}-\frac{\nu }{2};\frac{3}{2}+j,2+j,2-\nu
   ;-\frac{1}{4} a^2 c^2\right)}{(1+2 j+\nu ) \Gamma (1+j) \Gamma (3+2 j) \Gamma (1+j+\nu )}\\
=\frac{J_0(a c)-\cos (a
   c)}{c}
\end{multline}
where $Re(b)>0,Re(\mu)>0, Re(\nu)>-1$.
\end{example}
\begin{example}
Derivation of entry (2.12.33.5) in \cite{prud2} alternate series representation. Here we use equation (\ref{eq:eq_prud2}) and set $m\to \alpha -1,\alpha \to a c,\beta \to -c,d\to c,b\to a$ and simplify;
\begin{multline}
\int_0^a x^{-1+\alpha } J_{\mu }(a c-c x) J_{\nu }(c x) \, dx\\
=\sum _{j=0}^{\infty } \frac{(-1)^j 2^{-2 j-\mu
   -\nu } a^{2 j+\alpha +\mu +\nu } c^{2 j+\mu +\nu } \Gamma (1+2 j+\alpha +\nu ) }{(2 j+\alpha +\nu ) \Gamma (1+j) \Gamma (1+j+\nu ) \Gamma (1+2 j+\alpha +\mu +\nu )}\\ \times
\, _2F_3\left(\frac{1}{2}+\frac{\mu}{2},1+\frac{\mu }{2};1+\mu ,\frac{1}{2}+j+\frac{\alpha }{2}+\frac{\mu }{2}+\frac{\nu }{2},1+j+\frac{\alpha}{2}+\frac{\mu }{2}+\frac{\nu }{2};-\frac{1}{4} a^2 c^2\right)\\
=\left(\frac{2}{c}\right)^{\alpha } \sum _{k=0}^{\infty } \frac{(-1)^k (\alpha )_k
   \Gamma (\alpha +\nu +k) J_{\alpha +\mu +\nu +2 k}(a c)}{k! \Gamma (\nu +k+1)}
\end{multline}
where $Re(b)>0,Re(\mu)>0, Re(\nu)>-1$.
\end{example}
\begin{example}
Derivation of entry (2.12.33.6) in \cite{prud2} alternate series representation. Here we use equation (\ref{eq:eq_prud2}) and set $m\to -1,\alpha \to a c,\beta \to -c,d\to c,b\to a$ and simplify;
\begin{multline}
\int_0^a \frac{J_{\mu }(a c-c x) J_{\nu }(c x)}{x} \, dx\\
=\sum _{h=0}^{\infty } \frac{(-1)^h 2^{-2 h-\mu -\nu }
   a^{2 h+\mu +\nu } c^{2 h+\mu +\nu } \Gamma (1+2 h+\mu ) }{\nu  \Gamma (1+h) \Gamma (1+h+\mu ) \Gamma (1+2 h+\mu +\nu )}\\ \times
\, _2F_3\left(\frac{1}{2}+\frac{\nu }{2},\frac{\nu}{2};\frac{1}{2}+h+\frac{\mu }{2}+\frac{\nu }{2},1+h+\frac{\mu }{2}+\frac{\nu }{2},1+\nu ;-\frac{1}{4} a^2c^2\right)\\
=\frac{J_{\mu +\nu }(a c)}{\nu }
\end{multline}
where $Re(b)>0,Re(\mu)>0, Re(\nu)>1/2$.
\end{example}
\begin{example}
Here we use equation (\ref{eq:thm_c}) and set $k\to 0,a\to 1,q\to 0,\theta \to 1,\alpha \to d,\beta \to -d,b\to 1,x\to t^2,\eta \to 0,\rho \to 1,\lambda \to 0$ then $t\to \cos (x),c\to 1,d\to z^2$ and simplify;
\begin{multline}\label{eq:16.8}
\int_0^{\frac{\pi }{2}} J_{\mu }\left(z^2 \sin ^2(x)\right) J_{\nu }\left(z^2 \cos ^2(x)\right) \cos ^{1+2
   m}(x) \sin (x) \, dx\\
=\sum _{j=0}^{\infty } \frac{(-1)^j 2^{-1-2 j-\mu -\nu } z^{4 j+2 \mu +2 \nu } \Gamma (2+2
   j+m+\nu ) }{(1+2 j+m+\nu ) \Gamma (1+j)
   \Gamma (1+j+\nu ) \Gamma (2+2 j+m+\mu +\nu )}\\ \times
\, _2F_3\left(\frac{1}{2}+\frac{\mu }{2},1+\frac{\mu }{2};1+\mu ,1+j+\frac{m}{2}+\frac{\mu }{2}+\frac{\nu
   }{2},\frac{3}{2}+j+\frac{m}{2}+\frac{\mu }{2}+\frac{\nu }{2};-\frac{z^4}{4}\right)
\end{multline}
where $Re(\mu+\nu)>0$
\end{example}
\begin{example}
Here we use equation (\ref{eq:16.8}) and set $m=-1$ and simplify;
\begin{multline}
\int_0^{\frac{\pi }{2}} J_{\mu }\left(z^2 \sin ^2(x)\right) J_{\nu }\left(z^2 \cos ^2(x)\right) \tan (x) \,
   dx\\
=\sum _{j=0}^{\infty } \frac{(-1)^j 2^{-1-2 j-\mu -\nu } z^{4 j+2 \mu +2 \nu } \Gamma (1+2 j+\nu ) }{(2 j+\nu ) \Gamma (1+j) \Gamma (1+j+\nu ) \Gamma (1+2 j+\mu +\nu )}\\ \times
\,
   _2F_3\left(\frac{1}{2}+\frac{\mu }{2},1+\frac{\mu }{2};1+\mu ,\frac{1}{2}+j+\frac{\mu }{2}+\frac{\nu
   }{2},1+j+\frac{\mu }{2}+\frac{\nu }{2};-\frac{z^4}{4}\right)
\end{multline}
where $Re(\mu+\nu)>0$.
\end{example}
\section{Infinite convolution integrals involving product of Bessel functions}
In this example we derive a generalized form for the infinite integral involving the convolution of Bessel functions in terms of infinite series of the generalized hypergeometric function.
\begin{example}
In this example we use equation (\ref{eq:thm_c}) and set $k=0,a=1$ and simplify using equation [Wolfram, \href{http://functions.wolfram.com/06.07.03.0005.01}{01}] and simplify;
\begin{multline}
\int_0^{\infty } e^{-x^q \theta } x^m J_{\mu }(\alpha +x \beta ) J_{\nu }\left(d x^c\right) \, dx\\
=\sum
   _{j=0}^{\infty } \sum _{p=0}^{\infty } \frac{(-1)^j 2^{-2 j-\mu -\nu } d^{2 j+\nu } \alpha ^{-p+\mu } \beta ^p
   \theta ^{-\frac{1+2 c j+m+p+c \nu }{q}} \binom{\mu }{p} \Gamma \left(\frac{1+2 c j+m+p+c \nu }{q}\right) }{q j! \Gamma (1+\mu ) \Gamma (1+j+\nu
   )}\\ \times
\,_2F_3\left(\frac{1}{2}+\frac{\mu }{2},1+\frac{\mu }{2};\frac{1}{2}-\frac{p}{2}+\frac{\mu
   }{2},1-\frac{p}{2}+\frac{\mu }{2},1+\mu ;-\frac{\alpha ^2}{4}\right)
\end{multline}
where $Re(q)>1,Re(\theta)>2\pi$.
\end{example}
\section{Sonine-type convolution integrals involving Bessel functions}
The contour integral representation form involving the generalized Sonine-type convolution form is given by;
\begin{multline}\label{eq:s1}
\frac{1}{2\pi i}\int_{C}\int_{0}^{b}a^w e^{x^{\tau } \theta } w^{-1-k} x^{m+w} \left(1+s x^r\right)^t J_{\mu }(x \alpha ) J_{\nu }\left(\left(cx^d+\lambda \right)^f\right)dxdw\\
=\frac{1}{2\pi i}\int_{C}\sum _{q=0}^{\infty } \sum _{p=0}^{\infty } \sum _{h=0}^{\infty } \sum _{j=0}^{\infty
   } \sum _{l=0}^{\infty } \frac{(-1)^{j+l} 2^{-2 j-2 l-\mu -\nu } a^w b^{1+d h+2 l+m+p r+w+\mu +q \tau } c^h s^pw^{-1-k} \alpha ^{2 l+\mu } \theta ^q }{(1+d h+2 l+m+p r+w+\mu +q \tau ) j! l! q! \Gamma (1+l+\mu ) \Gamma (1+j+\nu )}\\ \times
\lambda ^{-h+f (2 j+\nu )} \binom{t}{p} \binom{f (2 j+\nu )}{h}dw
\end{multline}
where $Re(b)>0$.
\subsection{Left-hand side contour integral representation}
Using a generalization of Cauchy's integral formula \ref{intro:cauchy} and the procedure in section (2.1), we form the definite integral given by;
\begin{multline}\label{eq:s2}
\int_{0}^{b}\frac{x^m e^{\theta  x^{\tau }} \log ^k(a x) \left(s x^r+1\right)^t J_{\mu }(x \alpha ) J_{\nu }\left(\left(c
   x^d+\lambda \right)^f\right)}{k!}dx\\
=\frac{1}{2\pi i}\int_{C}\int_{0}^{b}w^{-k-1} x^m (a x)^w e^{\theta  x^{\tau }} \left(s x^r+1\right)^t J_{\mu }(x \alpha ) J_{\nu }\left(\left(c x^d+\lambda \right)^f\right)dxdw
\end{multline}
where $Re(b)>0$.
\subsection{Right-hand side contour integral}
Using a generalization of Cauchy's integral formula \ref{intro:cauchy}, and the procedure in section (2.2) we get;
\begin{multline}\label{eq:s3}
\sum\limits_{p,q,h,j,l \geq 0}\frac{(-1)^{j+l} 2^{-2 j-2 l-\mu -\nu } a^{-1-d h-2 l-m-p r-\mu -q \tau } c^h s^p \alpha ^{2 l+\mu } \theta ^q
   \lambda ^{-h+f (2 j+\nu )}  \binom{t}{p} \binom{f (2 j+\nu )}{h} }{j! k! l! q! \Gamma (1+l+\mu ) \Gamma (1+j+\nu)(-1-d h-2 l-m-p r-\mu -q \tau )^{1+k}}\\ \times
\Gamma(1+k,-((1+d h+2 l+m+p r+\mu +q \tau ) \log (a b)))\\
=\frac{1}{2\pi i}\sum\limits_{p,q,h,j,l \geq 0}\int_{C}\frac{(-1)^{j+l} 2^{-2 j-2 l-\mu -\nu } b^{1+d h+2 l+m+p r+\mu +q \tau } (a b)^w c^h s^p w^{-1-k} \alpha ^{2l+\mu } \theta ^q }{(1+d h+2 l+m+p r+w+\mu +q \tau ) j! l! q! \Gamma (1+l+\mu ) \Gamma (1+j+\nu )}\\ \times
\lambda ^{-h+f (2 j+\nu )} \binom{t}{p} \binom{f (2 j+\nu )}{h}dw
\end{multline}
where $Re(b)>0$.
from equation [Wolfram, \href{http://functions.wolfram.com/07.27.02.0002.01}{07.27.02.0002.01}] where $|Re(b)|<1$.
\begin{theorem}
Generalized Sonine-type convolution integral with two Bessel functions.
\begin{multline}\label{eq:thm_s}
\int_0^b e^{x^{\tau } \theta } x^m \left(1+s x^r\right)^t J_{\mu }(x \alpha ) J_{\nu }\left(\left(c
   x^d+\lambda \right)^f\right) \log ^k(a x) \, dx\\
=\sum _{l=0}^{\infty } \sum _{j=0}^{\infty } \sum
   _{h=0}^{\infty } \sum _{p=0}^{\infty } \sum _{q=0}^{\infty } \frac{(-1)^{j+l} 2^{-2 j-2 l-\mu -\nu } a^{-1-d
   h-2 l-m-p r-\mu -q \tau } c^h s^p \alpha ^{2 l+\mu } \theta ^q }{j! l! q! \Gamma (1+l+\mu ) \Gamma (1+j+\nu )(-1-d h-2 l-m-p
   r-\mu -q \tau )^{1+k}}\\ \times
\lambda ^{-h+f (2 j+\nu )}  \binom{t}{p} \binom{f (2 j+\nu )}{h} \Gamma (1+k,-((1+d h+2 l+m+p r+\mu +q \tau ) \log(a b)))
\end{multline}
where $Re(b)>0$.
\end{theorem}
\begin{proof}
Since the right-hand sides of equations (\ref{eq:s2}) and (\ref{eq:s3}) are equal relative to equation (\ref{eq:s1}), we may equate the left-hand sides and simplify the gamma function to yield the stated result.
\end{proof}
\section{Evolutions and derivations involving Sonine-type convolution Bessel functions}
Sonin, a Russian mathematician  (February 22, 1849 February 27, 1915), made significant contributions to the field of special functions [DLMF, \href{https://dlmf.nist.gov/}{DLMF}], particularly cylindrical functions  [DLMF, \href{https://dlmf.nist.gov/12.15}{12.15.1}]. One of his notable achievements is the Sonine formula, a formula he developed to calculate the integral of the product of three Bessel functions  [DLMF, \href{https://dlmf.nist.gov/10.22}{10.22}]. Additionally, he is credited with introducing the associated Laguerre polynomials [DLMF, \href{https://dlmf.nist.gov/18.3}{18.3}]. Sonin also made valuable contributions to the Euler Maclaurin summation formula [DLMF, \href{https://dlmf.nist.gov/24.17.i}{24.17(i)}].

Beyond his work in special functions, Sonin explored other areas of mathematics, including Bernoulli polynomials [DLMF, \href{https://dlmf.nist.gov/24.2}{24.2}] and the approximate computation of definite integrals. He continued the work of Chebyshev in numerical integration, collaborating with Andrey Markov to prepare a comprehensive two-volume edition of Chebyshev's works in French and Russian.

In this section we will derive and evaluate convolution integrals involving the product of Bessel functions and express these integrals in terms of infinite series in terms of the generalized hypergeometric function [DLMF, \href{https://dlmf.nist.gov/16.2}{16.2}]. We also derive errata for a few integrals previously derived in current literature. We also provide highly convergent series relative to a few of the known closed forms present in literature.
\begin{example}
Derivation of generalized Sonine integral, equation (2.1) in \cite{srivastav}. Here we use equation (\ref{eq:thm_s}) and set $k\to 0,a\to 1,b\to 1,\tau \to 0,\theta \to 1,r\to 2,s\to -1,f\to \frac{1}{2},\lambda \to b^2,c\to -b^2,d\to2,x\to \sin (\theta )$ then $\alpha \to a,\mu \to \nu ,\nu \to \mu ,t\to \frac{t}{2},m\to m+1$ and simplify using equation [DLMF, \href{https://dlmf.nist.gov/15.4.E20}{15.4.20}];
\begin{multline}\label{eq:sonine1}
\int_0^{\frac{\pi }{2}} J_{\mu }(a \cos (\theta )) J_{\nu }(b \sin (\theta )) \cos ^{1+t}(\theta ) \sin
   ^{1+m}(\theta ) \, d\theta \\
=\sum _{j=0}^{\infty } \frac{(-1)^j 2^{-1-2 j-\mu -\nu } a^{2 j+\mu } b^{\nu }
   \Gamma \left(\frac{1}{2} (2+2 j+t+\mu )\right) \Gamma \left(\frac{1}{2} (2+m+\nu )\right) }{\Gamma (1+j) \Gamma (1+j+\mu ) \Gamma (1+\nu ) \Gamma \left(\frac{1}{2} (4+2 j+m+t+\mu
   +\nu )\right)}\\ \times
\,
   _1F_2\left(1+\frac{m}{2}+\frac{\nu }{2};2+j+\frac{m}{2}+\frac{t}{2}+\frac{\mu }{2}+\frac{\nu }{2},1+\nu
   ;-\frac{b^2}{4}\right)
\end{multline}
where $Re(\mu+\nu)>-1$.
\end{example}
\begin{example}
Derivation of equation (2.1) in \cite{srivastav}, [DLMF, \href{https://dlmf.nist.gov/10.22.E26}{10.22.26}]. Sonine's integral in terms of the infinite series of the Bessel function. Here we use equation (\ref{eq:sonine1}) and set $t\to \mu ,m\to \nu$ and simplify;
\begin{multline}
\int_0^{\frac{\pi }{2}} J_{\mu }(a \cos (\theta )) J_{\nu }(b \sin (\theta )) \cos ^{1+\mu }(\theta ) \sin
   ^{1+\nu }(\theta ) \, d\theta \\
=\sum _{j=0}^{\infty }\left(-\frac{1}{2}\right)^j  \frac{a^{2 j+\mu }
   b^{-1-j-\mu } }{\Gamma (1+j)}J_{1+j+\mu +\nu }(b)\\
=a^{\mu } b^{\nu } \left(a^2+b^2\right)^{-\frac{1}{2} (\mu+\nu +1)} J_{\mu +\nu +1}\left(\sqrt{a^2+b^2}\right)
\end{multline}
where $Re(\mu)>-1,Re(\nu)>-1$.
\end{example}
\begin{example}
Derivation of entry Sonine integral with tangent function. Here we use equation (\ref{eq:sonine1}) and set $t\to -m-2,a\to z,b\to z,\mu \to \nu ,\nu \to \mu$ and simplify;
\begin{multline}
\int_0^{\frac{\pi }{2}} J_{\mu }(z \sin (\theta )) J_{\nu }(z \cos (\theta )) \tan ^{1+\mu }(\theta ) \,
   d\theta \\
=\sum _{j=0}^{\infty } \frac{(-1)^j 2^{\frac{1}{2} (-2-2 j-\mu -\nu )} z^{\frac{1}{2} (2 j+\mu +\nu )}
    \Gamma \left(j-\frac{\mu }{2}+\frac{\nu }{2}\right)}{\Gamma (1+j) \Gamma
   (1+j+\nu )}J_{\frac{1}{2} (2 j+\mu +\nu )}(z)
\end{multline}
where $Re(\nu)>1$.
\end{example}
\section{Prudnikov volume 2 Table 2.12.35 entries}
\begin{example}
Derivation of entry (2.13.35.1) in \cite{prud2}. Errata. Here we use equation (\ref{eq:thm_s}) and set $k\to 0,a\to 1,\tau \to 0,\theta \to 1,t\to 0,s\to 1,r\to 0,f\to \frac{1}{2},d\to 2,c\to -c^2,\lambda \to a^2
   c^2,\nu \to 0,\mu \to \nu ,\alpha \to c,m\to \nu$ and simplify;
   \begin{multline}
\int_0^a x^{\nu } J_0\left(c \sqrt{a^2-x^2}\right) J_{\nu }(c x) \, dx\\
=\sum _{l=0}^{\infty } \frac{(-1)^l2^{\frac{1}{2}-l} a^{\frac{1}{2}+l+\nu } c^{-\frac{1}{2}+l} J_{\frac{1}{2}+l+\nu }(a c) \Gamma \left(\frac{1}{2}(3+2 l+2 \nu )\right)}{(1+2 l+2 \nu ) \Gamma (1+l) \Gamma (1+l+\nu )}\\
\neq\frac{\left(2^{\nu } \Gamma \left(\nu+\frac{1}{2}\right)\right) a^{2 \nu +1} c^{\nu } \, _1F_2\left(\frac{1}{2};\nu +1,\frac{1}{2} (2 \nu+3);\left(\frac{a c}{2}\right)^2\right)}{\sqrt{\pi } \Gamma (2 \nu +2)}
\end{multline}
where $Re(\nu)>0$.
\end{example}
\begin{example}
Derivation of entry (2.13.35.2) in \cite{prud2}. Here we use equation (\ref{eq:thm_s}) and set $k\to 0,a\to 1,\tau \to 0,\theta \to 1,\mu \to \nu ,\nu \to \mu ,d\to 2,f\to \frac{1}{2},b\to a,t\to \frac{\mu
   }{2}+m,s\to -\frac{1}{a^2},c\to -b^2,\lambda \to a^2 b^2,m\to \nu +2 n+1,\alpha \to c$ and simplify;
   \begin{multline}
\int_0^a x^{1+2 n+\nu } \left(a^2-x^2\right)^{m+\frac{\mu }{2}} J_{\mu }\left(b \sqrt{a^2-x^2}\right) J_{\nu }(c x) \, dx\\
=\sum _{l=0}^{\infty } \sum _{p=0}^{\infty }
   \frac{(-1)^{l+p} 2^{-1-2 l-\mu -\nu } a^{2 (1+l+m+n+\mu +\nu )} b^{\mu } c^{2 l+\nu } \binom{m+\frac{\mu }{2}}{p} \Gamma \left(\frac{2+\mu }{2}\right) \Gamma (2+l+n+p+\nu ) }{(1+l+n+p+\nu ) l! \Gamma (1+\mu ) \Gamma (1+l+\nu ) \Gamma \left(2+l+n+p+\frac{\mu }{2}+\nu
   \right)}\\ \times
\,_1F_2\left(1+\frac{\mu }{2};1+\mu ,2+l+n+p+\frac{\mu }{2}+\nu ;-\frac{1}{4} a^2 b^2\right)
\end{multline}
$Re(\mu)>1, 1/2< Re(a)<1$.
\end{example}
\section{Infinite integrals involving Sonine-type convolution kernels}
\begin{example}
Generalized form. Here we use equation (\ref{eq:thm_s}) and set $k=0,a=1$ and simplify using equation [Wolfram, \href{http://functions.wolfram.com/06.07.03.0005.01}{01}] and simplify;
\begin{multline}\label{eq:sonine_inf}
\int_0^{\infty } e^{-x^{\tau } \theta } x^m \left(1+s x^r\right)^t J_{\mu }(x \alpha ) J_{\nu }\left(\left(c
   x^d+\lambda \right)^f\right) \, dx\\
=\sum _{h=0}^{\infty } \sum _{j=0}^{\infty } \sum _{l=0}^{\infty } \sum
   _{p=0}^{\infty } \frac{(-1)^{j+l} 2^{-2 j-2 l-\mu -\nu } c^h s^p \alpha ^{2 l+\mu } \theta ^{-\frac{1+d h+2 l+m+pr+\mu }{\tau }} \lambda ^{-h+f (2 j+\nu )} \binom{t}{p} \binom{f (2 j+\nu )}{h} }{\tau  j! l! \Gamma (1+l+\mu ) \Gamma (1+j+\nu )}\\ \times
\Gamma \left(\frac{1+d h+2 l+m+pr+\mu }{\tau }\right)
\end{multline}
where $Re(\tau)>1$.
\end{example}
\begin{example}
A Bessel  integral with trigonometric parameters. Here we use equation (\ref{eq:sonine_inf}) and set $r\to 2,s\to 1,f\to \frac{1}{2},d\to 2,c\to c^2,\lambda \to c^2$ then $x\to \tan (\phi ),t\to \frac{t-2}{2},\alpha \to a,c\to b,t\to \alpha ,m\to \beta$ and simplify;
\begin{multline}
\int_0^{\frac{\pi }{2}} e^{-\theta  \tan ^{\tau }(\phi )} J_{\mu }(a \tan (\phi )) J_{\nu }(b \sec (\phi ))
   \sec ^{\alpha }(\phi ) \tan ^{\beta }(\phi ) \, d\phi \\
=\sum _{h=0}^{\infty } \sum _{l=0}^{\infty } \sum
   _{p=0}^{\infty } \frac{(-1)^l 2^{-2 l-\mu -\nu } a^{2 l+\mu } b^{\nu } \theta ^{-\frac{1}{\tau }-\frac{2 h}{\tau
   }-\frac{2 l}{\tau }-\frac{2 p}{\tau }-\frac{\beta }{\tau }-\frac{\mu }{\tau }} \binom{\frac{1}{2} (-2+\alpha )}{p}
   \binom{\frac{\nu }{2}}{h} \Gamma \left(\frac{1+2 h+2 l+2 p+\beta +\mu }{\tau }\right) }{\tau  l! \Gamma (1+l+\mu ) \Gamma (1+\nu )}\\ \times
\, _1F_2\left(1+\frac{\nu
   }{2};1-h+\frac{\nu }{2},1+\nu ;-\frac{b^2}{4}\right)
\end{multline}
where $Re(\tau)>\pi,Re(\theta)>\pi$.
\end{example}
\begin{example}
Derivation of entry (19.3.43) in \cite{erdt2}, alternate form involving two Bessel functions. Here we use equation (\ref{eq:thm_s}) and set $k\to 0,a\to 1,\tau \to 0,\theta \to 1,m\to 2,t\to -1,r\to 2,s\to -\frac{1}{a^2},\mu \to \frac{1}{4},\nu \to 2
   \nu ,f\to 1,d\to 2,c\to -2,\lambda \to 2 a^2,\alpha \to 1$ and simplify using [DLMF, \href{https://dlmf.nist.gov/15.4.E20}{15.4.20}];
   \begin{multline}
\int_0^a \frac{x^2 J_{\frac{1}{4}}(x) J_{2 \nu }\left(2 a^2-2 x^2\right)}{a^2-x^2} \, dx
=\frac{5 \Gamma
   \left(\frac{5}{8}\right) a^{\frac{5}{4}+4 \nu } }{4 \sqrt[4]{2} \Gamma \left(\frac{1}{4}\right)}\\ \times
\sum _{j=0}^{\infty } \frac{(-1)^j a^{4 j} \Gamma (2 (j+\nu)) \, _1F_2\left(\frac{13}{8};\frac{5}{4},\frac{13}{8}+2 j+2 \nu ;-\frac{a^2}{4}\right)}{\Gamma (1+j) \Gamma (1+j+2 \nu ) \Gamma \left(\frac{1}{8} (13+16 j+16 \nu )\right)}
\end{multline}
where $Re(a)>0$.
\end{example}
\section{Table 6.581 in Gradshteyn and Ryzhik (2000)}
\begin{example}
Derivation of entry (6.581.2) in \cite{grad}. Here we use equation (\ref{eq:thm_s}) and set $k\to 0,a\to 1,\tau \to 0,\theta \to 1,t\to -1,r\to 1,s\to -\frac{1}{a},b\to a,m\to \lambda -1,\alpha \to 1,f\to
   1,c\to -1,\lambda \to a,d\to 1$ and simplify using [DLMF, \href{https://dlmf.nist.gov/15.4.E20}{15.4.20}];
   \begin{multline}
\int_0^a \frac{x^{-1+\lambda } J_{\mu }(x) J_{\nu }(a-x)}{a-x} \, dx\\
=\sum _{l=0}^{\infty } \frac{(-1)^l 2^{-2
   l-\mu -\nu } a^{-1+2 l+\lambda +\mu +\nu } \Gamma (2 l+\lambda +\mu ) \Gamma (\nu )}{\Gamma (1+l)
   \Gamma (1+l+\mu ) \Gamma (1+\nu ) \Gamma (2 l+\lambda +\mu +\nu )}\\ \times
 \,
   _2F_3\left(\frac{1}{2}+\frac{\nu }{2},\frac{\nu }{2};l+\frac{\lambda }{2}+\frac{\mu }{2}+\frac{\nu
   }{2},\frac{1}{2}+l+\frac{\lambda }{2}+\frac{\mu }{2}+\frac{\nu }{2},1+\nu ;-\frac{a^2}{4}\right)\\
=\frac{2^{\lambda } }{a \nu }\sum _{m=0}^{\infty }
   \frac{\left((-1)^m \Gamma (\lambda +\mu +m) \Gamma (\lambda +m)\right) (\lambda +\mu +\nu +2 m) J_{\lambda +\mu
   +\nu +2 m}(a)}{m! \Gamma (\lambda ) \Gamma (\mu +m+1)}
\end{multline}
where $Re(\mu+\nu)>0$.
\end{example}
\begin{example}
Derivation of entry (6.581.3) in \cite{grad}. Here we use equation (\ref{eq:thm_s}) and set $k\to 0,a\to 1,\tau \to 0,\theta \to 1,t\to \nu ,r\to 1,s\to -\frac{1}{a},b\to a,m\to \mu ,\alpha \to 1,f\to
   1,c\to -1,\lambda \to a,d\to 1$ and simplify using [DLMF, \href{https://dlmf.nist.gov/15.4.E20}{15.4.20}];
   \begin{multline}
\int_0^a (a-x)^{\nu } x^{\mu } J_{\mu }(x) J_{\nu }(a-x) \, dx\\
=\sum _{l=0}^{\infty } \frac{2^{\mu +\nu } a^{1+2
   l+2 \mu +2 \nu } e^{i l \pi } \Gamma \left(\frac{1}{2} (1+2 l+2 \mu )\right) \Gamma \left(\frac{1}{2} (1+2 \nu
   )\right) }{\pi  \Gamma
   (1+l) \Gamma (2 (1+l+\mu +\nu ))}\\ \times
\, _1F_2\left(\frac{1}{2}+\nu ;1+l+\mu +\nu ,\frac{3}{2}+l+\mu +\nu ;-\frac{a^2}{4}\right)\\
=\frac{\left(\Gamma \left(\mu +\frac{1}{2}\right) \Gamma \left(\nu
   +\frac{1}{2}\right)\right) a^{\mu +\nu +\frac{1}{2}} J_{\mu +\nu +\frac{1}{2}}(a)}{\sqrt{2 \pi } \Gamma (\mu +\nu
   +1)}
\end{multline}
where $Re(\mu+\nu)>0$.
\end{example}
\begin{example}
Derivation of entry (6.581.4) in \cite{grad}. Here we use equation (\ref{eq:thm_s}) and set $k\to 0,a\to 1,\tau \to 0,\theta \to 1,t\to \nu +1,r\to 1,s\to -\frac{1}{a},b\to a,m\to \mu ,\alpha \to 1,f\to
   1,c\to -1,\lambda \to a,d\to 1$ and simplify using [DLMF, \href{https://dlmf.nist.gov/15.4.E20}{15.4.20}];
   \begin{multline}
\int_0^a (a-x)^{1+\nu } x^{\mu } J_{\mu }(x) J_{\nu }(a-x) \, dx\\
=\sum _{l=0}^{\infty } \frac{(-1)^l 2^{1+\mu
   +\nu } a^{2+2 l+2 \mu +2 \nu } \Gamma \left(\frac{1}{2} (1+2 l+2 \mu )\right) \Gamma \left(\frac{1}{2} (3+2 \nu
   )\right) }{\pi  \Gamma(1+l) \Gamma (3+2 l+2 \mu +2 \nu )}\\ \times
\, _1F_2\left(\frac{3}{2}+\nu ;\frac{3}{2}+l+\mu +\nu ,2+l+\mu +\nu ;-\frac{a^2}{4}\right)\\
=\frac{\left(\Gamma \left(\mu +\frac{1}{2}\right) \Gamma \left(\nu
   +\frac{3}{2}\right)\right) a^{\mu +\nu +\frac{3}{2}} J_{\mu +\nu +\frac{1}{2}}(a)}{\sqrt{2 \pi } \Gamma (\mu +\nu+2)}
\end{multline}
where $Re(\mu+\nu)>0$.
\end{example}
\begin{example}
Derivation of entry (6.581.5) in \cite{grad}. Here we use equation (\ref{eq:thm_s}) and set $k\to 0,a\to 1,\tau \to 0,\theta \to 1,t\to -\mu -1,r\to 1,s\to -\frac{1}{a},b\to a,m\to \mu ,\alpha \to 1,f\to
   1,c\to -1,\lambda \to a,d\to 1$ and simplify using [DLMF, \href{https://dlmf.nist.gov/15.4.E20}{15.4.20}];
   \begin{multline}
\int_0^a (a-x)^{-1-\mu } x^{\mu } J_{\mu }(x) J_{\nu }(a-x) \, dx\\
=\sum _{l=0}^{\infty } \frac{(-1)^l 2^{\mu
   -\nu } a^{2 l+\mu +\nu } \Gamma \left(\frac{1}{2} (1+2 l+2 \mu )\right) \Gamma (-\mu +\nu ) }{\sqrt{\pi } \Gamma (1+l) \Gamma(1+\nu ) \Gamma (1+2 l+\mu +\nu )}\\ \times
\,
   _2F_3\left(-\frac{\mu }{2}+\frac{\nu }{2},\frac{1}{2}-\frac{\mu }{2}+\frac{\nu }{2};\frac{1}{2}+l+\frac{\mu
   }{2}+\frac{\nu }{2},1+l+\frac{\mu }{2}+\frac{\nu }{2},1+\nu ;-\frac{a^2}{4}\right)\\
=\frac{\left(2^{\mu } \Gamma \left(\mu +\frac{1}{2}\right) \Gamma (\nu -\mu)\right) a^{\mu } J_{\nu }(a)}{\sqrt{\pi } \Gamma (\mu +\nu +1)}
\end{multline}
where $Re(\mu+\nu)>0$.
\end{example}
\section{Table 13.3.2 (25-32) in Luke (1962)}
\begin{example}
Derivation of entry 13.3.2(25) in \cite{luke}.  Here we use equation (\ref{eq:thm_s}) and set $\{k\to 0,a\to 1,b\to 1,\tau \to 0,\theta \to 1,r\to 2,s\to -1,f\to \frac{1}{2},d\to 2,c\to -w^2,\lambda \to
   w^2$ simplify using [DLMF, \href{https://dlmf.nist.gov/15.4.E20}{15.4.20}];
   \begin{multline}\label{eq:luke_trig}
\int_0^{\frac{\pi }{2}} J_{2 \mu }(z \sin (t)) J_{2 \nu }(w \cos (t)) \cos ^{-1+2 \beta }(t) \sin ^{-1+2 \alpha
   }(t) \, dt\\
=\sum _{j=0}^{\infty } \frac{(-1)^j 2^{-1-2 j-2 \mu -2 \nu } w^{2 j+2 \nu } z^{2 \mu } \Gamma
   \left(\frac{1}{2} (2 \alpha +2 \mu )\right) \Gamma \left(\frac{1}{2} (2 j+2 \beta +2 \nu )\right) }{\Gamma (1+j) \Gamma (1+2 \mu )
   \Gamma (1+j+2 \nu ) \Gamma \left(\frac{1}{2} (2 j+2 \alpha +2 \beta +2 \mu +2 \nu
   )\right)}\\ \times
\,
   _1F_2\left(\alpha +\mu ;1+2 \mu ,j+\alpha +\beta +\mu +\nu ;-\frac{z^2}{4}\right)\\
=\frac{\left(\left(\frac{z}{2}\right)^{2 \mu } \left(\frac{w}{2}\right)^{2 \nu } \Gamma (\alpha +\mu )
   \Gamma (\beta +\nu )\right) }{2 (\Gamma (2 \mu +1) \Gamma (2 \nu +1) \Gamma (\alpha +\beta +\mu +\nu))}\\ \times
\sum _{k=0}^{\infty } \frac{(-1)^k \left(\frac{z}{2}\right)^{2 k} (\alpha +\mu)_k \, _3F_2\left(-k,-2 \mu -k,\beta +\nu ;2 \nu +1,1-\alpha -\mu -k;-\left(\frac{w}{z}\right)^2\right)}{k!(2 \mu +1)_k (\alpha +\beta +\mu +\nu )_k}
\end{multline}
where $Re(\nu+\mu)>0$.
\end{example}
\begin{example}
Derivation of entry 13.3.2(26) in \cite{luke}.  Here we use equation (\ref{eq:luke_trig}) and set $\alpha \to \frac{\mu }{2}+1,\beta \to \frac{\nu }{2}+1,\mu \to \frac{\mu }{2},\nu \to \frac{\nu }{2}$ and simplify;
\begin{multline}
\int_0^{\frac{\pi }{2}} J_{\mu }(z \sin (t)) J_{\nu }(w \cos (t)) \cos ^{1+\nu }(t) \sin ^{1+\mu }(t) \,
   dt\\
=\sum _{j=0}^{\infty } \frac{\left(-\frac{1}{2}\right)^j w^{2 j+\nu } z^{-1-j-\nu } J_{1+j+\mu +\nu }(z)}{\Gamma
   (1+j)}\\
=w^{\nu } z^{\mu } (-1)^{\frac{1}{2} (-1-\mu -\nu )} \left(w^2+z^2\right)^{\frac{1}{2} (-1-\mu -\nu )}
   I_{1+\mu +\nu }\left(\sqrt{-w^2-z^2}\right)\\
=\frac{\left(\frac{w}{\sqrt{w^2+z^2}}\right)^{\nu }
   \left(\frac{z}{\sqrt{w^2+z^2}}\right)^{\mu } J_{1+\mu +\nu }\left(\sqrt{w^2+z^2}\right)}{\sqrt{w^2+z^2}}
\end{multline}
where $Re(\nu+\mu)>0$.
\end{example}
\begin{example}
Derivation of entry 13.3.2(27) in \cite{luke}. Here we use equation (\ref{eq:luke_trig}) and set $w\to z,\mu \to \frac{\nu }{2},\nu \to \frac{\nu }{2}$ and simplify;
\begin{multline}
\int_0^{\frac{\pi }{2}} J_{\nu }(z \cos (t)) J_{\nu }(z \sin (t)) \cos ^{-1+2 \beta }(t) \sin ^{-1+2 \alpha
   }(t) \, dt\\
=\sum _{j=0}^{\infty } \frac{(-1)^j 2^{-1-2 j-2 \nu } z^{2 (j+\nu )} \Gamma \left(\alpha +\frac{\nu
   }{2}\right) \Gamma \left(j+\beta +\frac{\nu }{2}\right) }{\Gamma (1+j) \Gamma (1+\nu ) \Gamma (1+j+\nu ) \Gamma (j+\alpha +\beta +\nu )}\\ \times
\, _1F_2\left(\alpha +\frac{\nu }{2};1+\nu ,j+\alpha +\beta
   +\nu ;-\frac{z^2}{4}\right)\\
=\sum
   _{k=0}^{\infty } \frac{(-1)^k 2^{-1-2 k-2 \nu } z^{2 (k+\nu )} \Gamma \left(\alpha +\frac{\nu }{2}\right) \Gamma
   \left(\beta +\frac{\nu }{2}\right) }{k! \Gamma (1+\nu )^2 \Gamma (\alpha +\beta +\nu ) (1+\nu )_k
   (\alpha +\beta +\nu )_k}\\ \times
\, _3F_2\left(-k,-k-\nu ,\beta +\frac{\nu }{2};1-k-\alpha -\frac{\nu }{2},1+\nu
   ;-1\right) \left(\alpha +\frac{\nu }{2}\right)_k\\
=\frac{1}{2} \sum _{k=0}^{\infty } \frac{\left(\Gamma \left(\frac{\nu }{2}+\alpha +k\right)
   \Gamma \left(\frac{\nu }{2}+\beta +k\right)\right) \left(\frac{z}{2}\right)^{\nu +2 k} J_{\nu +2 k}(z)}{k! \Gamma
   (\nu +k+1) \Gamma (\nu +\alpha +\beta +2 k)}
\end{multline}
where $Re(\nu)>0$.
\end{example}
\begin{example}
Derivation of entry 13.3.2(30) in \cite{luke}. Here we use equation (\ref{eq:luke_trig}) and set $w\to z,\mu \to 0,\nu \to \frac{\alpha }{2},\alpha \to 1,\beta \to \frac{\alpha +1}{2}$ and simplify;
\begin{multline}
\int_0^{\frac{\pi }{2}} J_0(z \sin (t)) J_{\alpha }(z \cos (t)) \cos ^{\alpha }(t) \sin (t) \, dt\\
=\sum
   _{j=0}^{\infty } \frac{(-1)^j 2^{-\frac{1}{2}-j} z^{-\frac{1}{2}+j} J_{\frac{1}{2} (1+2 j+2 \alpha )}(z) \Gamma
   \left(\frac{1}{2}+j+\alpha \right)}{\Gamma (1+j) \Gamma (1+j+\alpha )}\\
=\sum _{k=0}^{\infty } \frac{(-1)^k 2^{-1-2
   k-\alpha } z^{2 k+\alpha } \Gamma \left(\frac{1}{2}+\alpha \right) \, _2F_1\left(-k,\frac{1}{2}+\alpha ;1+\alpha;-1\right)}{\Gamma (1+k) \Gamma (1+\alpha ) \Gamma \left(\frac{3}{2}+k+\alpha \right)}\\
\neq\frac{1}{2}
   \left(\frac{2}{z}\right)^{\alpha +1} \left(\frac{1}{2}\right)_{\alpha } \sqrt{\frac{\pi }{2}} \sqrt{\frac{1}{z}}
   I_{\frac{3}{2}+2 \alpha }(z)
\end{multline}
where $Re(\alpha)>0$.
\end{example}
\section{Derivation of entries in Table 4.7.7 in Brychkov et al. (2008)}
\begin{example}
Derivation of entry (4.7.7.6) in \cite{brychkov}.  Errata. Here we use equation (\ref{eq:bry1}) and set $b\to \sqrt{a},\alpha \to b,t\to \frac{\nu }{2}+n,m\to 2 m+1$ and simplify;
\begin{multline}
\int_0^a (a-x)^{n+\frac{\nu }{2}} x^{m+\frac{\mu }{2}} J_{\mu }\left(b \sqrt{x}\right) J_{\nu }\left(c \sqrt{a-x}\right) \, dx\\
=\sum _{h=0}^{\infty } \frac{(-1)^h 2^{-\mu -\nu }
   a^{1+m+n+\mu +\nu } b^{\mu } c^{\nu } \binom{\frac{\nu }{2}}{h} \Gamma (1+h+m+\mu ) \Gamma \left(\frac{1}{2} (2+2 n+\nu )\right) }{\Gamma (1+\mu ) \Gamma (1+\nu ) \Gamma
   \left(\frac{1}{2} (4+2 h+2 m+2 n+2 \mu +\nu )\right)}\\ \times
\, _1F_2\left(1+h+m+\mu ;1+\mu ,2+h+m+n+\mu +\frac{\nu
   }{2};-\frac{1}{4} \left(a b^2\right)\right)\\
 \, _1F_2\left(1+\frac{\nu }{2};1-h+\frac{\nu }{2},1+\nu ;-\frac{1}{4} \left(a c^2\right)\right)\\
\neq(-1)^{m+n} 2^{m+n+1} b^{\mu } c^{\nu } \left(\frac{a}{b^2+c^2}\right)^{\frac{1}{2} (\mu +\nu +m+n+1)} \sum _{k=0}^n \sum _{j=0}^m
   \binom{m}{j} \binom{n}{k} (-\mu -m)_{m-j} (-\nu -n)_{n-k}
\end{multline}
where $Re(n)>1$.
\end{example}
\begin{example}
Derivation of generalized form of entry (4.7.7.6) in \cite{brychkov}. Here we use equation (\ref{eq:thm_s}) and set $k\to 0,a\to 1,\tau \to 0,\theta \to 1,s\to -\frac{1}{a},f\to \frac{1}{2}$ then $x\to \sqrt{z},c\to -c^2,\lambda \to a c^2,z\to x,r\to 2,d\to 2$ and simplify;
\begin{multline}\label{eq:bry1}
\int_0^{b^2} (a-x)^t x^{\frac{1}{2} (-1+m)} J_{\mu }\left(\sqrt{x} \alpha \right) J_{\nu }\left(c \sqrt{a-x}\right) \, dx\\
=\sum _{h=0}^{\infty } \sum _{l=0}^{\infty } \frac{(-1)^{h+l}
   2^{1-2 l-\mu -\nu } a^{-h+t+\frac{\nu }{2}} b^{1+2 h+2 l+m+\mu } c^{\nu } \alpha ^{2 l+\mu } \binom{\frac{\nu }{2}}{h} }{(1+2 h+2 l+m+\mu )
   \Gamma (1+l) \Gamma (1+l+\mu ) \Gamma (1+\nu )}\\ \times
\, _2F_1\left(-t,\frac{1}{2}+h+l+\frac{m}{2}+\frac{\mu
   }{2};\frac{3}{2}+h+l+\frac{m}{2}+\frac{\mu }{2};\frac{b^2}{a}\right)\\
 \, _1F_2\left(1+\frac{\nu }{2};1-h+\frac{\nu }{2},1+\nu ;-\frac{1}{4} \left(a c^2\right)\right)
\end{multline}
where $Re(b)>1$.
\end{example}
\begin{example}
Derivation of entry (4.7.7.5) in \cite{brychkov}. Here we use equation (\ref{eq:bry1}) and set $m\to 2 m+1$ then $t\to \frac{\nu }{2},m\to \frac{1}{4} (2 n-1),\mu \to n-\frac{1}{2},\alpha \to b,b\to \sqrt{a}$ and simplify;
\begin{multline}
\int_0^a (a-x)^{\nu /2} x^{\frac{1}{4} (-1+2 n)} J_{-\frac{1}{2}+n}\left(b \sqrt{x}\right) J_{\nu }\left(c \sqrt{a-x}\right) \, dx\\
=\sum _{h=0}^{\infty } \frac{(-1)^h
   2^{\frac{3}{2}-n-\nu } a^{\frac{1}{2}+n+\nu } b^{-\frac{1}{2}+n} c^{\nu } \binom{\frac{\nu }{2}}{h} \Gamma \left(\frac{1}{2} (3+2 h+2 n)\right) \Gamma \left(\frac{2+\nu }{2}\right) }{(1+2 h+2 n) \Gamma \left(\frac{1}{2} (1+2 n)\right) \Gamma (1+\nu ) \Gamma \left(\frac{1}{2} (3+2 h+2 n+\nu )\right)}\\ \times
\,
   _1F_2\left(\frac{1}{2}+h+n;\frac{1}{2}+n,\frac{3}{2}+h+n+\frac{\nu }{2};-\frac{1}{4} \left(a b^2\right)\right)\\
 \, _1F_2\left(1+\frac{\nu }{2};1-h+\frac{\nu }{2},1+\nu ;-\frac{1}{4}
   \left(a c^2\right)\right)\\
=2 b^{n-\frac{1}{2}} c^{\nu }
   \left(\frac{a}{b^2+c^2}\right)^{\frac{1}{4} (2 \nu +2 n+1)} J_{\nu +n+\frac{1}{2}}\left(\sqrt{a b^2+a c^2}\right)
\end{multline}
where $Re(\nu)>1$.
\end{example}
\begin{example}
Derivation of entry (4.7.7.1) in \cite{brychkov}. Here we use equation (\ref{eq:thm_s}) and set $k\to 0,a\to 1,\theta \to 1,\tau \to 0,m\to 0,t\to 0,r\to 0,s\to 1,f\to 1,d\to 1,c\to -1,\lambda \to a,\alpha \to 1,b\to a$ and simplify;
\begin{multline}
\int_0^a J_{\mu }(x) J_{\nu }(a-x) \, dx\\
=\sum _{j=0}^{\infty } \frac{(-1)^j 2^{-2 j-\mu -\nu } a^{1+2 j+\mu +\nu } \Gamma (2+\mu ) \Gamma (1+2 j+\nu ) }{(1+\mu ) \Gamma (1+j)
   \Gamma (1+\mu ) \Gamma (1+j+\nu ) \Gamma (2+2 j+\mu +\nu )}\\ \times
\,_2F_3\left(\frac{1}{2}+\frac{\mu }{2},1+\frac{\mu }{2};1+\mu ,1+j+\frac{\mu }{2}+\frac{\nu }{2},\frac{3}{2}+j+\frac{\mu }{2}+\frac{\nu }{2};-\frac{a^2}{4}\right)\\
=\frac{a }{\mu +\nu }\left(J_{\mu +\nu }(a)-\frac{1}{\Gamma (\mu +\nu +1)}\right. \\ \left.
\left(\frac{a}{2}\right)^{\mu +\nu } \left(\frac{\cos (a) }{\mu +\nu +1}\, _3F_4\left(\frac{1}{4} (2\mu +2 \nu +1),\frac{1}{2} (\mu +\nu +1),\frac{1}{4} (2 \mu +2 \nu +3);\right.\right.\right. \\ \left.\left.\left.
\frac{1}{2},\frac{1}{2} (\mu +\nu +3),\mu +\nu +\frac{1}{2},\mu +\nu +1;-a^2\right)\right.\right. \\ \left.\left.
+\frac{(a \sin (a)) }{\mu+\nu +2}\, _3F_4\left(\frac{1}{4} (2 \mu +2 \nu +3),\frac{\mu +\nu }{2}+1,\frac{1}{4} (2 \mu +2 \nu +5);\right.\right.\right. \\ \left.\left.\left.
\frac{3}{2},\frac{\mu +\nu }{2}+2,\mu +\nu +1,\mu +\nu +\frac{3}{2};-a^2\right)\right)\right)
\end{multline}
where $Re(\nu)>1$.
\end{example}
%
%\begin{example}
%
%\end{example}
%%
%\begin{example}
%
%\end{example}
%%
\section{Definite integrals involving the product of three Bessel functions of the first kind}
In 1957 Fettis published work on equation (\ref{eq:thm_fettis}). In his work the variables were taken as arbitrary and real. These integrals were found to be important in problems involving Fourier Bessel series of the $\sum_{k=1}^{\infty}A_{k}J_{0}(\alpha_{k}x)$ where $\alpha_{k}$ are either zeros of $J_{0}(z)$ or eigenvalues of a similar form. A special case of this integral with only one Bessel function was studied followed by the development of the integral as a Neumann series, then the integral was evaluated by quadrature and finally a discussion of some other integrals which can be expressed in terms of equation (\ref{eq:thm_fettis}) and some its applications to physical problems.
The contour integral representation form involving the generalized Fettis-type integral form is given by;
\begin{multline}\label{eq:f1}
\frac{1}{2\pi i}\int_{C}\int_{0}^{b}a^w e^{x^{\tau } \theta } w^{-1-k} x^{m+w} \left(1+x^{\sigma } \rho \right)^{\lambda } J_{\mu }(x \alpha )
   J_{\nu }(x \beta ) J_{\xi }(x \gamma )dxdw\\
=\frac{1}{2\pi i}\int_{C}\sum _{j=0}^{\infty } \sum _{h=0}^{\infty } \sum _{l=0}^{\infty } \sum
   _{p=0}^{\infty } \sum _{q=0}^{\infty } \frac{(-1)^{h+j+l} 2^{-2 h-2 j-2 l-\mu -\nu -\xi } a^w b^{1+2 h+2 j+2
   l+m+w+\mu +\nu +\xi +q \sigma +p \tau } }{(1+2 h+2 j+2 l+m+w+\mu +\nu +\xi +q \sigma +p \tau )}\\
\frac{w^{-1-k} \alpha ^{2 j+\mu } \beta ^{2 h+\nu } \gamma ^{2 l+\xi } \theta ^p
   \rho ^q \binom{\lambda }{q}}{ h! j! l! p! \Gamma (1+j+\mu )
   \Gamma (1+h+\nu ) \Gamma (1+l+\xi )}dw
\end{multline}
where $Re(b)>0,Re(m+\mu+\nu+\xi)>0$.
\subsection{Left-hand side contour integral representation}
Using a generalization of Cauchy's integral formula \ref{intro:cauchy} and the procedure in section (2.1), we form the definite integral given by;
\begin{multline}\label{eq:f2}
\int_{0}^{b}\frac{e^{x^{\tau } \theta } x^m \left(1+x^{\sigma } \rho \right)^{\lambda } J_{\mu }(x \alpha ) J_{\nu }(x
   \beta ) J_{\xi }(x \gamma ) \log ^k(a x)}{k!}dx\\
=\frac{1}{2\pi i}\int_{C}\int_{0}^{b}e^{x^{\tau } \theta } w^{-1-k} x^m (a x)^w \left(1+x^{\sigma } \rho
   \right)^{\lambda } J_{\mu }(x \alpha ) J_{\nu }(x \beta ) J_{\xi }(x \gamma )dxdw
\end{multline}
where $Re(b)>0,Re(m+\mu+\nu+\xi)>0$.
\subsection{Right-hand side contour integral}
Using a generalization of Cauchy's integral formula \ref{intro:cauchy}, and the procedure in section (2.2) we get;
\begin{multline}\label{eq:f3}
\sum _{j=0}^{\infty } \sum _{h=0}^{\infty } \sum _{l=0}^{\infty } \sum _{p=0}^{\infty } \sum _{q=0}^{\infty }
   \frac{(-1)^{h+j+l} 2^{-2 h-2 j-2 l-\mu -\nu -\xi } a^{-1-2 h-2 j-2 l-m-\mu -\nu -\xi -q \sigma -p \tau } \alpha ^{2j+\mu } \beta ^{2 h+\nu }   }{h! j!
   k! l! p! \Gamma (1+j+\mu ) \Gamma (1+h+\nu ) \Gamma (1+l+\xi )}\\
\times\frac{\gamma ^{2 l+\xi } \theta ^p \rho ^q}{(-1-2 h-2 j-2 l-m-\mu -\nu -\xi -q \sigma -p \tau)^{-1-k}}\\ \times
\binom{\lambda }{q} \Gamma (1+k,-((1+2 h+2 j+2 l+m+\mu +\nu +\xi +q \sigma +p \tau ) \log (a b)))\\
=-\frac{1}{2\pi i}\int_{C}\sum _{j=0}^{\infty } \sum _{h=0}^{\infty } \sum
   _{l=0}^{\infty } \sum _{p=0}^{\infty } \sum _{q=0}^{\infty } \frac{(-1)^{h+j+l} 2^{-2 h-2 j-2 l-\mu -\nu -\xi } a^w
   b^{1+2 h+2 j+2 l+m+w+\mu +\nu +\xi +q \sigma +p \tau } }{(1+2 h+2 j+2 l+m+w+\mu +\nu +\xi +q \sigma +p \tau ) }\\
\times\frac{w^{-1-k} \alpha ^{2 j+\mu } \beta ^{2 h+\nu } \gamma ^{2
   l+\xi } \theta ^p \rho ^q \binom{\lambda }{q}}{h! j! l! p!
   \Gamma (1+j+\mu ) \Gamma (1+h+\nu ) \Gamma (1+l+\xi )}dw
\end{multline}
where $Re(b)>0,Re(m+\mu+\nu+\xi)>0$.
from equation [Wolfram, \href{http://functions.wolfram.com/07.27.02.0002.01}{07.27.02.0002.01}] where $|Re(b)|<1$.
\begin{theorem}
Generalized Fettis-type integral with three Bessel functions.
\begin{multline}\label{eq:thm_fettis}
\int_0^b e^{x^{\tau } \theta } x^m \left(1+x^{\sigma } \rho \right)^{\lambda } J_{\mu }(x \alpha ) J_{\nu
   }(x \beta ) J_{\xi }(x \gamma ) \log ^k\left(\frac{1}{a x}\right) \, dx\\
=\sum _{j=0}^{\infty } \sum
   _{h=0}^{\infty } \sum _{l=0}^{\infty } \sum _{p=0}^{\infty } \sum _{q=0}^{\infty } \frac{(-1)^{h+j+l} 2^{-2
   h-2 j-2 l-\mu -\nu -\xi } a^{-1-2 h-2 j-2 l-m-\mu -\nu -\xi -q \sigma -p \tau } \alpha ^{2 j+\mu } \beta ^{2
   h+\nu }
   }{h! j! l! p! \Gamma (1+j+\mu ) \Gamma (1+h+\nu ) \Gamma (1+l+\xi )}\\ \times
\frac{ \gamma ^{2 l+\xi } \theta ^p \rho ^q }{(1+2 h+2 j+2 l+m+\mu +\nu +\xi +q \sigma +p \tau )^{k+1}}\\
\times\binom{\lambda }{q} \Gamma \left(1+k,(1+2 h+2 j+2 l+m+\mu +\nu +\xi +q \sigma +p \tau ) \log \left(\frac{1}{ab}\right)\right)
\end{multline}
where $Re(b)>0$.
\end{theorem}
\begin{proof}
Since the right-hand sides of equations (\ref{eq:f2}) and (\ref{eq:f3}) are equal relative to equation (\ref{eq:f1}), we may equate the left-hand sides and simplify the gamma function to yield the stated result.
\end{proof}
\section{Evaluations and derivations}
\begin{example}
Generalized Fettis integral see \cite{fettis}. Here we use equation (\ref{eq:thm_fettis}) and set $k\to 0,a\to 1,\tau \to 0,\theta \to 1,\lambda \to 0,\sigma \to 0,\rho \to 1,m\to \lambda -1$ and simplify;
\begin{multline}\label{eq:30.1}
\int_0^b x^{-1+\lambda } J_{\mu }(x \alpha ) J_{\nu }(x \beta ) J_{\xi }(x \gamma ) \, dx\\
=\sum _{h=0}^{\infty }
   \sum _{j=0}^{\infty } \frac{(-1)^{h+j} 2^{-2 h-2 j-\mu -\nu -\xi } b^{2 h+2 j+\lambda +\mu +\nu +\xi } \alpha ^{2
   j+\mu } \beta ^{2 h+\nu } \gamma ^{\xi } }{(2 h+2 j+\lambda +\mu +\nu +\xi ) \Gamma (1+h) \Gamma (1+j) \Gamma (1+j+\mu ) \Gamma (1+h+\nu )
   \Gamma (1+\xi )}\\ \times
\, _1F_2\left(h+j+\frac{\lambda }{2}+\frac{\mu }{2}+\frac{\nu
   }{2}+\frac{\xi }{2};1+h+j+\frac{\lambda }{2}+\frac{\mu }{2}+\frac{\nu }{2}+\frac{\xi }{2},1+\xi ;-\frac{1}{4} b^2
   \gamma ^2\right)
\end{multline}
where $Re(b)>0,Re(\lambda+\mu+\nu+\xi)>-2$.
\end{example}
\begin{example}
Derivation of Fettis (1957) equation (1) in \cite{fettis} in terms of double series involving the 1F2 hypergeometric function [Wolfram, \href{https://functions.wolfram.com/HypergeometricFunctions/Hypergeometric1F2/}{01}]. The series involving the $1F2$ and its properties are discussed in \cite{straton}. Here we use equation (\ref{eq:thm_fettis}) and set $k\to 0,a\to 1,\tau \to 0,\theta \to 1,\lambda \to 0,\sigma \to 0,\rho \to 1,\mu \to 0,\nu \to 0,\xi \to 0,m\to 1$ and simplify;
\begin{multline}
\int_0^b x J_0(x \alpha ) J_0(x \beta ) J_0(x \gamma ) \, dx\\
=\sum _{h=0}^{\infty } \sum _{l=0}^{\infty }
   \frac{(-1)^{h+l} 2^{-1-2 h-2 l} b^{2+2 h+2 l} \beta ^{2 h} \gamma ^{2 l} \, _1F_2\left(1+h+l;1,2+h+l;-\frac{1}{4}
   b^2 \alpha ^2\right)}{(1+h+l) \Gamma (1+h)^2 \Gamma (1+l)^2}
\end{multline}
where $Re(b)>0$.
\end{example}
\section{Infinite Fettis integral forms and evaluations}
In this section we derive and evaluate the infinite form of integrals with three Bessel functions. Work found in current literature featuring this type of kernel are listed in the articles by Auluck \cite{auluck} where the evaluations of these types of integrals was done using Wolfram Mathematica, Lovatt et al. \cite{lovatt}, worked on a method to efficiently evaluate integrals containing the product of three Bessel functions of the first kind and of any non negative real order is presented, Glasser et al. \cite{glasser}, produced a number of new definite integrals involving Bessel functions derived by  finding new integral representations for the product of two Bessel functions of different order and argument in terms of the generalized hypergeometric function with subsequent reduction to special cases. Lin \cite{lin} used contour integration and the residue theorem to systematically develop and evaluate infinite integrals involving Bessel functions or closely related ones.
\begin{example}
A generalized infinite integral involving the product of three Bessel functions. Here we use equation (\ref{eq:thm_fettis}) and set $k\to 0,a\to 1,\theta \to -\theta $ and simplify the gamma function using equation [Wolfram, \href{http://functions.wolfram.com/06.07.03.0005.01}{01}];
\begin{multline}\label{eq:fettis_inf}
\int_0^{\infty } e^{-x^{\tau } \theta } x^m \left(1+x^{\sigma } \rho \right)^{\lambda } J_{\mu }(x \alpha )
   J_{\nu }(x \beta ) J_{\xi }(x \gamma ) \, dx\\
=\sum _{j=0}^{\infty } \sum _{h=0}^{\infty } \sum _{l=0}^{\infty } \sum
   _{q=0}^{\infty } \frac{(-1)^{h+j+l} 2^{-2 h-2 j-2 l-\mu -\nu -\xi } \alpha ^{2 j+\mu } \beta ^{2 h+\nu } \gamma ^{2
   l+\xi } }{\tau  h! j! l! \Gamma (1+j+\mu ) \Gamma
   (1+h+\nu ) \Gamma (1+l+\xi )}\\
    \theta ^{-\frac{1+2 h+2 j+2 l+m+\mu +\nu +\xi +q \sigma }{\tau }} \rho ^q \binom{\lambda }{q} \\ \times
\Gamma
   \left(\frac{1+2 h+2 j+2 l+m+\mu +\nu +\xi +q \sigma }{\tau }\right)
\end{multline}
where $Re(m+\mu+\nu+\xi)>0$.
\end{example}
\begin{example}
Derivation of equation [Glasser, \href{https://arxiv.org/abs/math/9307213}{3.8}]. Here we use equation (\ref{eq:fettis_inf}) and set $\lambda \to 0,\rho \to 1,\sigma \to 0,x\to \sqrt{t},\tau \to 2 \tau $ and simplify;
\begin{multline}
\int_0^{\infty } e^{-x \alpha } J_0\left(a \sqrt{x}\right) J_0\left(b \sqrt{x}\right) J_0\left(c
   \sqrt{x}\right) \, dx\\
=\begin{cases}
			\sum\limits_{h=0}^{\infty } \sum\limits _{l=0}^{\infty } \frac{(-1)^{h+l} 4^{-h-l} b^{2 h} c^{2 l} \alpha
   ^{-1-h-l} \Gamma (1+h+l) \, _1F_1\left(1+h+l;1;-\frac{a^2}{4 \alpha }\right)}{\Gamma (1+h)^2 \Gamma
   (1+l)^2}, & \text{if $Re(\alpha)>0$}\\
            \frac{1}{\alpha }\exp \left(-\frac{a^2+b^2+c^2}{4 \alpha }\right) \sum\limits_{n=0}^{\infty } (2-\delta _{0,n})I_n\left(\frac{a b}{2 \alpha }\right) I_n\left(\frac{a c}{2 \alpha }\right) I_n\left(\frac{b c}{2 \alpha
   }\right), & \text{if $Re(\alpha)>\pi$}
		 \end{cases}
\end{multline}
where $Re(a+b+c)>0$.
\end{example}
\begin{example}
Laplace transform. Here we use equation (\ref{eq:fettis_inf}) and set $\tau \to 1,\theta \to p,x\to t$ and simplify;
\begin{multline}
\int_0^{\infty } e^{-p t} J_z(t \alpha ) J_z(t \beta ) J_z(t \gamma ) \, dt\\
=\sum _{h=0}^{\infty } \sum
   _{l=0}^{\infty } \frac{(-1)^{h+l} 2^{-2 h-2 l-3 z} p^{-1-2 h-2 l-3 z} \alpha ^z \beta ^{2 h+z} \gamma ^{2 l+z}
   \Gamma (1+2 h+2 l+3 z) }{\Gamma (1+h) \Gamma (1+l) \Gamma (1+z) \Gamma (1+h+z) \Gamma (1+l+z)}\\ \times
\, _2F_1\left(\frac{1}{2}+h+l+\frac{3 z}{2},1+h+l+\frac{3 z}{2};1+z;-\frac{\alpha
   ^2}{p^2}\right)
\end{multline}
where $Re(p)>0$.
\end{example}
\begin{example}
 Gaussian factor. Extended Erdeyli form equation (6.633.1) in \cite{grad}. Here we use equation (\ref{eq:fettis_inf}) and set $\lambda \to 0,\sigma \to 0,\rho \to 1,\tau \to 2,x\to t,\theta \to p$ and simplify;
 \begin{multline}
\int_0^{\infty } e^{-p t^2} t^{-1+\lambda } J_{\mu }(t \alpha ) J_{\nu }(t \beta ) J_{\xi }(t \gamma ) \,
   dt\\
=\sum _{h=0}^{\infty } \sum _{j=0}^{\infty } \frac{(-1)^{h+j} 2^{-1-2 h-2 j-\mu -\nu -\xi } p^{\frac{1}{2} (-2
   h-2 j-\lambda -\mu -\nu -\xi )} \alpha ^{2 j+\mu } \beta ^{2 h+\nu } \gamma ^{\xi } }{h! j! \Gamma (1+j+\mu ) \Gamma (1+h+\nu ) \Gamma (1+\xi )}\\ \times
\Gamma \left(\frac{1}{2} (2 h+2
   j+\lambda +\mu +\nu +\xi )\right) \, _1F_1\left(\frac{1}{2} (2 h+2 j+\lambda +\mu +\nu +\xi );1+\xi ;-\frac{\gamma
   ^2}{4 p}\right)
\end{multline}
where $Re(p)>0,Re(\lambda+\mu+\nu+\xi)>-2$.
\end{example}
\begin{example}
Derivation of an infinite integral using Watson (1922) \cite{watson}, page 413 equation (13.47(7)). Here we use equation (\ref{eq:thm_fettis}) and set $k\to 0,a\to 1,b\to 1,\theta \to 1,\tau \to 0,\lambda \to 0,\rho \to 1,\sigma \to 0,\xi \to \rho ,\alpha \to \cos
   (\Phi ) \cos (\phi ),\beta \to \sin (\Phi ) \sin (\phi ),\gamma \to \cos (\theta ),m\to -\lambda -1$ then $\lambda \to \lambda -2$ and simplify;
   \begin{multline}
\int_1^{\infty } x^{1-\lambda } J_{\mu }(x \cos (\phi ) \cos (\Phi )) J_{\nu }(x \sin (\phi ) \sin (\Phi ))
   J_{\rho }(x \cos (\theta )) \, dx\\
=\frac{\left(\cos ^{\mu }(\phi ) \cos ^{\mu }(\Phi ) \sin ^{\nu }(\phi ) \sin ^{\nu
   }(\Phi ) \cos ^{\rho }(\theta )\right) }{2^{\lambda -1}\Gamma (\rho +1) \Gamma (\nu +1)^2}\\ \times
\sum _{n=0}^{\infty } \frac{(-1)^n (\mu +\nu +2 n+1) (\Gamma (\mu +\nu +n+1)
   \Gamma (\nu +n+1)) \Gamma \left(\frac{1}{2} (\mu +\nu +\rho -\lambda )+n+1\right) }{(n! \Gamma (\mu +n+1)) \Gamma \left(\frac{1}{2} (\mu +\nu -\rho +\lambda )+n+1\right)}\\ \times
\, _2F_1\left(\frac{1}{2} (\mu+\nu +\rho -\lambda )+n+1,\frac{1}{2} (\rho -\lambda -\mu -\nu )-n;\rho +1;\cos ^2(\theta )\right)\\ \,_2F_1\left(-n,\mu +\nu +n+1;\nu +1;\sin ^2(\phi )\right)\\
 \, _2F_1\left(-n,\mu +\nu +n+1;\nu +1;\sin ^2(\Phi)\right)\\
-\sum _{h=0}^{\infty } \sum _{j=0}^{\infty } \frac{(-1)^{h+j} 2^{-2 h-2 j-\mu
   -\nu -\rho } \cos ^{\rho }(\theta ) (\cos (\phi ) \cos (\Phi ))^{2 j+\mu } }{(2+2 h+2 j-\lambda +\mu+\nu +\rho ) \Gamma (1+h) \Gamma (1+j) \Gamma (1+j+\mu ) \Gamma (1+h+\nu ) \Gamma (1+\rho )}\\ \times
\, _1F_2\left(1+h+j-\frac{\lambda
   }{2}+\frac{\mu }{2}+\frac{\nu }{2}+\frac{\rho }{2};2+h+j-\frac{\lambda }{2}+\frac{\mu }{2}+\frac{\nu }{2}+\frac{\rho}{2},1+\rho ;-\frac{1}{4} \cos ^2(\theta )\right) \\
(\sin (\phi ) \sin (\Phi ))^{2 h+\nu }
\end{multline}
where $Re(\nu+\mu+\rho+2)>0,Re(\lambda)>-1/2$.
\end{example}
\begin{example}
A generalized form over the interval $[0,1]$. Here we use equation (\ref{eq:thm_fettis}) and set $k\to 0,a\to 1,b\to 1,\tau \to 0,\theta \to 1,m\to 0,\lambda \to 0,\rho \to 1,\sigma \to 0,\xi \to \lambda ,\alpha
   \to b,\beta \to c,\gamma \to a$ and simplify;
   \begin{multline}\label{eq:31.6}
\int_0^1 J_{\lambda }(a x) J_{\mu }(b x) J_{\nu }(c x) \, dx\\
=\sum _{h=0}^{\infty } \sum _{j=0}^{\infty }
   \frac{(-1)^{h+j} 2^{-2 h-2 j-\lambda -\mu -\nu } a^{\lambda } b^{2 j+\mu } c^{2 h+\nu } }{(1+2 h+2 j+\lambda +\mu
   +\nu ) \Gamma (1+h) \Gamma (1+j) \Gamma (1+\lambda ) \Gamma (1+j+\mu ) \Gamma (1+h+\nu )}\\ \times
\,
   _1F_2\left(\frac{1}{2}+h+j+\frac{\lambda }{2}+\frac{\mu }{2}+\frac{\nu }{2};1+\lambda
   ,\frac{3}{2}+h+j+\frac{\lambda }{2}+\frac{\mu }{2}+\frac{\nu }{2};-\frac{a^2}{4}\right)
\end{multline}
where $Re(\mu+\nu+\lambda)>0$.
\end{example}
\begin{example}
Derivation of equation (2.12.42.3) in \cite{prud2} over $[1,\infty)$. Here we use equation (\ref{eq:31.6}) and set $\lambda \to 0,\mu \to 0,\nu \to 1$ then expand the infinite integral and simplify;
\begin{multline}
\int_1^{\infty } J_0(a x) J_0(b x) J_1(c x) \, dx
=\frac{\cos ^{-1}\left(\frac{a^2+b^2-c^2}{2 a b}\right)}{\pi 
   c}\\
-\sum _{h=0}^{\infty } \sum _{j=0}^{\infty } \frac{(-1)^{h+j} 2^{-1-2 h-2 j} b^{2 j} c^{1+2 h} \,
   _1F_2\left(1+h+j;1,2+h+j;-\frac{a^2}{4}\right)}{(2+2 h+2 j) \Gamma (1+h) \Gamma (2+h) \Gamma (1+j)^2}
\end{multline}
where $|a-b|< Re(c) < Re(a+b), Re(a)>0, Re(b)>0$. 
\end{example}
\begin{example}
Derivation of equation (2.12.42.18) in \cite{prud2} over $[1,\infty)$. Errata. Here we use equation (\ref{eq:thm_fettis}) and set $k\to 0,a\to 1,b\to 1,\theta \to 1,\tau \to 0,\lambda \to 0,\rho \to 1,\sigma \to 0,m\to 0,\mu \to \nu ,\xi \to \nu
   ,\alpha \to c,\beta \to c,\gamma \to c$ then expand the infinite integral and simplify;
   \begin{multline}
\int_1^{\infty } J_{\nu }(c x){}^3 \, dx=\frac{\Gamma \left(\frac{1}{6} (3 \nu +1)\right)}{(3 c) \left(\Gamma
   \left(\frac{1}{6} (3 \nu +5)\right) \Gamma \left(\frac{2}{3}\right) \Gamma \left(\frac{2}{3}\right)\right)}-\\
\sum
   _{h=0}^{\infty } \sum _{j=0}^{\infty } \frac{(-1)^{h+j} 2^{-2 h-2 j-3 \nu } c^{2 h+2 j+3 \nu } \,
   _1F_2\left(\frac{1}{2}+h+j+\frac{3 \nu }{2};1+\nu ,\frac{3}{2}+h+j+\frac{3 \nu }{2};-\frac{c^2}{4}\right)}{(1+2 h+2
   j+3 \nu ) \Gamma (1+h) \Gamma (1+j) \Gamma (1+\nu ) \Gamma (1+h+\nu ) \Gamma (1+j+\nu )}
\end{multline}
where $Re(c)>0, Re(\nu)>-1/3$.
\end{example}
\begin{example}
Derivation of equation (2.12.42.19) in \cite{prud2} over $[1,\infty)$. Errata. Here we use equation (\ref{eq:thm_fettis}) and set $k\to 0,a\to 1,b\to 1,\theta \to 1,\tau \to 0,\lambda \to 0,\rho \to 1,\sigma \to 0,m\to 0,\mu \to -\nu ,\xi \to -\nu
   ,\alpha \to c,\beta \to c,\gamma \to c$  then expand the infinite integral and simplify;
   \begin{multline}
\int_1^{\infty } J_{-\nu }(c x){}^2 J_{\nu }(c x) \, dx=\frac{2 \sin \left(\frac{1}{3} (1-3 \nu ) \pi \right)
   \Gamma \left(\frac{1}{6} (1-3 \nu )\right)}{\left(3^{3/2} c\right) \left(\Gamma \left(\frac{1}{6} (5-3 \nu )\right)
   \Gamma \left(\frac{2}{3}\right) \Gamma \left(\frac{2}{3}\right)\right)}\\
-\sum _{h=0}^{\infty } \sum _{j=0}^{\infty }
   \frac{(-1)^{h+j} 2^{-2 h-2 j+\nu } c^{2 h+2 j-\nu } \, _1F_2\left(\frac{1}{2}+h+j-\frac{\nu }{2};1-\nu
   ,\frac{3}{2}+h+j-\frac{\nu }{2};-\frac{c^2}{4}\right)}{(1+2 h+2 j-\nu ) \Gamma (1+h) \Gamma (1+j) \Gamma (1-\nu )
   \Gamma (1+j-\nu ) \Gamma (1+h+\nu )}
\end{multline}
where $Re(c)>0, Re(\nu)>-1/3$.
\end{example}
\begin{example}
Derivation of equation (2.12.42.22) in \cite{prud2} over $[1,\infty)$. Errata. Here we use equation (\ref{eq:thm_fettis}) and set $k\to 0,a\to 1,b\to 1,\theta \to 1,\tau \to 0,\lambda \to 0,\rho \to 1,\sigma \to 0,m\to 1-\mu ,\alpha \to b,\beta
   \to b,\xi \to \nu ,\gamma \to c$  then expand the infinite integral and simplify;
   \begin{multline}
\int_1^{\infty } x^{1-\mu } J_{\mu }(b x) J_{\nu }(b x) J_{\nu }(c x) \, dx=\frac{c^{\mu -1} \left(4
   b^2-c^2\right)^{\frac{1}{4} (2 \mu -1)} }{2^{\mu } \sqrt{\pi } b^{\mu +\frac{1}{2}}}P_{\nu -\frac{1}{2}}^{\frac{1}{2}-\mu }\left(\frac{c}{2
   b}\right)\\
-\sum _{h=0}^{\infty } \sum _{j=0}^{\infty } \frac{(-1)^{h+j}
   2^{-1-2 h-2 j-\mu -2 \nu } b^{2 h+2 j+\mu +\nu } c^{\nu } }{(1+h+j+\nu ) \Gamma (1+h) \Gamma (1+j) \Gamma (1+j+\mu ) \Gamma (1+\nu ) \Gamma (1+h+\nu
   )}\\ \times
\, _1F_2\left(1+h+j+\nu ;1+\nu ,2+h+j+\nu
   ;-\frac{c^2}{4}\right)
\end{multline}
where $Re(a)>0,Re(c)>0,c \neq 2b, Re(\mu)>-1/2, Re(\mu+\nu)>-1$.
\end{example}
\begin{example}
Integrating triple product of three Bessel functions over a finite domain in StackExchange [StackExchange, \href{https://math.stackexchange.com/questions/3000535/integrating-triple-product-of-bessel-functions-over-a-finite-domain}{01}]. Here we use equation (\ref{eq:thm_fettis}) and set $\tau \to 0,\theta \to 1,k\to 0,a\to 1,b\to 1,\lambda \to 0,\rho \to 1,\sigma \to 0,\mu \to 0,\xi \to 0,\nu \to
   0,x\to r,\alpha \to a,\beta \to b,\gamma \to c,m\to 1$ and simplify;
   \begin{multline}
\int_0^1 r J_0(a r) J_0(b r) J_0(c r) \, dr\\
=\sum _{h=0}^{\infty } \sum _{j=0}^{\infty } \frac{(-1)^{h+j}
   2^{-1-2 h-2 j} a^{2 j} b^{2 h} \, _1F_2\left(1+h+j;1,2+h+j;-\frac{c^2}{4}\right)}{(1+h+j) \Gamma (1+h)^2 \Gamma
   (1+j)^2}
\end{multline}
where $a,b,c\in\mathbb{C}$.
\end{example}
\begin{example}
 The  Hurwitz-Lerch transcendent. Here we use equation (\ref{eq:thm_fettis}) and set $a\to 1,b\to 1,\lambda \to -1$ and simplify using equation [DLMF, \href{https://dlmf.nist.gov/25.14.i}{25.14.1}] and simplify;
 \begin{multline}\label{eq:fettis_lerch}
\int_0^1 \frac{e^{x^{\tau } \theta } x^m J_{\mu }(x \alpha ) J_{\nu }(x \beta ) J_{\xi }(x \gamma ) \log
   ^k\left(\frac{1}{x}\right)}{1+x^{\sigma } \rho } \, dx\\
=\sum _{h=0}^{\infty } \sum _{j=0}^{\infty } \sum
   _{l=0}^{\infty } \sum _{p=0}^{\infty } \frac{(-1)^{h+j+l} 2^{-2 h-2 j-2 l-\mu -\nu -\xi } \alpha ^{2 j+\mu } \beta^{2 h+\nu } \gamma ^{2 l+\xi } \theta ^p \sigma ^{-1-k} \Gamma (1+k) }{h! j! l! p! \Gamma (1+j+\mu ) \Gamma (1+h+\nu ) \Gamma (1+l+\xi)}\\ \times
\Phi \left(-\rho ,1+k,\frac{1+2 h+2 j+2 l+m+\mu +\nu +\xi +p \tau }{\sigma }\right)
\end{multline}
where $Re(\alpha+\beta+\gamma)>0,Re(\sigma)>1$.
\end{example}
\begin{example}
Equation (2.12.43.1) in \cite{prud2} over the interval $[0,1]$. Note equation (2.12.43.1) in \cite{prud2} is in error. Here we use equation (\ref{eq:fettis_lerch}) and set $k\to 0,\theta \to 1,\tau \to 0,\sigma \to 2,\rho \to \frac{1}{z^2}$ and simplify;
\begin{multline}
\int_0^1 \frac{x^m J_{\mu }(x \alpha ) J_{\nu }(x \beta ) J_{\xi }(x \gamma )}{x^2+z^2} \, dx\\
=\sum
   _{h=0}^{\infty } \sum _{j=0}^{\infty } \sum _{l=0}^{\infty } \frac{(-1)^{h+j+l} 2^{-1-2 h-2 j-2 l-\mu -\nu -\xi }
   \alpha ^{2 j+\mu } \beta ^{2 h+\nu } \gamma ^{2 l+\xi } }{z^2 \Gamma (1+h) \Gamma (1+j) \Gamma (1+l) \Gamma (1+j+\mu ) \Gamma (1+h+\nu ) \Gamma
   (1+l+\xi )}\\ \times
\Phi \left(-\frac{1}{z^2},1,\frac{1}{2} (1+2 h+2 j+2
   l+m+\mu +\nu +\xi )\right)
\end{multline}
where $|Re(m)|<1$.
\end{example}
\begin{example}
Derivation of equation (2.12.42.20) in \cite{prud2} over $[1,\infty)$. Here we use equation (\ref{eq:30.1}) we set $b\to 1,\lambda \to \alpha ,\mu \to \lambda ,\alpha \to b,\nu \to \mu ,\beta \to b,\xi \to \nu ,\gamma \to c$ and simplify;
\begin{multline}
\int_1^{\infty } x^{-1+\alpha } J_{\lambda }(b x) J_{\mu }(b x) J_{\nu }(c x) \, dx
=\frac{2^{\alpha -1} b^{\mu
   +\lambda } c^{-\alpha -\mu -\lambda } \Gamma \left(\frac{1}{2} (\alpha +\mu +\lambda +\nu )\right)}{\Gamma (\mu +1) \Gamma (\lambda +1) \Gamma \left(1+\frac{1}{2} (\nu -\alpha -\mu -\lambda )\right)}\\ \times \,_4F_3\left(\frac{1}{2} (\alpha +\mu +\lambda -\nu ),\frac{1}{2} (\alpha +\mu +\lambda +\nu ),\frac{1}{2} (1+\mu+\lambda ),1+\frac{\mu +\lambda }{2};\right. \\ \left.
1+\mu +\lambda ,1+\mu ,1+\lambda ;\left(\frac{2 b}{c}\right)^2\right)\\
-\sum _{h=0}^{\infty }\sum _{j=0}^{\infty } \frac{(-1)^{h+j} 2^{-2 h-2 j-\lambda -\mu -\nu } b^{2 h+2 j+\lambda +\mu } c^{\nu } }{(2 h+2 j+\alpha +\lambda +\mu+\nu ) \Gamma (1+h) \Gamma (1+j) \Gamma (1+j+\lambda ) \Gamma (1+h+\mu ) \Gamma (1+\nu )}\\ \times
\, _1F_2\left(h+j+\frac{\alpha }{2}+\frac{\lambda }{2}+\frac{\mu }{2}+\frac{\nu }{2};1+h+j+\frac{\alpha}{2}+\frac{\lambda }{2}+\frac{\mu }{2}+\frac{\nu }{2},1+\nu ;-\frac{c^2}{4}\right)
\end{multline}
where $0< Re(2b) < Re(c), -Re(\mu+\nu+\lambda) < Re(\alpha)<5/2$.
\end{example}
\begin{example}
Derivation of equation (2.12.42.22) in \cite{prud2} over $[1,\infty)$. Here we use equation (\ref{eq:30.1}) we set $b\to 1,\lambda \to 2-\mu ,\alpha \to b,\beta \to b,\xi \to \nu ,\gamma \to c$ and simplify;
\begin{multline}
\int_1^{\infty } x^{1-\mu } J_{\mu }(b x) J_{\nu }(b x) J_{\nu }(c x) \, dx
=\frac{\left(c^{\mu -1} \left(4
   b^2-c^2\right)^{\frac{1}{4} (2 \mu -1)}\right) P_{\nu -\frac{1}{2}}^{\frac{1}{2}-\mu }\left(\frac{c}{2
   b}\right)}{2^{\mu } \sqrt{\pi } b^{\mu +\frac{1}{2}}}\\
-\sum _{h=0}^{\infty } \sum _{j=0}^{\infty } \frac{(-1)^{h+j}
   2^{-2 h-2 j-\mu -2 \nu } b^{2 h+2 j+\mu +\nu } c^{\nu } \, _1F_2\left(1+h+j+\nu ;1+\nu ,2+h+j+\nu
   ;-\frac{c^2}{4}\right)}{(2+2 h+2 j+2 \nu ) \Gamma (1+h) \Gamma (1+j) \Gamma (1+j+\mu ) \Gamma (1+\nu ) \Gamma
   (1+h+\nu )}
\end{multline}
where $Re(b)>0,Re(c)>0, Re(c) \neq Re(2b), Re(\mu)>-1/2, Re(\nu+\mu)>-1$.
\end{example}
\begin{example}
Derivation of equation (2.12.42.23) in \cite{prud2} over $[1,\infty)$. Here we use equation (\ref{eq:30.1}) we set $b\to 1,\lambda \to \nu +2,\alpha \to b,\nu \to -\mu ,\beta \to b,\xi \to \nu ,\gamma \to c$ and simplify;
\begin{multline}
\int_1^{\infty } x^{1+\nu } J_{-\mu }(b x) J_{\mu }(b x) J_{\nu }(c x) \, dx=\frac{2^{\nu } \left(4
   b^2-c^2\right)^{-\frac{1}{4} (2 \nu +1)} P_{\mu -\frac{1}{2}}^{\nu +\frac{1}{2}}\left(\left(\frac{c}{b
   \sqrt{2}}\right)^2-1\right)}{b \sqrt{\pi  c}}\\
-\sum _{h=0}^{\infty } \sum _{j=0}^{\infty } \frac{(-1)^{h+j} 2^{-2
   h-2 j-\nu } b^{2 h+2 j} c^{\nu } \, _1F_2\left(1+h+j+\nu ;1+\nu ,2+h+j+\nu ;-\frac{c^2}{4}\right)}{(2+2 h+2 j+2 \nu
   ) \Gamma (1+h) \Gamma (1+j) \Gamma (1+h-\mu ) \Gamma (1+j+\mu ) \Gamma (1+\nu )}
\end{multline}
where $Re(b)>0,Re(c)>0, Re(c) \neq Re(2b), Re(\mu)>-1/2, Re(\nu+\mu)>-1$.
\end{example}
\begin{example}
Derivation of equation (2.12.42.28) in \cite{prud2} over $[1,\infty)$. Here we use equation (\ref{eq:30.1}) we set $b\to 1,\lambda \to 2,\alpha \to b,\nu \to \mu ,\beta \to b,\xi \to 0,\gamma \to c$ and simplify; 
\begin{multline}
\int_1^{\infty } x J_0(c x) J_{\mu }(b x){}^2 \, dx=\frac{2 \frac{1}{\sqrt{4 b^2-c^2}} \cos \left(\mu  \cos
   ^{-1}\left(1-\left(\frac{c}{b \sqrt{2}}\right)^2\right)\right)}{\pi  c}\\
-\sum _{h=0}^{\infty } \sum _{j=0}^{\infty }
   \frac{(-1)^{h+j} 2^{-2 h-2 j-2 \mu } b^{2 h+2 j+2 \mu } \, _1F_2\left(1+h+j+\mu ;1,2+h+j+\mu
   ;-\frac{c^2}{4}\right)}{(2+2 h+2 j+2 \mu ) \Gamma (1+h) \Gamma (1+j) \Gamma (1+h+\mu ) \Gamma (1+j+\mu )}
\end{multline}
where $Re(b)>0,Re(c)>0, Re(\mu)>-1$.
\end{example}
\section{Table 3.10.21 in  Brychkov et al. (2018) on page 170}
\begin{example}
Derivation of generalized form of entry (3.10.21.1) in \cite{brychkov_m} over $x\in[0,1]$. Here we use equation (\ref{eq:thm_fettis}) and set $k\to 0,a\to 1,b\to 1,\tau \to 0,\theta \to 1,\lambda \to 0,\rho \to 1,\sigma \to 0,m\to s-1$ and simplify;
 \begin{multline}\label{eq:bry_1}
\int_0^1 x^{-1+s} J_{\mu }(x \alpha ) J_{\nu }(x \beta ) J_{\xi }(x \gamma ) \, dx\\
=\sum _{h=0}^{\infty } \sum
   _{j=0}^{\infty } \frac{(-1)^{h+j} 2^{-2 h-2 j-\mu -\nu -\xi } \alpha ^{2 j+\mu } \beta ^{2 h+\nu } \gamma ^{\xi } }{(2 h+2 j+s+\mu +\nu +\xi ) \Gamma (1+h) \Gamma (1+j) \Gamma(1+j+\mu ) \Gamma (1+h+\nu ) \Gamma (1+\xi )}\\ \times
\,_1F_2\left(h+j+\frac{s}{2}+\frac{\mu }{2}+\frac{\nu }{2}+\frac{\xi }{2};1+h+j+\frac{s}{2}+\frac{\mu }{2}+\frac{\nu}{2}+\frac{\xi }{2},1+\xi ;-\frac{\gamma ^2}{4}\right)
\end{multline}
where $Re(\nu)>0,|Re(s)|<1$.
\end{example}
\begin{example}
Derivation of entry (3.10.21.1) in \cite{brychkov_m} over $x\in[1,\infty)$. Here we use equation (\ref{eq:bry_1}) and set $\xi \to \lambda ,\gamma \to a,\beta \to b,\alpha \to a$ then expand and simplify the infinite integral;
\begin{multline}
\int_1^{\infty } x^{-1+s} J_{\lambda }(a x) J_{\mu }(a x) J_{\nu }(b x) \, dx\\
=\frac{\left(2^{s-1} a^{\mu+\lambda }\right) \Gamma \left(\frac{1}{2} (s+\lambda +\mu +\nu )\right) }{b^{s+\mu +\lambda } \left(\Gamma (\mu +1) \Gamma(\lambda +1) \Gamma \left(\frac{1}{2} (\nu -\mu -\lambda -s+2)\right)\right)}\\ \times
\, _4F_3\left(\frac{1}{2} (\lambda +\mu+1),\frac{1}{2} (\lambda +\mu +2),\frac{1}{2} (s+\lambda +\mu -\nu ),\frac{1}{2} (s+\lambda +\mu +\nu );\right. \\ \left.
\lambda+1,\mu +1,\lambda +\mu +1;\left(\frac{2 a}{b}\right)^2\right)\\
-\sum _{h=0}^{\infty } \sum
   _{j=0}^{\infty } \frac{(-1)^{h+j} 2^{-2 h-2 j-\lambda -\mu -\nu } a^{2 j+\lambda +\mu } b^{2 h+\nu } 
}{(2 h+2 j+s+\lambda +\mu+\nu ) \Gamma (1+h) \Gamma (1+j) \Gamma (1+\lambda ) \Gamma (1+j+\mu ) \Gamma (1+h+\nu )}\\ \times
\,_1F_2\left(h+j+\frac{s}{2}+\frac{\lambda }{2}+\frac{\mu }{2}+\frac{\nu }{2};1+\lambda
   ,1+h+j+\frac{s}{2}+\frac{\lambda }{2}+\frac{\mu }{2}+\frac{\nu }{2};-\frac{a^2}{4}\right)
\end{multline}
where $Re(\nu)>0,|Re(s)|<1,0< Re(2a)< b$.
\end{example}
%
%\begin{example}
%
%\end{example}
%%
\section{Generalized Hyperharmonic integrals}
These types of integrals were studied by Professor Roy Hughes \cite{hughes1,hughes2}.
The contour integral representation form involving the generalized Hyperharmonic-type integral form is given by;
\begin{multline}\label{eq:h1}
\frac{1}{2\pi i}\int_{C}\int_{0}^{b}a^w e^{-x^{\gamma } \theta +\beta  \left(1+x^{\rho } \xi \right)^{\lambda }} w^{-1-k} x^{m+w} J_{\nu }(x
   \alpha )dxdw\\
=\frac{1}{2\pi i}\int_{C}\sum _{j=0}^{\infty } \sum _{h=0}^{\infty } \sum _{l=0}^{\infty } \sum _{p=0}^{\infty } \frac{(-1)^j
   2^{-2 j-\nu } a^w b^{1+2 j+m+w+p \gamma +\nu +h \rho } w^{-1-k} \alpha ^{2 j+\nu } \beta ^l (-\theta )^p \xi
   ^h \binom{l \lambda }{h}}{(1+2 j+m+w+p \gamma +\nu +h \rho ) j! l! p! \Gamma (1+j+\nu )}dw
\end{multline}
where $Re(b)>0, Re(\mu+\nu)>0$.
\subsection{Left-hand side contour integral representation}
Using a generalization of Cauchy's integral formula \ref{intro:cauchy} and the procedure in section (2.1), we form the definite integral given by;
\begin{multline}\label{eq:h2}
\int_{0}^{b}\frac{x^m \log ^k(a x) J_{\nu }(x \alpha ) e^{\beta  \left(\xi  x^{\rho }+1\right)^{\lambda }-\theta 
   x^{\gamma }}}{k!}dx\\
=\frac{1}{2\pi i}\int_{C}\int_{0}^{b}w^{-k-1} x^m (a x)^w J_{\nu }(x \alpha ) e^{\beta  \left(\xi  x^{\rho }+1\right)^{\lambda
   }-\theta  x^{\gamma }}dxdw
\end{multline}
where $Re(b)>0, Re(m+\nu)>0$.
\subsection{Right-hand side contour integral}
Using a generalization of Cauchy's integral formula \ref{intro:cauchy}, and the procedure in section (2.2) we get;
\begin{multline}\label{eq:h3}
\sum _{j=0}^{\infty } \sum _{h=0}^{\infty } \sum _{l=0}^{\infty } \sum _{p=0}^{\infty } \frac{(-1)^j 2^{-2
   j-\nu } a^{-1-2 j-m-p \gamma -\nu -h \rho } \alpha ^{2 j+\nu } \beta ^l (-\theta )^p \xi ^h \binom{l \lambda }{h} }{j! k!
   l! p! \Gamma (1+j+\nu )(-1-2 j-m-p \gamma-\nu -h \rho )^{1+k} }\\ \times
\Gamma (1+k,-((1+2 j+m+p \gamma +\nu +h \rho ) \log (a b)))\\
=-\frac{1}{2\pi i}\int_{C}\sum _{j=0}^{\infty } \sum _{h=0}^{\infty } \sum _{l=0}^{\infty } \sum
   _{p=0}^{\infty } \frac{(-1)^j 2^{-2 j-\nu } a^w b^{1+2 j+m+w+p \gamma +\nu +h \rho } w^{-1-k} \alpha ^{2 j+\nu} \beta ^l (-\theta )^p \xi ^h \binom{l \lambda }{h}}{(1+2 j+m+w+p \gamma +\nu +h \rho ) j! l! p! \Gamma(1+j+\nu )}dw
\end{multline}
where $Re(b)>0, Re(m+\nu)>0$.
from equation [Wolfram, \href{http://functions.wolfram.com/07.27.02.0002.01}{07.27.02.0002.01}] where $|Re(b)|<1$.
\begin{theorem}
Generalized Hyperharmonic-type integral.
\begin{multline}\label{eq:thm_luke}
\int_0^b e^{-x^{\gamma } \theta +\beta  \left(1+x^{\rho } \xi \right)^{\lambda }} x^m J_{\nu }(x \alpha )
   \log ^k\left(\frac{1}{a x}\right) \, dx\\
=\sum _{j=0}^{\infty } \sum _{h=0}^{\infty } \sum _{l=0}^{\infty } \sum
   _{p=0}^{\infty } \frac{(-1)^j 2^{-2 j-\nu } a^{-1-2 j-m-p \gamma -\nu -h \rho } \alpha ^{2 j+\nu } \beta ^l
   (-\theta )^p \xi ^h \binom{l \lambda }{h} }{j! l! p! \Gamma (1+j+\nu ) (1+2 j+m+p \gamma +\nu +h \rho )^{k+1}}\\ \times
\Gamma \left(1+k,(1+2 j+m+p \gamma +\nu +h \rho ) \log
   \left(\frac{1}{a b}\right)\right)
\end{multline}
where $Re(b)>0, Re(m+\nu)>0$.
\end{theorem}
\begin{proof}
Since the right-hand sides of equations (\ref{eq:h2}) and (\ref{eq:h3}) are equal relative to equation (\ref{eq:h1}), we may equate the left-hand sides and simplify the gamma function to yield the stated result.
\end{proof}
\section{Schwarz functions and generalizations in Luke (1962)}
\begin{example}
Derivation of the Hyperharmonic integral by Professor Hughes. Here we use equation (\ref{eq:thm_luke}) and set $k\to 0,a\to 1,\gamma \to 1,\lambda \to \frac{1}{2},\xi \to -\frac{1}{\xi ^2},\rho \to 2$ then $\beta \to i \beta  \xi ,\theta \to i \theta ,m\to 1,\nu \to 0,\alpha \to \gamma ,x\to k$ finally $\theta \to \alpha ,\xi \to \sigma ,b\to \sigma$ and simplify;
\begin{multline}
\int_0^{\sigma } e^{-i k \alpha +i \beta  \sqrt{-k^2+\sigma ^2}} k J_0(k \gamma ) \, dk\\
=\sum
   _{l=0}^{\infty } \sum _{p=0}^{\infty } \frac{(-1)^p i^{l+p} 2^{-1-l} \sqrt{\pi } \alpha ^p \beta ^l \sigma
   ^{2+l+p} \Gamma \left(\frac{2+p}{2}\right) }{\Gamma
   \left(\frac{1+l}{2}\right) \Gamma (1+p) \Gamma \left(\frac{1}{2} (4+l+p)\right)}\,
   _1F_2\left(1+\frac{p}{2};1,2+\frac{l}{2}+\frac{p}{2};-\frac{1}{4} \gamma ^2 \sigma ^2\right)
\end{multline}
where $Re(\sigma)>0$.
\end{example}
\begin{example}
Derivation of equation 10.2.15 in \cite{luke}. Here we use equation (\ref{eq:thm_luke}) and set $k\to 0,a\to 1,b\to z,x\to t,\lambda \to 0,\rho \to 0,\xi \to 1,\beta \to 1,m\to \mu ,\gamma \to 1,\theta \to-i,\alpha \to \lambda$. Next simplify the series two ways, first by $l\in[0,\infty),j\in[0,\infty)$ next $l\in[0,\infty),p\in[0,\infty)$ and simplify;
\begin{multline}
\int_0^z e^{i t} t^{\mu } J_{\nu }(t \lambda ) \, dt\\
=\sum _{p=0}^{\infty } \frac{i^p 2^{-\nu } z^{1+p+\mu
   +\nu } \lambda ^{\nu } \, _1F_2\left(\frac{1}{2}+\frac{p}{2}+\frac{\mu }{2}+\frac{\nu
   }{2};\frac{3}{2}+\frac{p}{2}+\frac{\mu }{2}+\frac{\nu }{2},1+\nu ;-\frac{1}{4} z^2 \lambda ^2\right)}{(1+p+\mu
   +\nu ) \Gamma (1+p) \Gamma (1+\nu )}\\
=\sum _{j=0}^{\infty } \frac{i (-1)^j 2^{-2 j-\nu } (-i z)^{-2 j-\mu -\nu
   } z^{2 j+\mu +\nu } \lambda ^{2 j+\nu } \Gamma (1+2 j+\mu +\nu ,0,-i z)}{\Gamma (1+j) \Gamma (1+j+\nu
   )}\\
=z^{\mu +1} \exp (i z) \sum _{k=0}^{\infty } \frac{\left((-1)^k \left(\frac{\lambda  z}{2}\right)^{\nu +2
   k}\right) \, _1F_1(1;\mu +\nu +2 k+2;-i z)}{k! (\mu +\nu +2 k+1) \Gamma (\nu +k+1)}
\end{multline}
where $Re(z)>0, Re(\mu+\nu)>0$.
\end{example}
\begin{example}
Derivation of equation 10.3.(3) in \cite{luke}. Here we use equation (\ref{eq:thm_luke}) and set $k\to 0,a\to 1,\lambda \to 0,\xi \to 1,\rho \to 0,\beta \to 1,\gamma \to 1,x\to t,\theta \to i \cos (\theta ),b\to z,\alpha \to 1,m\to 0$ and simplify;
\begin{multline}
\int_0^z e^{-i t \cos (\theta )} J_{\nu }(t) \, dt\\
=\sum _{p=0}^{\infty } \frac{2^{-\nu } z^{1+p+\nu } (-i
   \cos (\theta ))^p \, _1F_2\left(\frac{1}{2}+\frac{p}{2}+\frac{\nu }{2};\frac{3}{2}+\frac{p}{2}+\frac{\nu
   }{2},1+\nu ;-\frac{z^2}{4}\right)}{(1+p+\nu ) \Gamma (1+p) \Gamma (1+\nu )}\\
=2 \exp (-i z \cos (\theta )) \sum
   _{k=0}^{\infty } \frac{i^k \sin ((k+1) \theta ) J_{k+1+\nu }(z)}{\sin (\theta )}
\end{multline}
where $Re(\nu)>-1$.
\end{example}
\section{Representation in Series of Circular Functions in Luke (1962) in terms of a series involving a hypergeometric function}
\begin{example}
Derivation of equation 10.4.(1) in \cite{luke}. Here we use equation (\ref{eq:thm_luke}) and set $k\to 0,a\to 1,\lambda \to 0,\xi \to 1,\rho \to 0,\beta \to 1,\gamma \to 1,x\to t,\theta \to i \cos (\theta ),b\to z,\alpha \to 1,m\to 0$ and simplify. This evaluation represents a more easy and accurate closed form for the definite integral listed in Luke (1964).
\begin{multline}
\int_0^z e^{-i t} t^m J_{2 n}(t \lambda ) \, dt\\
=\sum _{p=0}^{\infty } \frac{(-i)^p 2^{-2 n} z^{1+m+2 n+p}
   \lambda ^{2 n} \, _1F_2\left(\frac{1}{2}+\frac{m}{2}+n+\frac{p}{2};1+2
   n,\frac{3}{2}+\frac{m}{2}+n+\frac{p}{2};-\frac{1}{4} z^2 \lambda ^2\right)}{(1+m+2 n+p) \Gamma (1+2 n) \Gamma
   (1+p)}
\end{multline}
where $Re(m+2n)>0$.
\end{example}
\begin{example}
Derivation of equation 10.4.(1) in \cite{luke}. Here we use equation (\ref{eq:thm_luke}) and set $k\to 0,a\to 1,\lambda \to 0,\xi \to 1,\rho \to 0,\beta \to 1$ and simplify using [Wolfram, \href{http://functions.wolfram.com/06.07.03.0005.01}{01}];
\begin{multline}
\int_0^{\infty } e^{-p t} t^{\mu } J_{\nu }(t \lambda ) \, dt\\
=\frac{2^{-\nu } p^{-1-\mu -\nu } \lambda
   ^{\nu } \Gamma (1+\mu +\nu ) \, _2F_1\left(\frac{1}{2}+\frac{\mu }{2}+\frac{\nu }{2},1+\frac{\mu
   }{2}+\frac{\nu }{2};1+\nu ;-\frac{\lambda ^2}{p^2}\right)}{\Gamma (1+\nu )}
\end{multline}
where $Re(p)>0, Re(\nu+\mu)>0$.
\end{example}
\begin{example}
Use the equation (\ref{eq:thm_luke}) and set $k\to 0,a\to 1,\gamma \to 0,\theta \to 1,\lambda \to \frac{1}{2},\beta \to i \beta ,\rho \to 2,\xi \to
   -\frac{1}{\xi ^2},m\to 0$, then $\beta \to \beta  \xi$ then take the infinite sums over $p\in[0,\infty)$ and $l\in[0,\infty)$ and simplify.
   \begin{multline}
\int_0^b e^{i \beta  \sqrt{-x^2+\xi ^2}} J_{\nu }(x \alpha ) \, dx\\
=\sum _{j=0}^{\infty } \sum
   _{l=0}^{\infty } \frac{(-1)^j i^l 2^{-2 j-\nu } b^{1+2 j+\nu } \alpha ^{2 j+\nu } \beta ^l \xi ^l \,
   _2F_1\left(-\frac{l}{2},\frac{1}{2}+j+\frac{\nu }{2};\frac{3}{2}+j+\frac{\nu }{2};\frac{b^2}{\xi
   ^2}\right)}{(1+2 j+\nu ) \Gamma (1+j) \Gamma (1+l) \Gamma (1+j+\nu )}
\end{multline}
where $Re(\nu)>-2$.
\end{example}
\begin{example}
Here we use equation (\ref{eq:thm_luke}) and set $k\to 0,a\to 1,\lambda \to 0,\xi \to 1,\rho \to 0,\beta \to 1$ and simplify using [Wolfram, \href{http://functions.wolfram.com/06.07.03.0005.01}{01}];
\begin{equation}
\int_0^{\infty } e^{-x^{\gamma } \theta } x^m J_{\nu }(x \alpha ) \, dx=\sum _{j=0}^{\infty } \frac{(-1)^j
   2^{-2 j-\nu } \alpha ^{2 j+\nu } \theta ^{-\frac{1+2 j+m+\nu }{\gamma }} \Gamma \left(\frac{1+2 j+m+\nu }{\gamma
   }\right)}{\gamma  \Gamma (1+j) \Gamma (1+j+\nu )}
\end{equation}
where $Re(\gamma)>1,Re(\theta)>1,Im(\theta)>0$.
\end{example}
\section{Table 2.12.10 in Prudnikov et al volume II, (1986b)}
In this section we derive entries for integrals over a finite domain involving the product of the Bessel function of the first kind, exponential function of a rational function and polynomial function. This integral and its evaluation are expressed in terms of infinite series involving special functions. The integrals in this Table of Prudnikov et al are used in the work by Eingorn et al. [Eingorn, \href{https://doi.org/10.3390/universe7040101}{12}] where they investigated the influence of the chimney topology $T \times T \times R$ of the Universe on the gravitational potential and force that are generated by point-like massive bodies. The work by Jos\'{e} Manuel M\'{e}ndez P\'{e}rez on page 35 in [P\'{e}rez, \href{https://riull.ull.es/xmlui/bitstream/handle/915/10177/cp290.pdf?sequence=1&isAllowed=y}{1.3.21}] where the integral distributional transformation Kontorovich-Lebedev and their applications was investigated. The work by Kacha Dzhaparidze on page 26 in [Dzhaparidze, \href{https://ir.cwi.nl/pub/10902/10902D.pdf}{7.2}], looks at representations of isotropic random fields with
homogeneous increments.
\\\\ 
The contour integral representation form involving the generalized integral form is given by;
\begin{multline}\label{eq:j1}
\frac{1}{2\pi i}\int_{C}\int_{0}^{b}a^w e^{-x^{\gamma } \theta +\beta  \left(1+x^{\rho } \xi \right)^{\lambda }} w^{-1-k} x^{m+w} \left(1+x^{\zeta } \eta \right)^{\mu } J_{\nu }(x \alpha )dxdw\\
=\frac{1}{2\pi i}\int_{C}\sum _{j=0}^{\infty } \sum
   _{h=0}^{\infty } \sum _{l=0}^{\infty } \sum _{p=0}^{\infty } \sum _{f=0}^{\infty } \frac{(-1)^j 2^{-2 j-\nu } a^w
   b^{1+2 j+m+w+p \gamma +f \zeta +\nu +h \rho } w^{-1-k}  }{(1+2 j+m+w+p \gamma +f \zeta +\nu +h \rho ) j! l! p! }\\ \times
\frac{\alpha ^{2 j+\nu } \beta ^l \eta ^f (-\theta )^p\xi ^h\binom{l \lambda }{h} \binom{\mu }{f}}{\Gamma (1+j+\nu
   )}dw
\end{multline}
where $Re(b)>0,Re(\zeta)>1/2$.
\subsection{Left-hand side contour integral representation}
Using a generalization of Cauchy's integral formula \ref{intro:cauchy} and the procedure in section (2.1), we form the definite integral given by;
\begin{multline}\label{eq:j2}
\int_{0}^{b}\frac{x^m \log ^k(a x) J_{\nu }(x \alpha ) \left(\eta  x^{\zeta }+1\right)^{\mu } e^{\beta  \left(\xi  x^{\rho
   }+1\right)^{\lambda }-\theta  x^{\gamma }}}{k!}dx\\
=\frac{1}{2\pi i}\int_{C}\int_{0}^{b}w^{-k-1} x^m (a x)^w J_{\nu }(x \alpha ) \left(\eta  x^{\zeta
   }+1\right)^{\mu } e^{\beta  \left(\xi  x^{\rho }+1\right)^{\lambda }-\theta  x^{\gamma }}dxdw
\end{multline}
where $Re(b)>0,Re(\zeta)>1/2$.
\subsection{Right-hand side contour integral}
Using a generalization of Cauchy's integral formula \ref{intro:cauchy}, and the procedure in section (2.2) we get;
\begin{multline}\label{eq:j3}
\sum _{j=0}^{\infty } \sum _{h=0}^{\infty } \sum _{l=0}^{\infty } \sum _{p=0}^{\infty } \sum _{f=0}^{\infty }
   \frac{(-1)^j 2^{-2 j-\nu } a^{-1-2 j-m-p \gamma -f \zeta -\nu -h \rho } \alpha ^{2 j+\nu } \beta ^l \eta ^f
   (-\theta )^p \xi ^h  \binom{l \lambda }{h} \binom{\mu }{f} }{j! k! l! p! \Gamma (1+j+\nu )(-1-2 j-m-p \gamma -f \zeta -\nu -h \rho )^{1+k}}\\ \times
\Gamma
   (1+k,-((1+2 j+m+p \gamma +f \zeta +\nu +h \rho ) \log (a b)))\\
=-\frac{1}{2\pi i}\int_{C}\sum _{j=0}^{\infty
   } \sum _{h=0}^{\infty } \sum _{l=0}^{\infty } \sum _{p=0}^{\infty } \sum _{f=0}^{\infty } \frac{(-1)^j 2^{-2 j-\nu
   } b^{1+2 j+m+p \gamma +f \zeta +\nu +h \rho } (a b)^w w^{-1-k} }{(1+2 j+m+w+p \gamma +f \zeta +\nu +h \rho ) j! l! p! }\\ \times 
\frac{ \alpha ^{2 j+\nu } \beta ^l\eta ^f (-\theta )^p \xi
   ^h \binom{l \lambda }{h} \binom{\mu }{f}}{\Gamma (1+j+\nu
   )}dw
\end{multline}
where $Re(b)>0,Re(\zeta)>1/2$.
from equation [Wolfram, \href{http://functions.wolfram.com/07.27.02.0002.01}{07.27.02.0002.01}] where $|Re(b)|<1$.
\begin{theorem}
Generalized integral.
\begin{multline}\label{eq:thm_prud}
\int_0^b e^{-x^{\gamma } \theta +\beta  \left(1+x^{\rho } \xi \right)^{\lambda }} x^m \left(1+x^{\zeta } \eta
   \right)^{\mu } J_{\nu }(x \alpha ) \log ^k\left(\frac{1}{a x}\right) \, dx\\
=\sum _{j=0}^{\infty } \sum
   _{h=0}^{\infty } \sum _{l=0}^{\infty } \sum _{p=0}^{\infty } \sum _{f=0}^{\infty } \frac{(-1)^j \alpha ^{2 j+\nu }
   \beta ^l \eta ^f (-\theta )^p \xi ^h \binom{l \lambda }{h} \binom{\mu }{f} }{j! l! p! \Gamma (1+j+\nu ) 2^{2 j+\nu } a^{1+2 j+m+p
   \gamma +f \zeta +\nu +h \rho } }\\ \times \frac{\Gamma \left(1+k,(1+2 j+m+p \gamma +f
   \zeta +\nu +h \rho ) \log \left(\frac{1}{a b}\right)\right)}{(1+2 j+m+p \gamma +f \zeta +\nu +h \rho )^{k+1}}
\end{multline}
where $Re(b)>0,Re(\zeta)>1/2$.
\end{theorem}
\begin{proof}
Since the right-hand sides of equations (\ref{eq:j2}) and (\ref{eq:j3}) are equal relative to equation (\ref{eq:j1}), we may equate the left-hand sides and simplify the gamma function to yield the stated result.
\end{proof}
\begin{example}
Derivation of generalized form of entry (2.12.10.3) in \cite{prud2}. Here we use equation (\ref{eq:thm_prud}) and set $k\to 0,a\to 1,\gamma \to 0,\theta \to 1,\lambda \to \frac{1}{2},\rho \to 2,\xi \to -\frac{1}{\xi ^2},\mu \to
   -\frac{1}{2},\zeta \to 2,\eta \to -\frac{1}{\eta ^2}$ and simplify;
   \begin{multline}\label{eq:exp_prud1}
\int_0^b \frac{e^{i \beta  \sqrt{-x^2+\xi ^2}} x^m J_{\nu }(x \alpha )}{\sqrt{-x^2+\eta ^2}} \, dx\\
=\sum
   _{f=0}^{\infty } \sum _{h=0}^{\infty } \sum _{l=0}^{\infty } \frac{(-1)^{f+h} i^l 2^{-\nu } b^{1+2 f+2 h+m+\nu
   } \alpha ^{\nu } \beta ^l \eta ^{-1-2 f} \xi ^{-2 h+l} \binom{-\frac{1}{2}}{f} \binom{\frac{l}{2}}{h} }{(1+2 f+2 h+m+\nu ) \Gamma (1+l) \Gamma (1+\nu )}\\ \times
\,
   _1F_2\left(\frac{1}{2}+f+h+\frac{m}{2}+\frac{\nu }{2};\frac{3}{2}+f+h+\frac{m}{2}+\frac{\nu }{2},1+\nu
   ;-\frac{1}{4} b^2 \alpha ^2\right)\\
=\sum _{f=0}^{\infty }
   \sum _{j=0}^{\infty } \sum _{l=0}^{\infty } \frac{(-1)^{f+j} i^l 2^{-2 j-\nu } b^{1+2 f+2 j+m+\nu } \alpha ^{2
   j+\nu } \beta ^l \eta ^{-1-2 f} \xi ^l \binom{-\frac{1}{2}}{f} }{(1+2 f+2 j+m+\nu ) \Gamma (1+j) \Gamma (1+l) \Gamma (1+j+\nu )}\\ \times
\,
   _2F_1\left(-\frac{l}{2},\frac{1}{2}+f+j+\frac{m}{2}+\frac{\nu }{2};\frac{3}{2}+f+j+\frac{m}{2}+\frac{\nu
   }{2};\frac{b^2}{\xi ^2}\right)
\end{multline}
where $Re(b)>0$.
\end{example}
\begin{example}
Derivation of entry (2.12.10.3) in \cite{prud2} page 188. Errata. In this example the Lommel function of two variables is employed [Wolfram, \href{https://mathworld.wolfram.com/LommelFunction.html}{6}]. Here we use equation (\ref{eq:exp_prud1}) and set $m\to 1,\beta \to p,\xi \to a,\nu \to 0,\alpha \to c,\eta \to a,b\to a$ and simplify;
\begin{multline}
\int_0^a \frac{x e^{i p \sqrt{a^2-x^2}} J_0(c x)}{\sqrt{a^2-x^2}} \, dx
=\sum _{l=0}^{\infty } \frac{i^l
   2^{-\frac{1}{2}+\frac{l}{2}} a^{1+l} (a c)^{-\frac{1}{2}-\frac{l}{2}} p^l J_{\frac{1+l}{2}}(a c) \Gamma
   \left(\frac{1+l}{2}\right)}{\Gamma (1+l)}\\
\neq\frac{\exp \left(i a \sqrt{p^2+c^2}\right)+2 i \sum\limits _{k=0}^{\infty }
   (-1)^k \left(\frac{-a p+a \sqrt{c^2+p^2}}{a c}\right)^{2 k} J_{2 k}(a c)-i J_0(a c)}{\sqrt{p^2+c^2}}
\end{multline}
where $Re(a)>0$.
\end{example}
\begin{example}
Derivation of generalized form of entry (2.12.10.4) in \cite{prud2}. Here we use equation (\ref{eq:thm_prud}) and set $k\to 0,a\to 1,\gamma \to 0,\theta \to 1,\lambda \to \frac{1}{2},\rho \to 2,\xi \to -\frac{1}{\xi ^2},\zeta \to
   2,\eta \to -\frac{1}{\eta ^2}$ and simplify;
   \begin{multline}\label{eq:exp_prud2}
\int_0^b e^{i \beta  \sqrt{-x^2+\xi ^2}} x^m \left(-x^2+\eta ^2\right)^{\mu } J_{\nu }(x \alpha ) \,
   dx\\
=\sum _{f=0}^{\infty } \sum _{j=0}^{\infty } \sum _{l=0}^{\infty } \frac{(-1)^{f+j} i^l 2^{-2 j-\nu } b^{1+2
   f+2 j+m+\nu } \alpha ^{2 j+\nu } \beta ^l \eta ^{-2 f+2 \mu } \xi ^l \binom{\mu }{f} }{(1+2 f+2 j+m+\nu ) \Gamma (1+j) \Gamma (1+l) \Gamma (1+j+\nu )}\\ \times
\,
   _2F_1\left(-\frac{l}{2},\frac{1}{2}+f+j+\frac{m}{2}+\frac{\nu }{2};\frac{3}{2}+f+j+\frac{m}{2}+\frac{\nu
   }{2};\frac{b^2}{\xi ^2}\right)\\
=\sum
   _{f=0}^{\infty } \sum _{h=0}^{\infty } \sum _{l=0}^{\infty } \frac{(-1)^{f+h} i^l 2^{-\nu } b^{1+2 f+2 h+m+\nu
   } \alpha ^{\nu } \beta ^l \eta ^{-2 f+2 \mu } \xi ^{-2 h+l} \binom{\frac{l}{2}}{h} \binom{\mu }{f} }{(1+2 f+2 h+m+\nu ) \Gamma (1+l) \Gamma (1+\nu )}\\ \times
\,
   _1F_2\left(\frac{1}{2}+f+h+\frac{m}{2}+\frac{\nu }{2};\frac{3}{2}+f+h+\frac{m}{2}+\frac{\nu }{2},1+\nu
   ;-\frac{1}{4} b^2 \alpha ^2\right)\\
=\sum _{h=0}^{\infty } \sum
   _{j=0}^{\infty } \sum _{l=0}^{\infty } \frac{(-1)^{h+j} i^l 2^{-2 j-\nu } b^{1+2 h+2 j+m+\nu } \alpha ^{2 j+\nu
   } \beta ^l \eta ^{2 \mu } \xi ^{-2 h+l} \binom{\frac{l}{2}}{h} }{(1+2 h+2 j+m+\nu ) \Gamma (1+j) \Gamma (1+l) \Gamma (1+j+\nu )}\\ \times
\, _2F_1\left(-\mu
   ,\frac{1}{2}+h+j+\frac{m}{2}+\frac{\nu }{2};\frac{3}{2}+h+j+\frac{m}{2}+\frac{\nu }{2};\frac{b^2}{\eta
   ^2}\right)
\end{multline}
where $Re(b)>0$.
\end{example}
\begin{example}
Derivation of entry (2.12.10.4) in \cite{prud2} in terms of infinite series and complex parameters. Here we use equation (\ref{eq:exp_prud2}) and set $\beta \to p,b\to a,\xi \to a,m\to n+1,\eta \to a,\nu \to n,\mu \to \frac{m-1}{2},\alpha \to c$ and simplify;
\begin{multline}
\int_0^a e^{i p \sqrt{a^2-x^2}} x^{1+n} \left(a^2-x^2\right)^{\frac{1}{2} (-1+m)} J_n(c x) \,
   dx\\
=\frac{a^n}{2}  \left(\frac{2 a}{c}\right)^{\frac{1+m}{2}} \sum _{l=0}^{\infty } \left(-\frac{2 a
   p^2}{c}\right)^{l/2}\frac{  \Gamma \left(\frac{1}{2} (1+l+m)\right)}{\Gamma
   (1+l)}J_{\frac{1}{2} (1+l+m)+n}(a c)
\end{multline}
where $Re(a)>0$.
\end{example}
\section{Definite integrals involving the product of four Bessel functions of the first kind}
Integrals of the product of four Bessel functions was first recorded in the work by Nicholson in (1920) and listed in Watson (1922) pages 414-415. Other authors have also published work on these integrals, namely; equation (10.2) on page 229 in \cite{ramanujan}
The contour integral representation form involving the generalized product of four Bessel functions is given by;
\begin{multline}\label{eq:4b1}
\frac{1}{2\pi i}\int_{C}\int_{0}^{b}a^w e^{x^{\tau } \theta } w^{-1-k} x^{m+w} \left(1+x^{\sigma } \phi \right)^{\lambda } J_{\mu }(x \alpha )
   J_{\nu }(x \beta ) J_{\xi }(x \gamma ) J_{\rho }(x \delta )dxdw\\
=\frac{1}{2\pi i}\int_{C}\sum _{f=0}^{\infty } \sum _{g=0}^{\infty } \sum
   _{h=0}^{\infty } \sum _{j=0}^{\infty } \sum _{l=0}^{\infty } \sum _{p=0}^{\infty } \frac{(-1)^{f+g+h+j} 2^{-2 f-2g-2 h-2 j-\mu -\nu -\xi -\rho } a^w }{(1+2 f+2 g+2 h+2 j+m+w+\mu +\nu +\xi +\rho +p \sigma +l \tau )}\\
\frac{b^{1+2 f+2 g+2 h+2 j+m+w+\mu +\nu +\xi +\rho +p \sigma +l \tau } w^{-1-k}\alpha ^{2 f+\mu } \beta ^{2 g+\nu } \gamma ^{2 h+\xi } \delta ^{2 j+\rho } \theta ^l \phi ^p \binom{\lambda}{p}}{ f! g! h! j! l! \Gamma (1+f+\mu ) \Gamma
   (1+g+\nu ) \Gamma (1+h+\xi ) \Gamma (1+j+\rho )}dw
\end{multline}
where $Re(\mu+\nu+\xi+\rho+m)>0,Re(b)>0,Re(\tau)>0$.
\subsection{Left-hand side contour integral representation}
Using a generalization of Cauchy's integral formula \ref{intro:cauchy} and the procedure in section (2.1), we form the definite integral given by;
\begin{multline}\label{eq:4b2}
\int_{0}^{b}\frac{x^m e^{\theta  x^{\tau }} \log ^k(a x) J_{\mu }(x \alpha ) J_{\nu }(x \beta ) J_{\xi }(x \gamma )
   J_{\rho }(x \delta ) \left(\phi  x^{\sigma }+1\right)^{\lambda }}{k!}dx\\
=\frac{1}{2\pi i}\int_{C}\int_{0}^{b}w^{-k-1} x^m (a x)^w e^{\theta  x^{\tau }}
   J_{\mu }(x \alpha ) J_{\nu }(x \beta ) J_{\xi }(x \gamma ) J_{\rho }(x \delta ) \left(\phi  x^{\sigma
   }+1\right)^{\lambda }dxdw
\end{multline}
where $Re(\mu+\nu+\xi+\rho+m)>0,Re(b)>0,Re(\tau)>0$.
\subsection{Right-hand side contour integral}
Using a generalization of Cauchy's integral formula \ref{intro:cauchy}, and the procedure in section (2.2) we get;
\begin{multline}\label{eq:4b3}
\sum _{f=0}^{\infty } \sum _{g=0}^{\infty } \sum _{h=0}^{\infty } \sum _{j=0}^{\infty } \sum _{l=0}^{\infty }
   \sum _{p=0}^{\infty } \frac{(-1)^{f+g+h+j} 2^{-2 f-2 g-2 h-2 j-\mu -\nu -\xi -\rho } a^{-1-2 f-2 g-2 h-2 j-m-\mu-\nu -\xi -\rho -p \sigma -l \tau }  }{f! g! h! j! k! l! \Gamma(1+f+\mu ) \Gamma (1+g+\nu ) \Gamma (1+h+\xi ) \Gamma (1+j+\rho )}\\ \times
(-1-2 f-2 g-2 h-2 j-m-\mu -\nu -\xi -\rho -p \sigma -l \tau )^{-1-k} \\
\Gamma(1+k,-((1+2 f+2 g+2 h+2 j+m+\mu +\nu +\xi +\rho +p \sigma +l \tau ) \log (a b)))\\ \times
\alpha ^{2 f+\mu } \beta ^{2 g+\nu } \gamma ^{2 h+\xi } \delta ^{2 j+\rho }\theta ^l \phi ^p \binom{\lambda }{p}\\
=-\frac{1}{2\pi i}\int_{C}\sum _{f=0}^{\infty } \sum _{g=0}^{\infty }
   \sum _{h=0}^{\infty } \sum _{j=0}^{\infty } \sum _{l=0}^{\infty } \sum _{p=0}^{\infty } \frac{(-1)^{f+g+h+j}
   2^{-2 f-2 g-2 h-2 j-\mu -\nu -\xi -\rho } }{(1+2 f+2 g+2 h+2 j+m+w+\mu +\nu +\xi +\rho +p \sigma +l \tau ) }\\ \times
\frac{a^w b^{1+2 f+2 g+2 h+2 j+m+w+\mu +\nu +\xi +\rho +p \sigma +l \tau }w^{-1-k} \alpha ^{2 f+\mu } \beta ^{2 g+\nu } \gamma ^{2 h+\xi } \delta ^{2 j+\rho } \theta ^l \phi ^p\binom{\lambda }{p}}{f! g! h! j! l! \Gamma (1+f+\mu) \Gamma (1+g+\nu ) \Gamma (1+h+\xi ) \Gamma (1+j+\rho )}dw
\end{multline}
where $Re(\mu+\nu+\xi+\rho+m)>0,Re(b)>0,Re(\tau)>0$.
from equation [Wolfram, \href{http://functions.wolfram.com/07.27.02.0002.01}{07.27.02.0002.01}] where $|Re(b)|<1$.
\begin{theorem}
Generalized integral with four Bessel functions.
\begin{multline}\label{eq:thm_4bessel}
\int_0^b e^{x^{\tau } \theta } x^m \left(1+x^{\sigma } \phi \right)^{\lambda } J_{\mu }(x \alpha ) J_{\nu }(x
   \beta ) J_{\xi }(x \gamma ) J_{\rho }(x \delta ) \log ^k\left(\frac{1}{a x}\right) \, dx\\
=\sum _{f=0}^{\infty }\sum _{g=0}^{\infty } \sum _{h=0}^{\infty } \sum _{j=0}^{\infty } \sum _{l=0}^{\infty } \sum _{p=0}^{\infty } \frac{(-1)^{f+g+h+j} 2^{-2 f-2 g-2 h-2 j-\mu -\nu -\xi -\rho }  }{f! g! h! j! l! \Gamma(1+f+\mu ) \Gamma (1+g+\nu ) \Gamma (1+h+\xi ) \Gamma (1+j+\rho )}\\ \times
\frac{ \Gamma \left(1+k,(1+2 f+2g+2 h+2 j+m+\mu +\nu +\xi +\rho +p \sigma +l \tau ) \log \left(\frac{1}{a b}\right)\right)}{(1+2 f+2g+2 h+2 j+m+\mu +\nu +\xi +\rho +p \sigma +l \tau )^{1+k}}\\ \times
\alpha ^{2 f+\mu } \beta ^{2 g+\nu } \gamma ^{2 h+\xi } \delta ^{2 j+\rho } \theta ^l  \phi ^p \binom{\lambda }{p}a^{-1-2 f-2 g-2 h-2 j-m-\mu -\nu -\xi -\rho -p\sigma -l \tau }
\end{multline}
where $Re(\mu+\nu+\xi+\rho+m)>0,Re(b)>0,Re(\tau)>0$.
\end{theorem}
\begin{proof}
Since the right-hand sides of equations (\ref{eq:4b2}) and (\ref{eq:4b3}) are equal relative to equation (\ref{eq:4b1}), we may equate the left-hand sides and simplify the gamma function to yield the stated result.
\end{proof}
\section{Evaluations of the infinite integral form}
\begin{example}
Generalized infinite form. Here we use equation (\ref{eq:thm_4bessel}) and set $k\to 0,a\to 1,\theta \to -\theta$ and take the sum over $l\in[0,\infty)$ and simplify the gamma function using equation [Wolfram, \href{http://functions.wolfram.com/06.07.03.0005.01}{01}];
\begin{multline}\label{eq:4b1}
\int_0^{\infty } e^{-x^{\tau } \theta } x^m \left(1+x^{\sigma } \phi \right)^{\lambda } J_{\mu }(x \alpha )
   J_{\nu }(x \beta ) J_{\xi }(x \gamma ) J_{\rho }(x \delta ) \, dx\\
=\sum _{f=0}^{\infty } \sum _{g=0}^{\infty } \sum _{h=0}^{\infty } \sum _{j=0}^{\infty } \sum _{p=0}^{\infty } \frac{(-1)^{f+g+h+j} 2^{-2 f-2 g-2 h-2 j-\mu -\nu-\xi -\rho } \alpha ^{2 f+\mu } \beta ^{2 g+\nu } \gamma ^{2 h+\xi } \delta ^{2 j+\rho } }{\tau  f! g! h! j! \Gamma (1+f+\mu ) \Gamma (1+g+\nu ) \Gamma(1+h+\xi ) \Gamma (1+j+\rho )}\\ \times
\theta ^{-\frac{1+2 f+2 g+2 h+2 j+m+\mu +\nu +\xi +\rho +p \sigma }{\tau }} \phi ^p \binom{\lambda }{p} \Gamma \left(\frac{1+2 f+2 g+2 h+2j+m+\mu +\nu +\xi +\rho +p \sigma }{\tau }\right)
\end{multline}
where $Re(\mu+\nu+\xi+\rho+m)>0,Re(\tau)>0$.
\end{example}
\begin{example}
An example from [StackExchange, \href{https://mathoverflow.net/questions/358552/integral-with-4-bessel-functions-and-an-exponential}{3}]. Here we use equation (\ref{eq:4b1}) and set $\tau \to 2,\theta \to a,m\to \lambda -1,\lambda \to 0,\sigma \to 0,\phi \to 1$ and simplify;
\begin{multline}
\int_0^{\infty } e^{-a x^2} x^{-1+\lambda } J_{\mu }(x \alpha ) J_{\nu }(x \beta ) J_{\xi }(x \gamma )
   J_{\rho }(x \delta ) \, dx\\
=\sum _{f=0}^{\infty } \sum _{g=0}^{\infty } \sum _{h=0}^{\infty }
   \frac{(-1)^{f+g+h} 2^{-1-2 f-2 g-2 h-\mu -\nu -\xi -\rho } a^{-f-g-h-\frac{\lambda }{2}-\frac{\mu
   }{2}-\frac{\nu }{2}-\frac{\xi }{2}-\frac{\rho }{2}} \alpha ^{2 f+\mu } \beta ^{2 g+\nu } \gamma ^{2 h+\xi }
   \delta ^{\rho } }{f! g! h! \Gamma (1+f+\mu ) \Gamma (1+g+\nu ) \Gamma (1+h+\xi ) \Gamma (1+\rho
   )}\\ \times
\Gamma \left(\frac{1}{2} (2 f+2 g+2 h+\lambda +\mu +\nu +\xi +\rho )\right)\\
 \,
   _1F_1\left(f+g+h+\frac{\lambda }{2}+\frac{\mu }{2}+\frac{\nu }{2}+\frac{\xi }{2}+\frac{\rho }{2};1+\rho
   ;-\frac{\delta ^2}{4 a}\right)
\end{multline}
where $Re(a)>0,Re(\nu+\nu+\xi+\rho)>-2,|Re(\lambda)|<1$.
\end{example}
\begin{example}
Laplace transform. Here we use equation (\ref{eq:4b1}) and set $\tau \to 1,\theta \to a,m\to 0,\lambda \to 0,\sigma \to 0,\phi \to 1$ and simplify;
\begin{multline}
\int_0^{\infty } e^{-a x} J_{\mu }(x \alpha ) J_{\nu }(x \beta ) J_{\xi }(x \gamma ) J_{\rho }(x \delta )
   \, dx\\
=\sum _{f=0}^{\infty } \sum _{g=0}^{\infty } \sum _{h=0}^{\infty } \frac{(-1)^{f+g+h} 2^{-2 f-2 g-2 h-\mu
   -\nu -\xi -\rho } a^{-1-2 f-2 g-2 h-\mu -\nu -\xi -\rho } \alpha ^{2 f+\mu } \beta ^{2 g+\nu } \gamma ^{2
   h+\xi } \delta ^{\rho } }{f! g! h! \Gamma (1+f+\mu ) \Gamma (1+g+\nu )
   \Gamma (1+h+\xi ) \Gamma (1+\rho )}\\ \times
\Gamma (1+2 f+2 g+2 h+\mu +\nu +\xi +\rho )\\
 \, _2F_1\left(\frac{1}{2}+f+g+h+\frac{\mu
   }{2}+\frac{\nu }{2}+\frac{\xi }{2}+\frac{\rho }{2},1+f+g+h+\frac{\mu }{2}+\frac{\nu }{2}+\frac{\xi
   }{2}+\frac{\rho }{2};1+\rho ;-\frac{\delta ^2}{a^2}\right)
\end{multline}
where $Re(a)>0,Re(\nu+\nu+\xi+\rho)>-2,|Re(\lambda)|<1$.
\end{example}
\section{Infinite integrals derived using integrals over the range $[0,1]$}
\begin{example}
Derivation of entry (6.579.3) over $[1,\infty)$. Errata. Here we use equation (\ref{eq:thm_4bessel}) and set $k\to 0,a\to 1,\tau \to 0,\theta \to 1,\lambda \to 0,\sigma \to 0,\phi \to 1,\mu \to \nu ,\xi \to \nu ,\rho \to \nu
   ,\alpha \to a,\beta \to a,\gamma \to a,\delta \to a,m\to 1-2 \nu ,b\to 1$ and expand integral over $[0,\infty)$ and simplify using equation [DLMF, \href{https://dlmf.nist.gov/15.4.E20}{15.4.20}];
   \begin{multline}
\int_1^{\infty } x^{1-2 \nu } J_{\nu }(a x){}^4 \, dx=\frac{a^{-2+2 \nu } \Gamma (\nu ) \Gamma (2 \nu )}{2 \pi\Gamma (3 \nu ) \Gamma \left(\frac{1}{2}+\nu \right)^2}\\
-\sum _{g=0}^{\infty } \sum _{h=0}^{\infty } \sum_{j=0}^{\infty } \frac{(-1)^{g+h+j} 2^{-1-2 g-2 h-2 j-4 \nu } a^{2 g+2 h+2 j+4 \nu } }{(1+g+h+j+\nu ) \Gamma (1+g) \Gamma (1+h) \Gamma (1+j) \Gamma (1+\nu ) }\\ \times
\frac{\, _1F_2\left(1+g+h+j+\nu;1+\nu ,2+g+h+j+\nu ;-\frac{a^2}{4}\right)}{ \Gamma (1+h+\nu ) \Gamma (1+j+\nu )\Gamma (1+g+\nu )}
\end{multline}
where $Re(a)>0,Re(b)>0,Re(\nu)>0$.
\end{example}
\begin{example}
Derivation of entry (6.579.4). Here we use equation (\ref{eq:thm_4bessel}) and set $k\to 0,a\to 1,\tau \to 0,\theta \to 1,\lambda \to 0,\sigma \to 0,\phi \to 1,\mu \to \nu ,\xi \to \nu ,\rho \to \nu
   ,\alpha \to a,\beta \to a,\gamma \to b,\delta \to b,m\to 1-2 \nu ,b\to 1$ and expand integral over $[0,\infty)$ and simplify;
   \begin{multline}
\int_1^{\infty } x^{1-2 \nu } J_{\nu }(a x){}^2 J_{\nu }(b x){}^2 \, dx
=\frac{\left(a^{2 \nu -1} \Gamma (\nu
   )\right) \, _2F_1\left(\nu ,\frac{1}{2}-\nu ;2 \nu +\frac{1}{2};\left(\frac{a}{b}\right)^2\right)}{2 \pi  b \Gamma
   \left(\nu +\frac{1}{2}\right) \Gamma \left(2 \nu +\frac{1}{2}\right)}\\
=-\sum _{f=0}^{\infty } \sum _{g=0}^{\infty }
   \sum _{h=0}^{\infty } \frac{(-1)^{f+g+h} 2^{-1-2 f-2 g-2 h-4 \nu } a^{2 f+2 g+2 \nu } b^{2 h+2 \nu } }{(1+f+g+h+\nu ) \Gamma (1+f) \Gamma (1+g) \Gamma(1+h) }\\
\frac{\,_1F_2\left(1+f+g+h+\nu ;1+\nu ,2+f+g+h+\nu ;-\frac{b^2}{4}\right)}{\Gamma (1+\nu ) \Gamma (1+f+\nu ) \Gamma (1+g+\nu ) \Gamma (1+h+\nu )}
\end{multline}
where $Re(a)>0,Re(b)>0,Re(\nu)>0$.
\end{example}
\section{Extended Lindsey integrals (1964)}
In 1964, Lindsey \cite{lindsey} evaluated an integral which arose in connection with a multiple-sample detection theory problem. The integral had been encountered in specifying the probability of error for a binary multilink signaling communication system which is connected to the Rician fading multichannel. In classical statistics the integral occurs when one attempts to find the probability distribution function for the ratio of the lengths of two correlated variates which, when two correlated, obey the noncentral and central chi-squared density functions of $m$ degrees of freedom. The main theorem in this section can be used to derive alternate series form for integrals listed in [Wolfram, \href{https://functions.wolfram.com/Bessel-TypeFunctions/BesselI/21/ShowAll.html}{BesselI}].
\\\\
The contour integral representation form involving the generalized Lindsey integral form is given by;
\begin{multline}\label{eq:k1}
\frac{1}{2\pi i}\int_{C}\int_{0}^{b}a^w e^{x^{\tau } \theta } w^{-1-k} x^{m+w} I_{\mu }\left(x^r \alpha \right) I_{\nu }\left(x^s \beta
   \right)dxdw\\
=\frac{1}{2\pi i}\int_{C}\sum _{f=0}^{\infty } \sum _{g=0}^{\infty } \sum _{l=0}^{\infty } \frac{2^{-2 f-2 g-\mu -\nu } a^wb^{1+m+w+r (2 f+\mu )+s (2 g+\nu )+l \tau } w^{-1-k} }{(1+m+w+r (2 f+\mu )+s (2 g+\nu )+l \tau ) f! g! l! }\\ \times
\frac{\alpha ^{2 f+\mu } \beta ^{2 g+\nu } \theta^l}{\Gamma (1+f+\mu ) \Gamma (1+g+\nu )}dw
\end{multline}
where $Re(b)>0$.
\subsection{Left-hand side contour integral representation}
Using a generalization of Cauchy's integral formula \ref{intro:cauchy} and the procedure in section (2.1), we form the definite integral given by;
\begin{multline}\label{eq:k2}
\int_{0}^{b}\frac{x^m e^{\theta  x^{\tau }} \log ^k(a x) I_{\mu }\left(x^r \alpha \right) I_{\nu }\left(x^s \beta
   \right)}{k!}dx\\
=\frac{1}{2\pi i}\int_{C}\int_{0}^{b}w^{-k-1} x^m (a x)^w e^{\theta  x^{\tau }} I_{\mu }\left(x^r \alpha \right) I_{\nu }\left(x^s
   \beta \right)dxdw
\end{multline}
where $Re(b)>0$.
\subsection{Right-hand side contour integral}
Using a generalization of Cauchy's integral formula \ref{intro:cauchy}, and the procedure in section (2.2) we get;
\begin{multline}\label{eq:k3}
\sum _{f=0}^{\infty } \sum _{g=0}^{\infty } \sum _{l=0}^{\infty } \frac{2^{-2 f-2 g-\mu -\nu } a^{-1-m-r
   (2 f+\mu )-s (2 g+\nu )-l \tau } \alpha ^{2 f+\mu } \beta ^{2 g+\nu } \theta ^l }{f! g! k! l! \Gamma (1+f+\mu ) \Gamma (1+g+\nu ) (-1-m-r (2 f+\mu)-s (2 g+\nu )-l \tau )^{k+1}}\\ \times
\Gamma (1+k,-((1+m+r (2
   f+\mu )+s (2 g+\nu )+l \tau ) \log (a b)))\\
=-\frac{1}{2\pi i}\int_{C}\sum _{f=0}^{\infty } \sum _{g=0}^{\infty } \sum _{l=0}^{\infty }
   \frac{2^{-2 f-2 g-\mu -\nu } b^{1+m+r (2 f+\mu )+s (2 g+\nu )+l \tau } (a b)^w w^{-1-k} \alpha ^{2 f+\mu }}{(1+m+w+r (2 f+\mu )+s (2 g+\nu )+l \tau ) f! g! l! }\\ \times
\frac{\beta ^{2 g+\nu } \theta ^l}{\Gamma (1+f+\mu ) \Gamma(1+g+\nu )}dw
\end{multline}
where $Re(b)>0$.
from equation [Wolfram, \href{http://functions.wolfram.com/07.27.02.0002.01}{07.27.02.0002.01}] where $|Re(b)|<1$.
\begin{theorem}
Generalized Lindsey integral.
\begin{multline}\label{eq:thm_lindsey}
\int_0^b e^{x^{\tau } \theta } x^m I_{\mu }\left(x^r \alpha \right) I_{\nu }\left(x^s \beta \right) \log^k(a x) \, dx\\
=\sum _{f=0}^{\infty } \sum _{g=0}^{\infty } \sum _{l=0}^{\infty } \frac{2^{-2 f-2 g-\mu -\nu }
   a^{-1-m-r (2 f+\mu )-s (2 g+\nu )-l \tau } \alpha ^{2 f+\mu } \beta ^{2 g+\nu } \theta ^l }{f! g! l! \Gamma (1+f+\mu ) \Gamma (1+g+\nu )(-1-m-r (2 f+\mu)-s (2 g+\nu )-l \tau )^{1+k} }\\ \times
\Gamma (1+k,-((1+m+r (2 f+\mu )+s (2 g+\nu )+l \tau ) \log (a b)))
\end{multline}
where $Re(b)>0,Re(\tau)>0$.
\end{theorem}
\begin{proof}
Since the right-hand sides of equations (\ref{eq:k2}) and (\ref{eq:k3}) are equal relative to equation (\ref{eq:k1}), we may equate the left-hand sides and simplify the gamma function to yield the stated result.
\end{proof}
\begin{example}
Here we use equation (\ref{eq:thm_lindsey}) and set $k\to 0,a\to 1,\theta \to -\theta$  and simplify using equation 
[Wolfram, \href{http://functions.wolfram.com/06.07.03.0005.01}{01}]
\begin{multline}\label{eq:thm_inf}
\int_0^{\infty } e^{-x^{\tau } \theta } x^m I_{\mu }\left(x^r \alpha \right) I_{\nu }\left(x^s \beta
   \right) \, dx\\
=\sum _{f=0}^{\infty } \sum _{g=0}^{\infty } \frac{2^{-2 f-2 g-\mu -\nu } \alpha ^{2 f+\mu }
   \beta ^{2 g+\nu } \theta ^{-\frac{1+m+2 f r+2 g s+r \mu +s \nu }{\tau }} \Gamma \left(\frac{1+m+2 f r+2 g
   s+r \mu +s \nu }{\tau }\right)}{\tau  f! g! \Gamma (1+f+\mu ) \Gamma (1+g+\nu )}
\end{multline}
where $Re(\tau)>0$.
\end{example}
\begin{example}
Derivation of equation (1) in \cite{lindsey}. Extended form involving infinite series of the hypergeometric function. Here we use equation (\ref{eq:thm_inf}) and set $\tau \to 2,\theta \to k,\mu \to m-1,r\to 1,\alpha \to b,\nu \to n,s\to 2,\beta \to \lambda$ and simplify.
\begin{multline}\label{eq:thm_inf1}
\int_0^{\infty } e^{-k x^2} x^m I_{-1+m}(b x) I_n\left(x^2 \lambda \right) \, dx\\
=\sum _{f=0}^{\infty } \frac{2^{-2 f-m-n} b^{-1+2 f+m} k^{-f-m-n} \lambda ^n \Gamma (f+m+n) }{\Gamma (1+f) \Gamma (f+m) \Gamma (1+n)}\\ \times
\,
   _2F_1\left(\frac{f}{2}+\frac{m}{2}+\frac{n}{2},\frac{1}{2}+\frac{f}{2}+\frac{m}{2}+\frac{n}{2};1+n;\frac{\lambda ^2}{k^2}\right)\\
=\sum _{g=0}^{\infty } \frac{2^{-2 g-m-n} b^{-1+m} k^{-2 g-m-n} \lambda ^{2 g+n}
   \Gamma (2 g+m+n) \, _1F_1\left(2 g+m+n;m;\frac{b^2}{4 k}\right)}{\Gamma (1+g) \Gamma (m) \Gamma (1+g+n)}
\end{multline}
where $Re(k)>0$.
\end{example}
\begin{example}
Derivation of equation (10(i)) in \cite{lindsey}. Extended form involving infinite series of the hypergeometric function. Here we use equation (\ref{eq:thm_inf1}) and set $n\to 0,m\to 1$ and simplify;
\begin{multline}
\int_0^{\infty } e^{-k x^2} x I_0(b x) I_0\left(x^2 \lambda \right) \, dx\\
=\sum _{f=0}^{\infty } \frac{2^{-1-2 f} b^{2 f} k^{-1-f} \, _2F_1\left(\frac{1}{2}+\frac{f}{2},1+\frac{f}{2};1;\frac{\lambda ^2}{k^2}\right)}{\Gamma (1+f)}\\
=\sum _{g=0}^{\infty
   } \frac{2^{-1-2 g} k^{-1-2 g} \lambda ^{2 g} \Gamma (1+2 g) \, _1F_1\left(1+2 g;1;\frac{b^2}{4 k}\right)}{\Gamma (1+g)^2}\\
=\frac{\sqrt{\frac{k^2}{k^2-\lambda ^2}} \exp \left(\frac{b^2 k}{4 \left(k^2-\lambda ^2\right)}\right) I_0\left(\frac{b^2 \lambda
   }{4 \left(k^2-\lambda ^2\right)}\right)}{2 k}
\end{multline}
where $Re(k)>0$.
\end{example}
\subsection{Table 3.15.16 in Prudnikov et al. volume  IV, Direct Laplace transforms, 1998, page 341 }
\begin{example}
Derivation of entry (3.15.16.1) in \cite{prud4}. Here we use equation (\ref{eq:thm_inf}) and set $\tau \to 1,\theta \to p,m\to 0,\mu \to \nu ,\alpha \to a,\beta \to b,r\to 1,s\to 1$ and simplify;
\begin{multline}
\int_0^{\infty } e^{-p x} I_{\nu }(a x) I_{\nu }(b x) \, dx\\
=\sum _{f=0}^{\infty } \frac{a^{2 f+\nu }
   b^{\nu } p^{-1-2 f-2 \nu } \Gamma \left(\frac{1}{2} (1+2 f+2 \nu )\right) \, _2F_1\left(\frac{1}{2}+f+\nu
   ,1+f+\nu ;1+\nu ;\frac{b^2}{p^2}\right)}{\sqrt{\pi } \Gamma (1+f) \Gamma (1+\nu )}\\
=\sum _{g=0}^{\infty }
   \frac{a^{\nu } b^{2 g+\nu } p^{-1-2 g-2 \nu } \Gamma \left(\frac{1}{2} (1+2 g+2 \nu )\right) \,
   _2F_1\left(\frac{1}{2}+g+\nu ,1+g+\nu ;1+\nu ;\frac{a^2}{p^2}\right)}{\sqrt{\pi } \Gamma (1+g) \Gamma (1+\nu
   )}\\
=\Re\left(\frac{Q_{\nu -\frac{1}{2}}\left(\frac{p^2-a^2-b^2}{2 a b}\right)}{\pi  \sqrt{a
   b}}\right)
\end{multline}
where $Re(\nu)>-1/2, Re(p)> |Re(a)|+|Re(b)|$.
\end{example}
\begin{example}
Derivation of entry (3.15.16.2) in \cite{prud4}. Here we use equation (\ref{eq:thm_inf}) and set $\tau \to 1,\theta \to p,m\to 0,\mu \to 0,\nu \to 0,\alpha \to a,\beta \to b,r\to 1,s\to 1$ and simplify; 
\begin{multline}
\int_0^{\infty } e^{-p x} I_0(a x) I_0(b x) \, dx=\sum _{f=0}^{\infty } \frac{a^{2 f} p^{-1-2 f} \Gamma
   \left(\frac{1}{2}+f\right) \, _2F_1\left(\frac{1}{2}+f,1+f;1;\frac{b^2}{p^2}\right)}{\sqrt{\pi } \Gamma
   (1+f)}\\
=\sum _{g=0}^{\infty } \frac{b^{2 g} p^{-1-2 g} \Gamma \left(\frac{1}{2}+g\right) \,
   _2F_1\left(\frac{1}{2}+g,1+g;1;\frac{a^2}{p^2}\right)}{\sqrt{\pi } \Gamma (1+g)}=\frac{2 K\left(\frac{4 a
   b}{-(a-b)^2+p^2}\right)}{\sqrt{-(a-b)^2+p^2} \pi }
\end{multline}
where $Re(p)> |Re(a)|+|Re(b)|$.
\end{example}
\begin{example}
Derivation of entry (3.15.16.3) in \cite{prud4}. Here we use equation (\ref{eq:thm_inf}) and set $\tau \to 1,\theta \to p,m\to -1,\mu \to \nu ,\alpha \to a,\beta \to a,r\to \frac{1}{2},s\to
   \frac{1}{2}$ and simplify; The series involving the hypergeometric function converges faster than the series involving the Bessel function listed in Prudnikov et al.
   \begin{multline}
\int_0^{\infty } \frac{e^{-p x} I_{\nu }\left(a \sqrt{x}\right){}^2}{x} \, dx\\
=\sum _{g=0}^{\infty }
   \frac{2^{-2 g-2 \nu } a^{2 g+2 \nu } p^{-g-\nu } \, _1F_1\left(g+\nu ;1+\nu ;\frac{a^2}{4 p}\right)}{(g+\nu
   ) \Gamma (1+g) \Gamma (1+\nu )}\\
=\frac{\exp \left(\frac{a^2}{2 p}\right) \left(I_{\nu }\left(\frac{a^2}{2
   p}\right)+2 \sum\limits _{k=1}^{\infty } (-1)^k I_{\nu +k}\left(\frac{a^2}{2 p}\right)\right)}{\nu }
\end{multline}
where $Re(p)>0$.
\end{example}
\begin{example}
Derivation of entry (3.15.16.4) in \cite{prud4}. Errata. Here we use equation (\ref{eq:thm_inf}) and set $\tau \to 1,\theta \to p,m\to 1,\mu \to n+2,\nu \to n,\alpha \to a,\beta \to b,r\to \frac{1}{2},s\to
   \frac{1}{2}$ and simplify;
   \begin{multline}
\int_0^{\infty } e^{-p x} x I_n\left(b \sqrt{x}\right) I_{2+n}\left(a \sqrt{x}\right) \, dx\\
=\sum_{g=0}^{\infty } \frac{2^{-2-2 g-2 n} a^{2+n} b^{2 g+n} (1+g+n) (2+g+n) p^{-2-g+\frac{1}{2}
   (-2-n)-\frac{n}{2}} }{\Gamma (1+g) \Gamma (3+n)}\\ \times
\, _1F_1\left(3+g+n;3+n;\frac{a^2}{4 p}\right)\\
=\frac{\exp\left(\frac{a^2+b^2}{4 p}\right) }{4 a^2 p^3}\\ \times
\left(a b \left(a^2+4 p\right) I_{n+1}\left(\frac{a b}{2 p}\right)-a b
   \left(4 p n+8 p-a^2\right) I_{n-1}\left(\frac{a b}{2 p}\right)\right. \\ \left.
+\left(16 p^2 n (n+2)-\left(4 p n-b^2\right)
   a^2+a^4\right) I_n\left(\frac{a b}{2 p}\right)\right)
\end{multline}
where $Re(p)>0$.
\end{example}
\begin{example}
Derivation of entry (3.15.16.5) in \cite{prud4}. Errata. Here we use equation (\ref{eq:thm_inf}) and set $\tau \to 1,\theta \to p,m\to 1,\mu \to 0,\nu \to 0,\alpha \to a,\beta \to a,r\to \frac{1}{2},s\to
   \frac{1}{2}$ and simplify;
   \begin{multline}
\int_0^{\infty } e^{-p x} x I_0\left(a \sqrt{x}\right){}^2 \, dx\\
=\sum _{g=0}^{\infty } \frac{2^{-2 g}
   a^{2 g} (1+g) p^{-2-g} \, _1F_1\left(2+g;1;\frac{a^2}{4 p}\right)}{\Gamma (1+g)}\\
=\frac{\exp
   \left(\frac{a^2}{2 p}\right) \left(\left(2+\frac{a^2}{p}\right) I_0\left(\frac{a^2}{2 p}\right)+\frac{a^2
   I_1\left(\frac{a^2}{2 p}\right)}{p}\right)}{2 p^2}
\end{multline}
where $Re(p)>0$.
\end{example}
\begin{example}
Derivation of entry (3.15.16.6) in \cite{prud4}. Errata. Here we use equation (\ref{eq:thm_inf}) and set $\tau \to 1,\theta \to p,m\to 2,\mu \to 0,\nu \to 0,\alpha \to a,\beta \to a,r\to \frac{1}{2},s\to
   \frac{1}{2}$ and simplify;
   \begin{multline}
\int_0^{\infty } e^{-p x} x^2 I_0\left(a \sqrt{x}\right){}^2 \, dx\\
=\sum _{g=0}^{\infty } \frac{2^{-2 g}
   a^{2 g} (1+g) (2+g) p^{-3-g} \, _1F_1\left(3+g;1;\frac{a^2}{4 p}\right)}{\Gamma (1+g)}\\
=\frac{\exp
   \left(\frac{a^2}{2 p}\right) \left(\left(4+\frac{4 a^2}{p}+\frac{a^4}{p^2}\right) I_0\left(\frac{a^2}{2
   p}\right)+\frac{a^2 \left(3+\frac{a^2}{p}\right) I_1\left(\frac{a^2}{2 p}\right)}{p}\right)}{2 p^3}
\end{multline}
where $Re(p)>0$.
\end{example}
\begin{example}
 Derivation of entry (3.15.16.7) in \cite{prud4}. Here we use equation (\ref{eq:thm_inf}) and set $\tau \to 1,\theta \to p,m\to -\frac{1}{2},\mu \to 0,\nu \to 1,r\to \frac{1}{2},s\to \frac{1}{2},\beta
   \to a,\alpha \to a$ and simplify;
   \begin{multline}
\int_0^{\infty } \frac{e^{-p x} I_0\left(a \sqrt{x}\right) I_1\left(a \sqrt{x}\right)}{\sqrt{x}} \,
   dx\\
=\sum _{g=0}^{\infty } \frac{2^{-1-2 g} a^{1+2 g} p^{-1-g} \, _1F_1\left(1+g;1;\frac{a^2}{4
   p}\right)}{\Gamma (2+g)}\\
=\frac{\exp \left(\frac{a^2}{2 p}\right) I_0\left(\frac{a^2}{2
   p}\right)-1}{a}
\end{multline}
where $Re(p)>0$.
\end{example}
\begin{example}
Derivation of entry (3.15.16.8) in \cite{prud4}. Here we use equation (\ref{eq:thm_inf}) and set $\tau \to 1,\theta \to p,m\to \frac{1}{2},\mu \to 0,\nu \to 1,r\to \frac{1}{2},s\to \frac{1}{2},\beta \to
   a,\alpha \to a$ and simplify;
   \begin{multline}
\int_0^{\infty } e^{-p x} \sqrt{x} I_0\left(a \sqrt{x}\right) I_1\left(a \sqrt{x}\right) \, dx\\
=\sum
   _{g=0}^{\infty } \frac{2^{-1-2 g} a^{1+2 g} p^{-2-g} \, _1F_1\left(2+g;1;\frac{a^2}{4 p}\right)}{g!}\\
=\frac{a
   \exp \left(\frac{a^2}{2 p}\right) \left(I_0\left(\frac{a^2}{2 p}\right)+I_1\left(\frac{a^2}{2
   p}\right)\right)}{2 p^2}
\end{multline}
where $Re(p)>0$.
\end{example}
\begin{example}
Derivation of entry (3.15.16.9) in \cite{prud4}. Here we use equation (\ref{eq:thm_inf}) and set $\tau \to 1,\theta \to p,m\to \frac{3}{2},\mu \to 0,\nu \to 1,r\to \frac{1}{2},s\to \frac{1}{2},\beta \to
   a,\alpha \to a$ and simplify;
   \begin{multline}
\int_0^{\infty } e^{-p x} x^{3/2} I_0\left(a \sqrt{x}\right) I_1\left(a \sqrt{x}\right) \, dx\\
=\sum
   _{g=0}^{\infty } \frac{2^{-1-2 g} a^{1+2 g} (2+g) p^{-3-g} \, _1F_1\left(3+g;1;\frac{a^2}{4
   p}\right)}{\Gamma (1+g)}\\
=\frac{a \exp \left(\frac{a^2}{2 p}\right) \left(\left(2+\frac{a^2}{p}\right)
   I_0\left(\frac{a^2}{2 p}\right)+\left(1+\frac{a^2}{p}\right) I_1\left(\frac{a^2}{2 p}\right)\right)}{2
   p^3}
\end{multline}
where $Re(p)>0$.
\end{example}
\begin{example}
Derivation of entry (3.15.16.10) in \cite{prud4}. Here we use equation (\ref{eq:thm_inf}) and set $\tau \to 1,\theta \to p,m\to -1,\mu \to 1,\nu \to 1,\alpha \to a,\beta \to a,r\to \frac{1}{2},s\to
   \frac{1}{2}$ and simplify;
   \begin{multline}
\int_0^{\infty } \frac{e^{-p x} I_1\left(a \sqrt{x}\right){}^2}{x} \, dx\\
=\sum _{g=0}^{\infty }
   \frac{2^{-2-2 g} a^{2+2 g} p^{-1-g} \, _1F_1\left(1+g;2;\frac{a^2}{4 p}\right)}{\Gamma (2+g)}\\
=\exp
   \left(\frac{a^2}{2 p}\right) \left(I_0\left(\frac{a^2}{2 p}\right)-I_1\left(\frac{a^2}{2
   p}\right)\right)-1
\end{multline}
where $Re(p)>0$.
\end{example}
\begin{example}
Derivation of entry (3.15.16.11) in \cite{prud4}. Here we use equation (\ref{eq:thm_inf}) and set $\tau \to 1,\theta \to p,m\to 1,\mu \to 1,\nu \to 1,\alpha \to a,\beta \to a,r\to \frac{1}{2},s\to
   \frac{1}{2}$ and simplify;
   \begin{multline}
\int_0^{\infty } e^{-p x} x I_1\left(a \sqrt{x}\right){}^2 \, dx\\
=\sum _{g=0}^{\infty } \frac{2^{-2-2 g}
   a^{2+2 g} (2+g) p^{-3-g} \, _1F_1\left(3+g;2;\frac{a^2}{4 p}\right)}{\Gamma (1+g)}\\
=\frac{a^2 \exp
   \left(\frac{a^2}{2 p}\right) \left(I_0\left(\frac{a^2}{2 p}\right)+I_1\left(\frac{a^2}{2 p}\right)\right)}{2
   p^3}
\end{multline}
where $Re(p)>0$.
\end{example}
\begin{example}
Derivation of entry (3.15.16.12) in \cite{prud4}. Errata. Here we use equation (\ref{eq:thm_inf}) and set $\tau \to 1,\theta \to p,m\to 2,\mu \to 1,\nu \to 1,\alpha \to a,\beta \to a,r\to \frac{1}{2},s\to
   \frac{1}{2}$ and simplify; 
   \begin{multline}
\int_0^{\infty } e^{-p x} x^2 I_1\left(a \sqrt{x}\right){}^2 \, dx\\
=\sum _{g=0}^{\infty } \frac{2^{-2-2
   g} a^{2+2 g} (2+g) (3+g) p^{-4-g} \, _1F_1\left(4+g;2;\frac{a^2}{4 p}\right)}{\Gamma (1+g)}\\
=\frac{a^2 \exp\left(\frac{a^2}{2 p}\right) \left(\left(3+\frac{a^2}{p}\right) I_0\left(\frac{a^2}{2
   p}\right)+\left(2+\frac{a^2}{p}\right) I_1\left(\frac{a^2}{2 p}\right)\right)}{2 p^4}
\end{multline}
where $Re(p)>0$.
\end{example}
\begin{example}
Derivation of entry (3.15.16.15) in \cite{prud4}. Here we use equation (\ref{eq:thm_inf}) and set $\tau \to 1,\theta \to p,m\to 0,\mu \to 2 \nu ,r\to \frac{1}{2},\alpha \to a,s\to 1,\beta \to b$ and simplify; 
\begin{multline}
\int_0^{\infty } e^{-p x} I_{\nu }(b x) I_{2 \nu }\left(a \sqrt{x}\right) \, dx\\
=\sum _{g=0}^{\infty }
   \frac{2^{-\nu } a^{2 \nu } b^{2 g+\nu } p^{-1-2 g-2 \nu } \Gamma \left(\frac{1}{2} (1+2 g+2 \nu )\right) \,
   _1F_1\left(1+2 g+2 \nu ;1+2 \nu ;\frac{a^2}{4 p}\right)}{\sqrt{\pi } \Gamma (1+g) \Gamma (1+2 \nu
   )}\\
=\frac{\exp \left(\frac{a^2 p}{4 p^2-4 b^2}\right) I_{\nu }\left(\frac{a^2 b}{4 p^2-4
   b^2}\right)}{\sqrt{p^2-b^2}}
\end{multline}
where $Re(p)>0$.
\end{example}
\subsection{Table 4.9.7 in Brychkov et al. (2008) page 214}
\begin{example}
Derivation of entry (4.9.7.12) in \cite{brychkov}. Here we use equation (\ref{eq:thm_lindsey}) and set $k\to 1,a\to 1,b\to 1,\tau \to 0,\theta \to 1,m\to 1,\mu \to 0,\nu \to 0,r\to 1,s\to 1,\alpha \to a,\beta \to
   a$ and simplify; 
   \begin{multline}
\int_0^1 x I_0(a x){}^2 \log (x) \, dx\\
=-\sum _{g=0}^{\infty } \frac{2^{-2-2 g} a^{2 g} \,
   _2F_3\left(1+g,1+g;1,2+g,2+g;\frac{a^2}{4}\right)}{(1+g)^2 \Gamma (1+g)^2}\\
=\frac{1}{2}
   \left(I_1(a){}^2-I_0(a){}^2+\frac{I_0(a) I_1(a)}{a}\right)
\end{multline}
where $a\in\mathbb{C}$.
\end{example}
\begin{example}
Derivation of entry (4.9.7.13) in \cite{brychkov}. Here we use equation (\ref{eq:thm_lindsey}) and set $k\to 1,a\to 1,b\to 1,\tau \to 0,\theta \to 1,m\to 1,\mu \to 1,\nu \to 1,r\to 1,s\to 1,\alpha \to a,\beta \to
   a$ and simplify; 
   \begin{multline}
\int_0^1 x I_1(a x){}^2 \log (x) \, dx\\
=-\sum _{g=0}^{\infty } \frac{2^{-4-2 g} a^{2+2 g} \,
   _2F_3\left(2+g,2+g;2,3+g,3+g;\frac{a^2}{4}\right)}{(2+g)^2 \Gamma (1+g) \Gamma
   (2+g)}\\
=\frac{1+\left(a^2-1\right) I_0(a){}^2-a I_0(a) I_1(a)-a^2 I_1(a){}^2}{2 a^2}
\end{multline}
where $a\in\mathbb{C}$.
\end{example}
\begin{example}
Derivation of entry [DLMF, \href{https://dlmf.nist.gov/10.43.E28}{10.43.28}]. Here we use equation (\ref{eq:thm_inf}) and set $\tau \to 2,\theta \to p^2,m\to 1,\mu \to \nu ,\alpha \to a,\beta \to b,r\to 1,s\to 1,x\to t$ and simplify;
\begin{multline}
\int_0^{\infty } e^{-p^2 t^2} t I_{\nu }(a t) I_{\nu }(b t) \, dt\\
=\sum _{g=0}^{\infty } \frac{2^{-1-2
   g-2 \nu } a^{\nu } b^{2 g+\nu } p^{-2-2 g-2 \nu } \, _1F_1\left(1+g+\nu ;1+\nu ;\frac{a^2}{4
   p^2}\right)}{\Gamma (1+g) \Gamma (1+\nu )}\\
=\frac{\exp \left(\frac{a^2+b^2}{4 p^2}\right) I_{\nu
   }\left(\frac{a b}{2 p^2}\right)}{2 p^2}
\end{multline}
where $Re(\nu)>-1,Re(p^2)>0$.
\end{example}
\subsection{Table 2.15.20 in Prudnikov et al volume II (1986b)}
In this section we derive entries for infinite integrals in terms of infinite series involving special functions and fundamental constants.
\begin{example}
Derivation of entry (2.15.20.2) in \cite{prud2}. Here we use equation (\ref{eq:thm_inf}) and set $\tau \to 1,\theta \to p,m\to \alpha -1,r\to 1,s\to 1,\alpha \to b,\beta \to c$ and simplify using [Wolfram, \href{https://reference.wolfram.com/language/ref/AppellF4.html}{2}];
\begin{multline}
\int_0^{\infty } e^{-p x} x^{-1+\alpha } I_{\mu }(b x) I_{\nu }(c x) \, dx\\
=\sum _{g=0}^{\infty }
   \frac{2^{-2 g-\mu -\nu } b^{\mu } c^{2 g+\nu } p^{-2 g-\alpha -\mu -\nu } \Gamma (2 g+\alpha +\mu +\nu ) }{\Gamma (1+g) \Gamma (1+\mu ) \Gamma (1+g+\nu)}\\ \times
\,_2F_1\left(g+\frac{\alpha }{2}+\frac{\mu }{2}+\frac{\nu }{2},\frac{1}{2}+g+\frac{\alpha }{2}+\frac{\mu }{2}+\frac{\nu }{2};1+\mu ;\frac{b^2}{p^2}\right)\\
=\frac{\left(b^{\mu } c^{\nu }\right) \Gamma (\alpha +\mu +\nu ) }{\left(2^{\mu +\nu } p^{\alpha +\mu +\nu }\right) (\Gamma (\mu +1) \Gamma (\nu +1))}\\ \times
\sum\limits_{m=0}^{\infty } \sum\limits_{n=0}^{\infty
   } \frac{\left(\frac{b^2}{p^2}\right)^m \left(\frac{c^2}{p^2}\right)^n \left(\frac{1}{2} (\alpha +\mu +\nu
   )\right)_{m+n} \left(\frac{1}{2} (1+\alpha +\mu +\nu )\right)_{m+n}}{m! n! (1+\mu )_m (1+\nu
   )_n}\\
=\frac{\left(b^{\mu
   } c^{\nu }\right) p^{-\alpha -\mu -\nu } }{2^{\mu +\nu } \Gamma (\nu +1)}\sum\limits_{k=0}^{\infty } \frac{\Gamma (\alpha +\mu +\nu +2 k) \,
   _2F_1\left(-k,-\mu -k;\nu +1;\left(\frac{c}{b}\right)^2\right) \left(\frac{b}{2 p}\right)^{2 k}}{k! \Gamma
   (\mu +k+1)}
\end{multline}
where $Re(\nu+\mu)>-1,Re(p)>0$.
\end{example}
\begin{example}
Derivation of entry (2.15.20.3) in \cite{prud2}. Errata. Here we use equation (\ref{eq:thm_inf}) and set $\tau \to 1,\theta \to p,m\to -\frac{1}{2},r\to 1,s\to 1,\alpha \to b,\beta \to c$ and simplify;
\begin{multline}
\int_0^{\infty } \frac{e^{-p x} I_{\mu }(b x) I_{\nu }(c x)}{\sqrt{x}} \, dx\\
=\sum _{g=0}^{\infty }
   \frac{2^{-2 g-\mu -\nu } b^{\mu } c^{2 g+\nu } p^{-\frac{1}{2}-2 g-\mu -\nu } \Gamma \left(\frac{1}{2} (1+4
   g+2 \mu +2 \nu )\right) }{\Gamma (1+g) \Gamma (1+\mu ) \Gamma (1+g+\nu
   )}\\ \times
\, _2F_1\left(\frac{1}{4}+g+\frac{\mu }{2}+\frac{\nu }{2},\frac{3}{4}+g+\frac{\mu
   }{2}+\frac{\nu }{2};1+\mu ;\frac{b^2}{p^2}\right)\\
\neq\frac{\sqrt{2} \Gamma \left(\frac{1}{2}+\mu +\nu \right) }{\sqrt{\sqrt{-(b-c)^2+p^2}+\sqrt{-(b+c)^2+p^2}}}\\ \times
P_{-\frac{1}{2}+\mu }^{-\nu }\left(\frac{2\sqrt{c^2+\frac{1}{4}
   \left(\sqrt{-(b-c)^2+p^2}+\sqrt{-(b+c)^2+p^2}\right)^2}}{\sqrt{-(b-c)^2+p^2}+\sqrt{-(b+c)^2+p^2}}\right)\\
\times   P_{-\frac{1}{2}+\nu }^{-\mu }\left(\frac{2 \sqrt{b^2+\frac{1}{4}
   \left(\sqrt{-(b-c)^2+p^2}+\sqrt{-(b+c)^2+p^2}\right)^2}}{\sqrt{-(b-c)^2+p^2}+\sqrt{-(b+c)^2+p^2}}\right)
\end{multline}
where $Re(\nu+\mu)>-1,Re(p)>0$.
\end{example}
\begin{example}
Derivation of entry (2.15.20.4) in \cite{prud2}. Generalized form. Here we use equation (\ref{eq:thm_inf}) and set $\tau \to 1,\theta \to p,m\to 0,r\to 1,s\to 1,\alpha \to c,\beta \to c$ and simplify;
\begin{multline}
\int_0^{\infty } e^{-p x} I_{\mu }(c x) I_{\nu }(c x) \, dx\\
=\sum _{g=0}^{\infty } \frac{2^{-2 g-\mu -\nu
   } c^{2 g+\mu +\nu } p^{-1-2 g-\mu -\nu } \Gamma (1+2 g+\mu +\nu ) }{\Gamma (1+g) \Gamma
   (1+\mu ) \Gamma (1+g+\nu )}\\ \times
\, _2F_1\left(\frac{1}{2}+g+\frac{\mu
   }{2}+\frac{\nu }{2},1+g+\frac{\mu }{2}+\frac{\nu }{2};1+\mu ;\frac{c^2}{p^2}\right)
\end{multline}
where $Re(\nu+\mu)>-1,Re(p)>0$.
\end{example}
\begin{example}
Derivation of entry (2.15.20.5) in \cite{prud2}. Here we use equation (\ref{eq:thm_inf}) and set $\tau \to 1,\theta \to p,m\to \alpha -1,r\to 1,s\to 1,\alpha \to c,\beta \to c$ and simplify;
\begin{multline}
\int_0^{\infty } e^{-p x} x^{-1+\alpha } I_{\mu }(c x) I_{\nu }(c x) \, dx\\
=\sum _{g=0}^{\infty }
   \frac{2^{-2 g-\mu -\nu } c^{2 g+\mu +\nu } p^{-2 g-\alpha -\mu -\nu } \Gamma (2 g+\alpha +\mu +\nu ) }{\Gamma (1+g) \Gamma (1+\mu ) \Gamma (1+g+\nu
   )}\\ \times
\,_2F_1\left(g+\frac{\alpha }{2}+\frac{\mu }{2}+\frac{\nu }{2},\frac{1}{2}+g+\frac{\alpha }{2}+\frac{\mu
   }{2}+\frac{\nu }{2};1+\mu ;\frac{c^2}{p^2}\right)\\
=\frac{\left(\frac{c}{2}\right)^{\mu +\nu } p^{-\alpha -\mu -\nu } \Gamma (\alpha +\mu +\nu ) 
}{\Gamma (\mu +1) \Gamma(\nu +1)}\\ \times
\,_4F_3\left(\frac{1}{2} (\mu +\nu +1),\frac{\mu +\nu }{2}+1,\frac{1}{2} (\alpha +\mu +\nu ),\frac{1}{2}
   (\alpha +\mu +\nu +1);\right. \\ \left.
\mu +\nu +1,\nu +1,\mu +1;\left(\frac{2 c}{p}\right)^2\right)
\end{multline}
where $Re(\nu+\mu)>-1,Re(p)>0$.
\end{example}
\begin{example}
Derivation of entry (2.15.20.7) in \cite{prud2}. Here we use equation (\ref{eq:thm_inf}) and set $\tau \to 2,\theta \to p,m\to \alpha -1,r\to 1,s\to 1,\alpha \to b,\beta \to c$ and simplify;
\begin{multline}
\int_0^{\infty } e^{-p x^2} x^{-1+\alpha } I_{\mu }(b x) I_{\nu }(c x) \, dx\\
=\sum _{g=0}^{\infty }
   \frac{2^{-1-2 g-\mu -\nu } b^{\mu } c^{2 g+\nu } p^{-g-\frac{\alpha }{2}-\frac{\mu }{2}-\frac{\nu }{2}}
   \Gamma \left(\frac{1}{2} (2 g+\alpha +\mu +\nu )\right) }{\Gamma (1+g) \Gamma (1+\mu ) \Gamma (1+g+\nu
   )}\\ \times
\, _1F_1\left(g+\frac{\alpha }{2}+\frac{\mu
   }{2}+\frac{\nu }{2};1+\mu ;\frac{b^2}{4 p}\right)\\
=\frac{\left(b^{\mu } c^{\nu } p^{-\frac{1}{2} (\alpha +\mu +\nu )}\right) }{2^{\mu +\nu +1}
   \Gamma (\nu +1)}\sum\limits_{k=0}^{\infty }
   \frac{\Gamma \left(k+\frac{1}{2} (\alpha +\mu +\nu )\right) \left(\frac{b}{2 p}\right)^{2 k} 
}{\Gamma (\mu +k+1) k!}\\ \times\,_2F_1\left(-k,-\mu -k;\nu +1;\left(\frac{c}{b}\right)^2\right)
\end{multline}
where $Re(\nu+\mu)>-1,Re(p)>0$.
\end{example}
\begin{example}
Derivation of entry (2.15.20.8) in \cite{prud2}. Here we use equation (\ref{eq:thm_inf}) and set $\tau \to 2,\theta \to p,m\to 1,r\to 1,s\to 1,\alpha \to b,\beta \to c,\mu \to \nu $ and simplify;
\begin{multline}
\int_0^{\infty } e^{-p x^2} x I_{\nu }(b x) I_{\nu }(c x) \, dx\\
=\sum _{g=0}^{\infty } \frac{2^{-1-2 g-2
   \nu } b^{\nu } c^{2 g+\nu } p^{-1-g-\nu } \, _1F_1\left(1+g+\nu ;1+\nu ;\frac{b^2}{4 p}\right)}{\Gamma (1+g)
   \Gamma (1+\nu )}\\
=\frac{\exp \left(\frac{b^2+c^2}{4 p}\right) I_{\nu }\left(\frac{b c}{2 p}\right)}{2
   p}
\end{multline}
where $Re(\nu+\mu)>-1,Re(p)>0$.
\end{example}
\begin{example}
Derivation of entry (2.15.20.9) in \cite{prud2}. Here we use equation (\ref{eq:thm_inf}) and set $\tau \to 2,\theta \to p,m\to 3,r\to 1,s\to 1,\alpha \to b,\beta \to c,\mu \to n+2,\nu \to n$ and simplify;
\begin{multline}
\int_0^{\infty } e^{-p x^2} x^3 I_n(c x) I_{2+n}(b x) \, dx\\
=\sum _{g=0}^{\infty } \frac{2^{-3-2 g-2 n}
   b^{2+n} c^{2 g+n} (1+g+n) (2+g+n) p^{-3-g-n} \, _1F_1\left(3+g+n;3+n;\frac{b^2}{4 p}\right)}{\Gamma (1+g)
   \Gamma (3+n)}\\
=\frac{\exp \left(\frac{b^2+c^2}{4 p}\right) }{8 p^3 b^2}\left(b c \left(4 p+b^2\right)
   I_{n+1}\left(\frac{b c}{2 p}\right)-b c \left(4 p n+8 p-b^2\right) I_{n-1}\left(\frac{b c}{2
   p}\right)\right. \\ \left.
+\left(16 p^2 n (n+2)-\left(4 p n-c^2\right) b^2+b^4\right) I_n\left(\frac{b c}{2
   p}\right)\right)
\end{multline}
where $Re(\nu+\mu)>-1,Re(p)>0$.
\end{example}
\begin{example}
 Derivation of entry (2.15.20.10) in \cite{prud2}. Here we use equation (\ref{eq:thm_inf}) and set $\tau \to 2,\theta \to p,m\to \alpha -1,r\to 1,s\to 1,\alpha \to c,\beta \to c$ and simplify;
 \begin{multline}
\int_0^{\infty } e^{-p x^2} x^{-1+\alpha } I_{\mu }(c x) I_{\nu }(c x) \, dx\\
=\sum _{g=0}^{\infty }
   \frac{2^{-1-2 g-\mu -\nu } c^{2 g+\mu +\nu } p^{-g-\frac{\alpha }{2}-\frac{\mu }{2}-\frac{\nu }{2}} \Gamma
   \left(\frac{1}{2} (2 g+\alpha +\mu +\nu )\right) }{\Gamma (1+g) \Gamma (1+\mu ) \Gamma (1+g+\nu )}\\ \times
\, _1F_1\left(g+\frac{\alpha }{2}+\frac{\mu }{2}+\frac{\nu
   }{2};1+\mu ;\frac{c^2}{4 p}\right)\\
=\frac{2^{-\mu -\nu -1}
   c^{\mu +\nu } p^{-\frac{1}{2} (\alpha +\mu +\nu )} \Gamma \left(\frac{1}{2} (\alpha +\mu +\nu )\right) }{\Gamma (\mu +1) \Gamma (\nu +1)}\\ \times
\,
   _3F_3\left(\frac{1}{2} (\mu +\nu +1),\frac{\mu +\nu }{2}+1,\frac{1}{2} (\alpha +\mu +\nu );\mu +1,\nu +1,\mu
   +\nu +1;\frac{c^2}{p}\right)
\end{multline}
where $Re(\nu+\mu)>-1,Re(p)>0$.
\end{example}
\begin{example}
Derivation of entry (2.15.20.11) in \cite{prud2}. Here we use equation (\ref{eq:thm_inf}) and set $\tau \to 2,\theta \to p,m\to -1,\mu \to \nu ,r\to 1,s\to 1,\alpha \to c,\beta \to c$ and simplify;
\begin{multline}
\int_0^{\infty } \frac{e^{-p x^2} I_{\nu }(c x){}^2}{x} \, dx\\
=\sum _{g=0}^{\infty } \frac{2^{-1-2 g-2
   \nu } c^{2 g+2 \nu } p^{-g-\nu } \, _1F_1\left(g+\nu ;1+\nu ;\frac{c^2}{4 p}\right)}{(g+\nu ) \Gamma (1+g)
   \Gamma (1+\nu )}\\
=\frac{\exp \left(\frac{c^2}{2 p}\right) \left(I_{\nu }\left(\frac{c^2}{2 p}\right)+2 \sum
   _{k=1}^{\infty } (-1)^k I_{\nu +k}\left(\frac{c^2}{2 p}\right)\right)}{2 \nu }
\end{multline}
where $Re(\nu+\mu)>-1,Re(p)>0$.
\end{example}
\begin{example}
 Derivation of entry (2.15.20.12) in \cite{prud2}. Here we use equation (\ref{eq:thm_inf}) and set $\tau \to 2,\theta \to p,m\to 1,\mu \to \nu ,r\to 1,s\to 1,\alpha \to c,\beta \to c$ and simplify;
 \begin{multline}
\int_0^{\infty } e^{-p x^2} x I_{\nu }(c x){}^2 \, dx\\
=\sum _{g=0}^{\infty } \frac{2^{-1-2 g-2 \nu } c^{2
   g+2 \nu } p^{-1-g-\nu } \, _1F_1\left(1+g+\nu ;1+\nu ;\frac{c^2}{4 p}\right)}{\Gamma (1+g) \Gamma (1+\nu
   )}=\frac{\exp \left(\frac{c^2}{2 p}\right) I_{\nu }\left(\frac{c^2}{2 p}\right)}{2 p}
\end{multline}
where $Re(\nu+\mu)>-1,Re(p)>0$.
\end{example}
\begin{example}
Derivation of entry (6.626.4) in \cite{grad}. Here we use equation (\ref{eq:thm_inf}) and set $\tau \to 1,\theta \to 2 \alpha ,m\to 1,\mu \to 0,\alpha \to \beta ,\nu \to 1,r\to 1,s\to 1$ and simplify;
\begin{multline}
\int_0^{\infty } e^{-2 x \alpha } x I_0(x \beta ) I_1(x \beta ) \, dx\\
=\sum _{g=0}^{\infty }
   \frac{2^{-2-2 g} \alpha ^{-3-2 g} \beta ^{1+2 g} \Gamma \left(\frac{1}{2} (3+2 g)\right) \,
   _2F_1\left(\frac{3}{2}+g,2+g;1;\frac{\beta ^2}{4 \alpha ^2}\right)}{\sqrt{\pi } \Gamma
   (1+g)}\\
=\frac{\frac{\alpha  E\left(\left(\frac{\beta }{\alpha }\right)^2\right)}{\alpha ^2-\beta
   ^2}-\frac{K\left(\left(\frac{\beta }{\alpha }\right)^2\right)}{\alpha }}{2 \pi  \beta }
\end{multline}
where $Re(\alpha)>0$.
\end{example}
\section{Conclusion}
In this paper, we have presented a method for deriving definite integrals involving the product of Bessel functions used in diffraction theory. These integral forms were first investigated by Carslaw with other researchers presenting extended methods and derivations of these formulae. The results presented were numerically verified for both real and imaginary and complex values of the parameters in the integrals using Mathematica by Wolfram.
\end{document}